\documentclass[11pt]{article}
\usepackage[a4paper,margin=1in]{geometry}
\usepackage{amsmath,amssymb,amsthm,bm,mathtools}
\usepackage{booktabs,longtable,array,graphicx,flafter}
\usepackage{setspace,natbib,needspace}
\usepackage[T1]{fontenc}
\usepackage{lmodern,microtype}
\usepackage{xcolor,caption,threeparttable}
\definecolor{linkblue}{RGB}{0,73,128}
\usepackage[colorlinks=true,linkcolor=linkblue,citecolor=linkblue,
  urlcolor=linkblue,bookmarks=false]{hyperref}
\hypersetup{pdftitle={Two-sample Behrens-Fisher problems for high-dimensional data: noncentral normal-reference inference},pdfauthor={Pengfei Wang}}
\allowdisplaybreaks
\theoremstyle{plain}
\newtheorem{theorem}{Theorem}[section]
\newtheorem{prop}[theorem]{Proposition}
\newtheorem{lem}[theorem]{Lemma}
\newtheorem{cor}[theorem]{Corollary}
\newtheorem{assumption}[theorem]{Assumption}
\theoremstyle{remark}
\newtheorem{remark}[theorem]{Remark}

\newcommand{\E}{\mathbb E}
\newcommand{\Var}{\operatorname{Var}}
\newcommand{\Cov}{\operatorname{Cov}}
\newcommand{\tr}{\operatorname{tr}}
\newcommand{\diag}{\operatorname{diag}}

\newcommand{\R}{\mathbb R}

\title{Two-sample Behrens--Fisher problems for high-dimensional data: noncentral normal-reference inference}
\author{Pengfei Wang\\[3pt]
\small Mathematics and Mathematics Education, National Institute of Education,\\
\small Nanyang Technological University, Singapore\\[2pt]
\small \texttt{nie23.wp8738@e.ntu.edu.sg}}
\date{}
\begin{document}
\maketitle
\begin{abstract}
We study noncentral normal-reference inference for high-dimensional two-sample Behrens--Fisher problems. Within the Gaussian reference, a shifted null law is accurate if and only if covariance-weighted signal energy is negligible relative to the trace of the squared covariance. We establish the sharp uniform fourth-root error order over unrestricted spectra and a linear bound under spectral spread. A separate spectral-range criterion determines when signal magnitude alone suffices. An augmented spherical reference gives finite-sample calibration for Gaussian samples and asymptotically recovers Gaussian-reference local power under projection fourth-moment bounds and control of diagonal fluctuations, without a dimension-to-sample-size restriction. For unknown Gaussian covariance, explicit trace and hollow-Gram confidence bounds turn the linear approximation bound into a simultaneous, computable certificate over local effect ranges. Eigenspace signal energies describe directional power, while a two-point argument limits estimation of unknown local power. Simulations with 10,000 evaluation replications per configuration examine unequal and noncommuting covariances, severe sample imbalance, non-Gaussian observations and the informativeness of the certificate.
\end{abstract}

\section{Introduction}
Testing the equality of two mean vectors without assuming equal covariance matrices is the two-sample Behrens--Fisher problem. For high-dimensional data, the dimension may be comparable to or much larger than the sample sizes, so that procedures involving covariance inversion can become unstable or undefined. Tests based on squared mean norms provide an alternative \citep{BaiSaranadasa1996,ChenQin2010}. Diagonal standardization leads to scale-invariant procedures \citep{SrivastavaDu2008,SrivastavaEtAl2013}. For these tests, the approximation to the sampling distribution is important for both size control and power assessment.

The normal-reference approach uses the distribution of a statistic under Gaussian observations with the same means and covariances. For a squared mean norm, the null reference is a weighted chi-squared distribution, which can accommodate both normal and nonnormal limits \citep{ZhangEtAl2020,ZhangEtAl2021,ZhangEtAl2022}. Subsequent developments include three-cumulant approximations \citep{ZhangZhuBS2022,ZhangZhuCQ2022}, scale-invariant tests \citep{ZhangZhuZhang2023,CaoEtAl2024}, and an F-type test \citep{ZhuWangZhang2024}. These methods and their implementations in the R package HDNRA are reviewed by \citet{WangZhuZhang2026}. The present paper studies the use of the normal reference under unequal population means, with particular attention to the information needed to approximate power.

To explain the problem, let $n_1,n_2$ be the sample sizes, $\bm{\mu}_1,\bm{\mu}_2$ the population means, and $\bar{\bm{Y}}_1,\bar{\bm{Y}}_2$ the sample means. Write
\[
 \widehat{\bm{\theta}}=\sqrt{\frac{n_1n_2}{n_1+n_2}}(\bar{\bm{Y}}_1-\bar{\bm{Y}}_2),
 \qquad
 \bm{\theta}=\sqrt{\frac{n_1n_2}{n_1+n_2}}(\bm{\mu}_1-\bm{\mu}_2).
\]
The normalized mean difference has the representation
\begin{equation}\label{eq:intro-representation}
 \widehat{\bm{\theta}}=\bm{\theta}+\bm{W},\qquad \E\bm{W}=\bm{0},\qquad
 \Cov(\bm{W})=\bm{\Omega},\qquad T=\|\widehat{\bm{\theta}}\|^2.
\end{equation}
Vectors and matrices are written in bold; scalar coordinates, eigenvalues and statistics are written in ordinary type. The statistic satisfies
\begin{equation}\label{eq:intro-decomposition}
 T-\tr(\bm{\Omega})
 =\{\|\bm{W}\|^2-\tr(\bm{\Omega})\}+2\bm{\theta}^\top\bm{W}+\|\bm{\theta}\|^2.
\end{equation}
Under the null hypothesis, $\bm{\theta}=\bm{0}$ and only the quadratic fluctuation remains. Under an alternative, the signal--error interaction $2\bm{\theta}^\top\bm{W}$ may have the same order as that fluctuation. Thus, shifting the null distribution by the squared signal size need not give an adequate power approximation.

Existing normal-reference power analyses control this interaction through local-alternative conditions; see \citet[][Section 2.3]{ZhangEtAl2021} and \citet[][Section 2.3]{ZhuWangZhang2024}. A corresponding condition appears in \citet[][Corollary 2]{WangXu2022}. In the normalization above, it takes the form $\bm{\theta}^\top\bm{\Omega}\bm{\theta}=o\{\tr(\bm{\Omega}^2)\}$. We establish that this condition is necessary and sufficient for the shifted-null approximation within the covariance-matched Gaussian reference, without requiring a limiting covariance spectrum. The statement concerns the whole distribution, or equivalently power over all significance levels; necessity does not imply failure at every fixed level. We also obtain the sharp uniform fourth-root error order over unrestricted spectra, a three-quarter-order upper bound when three spectral blocks retain nonnegligible variance, and a linear bound under stronger spectral spread. Uniform validity over bounded local signal balls is characterized by spectral diffuseness.

When the shifted-null approximation is inadequate, a natural question is whether another reference based only on signal magnitude can be used. At fixed positive local signal strength, we show that the answer is determined by the normalized spectral range. In growing dimension, the resulting condition again reduces to spectral diffuseness. In fixed dimension, it also permits isotropic noncentral laws which need not be shifted null laws. When signal direction matters, the relevant quantities are its energies in the distinct covariance eigenspaces. A joint quadratic--linear spectral limit describes the resulting local power. This connects with the centered spectral limit of \citet[][Supplement, S.1]{WangXu2022} and complements studies of strongly spiked spectra, heterogeneous factor models and structured precision matrices \citep{AoshimaYata2018,HyodoEtAl2023,YangEtAl2024}. A separate two-point argument shows that even known covariance does not generally permit consistent estimation of power at an unknown local alternative. Directional power ranges are therefore interpreted as sensitivity analyses at specified effects.

To conduct the test with unknown covariance, we use an augmented spherical reference constructed from the observed mean contrast and matched residual columns. The Gaussian residual reduction is related to \citet{Bennett1951} and its high-dimensional use by \citet{NishiyamaEtAl2013}; finite-sample validity follows from the conditional rotation principle of \citet{Langsrud2005}. Our analysis connects this feasible reference to the Gaussian population law and its local power, without an effective-dimension restriction. The same reference is asymptotically justified under projection fourth-moment bounds and control of diagonal fluctuations, although the non-Gaussian residuals need not be independent. For Gaussian observations, trace and hollow-Gram confidence bounds further make it possible to assess the shifted-null approximation with unknown covariance. Combined with the linear error bound, they give a simultaneous confidence certificate over local effect ranges.

It is useful to distinguish these conclusions ..." with this paragraph.
It is useful to distinguish these conclusions from the approximation tools used to establish them. The Gaussian comparison specializes the independent-polynomial theory of \citet{HuangAusternOrbanz2025}, which also treats imbalanced two-sample MMD; related quadratic--linear limits appear in \citet{HuangEtAl2023}. Gaussian-ball comparison and smooth quadratic-form perturbation bounds are developed by \citet{GoetzeEtAl2019} and \citet{Kosov2022}. The general bound of \citet{Kosov2021} supplies our unrestricted upper rate, while spectral refinements are related to the anti-concentration theory of \citet{ChoudharyKuchibhotla2026}. Here the sharper shifted-null bounds exploit the centered perturbation supplied by a classical Poisson-mixture identity. Existing methods already provide feasible nonnormal calibration \citep{WangXu2022,CaiHu2026}. Our emphasis is on the approximation boundaries, the role of signal direction, and their connection to feasible calibration and diagnosis. Appendix~\ref{app:literature-comparison} distinguishes these questions from the Gaussian-vector and Gaussian-chaos results of \citet{Imai2026} and \citet{CaiHu2026}.

Simulation studies examine the theoretical distinctions through covariance spectra, signal directions, unequal and noncommuting covariance matrices, and unbalanced sample sizes. They also assess non-Gaussian calibration and the practical informativeness of the confidence certificate. The distinction between assumptions is maintained throughout: both moment conditions are used for asymptotic reference comparison and spherical calibration, whereas Gaussian sampling is required for the exact test and the finite-sample certificate.

The rest of the paper is organized as follows. Section~\ref{sec:setup} introduces the two-sample statistic, its normal reference and the feasible test. Section~\ref{prof:section} presents the approximation boundaries, directional power results and unknown-covariance diagnostic. Simulation studies and a discussion are given in Sections~\ref{sec:numerics} and~\ref{sec:discussion}, respectively. Additional results, simulation details and technical proofs are collected in the appendices.

\section{Test statistic and calibration}\label{sec:setup}
\subsection{Test statistic and normal-reference approximation}
\label{sec:uniform-transfer}
Let $\bm{Y}_{ij}\in\R^p$, $j=1,\ldots,n_i$, $i=1,2$, be independent observations, identically distributed within each sample, with mean $\bm{\mu}_i$ and covariance $\bm{\Sigma}_i$. The two covariance matrices may differ and may be singular. We consider
\[
 H_0:\bm{\mu}_1=\bm{\mu}_2
 \qquad\text{against}\qquad
 H_1:\bm{\mu}_1\ne\bm{\mu}_2.
\]
Write $n=n_1+n_2$, $\bar{\bm{Y}}_i=n_i^{-1}\sum_{j=1}^{n_i}\bm{Y}_{ij}$, and define
\[
 \widehat{\bm{\theta}}=\sqrt{\frac{n_1n_2}{n}}
 (\bar{\bm{Y}}_1-\bar{\bm{Y}}_2),\qquad
 \bm{\theta}=\sqrt{\frac{n_1n_2}{n}}(\bm{\mu}_1-\bm{\mu}_2),
 \qquad \bm{W}=\widehat{\bm{\theta}}-\bm{\theta}.
\]
The covariance contributions from the two samples are
\begin{equation}\label{eq:group-covariance}
 \bm{\Omega}_1=\frac{n_2}{n}\bm{\Sigma}_1,\qquad
 \bm{\Omega}_2=\frac{n_1}{n}\bm{\Sigma}_2,\qquad
 \bm{\Omega}=\Cov(\bm{W})=\bm{\Omega}_1+\bm{\Omega}_2.
\end{equation}
We use the squared-norm statistic
\[
 T=\|\widehat{\bm{\theta}}\|^2
   =\frac{n_1n_2}{n}\|\bar{\bm{Y}}_1-\bar{\bm{Y}}_2\|^2.
\]
Its decomposition in \eqref{eq:intro-decomposition} separates the quadratic sampling fluctuation, the signal--error interaction and the squared signal size. Under $H_0$, the interaction term vanishes. Under $H_1$, it can remain of the same order as the quadratic fluctuation, which motivates retaining the full noncentral reference. The term \emph{noncentral} refers to the distribution under a nonzero mean difference, not to whether $T$ has been centered by its null expectation.

The Gaussian reference replaces the centered two-sample error by a Gaussian vector with the same covariance.
Let \(\bm{Z}\sim N(\bm{0},\bm{\Omega})\), and define
\begin{equation}\label{eq:gaussian-reference-setup}
 T^*=\|\bm{\theta}+\bm{Z}\|^2,\qquad
 \mu_0=\tr(\bm{\Omega}),\quad A=\tr(\bm{\Omega}^2)>0,\quad
 s_0^2=2A,\quad s_{\bm{\theta}}^2=2A+4\bm{\theta}^\top\bm{\Omega}\bm{\theta}.
\end{equation}
The null scale \(s_0\) is used for local-power limits. The full reference standard deviation \(s_{\bm{\theta}}\) is used for uniform distributional comparison. Neither scale is asserted to equal the exact non-Gaussian variance in finite samples.

If \(\lambda_j>0\) are the nonzero eigenvalues of \(\bm{\Omega}\), \(\bm{e}_j\) are associated orthonormal eigenvectors, and \(\bm{\theta}_\perp\) is the projection of \(\bm{\theta}\) onto the null space of \(\bm{\Omega}\), then
\begin{equation}\label{eq:noncentral-representation}
 T^*\ \stackrel{d}{=}\ \|\bm{\theta}_\perp\|^2+
 \sum_j\lambda_j\left(Z_j+\frac{\bm{e}_j^\top\bm{\theta}}{\sqrt{\lambda_j}}\right)^2.
\end{equation}
Here the $Z_j$ are independent standard normal variables. Thus, for \(r\ge1\), its cumulants are
\[
 K_r(T^*)=\mathbf1_{\{r=1\}}\|\bm{\theta}_\perp\|^2+
 2^{r-1}(r-1)!\sum_j\lambda_j^r
 \left\{1+r\frac{(\bm{e}_j^\top\bm{\theta})^2}{\lambda_j}\right\}.
\]
The representation and cumulants are standard properties of noncentral
chi-squared mixtures; see \citet{Zhang2005}.\footnote{The correction by
\citet{Duivesteijn2016} concerns constants in an approximation bound;
the distributional identity and cumulant formulas used here are unchanged.}

The oracle null critical value \(q_{\alpha}^{\mathrm{NR}}\) is the \((1-\alpha)\)-quantile of \(\|\bm{Z}\|^2\). It depends on \(\bm{\Omega}\) and does not depend on the unknown alternative \(\bm{\theta}\). Noncentral distributions describe power; they are not used to estimate an alternative-dependent rejection threshold.

The reference is useful for non-Gaussian observations only when its distribution approximates that of $T$. Write $\bm{\xi}_{ij}=\bm{Y}_{ij}-\bm{\mu}_i$. A sufficient two-sample condition is
\begin{equation}
 \E(\bm{a}^\top\bm{\xi}_{ij})^4\le K(\bm{a}^\top\bm{\Sigma}_i\bm{a})^2,
 \qquad
 \Var(\|\bm{\xi}_{ij}\|^2)\le K\tr(\bm{\Sigma}_i^2),
 \label{eq:raw-group-moments}
\end{equation}
for every deterministic $\bm{a}\in\R^p$, both samples and all observations, with a common finite constant $K$. The first inequality bounds projection kurtosis; the second controls diagonal fluctuations and does not follow from the first. Both hold for linear transformations of independent standardized coordinates with uniformly bounded fourth moments.

Put $n_{\min}=\min(n_1,n_2)$ and let $d_{\mathrm K}$ denote Kolmogorov distance. Corollary~\ref{cor:grouped-transfer} gives
\begin{equation}\label{eq:main-two-sample-comparison}
 \sup_{\bm{\theta}\in\R^p}d_{\mathrm K}(T,T^*)
 \le \min\{1,C_K(n_{\min}^{-1/5}+n_{\min}^{-1/14})\}.
\end{equation}
Uniformity concerns location shifts of fixed centered sample distributions. Thus the reference comparison is valid as both sample sizes increase, without requiring a fixed sample-size ratio or diffuse covariance eigenvalues. The bound is a specialization of Gaussian polynomial comparison; sharper bounds and the full argument are given in Appendices~\ref{app:transfer}--\ref{app:adaptive}. The next subsection gives an operational reference and distinguishes its exact Gaussian validity from its asymptotic validity under these moment conditions.

\subsection{Feasible reference calibration}
\label{sec:spectral-calibration}
Under Gaussian sampling, the residuals permit a calibration that avoids
estimating a population spectrum on an operator-norm scale. The construction uses the Gaussian Behrens--Fisher reduction of \citet{Bennett1951} and the conditional rotation principle of \citet{Langsrud2005}; see also \citet{NishiyamaEtAl2013}.
Its finite-sample validity is a randomization property, whereas its
asymptotic connection to the central normal reference is a separate statement.

Write $\bm Y_i=(\bm Y_{i1},\ldots,\bm Y_{in_i})$. For each sample, let $\bm H_i$ be an $n_i\times(n_i-1)$ matrix with
$\bm H_i^\top\bm H_i=\bm I$ and $\bm H_i^\top\bm1=\bm0$.
Choose these matrices deterministically before observing the data. Use the residual contrasts
$\bm V_{ir}=\sqrt{w_i}\bm Y_i\bm H_i\bm e_r$, where
$w_1=n_2/(n_1+n_2)$ and $w_2=n_1/(n_1+n_2)$.
Set $m=\min(n_1-1,n_2-1)$, $N=m+1$, and form
\begin{equation}
 \bm X=(\bm X_0,\ldots,\bm X_m),\qquad
 \bm X_0=\widehat{\bm\theta},\qquad
 \bm X_r=\bm V_{1r}+\bm V_{2r}\quad(1\le r\le m).
 \label{eq:spherical-columns}
\end{equation}
Under Gaussian sampling the columns are independent, with
$\bm X_0\sim N(\bm\theta,\bm\Omega)$ and
$\bm X_r\sim N(\bm0,\bm\Omega)$ for $r\ge1$.
In particular, all $N$ columns have the same centered law under $H_0$,
despite the two original populations having different covariance matrices.
All observations enter $\bm X_0$; only the reference construction uses
the first $m$ residual columns from each group.

Let $\bm U$ be uniform on the unit sphere in $\mathbb R^N$, independent
of the data. The conditional spherical reference is
\begin{equation}
 R^{\mathrm{sph}}=\|\bm X\bm U\|^2
 \ \stackrel{d\mid\bm X}{=}
 \frac{\sum_{j=1}^N d_j Z_j^2}{\sum_{j=1}^N Z_j^2},
 \label{eq:spherical-reference}
\end{equation}
where $d_1,\ldots,d_N$ are the eigenvalues of the Gram matrix
$\bm X^\top\bm X$, including zeros, and the $Z_j$ are independent
standard normal variables. Thus its spectral weights are normalized
jointly, rather than being used in an unnormalized chi-squared mixture.
No spectral threshold or variance projection is needed. For numerical
stability, put $\bar d=N^{-1}\sum_jd_j$ and compute
$R^{\mathrm{sph}}-\bar d=
\sum_j(d_j-\bar d)Z_j^2/\sum_jZ_j^2$.

With $B$ independent draws $R_b^{\mathrm{sph}}$ from this reference,
use the Monte Carlo $p$-value
\begin{equation}
 p_B=\frac{1+\sum_{b=1}^B
       \boldsymbol1\{R_b^{\mathrm{sph}}\ge T\}}{B+1}.
 \label{eq:spherical-monte-carlo-p}
\end{equation}

\begin{theorem}[Gaussian spherical calibration]
\label{thm:spherical-calibration}
Suppose both samples are Gaussian, $n_1,n_2\ge2$, and
$A=\tr(\bm\Omega^2)>0$. The covariance matrices may be unequal
or singular, and the dimension is unrestricted.
Under $H_0$, for every $B\ge1$ and $\alpha\in(0,1)$,
\begin{equation}
 \Pr\{p_B\le\alpha\mid\bm X\bm X^\top\}\le\alpha
 \quad\text{almost surely}.
 \label{eq:spherical-finite-level}
\end{equation}
If the conditional reference is continuous, the left side is
$\lfloor(B+1)\alpha\rfloor/(B+1)$.

For the asymptotic reference statement, let $F_{\bm\Omega}$ be the
CDF of $\|\bm G\|^2$, $\bm G\sim N(\bm0,\bm\Omega)$,
and let $\widehat F_{\bm X}$ be the conditional CDF in
\eqref{eq:spherical-reference}. For every fixed $K<\infty$,
there is $C_K<\infty$ such that, uniformly over all dimensions,
nonzero positive semidefinite $\bm\Omega$, and
$\|\bm\theta\|^2/\sqrt{2A}\le K$,
\begin{equation}
 \E_{\bm\theta}
 \|\widehat F_{\bm X}-F_{\bm\Omega}\|_\infty
 \le C_K N^{-1/15}.
 \label{eq:spherical-reference-consistency}
\end{equation}
The displayed rate is a sufficient bound from coupling and a quantile
grid, and is not asserted to be sharp.
Consequently, if $\min(n_1,n_2)\to\infty$ and $B\to\infty$,
then, for each fixed $\alpha\in(0,1)$,
\begin{equation}
 \sup_{\|\bm\theta\|^2\le K\sqrt{2A}}
 \left|\Pr_{\bm\theta}(p_B\le\alpha)
  -\Pr\{\|\bm G+\bm\theta\|^2>
                   q_{\alpha}^{\mathrm{NR}}\}\right|\longrightarrow0.
 \label{eq:spherical-local-power}
\end{equation}
No condition of the form $\tr(\bm\Omega)/(m\sqrt A)\to0$
is required. The convergence statement is along any covariance and
dimension sequence satisfying the displayed assumptions.
\end{theorem}
\begin{proof}See Appendix~\ref{app:spherical-proof}.\end{proof}

The observed contrast $\bm X_0$ must be included in the rotation matrix;
rotating residuals alone does not establish \eqref{eq:spherical-finite-level}.
The residual bases and matching are fixed before observing the data.
Discarding excess residual contrasts can have a finite-sample efficiency
cost, especially with unequal sample sizes.
The simulations use a fixed auxiliary draw count. Its binomial power
expression and an explicit comparison with the ideal conditional
$p$-value are given in Appendix~\ref{par:spherical-finite-B-power};
the local-power limit above requires $B\to\infty$.

The same calibration has an asymptotic justification beyond Gaussian
sampling. This does not follow from independence of the residual columns,
which generally fails. Instead, a two-reference-draw comparison controls
the conditional distribution directly.

\begin{theorem}[Spherical calibration under moment conditions]
\label{thm:nongaussian-spherical-calibration}
Suppose the observations are independent, are identically distributed
within each sample, and satisfy both inequalities in
\eqref{eq:raw-group-moments} with a common constant $K_0<\infty$.
Use the fixed residual bases, augmented matrix $\bm X$ and Monte Carlo
$p$-value in \eqref{eq:spherical-columns}--\eqref{eq:spherical-monte-carlo-p}.
Put $N=\min(n_1,n_2)$ and suppose $A=\tr(\bm\Omega^2)>0$.
For every fixed $K_1<\infty$, there is $C_{K_0,K_1}<\infty$ such that
\begin{equation}
 \sup_{\|\bm\theta\|^2\le K_1\sqrt{2A}}
 \E_{\bm\theta}\|\widehat F_{\bm X}-F_{\bm\Omega}\|_\infty
 \le \min\{1,C_{K_0,K_1}N^{-1/42}\}.
 \label{eq:nongaussian-spherical-cdf}
\end{equation}
The constant and rate are uniform over all dimensions, unequal or singular
covariances, and fixed choices of the residual bases. As $N\to\infty$ and
$B\to\infty$, for every fixed $\alpha\in(0,1)$,
\begin{equation}
 \sup_{\|\bm\theta\|^2\le K_1\sqrt{2A}}
 \left|\Pr_{\bm\theta}\{p_B\le\alpha\}
       -\Pr\{\|\bm G+\bm\theta\|^2>q_\alpha^{\mathrm{NR}}\}\right|
 \longrightarrow0,
 \qquad \bm G\sim N(\bm0,\bm\Omega).
 \label{eq:nongaussian-spherical-power}
\end{equation}
In particular, the null rejection probability tends to $\alpha$.
There is no restriction relating $p$ to $N$, and no sample-size balance
condition. Uniformity in $\bm\theta$ concerns location shifts of fixed
centered sample distributions satisfying the displayed assumptions.
\end{theorem}
\begin{proof}See Appendix~\ref{app:nongaussian-spherical-proof}.\end{proof}

The theorem retains both moment conditions: projection kurtosis alone is
insufficient to control fluctuations of the diagonal quadratic terms.
The finite-sample rotation identity in
\eqref{eq:spherical-finite-level} remains a Gaussian assertion.
The bound \eqref{eq:nongaussian-spherical-cdf} is a theoretical bound with
an unspecified constant, not a numerically calibrated finite-sample
certificate. Its exponent is a consequence of smoothing and a quantile
grid and is not asserted to be optimal. The construction uses the original
observations and the same residual matrix as in the Gaussian case;
block lengths, a separate estimation sample and spectral thresholds are
unnecessary for this result.

The projected Gaussian and gamma calibrations in
Appendix~\ref{app:moment-calibration-summary} provide useful comparisons.
Their sufficient estimation condition
$\tr(\bm\Omega)/(m_{\min}\sqrt A)\to0$ is not required by either spherical
calibration theorem. The non-Gaussian result above is asymptotic, whereas
the finite-sample assertion of Theorem~\ref{thm:spherical-calibration}
uses Gaussian rotation invariance. Approximate randomization under
general covariance spectra is already developed by \citet{WangXu2022};
the present statement connects the particular augmented spherical
reference to the full quadratic normal reference and its noncentral
power under the stated moment class.

\section{Noncentral normal-reference inference}
\label{prof:section}

We now study the distribution of the normal reference under the alternative.
Two questions arise: when is a shift of the null distribution sufficient,
and when must the signal direction also be retained? We first characterize
these approximation boundaries and then derive the corresponding power
statements. The centered Gaussian-plus-chi-squared limit is given by
\citet{WangXu2022}; Theorem~\ref{prof:theorem} describes the contribution
of the signal coordinates. Corollary~\ref{cor:grouped-transfer} transfers
the reference to non-Gaussian observations under the stated moment conditions.
The final part of this section considers how the approximation error can
be assessed when the covariance is unknown.

\subsection{Approximation boundaries}

\begin{theorem}[Exact boundary for a shifted null reference]
\label{thm:shift-boundary}
Let $\bm G\sim N(\bm0,\bm\Omega)$, where
$\bm\Omega\succeq0$ and $A=\tr(\bm\Omega^2)>0$. Define
\[
 T_{\bm\theta}^*=\|\bm G+\bm\theta\|^2,\qquad
 T_0^*=\|\bm G\|^2,\qquad
 \Delta_{\bm\theta}=\frac{\bm\theta^\top\bm\Omega\bm\theta}{A},\qquad
 d_{\bm\theta}=d_{\mathrm K}
       (T_{\bm\theta}^*,T_0^*+\|\bm\theta\|^2).
\]
There is an absolute constant $C$ such that, for every dimension,
nonzero covariance matrix and deterministic signal,
\begin{equation}
 \frac{\Delta_{\bm\theta}^2}
 {60+120\Delta_{\bm\theta}+24\Delta_{\bm\theta}^2}
 \ \leq\ d_{\bm\theta}\ \leq\
 \min\{1/2,C\Delta_{\bm\theta}^{1/4}\}.
 \label{eq:shift-boundary-bounds}
\end{equation}
If the eigenvalues admit three disjoint blocks, each containing at least
$cA$ of the total squared-eigenvalue mass, for a fixed $c>0$, then
\begin{equation}
 d_{\bm\theta}\leq\min\{1/2,C_c\Delta_{\bm\theta}^{3/4}\}.
 \label{eq:shift-boundary-spread}
\end{equation}
One may take $C_c=8\sqrt2/(3\pi c^{3/4})$.
Consequently, along any sequence of covariance matrices and signals,
\begin{equation}
 d_{\bm\theta_n}\longrightarrow0
 \quad\Longleftrightarrow\quad
 \frac{\bm\theta_n^\top\bm\Omega_n\bm\theta_n}
      {\tr(\bm\Omega_n^2)}\longrightarrow0.
 \label{eq:shift-boundary-equivalence}
\end{equation}
No spectral limit or bound on
$\|\bm\theta_n\|^2/\sqrt{2\tr(\bm\Omega_n^2)}$ is required.
The upper bound with exponent $1/4$ follows from the polynomial comparison
bound of \citet{Kosov2021}. Proposition~\ref{prop:shift-quarter-sharp} establishes that its exponent is optimal uniformly for the present family. The three-block exponent is a sufficient bound and is not asserted to be optimal.
\end{theorem}
\begin{proof}See Appendix~\ref{app:shift-boundary-proofs}.\end{proof}

\begin{prop}[Sharp uniform order of the shifted-null error]
\label{prop:shift-quarter-sharp}
Fix $\kappa_0>0$. There is a family of positive definite
two-dimensional covariance matrices $\bm\Omega_\delta$ and signals
$\bm\theta_\delta$ for which
\[
 \tr(\bm\Omega_\delta^2)=1,\qquad
 \frac{\|\bm\theta_\delta\|^2}{\sqrt2}=\kappa_0,\qquad
 \Delta_{\bm\theta_\delta}=\delta\longrightarrow0,
\]
and
\begin{equation}
 \liminf_{\delta\downarrow0}
 \frac{d_{\mathrm K}(T_{\bm\theta_\delta}^*,
              T_0^*+\|\bm\theta_\delta\|^2)}{\delta^{1/4}}
 \ge c_*:=\frac{2^{1/4}\Gamma(3/4)}{\pi}>0.
 \label{eq:shift-quarter-lower}
\end{equation}
Together with Theorem~\ref{thm:shift-boundary}, this establishes the
optimal uniform order $\Delta^{1/4}$ over unrestricted covariance
spectra and signal directions, even at a fixed positive local effect.
In particular, no bound $C_{\kappa_0}\Delta^{1/2}$, or more generally
$o(\Delta^{1/4})$, can hold uniformly on this class.
\end{prop}\begin{proof}See Appendix~\ref{app:sharp-quarter-proof}.\end{proof}

The signal in this construction lies in a shrinking-eigenvalue direction.
Its linear fluctuation perturbs the support boundary of a dominant
squared-normal term. Corollary~\ref{cor:shift-quarter-fixed-level}
shows the same fourth-root obstruction at any prescribed significance
level for a particular fixed effect magnitude; both powers tend to one
in that construction. The obstruction does not apply to every fixed
effect or to spectra bounded away from this boundary.

The condition that removes the signal--error interaction is therefore
necessary and sufficient for the entire shifted-null approximation,
without requiring a limiting spectral profile. This is stronger than
requiring that the interaction be negligible as a convenient sufficient
assumption. It does not assert necessity for power approximation at a
single prescribed significance level.

\begin{cor}[Uniformity over directions and significance levels]
\label{cor:shift-uniform-directions}
Let $q_{\alpha}^{\mathrm{NR}}$ be the exact $(1-\alpha)$ quantile of $T_0^*$.
For every $\bm\theta$,
\begin{equation}
 \sup_{0<\alpha<1}\left|
  \Pr\{T_{\bm\theta}^*>q_{\alpha}^{\mathrm{NR}}\}
 -\Pr\{T_0^*+\|\bm\theta\|^2>q_{\alpha}^{\mathrm{NR}}\}
 \right|=d_{\bm\theta}.
 \label{eq:shift-power-supremum}
\end{equation}
For fixed $K>0$, put
\[
 \Theta_n(K)=\left\{\bm\theta:
       \|\bm\theta\|^2\leq K\sqrt{2A_n}\right\},\qquad
 D_n(K)=\sup_{\bm\theta\in\Theta_n(K)}d_{\bm\theta},\qquad
 r_n=\sqrt2K\frac{\lambda_{\max}(\bm\Omega_n)}{\sqrt{A_n}}.
\]
Then
\begin{equation}
 \frac{r_n^2}{60+120r_n+24r_n^2}
 \leq D_n(K)\leq\min\{1/2,Cr_n^{1/4}\},
 \label{eq:shift-uniform-bounds}
\end{equation}
and hence
\begin{equation}
 D_n(K)\longrightarrow0
 \quad\Longleftrightarrow\quad
 \frac{\lambda_{\max}(\bm\Omega_n)^2}
      {\tr(\bm\Omega_n^2)}\longrightarrow0.
 \label{eq:shift-uniform-equivalence}
\end{equation}
Under uniform Gaussian-reference transfer on $\Theta_n(K)$, the same
if-and-only-if conclusion holds when the first statistic in the
Kolmogorov distance or in \eqref{eq:shift-power-supremum} is replaced by
the original-data statistic, with exact central-reference critical values.
\end{cor}
\begin{proof}See Appendix~\ref{app:shift-boundary-proofs}.\end{proof}

Thus spectral diffuseness has a precise role: it is equivalent to
validity of the shifted-null approximation uniformly over bounded local
signal directions and all significance levels. A persistent eigenvalue
precludes that uniform approximation. For an individual signal, however,
only its covariance-weighted energy
$\bm\theta^\top\bm\Omega\bm\theta$ matters; a persistent spike need not
invalidate the approximation when the signal avoids its eigenspace.
The supremum over significance levels is essential to this conclusion.
At a fixed level, two different distributions can have the same tail
probability, and an accidental equality does not validate the full
shifted-null reference.

The distributional bound takes a supremum over all significance levels.
Write $\Phi$ and $\phi$ for the standard normal CDF and density.
The next proposition gives an exact benchmark in which that supremum and
the error at a fixed level have different orders. It also identifies the
sign of the fixed-level error for a local signal.

\begin{prop}[Rank-one approximation and fixed-level power]
\label{prop:rankone-sharp-error}
Let $Z\sim N(0,1)$ and, for $a\in\mathbb R$, put
$X_a=(Z+a)^2$ and $Y_a=Z^2+a^2$. Then
\begin{equation}
 d_{\mathrm K}(X_a,Y_a)
 =\Phi(2|a|)-\frac12
 =\sqrt{\frac2\pi}\,|a|+O(|a|^3)
 \qquad(a\longrightarrow0).
 \label{eq:rankone-exact-distance}
\end{equation}
For $\alpha\in(0,1)$, let
$z_\alpha=\Phi^{-1}(1-\alpha/2)$,
$q_\alpha=z_\alpha^2$ and
\[
 \pi_\alpha(a)=\Pr(X_a>q_\alpha),\qquad
 \pi_\alpha^{\mathrm{shift}}(a)=\Pr(Y_a>q_\alpha).
\]
For any fixed $0<\alpha_0<\alpha_1<1$,
\begin{equation}
 \pi_\alpha(a)-\pi_\alpha^{\mathrm{shift}}(a)
 =a^2\phi(z_\alpha)
       \bigl(z_\alpha-z_\alpha^{-1}\bigr)+O(a^4),
 \label{eq:rankone-fixed-level}
\end{equation}
uniformly over $\alpha\in[\alpha_0,\alpha_1]$ as $a\to0$.
\end{prop}
\begin{proof}See Appendix~\ref{app:rankone-sharp-proof}.\end{proof}

This is the Gaussian-reference model with
$\bm\Omega=\lambda\bm e\bm e^\top$ and
$\bm\theta=\sqrt\lambda\,a\bm e$, where $\lambda>0$ and
$\|\bm e\|=1$; scaling by $\lambda$ gives
$\Delta_{\bm\theta}=a^2$.
Consequently, the rank-one Kolmogorov error has exact order
$\Delta_{\bm\theta}^{1/2}$, while the fixed-level error in
\eqref{eq:rankone-fixed-level} is of order
$\Delta_{\bm\theta}$ unless its leading coefficient vanishes.
For $a\ne0$, the distance in \eqref{eq:rankone-exact-distance} is attained
at $x=a^2$, corresponding to the null-reference level
$\alpha_a=2\{1-\Phi(|a|)\}$, which tends to one as $a\to0$.
At any fixed $\alpha<2\{1-\Phi(1)\}$, including the usual levels
$0.01$, $0.05$ and $0.10$, the shifted-null formula underestimates power
for all sufficiently small nonzero $a$. The sign is reversed for fixed
$\alpha>2\{1-\Phi(1)\}$.
These statements concern the specified rank-one model and do not assert
a universal ordering across spectra or signal directions. Proposition~\ref{prop:shift-quarter-sharp} gives the distinct shrinking-eigenvalue construction needed for sharpness of the general upper bound.
In the parallel rank-one model, $a\to0$ also implies $\kappa=a^2/\sqrt2\to0$; the small-signal expansion is not a sharp-rate assertion at fixed positive $\kappa$.
As $|a|\to\infty$, the exact distance tends to $1/2$, so the universal cap in Theorem~\ref{thm:shift-boundary} is sharp.

The general fourth-root rate reflects covariance sequences approaching a
low-rank boundary. A more informative bound is available when the normalized
leading eigenvalue is sufficiently small. Put
$A=\tr(\bm\Omega^2)>0$,
$r=\lambda_{\max}(\bm\Omega)/\sqrt A$, and
$\Delta_{\bm\theta}=\bm\theta^\top\bm\Omega\bm\theta/A$.
If $r<1/2$, then
\begin{equation}
 d_{\bm\theta}
 :=d_K\{\|\bm G+\bm\theta\|^2,
              \|\bm G\|^2+\|\bm\theta\|^2\}
 \le \min\left\{\frac12,
          \frac{\Delta_{\bm\theta}}{\pi(1-4r^2)}\right\},
 \qquad \bm G\sim N(\bm0,\bm\Omega).
 \label{eq:smooth-shift-linear-bound}
\end{equation}
The improvement follows from the classical Poisson-mixture representation
of a noncentral chi-squared variable; see, for example,
\citet{JonesMarchand2021}. It gives an independent, centered
perturbation of the central quadratic form, so a Taylor expansion of the
null CDF cancels the first-order term. Fourier inversion provides an
explicit bound for the derivative of its density. This distributional
representation does not make the quadratic and linear terms in the
original Gaussian decomposition independent.
 The full argument is given in Appendix~\ref{app:smooth-shift-certificate}.

The preceding result concerns a particular shifted reference. We next ask whether any reference given the covariance and signal norm, but not the signal direction, can approximate all alternatives of that norm.
\begin{theorem}[When a signal-norm-only reference is possible]
\label{thm:directional-diameter}
Let $\bm G\sim N(\bm0,\bm\Omega)$, where
$\bm\Omega\succeq0$ and $A=\tr(\bm\Omega^2)>0$. Write
\[
 U_{\bm\theta}
   =\frac{\|\bm G+\bm\theta\|^2-\tr(\bm\Omega)}{\sqrt{2A}},
 \qquad
 \mathcal S_\kappa
   =\{\bm\theta:\|\bm\theta\|^2=\kappa\sqrt{2A}\},
 \qquad \kappa>0.
\]
Include zero eigenvalues when defining
\[
 d(\bm\Omega)=
 \frac{\lambda_{\max}(\bm\Omega)-\lambda_{\min}(\bm\Omega)}{\sqrt A},
 \qquad
 D_\kappa(\bm\Omega)=
 \sup_{\bm\theta,\bm\vartheta\in\mathcal S_\kappa}
 d_{\mathrm K}(U_{\bm\theta},U_{\bm\vartheta}).
\]
There are constants $0<c_\kappa\leq C_\kappa<\infty$, independent of
$p$ and $\bm\Omega$, such that
\begin{equation}
 c_\kappa d(\bm\Omega)^2
 \leq D_\kappa(\bm\Omega)
 \leq \min\{1,C_\kappa d(\bm\Omega)^{1/4}\}.
 \label{eq:directional-diameter-bounds}
\end{equation}
Consequently, along any covariance sequence and for fixed $\kappa>0$,
\begin{equation}
 D_\kappa(\bm\Omega_n)\longrightarrow0
 \quad\Longleftrightarrow\quad
 d(\bm\Omega_n)\longrightarrow0.
 \label{eq:directional-diameter-boundary}
\end{equation}
The constants can be chosen uniformly for
$\kappa\in[\kappa_0,K]$, where $0<\kappa_0\leq K<\infty$.
If $p_n\to\infty$, the right-hand condition is equivalent to
$\lambda_{\max}(\bm\Omega_n)/\sqrt{A_n}\to0$.
No limiting spectrum is required.
\end{theorem}
\begin{proof}See Appendix~\ref{app:directional-boundary-proofs}.\end{proof}

\begin{cor}[An unavoidable error for direction-blind approximation]
\label{cor:direction-blind-reference}
In the setting of Theorem~\ref{thm:directional-diameter}, let
$F_{\bm\theta}$ denote the distribution function of $U_{\bm\theta}$.
Among all distribution functions $H$ allowed to depend on
$(\bm\Omega,\kappa)$ but not on the signal direction, define
\[
 R_\kappa(\bm\Omega)=
 \inf_H\sup_{\bm\theta\in\mathcal S_\kappa}
       \sup_{x\in\mathbb R}|F_{\bm\theta}(x)-H(x)|.
\]
Then
\begin{equation}
 R_\kappa(\bm\Omega)=\frac12D_\kappa(\bm\Omega).
 \label{eq:direction-blind-radius}
\end{equation}
In particular, uniformly accurate approximation using the signal norm
alone is possible if and only if the normalized spectral range tends to
zero. This statement allows the approximation to know the entire
population covariance matrix.

For the exact null-reference power functions
$\beta_{\bm\theta}(\alpha)
 =\Pr\{U_{\bm\theta}>c_\alpha\}$, where $c_\alpha$ is the
$(1-\alpha)$ quantile of $U_{\bm0}$,
\begin{equation}
 D_\kappa(\bm\Omega)=
 \sup_{\bm\theta,\bm\vartheta\in\mathcal S_\kappa}
 \sup_{0<\alpha<1}
 |\beta_{\bm\theta}(\alpha)-\beta_{\bm\vartheta}(\alpha)|.
 \label{eq:directional-power-diameter}
\end{equation}
The result concerns approximation of the entire distribution or power
curve. It does not assert a positive separation at every fixed
significance level or a common ordering of directions.
\end{cor}
\begin{proof}See Appendix~\ref{app:directional-boundary-proofs}.\end{proof}

\begin{remark}
The normalized spectral range and the normalized largest eigenvalue
answer different questions in finite dimension. If
$\bm\Omega=\lambda\bm I_p$ with fixed $p$, the exact noncentral law is
independent of direction even though the signal--error interaction need
not vanish. Thus a signal-norm-only reference can be valid while a
shifted-null reference is not. When $p\to\infty$, both uniform questions
reduce to spectral diffuseness. The exponents in
\eqref{eq:directional-diameter-bounds} are bounds obtained from moment
comparison and Gaussian polynomial stability; their optimality is not
claimed. The approximation criterion in
\eqref{eq:direction-blind-radius} is not a minimax detection result.
\end{remark}

\subsection{Signal information and directional power}
The boundary parameter $\Delta_{\bm\theta}$ determines whether the
signal--error interaction may be omitted, but it does not generally
determine the noncentral reference. The following consequence of the
Gaussian quadratic-form identity \citep{MathaiProvost1992} makes the required information precise.

\begin{prop}[Signal information in the Gaussian reference]
\label{prop:spectral-signal-information}
Fix a positive semidefinite $p\times p$ matrix $\bm\Omega$. Let
$\mathcal L$ be its set of distinct eigenvalues, including zero if
$\bm\Omega$ is singular, and let $\bm P_\lambda$ be the orthogonal
projector onto the eigenspace corresponding to $\lambda$. Put
\[
 e_\lambda(\bm\theta)=\|\bm P_\lambda\bm\theta\|^2,
 \qquad T_{\bm\theta}^*=\|\bm G+\bm\theta\|^2,
 \qquad \bm G\sim N(\bm0,\bm\Omega).
\]
For any deterministic $\bm\theta,\bm\vartheta$, the following are
equivalent:
\begin{enumerate}\setlength{\itemsep}{0pt}
 \item $T_{\bm\theta}^*$ and $T_{\bm\vartheta}^*$ have the same law;
 \item $e_\lambda(\bm\theta)=e_\lambda(\bm\vartheta)$ for every
       $\lambda\in\mathcal L$;
 \item $\bm\theta^\top\bm\Omega^j\bm\theta
       =\bm\vartheta^\top\bm\Omega^j\bm\vartheta$ for
       $j=0,\ldots,|\mathcal L|-1$, where $\bm\Omega^0=\bm I_p$.
\end{enumerate}
In particular, the reference depends on the signal energy in each
distinct eigenspace, but not on a choice of basis or orientation within
that eigenspace. For every integer $r\geq1$, its $r$th cumulant is
\begin{equation}
 \operatorname{cum}_r(T_{\bm\theta}^*)
 =2^{r-1}(r-1)!\left\{
   \tr(\bm\Omega^r)
   +r\bm\theta^\top\bm\Omega^{r-1}\bm\theta\right\}.
 \label{eq:noncentral-cumulants-information}
\end{equation}
\end{prop}
\begin{proof}See Appendix~\ref{app:spectral-information-proof}.\end{proof}

For example, take $\bm\Omega=\diag(1,2,3)$ and
\[
 \bm\theta_a=(1/\sqrt2,0,1/\sqrt2)^\top,
 \qquad \bm\theta_b=(0,1,0)^\top.
\]
Both signals have squared norm $1$ and covariance-weighted energy $2$,
so they give the same $\kappa$, $\Delta_{\bm\theta}$, mean and variance.
However, their values of $\bm\theta^\top\bm\Omega^2\bm\theta$ are
$5$ and $4$, respectively. Their third cumulants therefore differ by
$24$, and their reference distributions differ. Thus even knowledge of
the signal magnitude and $\Delta_{\bm\theta}$ does not generally
recover the noncentral law. More generally, a fixed number of
cumulants can discard signal information when the number of distinct
eigenvalues is unrestricted; Proposition~\ref{prop:finite-cumulant-information}
gives an explicit construction.

Along a sequence with $n=n_1+n_2$ and dimension $p=p_n$, write
\[
 \bm{G}_n\sim N(\bm{0},\bm{\Omega}_n),\qquad
 T_n^*=\|\bm{G}_n+\bm{\theta}_n\|^2,\qquad
 A_n=\operatorname{tr}(\bm{\Omega}_n^2)>0,\qquad
 s_{0,n}=\sqrt{2A_n}.
\]
Let $\lambda_{n,1}\geq\cdots\geq\lambda_{n,p_n}\geq0$ be the
eigenvalues of $\bm{\Omega}_n$, with a corresponding orthonormal eigenbasis
$\bm{e}_{n,1},\ldots,\bm{e}_{n,p_n}$.  Define
\begin{equation}
 \rho_{n,j}=\frac{\lambda_{n,j}}{\sqrt{A_n}},\qquad
 \eta_{n,j}=\frac{\bm{e}_{n,j}^{\top}\bm{\theta}_n}{A_n^{1/4}},\qquad
 \kappa_n=\frac{\|\bm{\theta}_n\|^2}{s_{0,n}}.
 \label{prof:coordinates}
\end{equation}
Coordinates are extended by zero beyond $p_n$.  In particular,
$\sum_j\rho_{n,j}^2=1$ and
$\sum_j\eta_{n,j}^2=\sqrt{2}\kappa_n$.  With independent standard
normal variables $Z_1,Z_2,\ldots$, the exact decomposition is
\begin{align}
 \frac{T_n^*-\operatorname{tr}(\bm{\Omega}_n)}{s_{0,n}}
 &=Q_n^*+L_n^*+\kappa_n, \label{prof:decomposition}\\
 Q_n^*&=\frac{1}{\sqrt{2}}\sum_j\rho_{n,j}(Z_j^2-1),\qquad
 L_n^*=\sqrt{2}\sum_j\sqrt{\rho_{n,j}}\eta_{n,j}Z_j.
 \label{prof:QL}
\end{align}
The quadratic and linear sums share the same Gaussian coordinates.
They are uncorrelated, but need not be independent.

\begin{theorem}[Explicit local directional profile]
\label{prof:theorem}
Suppose $\kappa_n\to\kappa<\infty$ and, for every fixed $j$,
$\rho_{n,j}\to\rho_j$.  For each $j$ such that $\rho_j>0$, suppose
also that $\eta_{n,j}\to\eta_j$, after choosing an orthonormal eigenbasis
as in \eqref{prof:coordinates}.  Set $\eta_j=0$ when $\rho_j=0$ and
\[
 \sigma^2=1-\sum_{j=1}^{\infty}\rho_j^2.
\]
Then $\sigma^2\geq0$,
$\sum_j\eta_j^2\leq\sqrt{2}\kappa$, and
\begin{equation}
 (Q_n^*,L_n^*)\ \xrightarrow{d}\ (\zeta,L),
 \label{prof:joint-limit}
\end{equation}
where, for independent standard normal variables $Z_0,Z_1,\ldots$,
\begin{equation}
 \zeta=\sigma Z_0+\frac{1}{\sqrt{2}}
              \sum_{j=1}^{\infty}\rho_j(Z_j^2-1),
 \qquad
 L=\sqrt{2}\sum_{j=1}^{\infty}\sqrt{\rho_j}\eta_jZ_j.
 \label{prof:limit}
\end{equation}
Both series converge in $L^2$, and
\begin{equation}
 \frac{T_n^*-\operatorname{tr}(\bm{\Omega}_n)}{s_{0,n}}
 \ \xrightarrow{d}\ V=\zeta+L+\kappa.
 \label{prof:alternative-limit}
\end{equation}
Moreover,
\begin{equation}
 \operatorname{Var}(L)=2\sum_j\rho_j\eta_j^2,
 \qquad
 \operatorname{Var}(V)=1+2\sum_j\rho_j\eta_j^2.
 \label{prof:limit-variance}
\end{equation}
In particular, the random linear contribution can survive only in
eigendirections whose normalized eigenvalues have positive limits.
\end{theorem}

\begin{proof}See Appendix~\ref{app:profiles}.\end{proof}

\begin{cor}[Power at the exact central-reference critical value]
\label{prof:power}
Let the assumptions of Theorem~\ref{prof:theorem} hold.  Fix
$\alpha\in(0,1)$, and let $c_{n,\alpha}$ be the $(1-\alpha)$ quantile
of $Q_n^*$, and let $c_{\alpha}$ be the $(1-\alpha)$ quantile of
$\zeta$.  This limiting quantile is unique, and
\begin{equation}
 \Pr\{T_n^*>\operatorname{tr}(\bm{\Omega}_n)+s_{0,n}c_{n,\alpha}\}
 \longrightarrow\Pr\{\zeta+L+\kappa>c_{\alpha}\}.
 \label{prof:power-limit}
\end{equation}
The same limit holds for the original-data statistic whenever its
standardized distribution has vanishing Kolmogorov distance from that
of $T_n^*$ and its standardized critical value converges in probability
to $c_{\alpha}$.
\end{cor}

\begin{proof}See Appendix~\ref{app:profiles}.\end{proof}

Theorem~\ref{prof:theorem} allows singular covariance matrices: a signal component in the null space contributes only to the deterministic shift. For repeated eigenvalues, the relevant quantity is the signal energy in the corresponding eigenspace; Remark~\ref{prof:multiplicity} gives the basis-invariant formulation. The limiting formula is sequential. Uniform finite-sequence comparison does not, by itself, imply a uniform limiting profile; this distinction is detailed in Remark~\ref{prof:uniformity}.

Under the estimation and moment conditions of Corollary~\ref{cor:feasible-final}, the feasible rejection rule has the same limit in \eqref{prof:power-limit}. Thus the noncentral law describes power while calibration continues to use the central reference.

The next example compares two signals with the same Euclidean norm and makes the directional effect explicit.

\begin{prop}[A single persistent spike]
\label{prof:spike}
Let
\[
 \bm{\Omega}_p=\operatorname{diag}(\sqrt{p},1,\ldots,1),\qquad
 A_p=2p-1,\qquad s_{0,p}=\sqrt{4p-2},
\]
and fix $\kappa>0$.  Consider the two signals
\[
 \bm{\theta}_p^{\parallel}=\sqrt{\kappa s_{0,p}}\,\bm{e}_1,
 \qquad
 \bm{\theta}_p^{\perp}=\sqrt{\kappa s_{0,p}}\,\bm{e}_2.
\]
They have the same Euclidean norm and the same local signal size
$\kappa$.  Nevertheless, their standardized Gaussian-reference limits
are respectively
\begin{align}
 V_{\parallel}
 &=\kappa+\frac{Z_1^2-1}{2}+\frac{Z_0}{\sqrt{2}}
                 +\sqrt{2\kappa}\,Z_1,\label{prof:parallel}\\
 V_{\perp}
 &=\kappa+\frac{Z_1^2-1}{2}+\frac{Z_0}{\sqrt{2}}.
 \label{prof:perpendicular}
\end{align}
Their variances are $1+2\kappa$ and $1$, respectively, so their
distributions differ.  Their rejection probabilities at the exact
central-reference thresholds differ for at least one significance
level.  This assertion does not impose a universal power ordering.
\end{prop}

\begin{proof}See Appendix~\ref{app:profiles}.\end{proof}

The explicit profile recovers the familiar normal local-power formula under
\(\lambda_{\max}(\bm{\Omega}_n)^2/A_n\to0\). In this case every
\(\rho_j\) vanishes, the linear contribution vanishes for bounded
\(\kappa_n\), and \(\zeta\sim N(0,1)\). Consequently the oracle power
tends to \(\Phi(\kappa-z_{1-\alpha})\). A persistent spike does not by
itself force a directional contribution: the signal must have nonvanishing
energy in its eigenspace. Thus covariance shape and signal orientation
play distinct roles.

\begin{prop}[Consistency beyond the local scale]\label{prop:consistency}
Suppose the Gaussian comparison error for \(T_n\) and \(T_n^*\) tends
to zero, and \(\kappa_n=\|\bm{\theta}_n\|^2/\sqrt{2A_n}\to\infty\).
Then rejection at the oracle null critical value has power tending to one.
The same conclusion holds for any estimated critical value satisfying
\(\widehat q_{n,\alpha}-q_{n,\alpha}^{\mathrm{NR}}=o_p(s_{0,n})\).
No spectral limit is required.
\end{prop}
\begin{proof}See Appendix~\ref{app:profiles}.\end{proof}
The proposition establishes consistency of the specified test; minimax
optimality over signal directions or covariance classes is not considered.

Two approximation boundaries are useful when interpreting these power statements. First, matching two or three cumulants need not recover the exact reference under a persistent spectrum (Proposition~\ref{prof:matching}). Second, bounded projection fourth moments alone do not control the diagonal fluctuations of the full squared-mean statistic (Propositions~\ref{prof:radial}--\ref{prof:diagonal}). These examples, together with the qualifications on eigenspaces and uniformity, are given in Appendix~\ref{app:power-details}.

The noncentral reference describes power at a specified population signal.
It does not follow that this power can be consistently estimated from the
same observations when the signal is unknown. The following Gaussian
subexperiment gives a two-point estimation bound \citep[Chapter~2]{Tsybakov2009}. Its covariance is known,
so its limitation cannot be attributed to spectral estimation.

\begin{prop}[A limitation on estimating local power]
\label{prop:local-power-estimation}
For each integer $m\geq1$, observe two independent samples of size $m$,
\[
 \bm Y_{1j}\sim N\!\left(\frac{\bm\theta}{\sqrt{2m}},\bm\Omega\right),
 \qquad
 \bm Y_{2j}\sim N\!\left(-\frac{\bm\theta}{\sqrt{2m}},\bm\Omega\right),
 \qquad j=1,\ldots,m,
\]
where $\bm\Omega\succeq0$ is known and
$\tr(\bm\Omega^2)>0$. Let $P_{\bm\theta,m}$ denote the law of all
observations. Fix $\alpha\in(0,1)$ and write
\[
 \beta_\alpha(\bm\theta)
 =\Pr\{\|\bm G+\bm\theta\|^2>q_\alpha^{\mathrm{NR}}\},
 \qquad \bm G\sim N(\bm0,\bm\Omega),
\]
where $q_\alpha^{\mathrm{NR}}$ is the exact null critical value.
For two signals $\bm\theta_0,\bm\theta_1\in
\operatorname{range}(\bm\Omega)$, put
\[
 h^2=(\bm\theta_1-\bm\theta_0)^\top
       \bm\Omega^+(\bm\theta_1-\bm\theta_0),
 \qquad
 \delta=|\beta_\alpha(\bm\theta_1)
               -\beta_\alpha(\bm\theta_0)|,
\]
where $\bm\Omega^+$ is the Moore--Penrose inverse. If $\delta>0$,
every estimator $\widehat\beta$ based on all observations, with or without
independent auxiliary randomization, satisfies
\begin{align}
 \max_{\ell\in\{0,1\}}
 P_{\bm\theta_\ell,m}
 \left\{|\widehat\beta-\beta_\alpha(\bm\theta_\ell)|
                                      \geq\frac\delta2\right\}
 &\geq\Phi(-h/2),
 \label{eq:power-estimation-probability-bound}\\
 \max_{\ell\in\{0,1\}}
 \E_{\bm\theta_\ell,m}
 |\widehat\beta-\beta_\alpha(\bm\theta_\ell)|
 &\geq\delta\Phi(-h/2).
 \label{eq:power-estimation-absolute-bound}
\end{align}
Consequently, uniform consistency is impossible along any sequence of
such pairs for which $\liminf\delta>0$ and $\limsup h<\infty$.
\end{prop}
\begin{proof}See Appendix~\ref{app:power-estimation-proofs}.\end{proof}

\begin{cor}[Persistent eigenvalues and unknown local effects]
\label{cor:persistent-power-estimation}
Consider the experiment of Proposition~\ref{prop:local-power-estimation}
along any covariance sequence $\bm\Omega_m$ with
$A_m=\tr(\bm\Omega_m^2)>0$. Suppose
$\bm\Omega_m\bm e_m=\lambda_m\bm e_m$, $\|\bm e_m\|=1$, and
$\lambda_m/\sqrt{A_m}\geq c>0$. Fix
$0<a_0<a_1$ and $\alpha\in(0,1)$, and set
$\bm\theta_{\ell,m}=a_\ell\sqrt{\lambda_m}\bm e_m$.
There is a constant $b=b(c,\alpha,a_0,a_1)>0$ such that
\[
 \beta_{\alpha,m}(\bm\theta_{1,m})
 -\beta_{\alpha,m}(\bm\theta_{0,m})\geq b
 \quad\text{for every }m.
\]
Consequently, every estimator based on the two full samples satisfies
\begin{equation}
 \max_{\ell\in\{0,1\}}P_{\bm\theta_{\ell,m},m}
 \left\{|\widehat\beta_m-
       \beta_{\alpha,m}(\bm\theta_{\ell,m})|\geq b/2\right\}
 \geq\Phi\{-(a_1-a_0)/2\}.
 \label{eq:persistent-power-estimation-bound}
\end{equation}
Both alternatives have local strengths bounded above and away from zero:
$a_\ell^2c/\sqrt2\leq\kappa_{\ell,m}
\leq a_\ell^2/\sqrt2$.
\end{cor}
\begin{proof}See Appendix~\ref{app:power-estimation-proofs}.\end{proof}

Knowing the signal norm need not resolve this difficulty. For every
fixed $\alpha\in(0,1)$, Appendix~\ref{app:power-estimation-proofs}
constructs a fixed $\varepsilon\in(0,1)$ such that
$\bm\Omega=\diag(1,\varepsilon,0,\ldots,0)$ and the equal-norm
signals $z_\alpha\bm e_1,z_\alpha\bm e_2$ have distinct powers and
a finite $h$. Thus a direction-sensitive power can remain inestimable
uniformly over a known local signal sphere, even with known covariance.
This is a limitation on estimating the power of the specified test; it
is not a lower bound on the detection risk of all tests.

A useful operational target is therefore power under specified effect
assumptions. For a nonempty compact signal class $\Theta$, define
\begin{equation}
 \beta^-_\alpha(\Theta)
    =\min_{\bm\theta\in\Theta}\beta_\alpha(\bm\theta),
 \qquad
 \beta^+_\alpha(\Theta)
    =\max_{\bm\theta\in\Theta}\beta_\alpha(\bm\theta).
 \label{eq:specified-effect-power-envelope}
\end{equation}
The extrema exist because the reference power is continuous in
$\bm\theta$. These are population sensitivity bounds for the specified
class, rather than confidence limits inferred from an unknown signal.
They allow a proposed scientific effect, or a set of plausible spectral
energy allocations, to be assessed for study planning. Replacing the
covariance by an estimate requires separate error control. Consistent
null calibration and the estimation limitation above are compatible:
the former concerns the centered distribution, whereas the latter
concerns an unknown local alternative.

The boundary criterion can guide an assessment over a specified effect
range. For $\|\bm\theta\|^2\le K\sqrt{2A}$,
$\Delta_{\bm\theta}\le\sqrt2Kr$. Thus a valid upper bound for
$r$ controls approximation over the entire signal ball, without
estimating an unknown signal direction.

Use the $m\ge2$ matched Gaussian residual columns in
\eqref{eq:spherical-columns} and compute
\[
 U_m=\frac{\sum_{r\ne s}(\bm X_r^\top\bm X_s)^2}{m(m-1)},
 \qquad H_{rs}=m^{-1}\bm X_r^\top\bm X_s\ (r\ne s),\quad H_{rr}=0.
\]
Proposition~\ref{prop:v104-trace-convex-order} gives an explicit lower
confidence bound $A_L$ from $U_m$. Proposition~\ref{prop:v104-exponential-hollow-ratio}
then obtains $r_U$ by a scalar inversion using
$\lambda_{\max}(\bm H)/\sqrt{A_L}$. Both constructions allow arbitrary
dimension and singular covariance. Their constants depend only on $m$
and prespecified error allocations. The four events controlling the
lower trace bound, a Gaussian norm, its direction and the hollow
remainder can share a total error probability $\gamma_r$.

In particular, let $r_U$ be any upper confidence bound satisfying
$\Pr(r\le r_U)\ge1-\gamma_r$. Define
\begin{equation}
 B_{\mathrm{smooth}}(K;r_U)=
 \begin{cases}
 \displaystyle\min\left\{\frac12,
       \frac{\sqrt2K r_U}{\pi(1-4r_U^2)}\right\},&0\le r_U<1/2,\\[6pt]
 \displaystyle\frac12,&1/2\le r_U\le1.
 \end{cases}
 \label{eq:smooth-finite-certificate}
\end{equation}
With probability at least $1-\gamma_r$, simultaneously for every $K>0$,
\begin{equation}
 \sup_{\|\bm\theta\|^2\le K\sqrt{2A}}d_{\bm\theta}
       \le B_{\mathrm{smooth}}(K;r_U).
 \label{eq:smooth-finite-coverage}
\end{equation}
Thus the Gaussian residual-based hollow-Gram bound supplies a finite-sample
certificate with unknown covariance using only $r_U$. It does not require
a separate upper trace bound $A_U$ or estimating the entire null CDF. The condition
$r_U<1/2$ is sufficient for this certificate; failure to certify it is not
evidence that the approximation is inaccurate. The empirical-interval
certificate remains applicable outside this range.

The bound can also be inverted. For a tolerance $0<\varepsilon<1/2$
and $0<r_U<1/2$, the largest local effect range certified by the smooth
bound is
$\widehat K_{\varepsilon}=
\varepsilon\pi(1-4r_U^2)/(\sqrt2r_U)$.
Set $\widehat K_{\varepsilon}=0$ when $r_U\ge1/2$.
The same confidence event covers every tolerance, so this inversion
allows the effect range to be selected after seeing the bound;
Corollary~\ref{cor:certified-effect-range} gives the precise statement.
This is a sufficient certified range, not an estimate of the maximal
range of accurate approximation.

For the empirical-interval alternative, let $A_U$ cover $A$ with error
probability $\gamma_A$. Simultaneous beta-spacing and DKW bounds for
the residual squared norms yield $B_{\mathrm{interval}}(K)$ in
Proposition~\ref{prop:finite-unknown-covariance-diagnostic}, with an
additional error probability $\gamma_F$. On the joint event, using
the same $r_U$ in both bounds,
\begin{equation}
 \sup_{\|\bm\theta\|^2\le K\sqrt{2A}}d_{\bm\theta}
 \le\min\{B_{\mathrm{smooth}}(K;r_U),B_{\mathrm{interval}}(K)\}
 \quad\text{for every }K>0,
 \label{eq:main-finite-diagnostic}
\end{equation}
with probability at least $1-\gamma_r-\gamma_A-\gamma_F$.
The common event handles dependence among estimates from the same
residuals. A standalone smooth certificate allocates the full error
budget to $r_U$; combining it with the interval certificate requires
the joint allocation just stated. The conclusion bounds a population
approximation error, not the type-I error of a selected test. These
finite-sample confidence statements remain Gaussian-model assertions.

\section{Simulation studies}\label{sec:numerics}
We conduct three simulation studies. Simulation 1 examines the accuracy
of the noncentral reference and the effect of signal direction on power.
Simulation 2 examines feasible calibration under unequal covariances,
including sample imbalance and non-Gaussian observations. Simulation 3
compares the augmented spherical and moment-based spectral procedures
with existing tests on the same datasets.
Additional experiments in Appendix~\ref{app:numerical-details} assess
the finite-sample covariance diagnostic and the approximation bounds.
Throughout, the nominal significance level is $\alpha=0.05$ unless
otherwise stated, and the normalized signal strength is
$\kappa=\|\bm{\theta}\|^2/s_0$.

Each Monte Carlo configuration uses $R=10{,}000$ independent evaluation samples.
Known-covariance experiments generate the normalized contrast directly.
Feasible experiments generate either the original samples or their exact
Gaussian contrast--residual Gram representation, as specified below, and
estimate their critical values. Diagnostic coverage experiments generate
independent residual datasets and evaluate confidence events. Deterministic calculations and auxiliary reference
draws are identified separately. In each evaluation replication we record
rejection. The empirical size or power
is the proportion of rejections, reported as a percentage in the tables. For an estimate $\widehat\pi$, the
binomial Monte Carlo standard error is
$\{\widehat\pi(1-\widehat\pi)/R\}^{1/2}$, approximately $0.00218$
at $\pi=0.05$ and at most $0.005$, or $0.22$ and $0.50$ percentage points, respectively. Common simulated errors are used across
methods and mean shifts within a configuration, so differences between
methods can also be assessed using paired outcomes.

The moment-based calibration and external-comparison experiments use an additional $R_0=10{,}000$ independent
null samples for the size-adjusted power diagnostic. These samples are
separate from the evaluation samples. Reference draws within a dataset
are distinct from both outer replication counts; their budgets are given
in Table~\ref{tab:simulation-budgets}. Population-oracle computations and
the bootstrap calculations for Monte Carlo uncertainty are described
separately below. Models, seeds and numerical conventions are given in
Appendix~\ref{app:numerical-details} and the separately maintained computational code.

\begin{table}[!htbp]
\centering\small
\begin{threeparttable}
\caption{Simulation budgets.}
\label{tab:simulation-budgets}
\begin{tabular}{@{}lrrp{0.38\textwidth}@{}}
\toprule
Study & $R$ & $R_0$ & Conditional reference computation\\
\midrule
Population reference & 10,000 & --- & Deterministic quadrature or exact CDFs\\
Feasible two-sample calibration & 10,000 & 10,000 & 3,000 draws for each spectral reference\\
External comparison & 10,000 & 10,000 & 999 signs for Wang--Xu; 1,999 draws for each spectral reference\\
Gaussian spherical calibration & 10,000 & --- & 999 draws per feasible method\\
Rank-one validation & 10,000 & --- & 999 spherical draws; exact finite-$B$ power\\
Unknown-covariance diagnostic & 10,000 & --- & Deterministic bounds; refinements use paired replay\\
Noncommuting/imbalanced Gaussian & 10,000 & --- & 999 draws per feasible method\\
Non-Gaussian spherical calibration & 10,000 & --- & 999 draws per feasible method\\
\bottomrule
\end{tabular}
\begin{tablenotes}[flushleft]
\footnotesize
\item[] \textit{Note.} $R$ denotes evaluation replications and $R_0$ independent null-training replications. Inner draws are generated separately for each observed dataset.
\end{tablenotes}
\end{threeparttable}
\end{table}

\subsection{Simulation 1: Reference approximation and directional power}
We examine the boundaries in Theorems~\ref{thm:shift-boundary} and~\ref{thm:directional-diameter} using known-population Gaussian references. Take $p\in\{50,200,800,3200\}$ and either
$\bm\Omega=\bm I_p$ or
$\bm\Omega=\operatorname{diag}(\sqrt p,1,\ldots,1)$. Set $\kappa=1$.
In the latter model, use
\[
 \bm\theta=\sqrt{s_0}\{\sqrt h\,\bm e_1+\sqrt{1-h}\,\bm e_2\},
 \qquad h\in\{0,0.25,0.5,0.75,1\}.
\]
Thus $h$ is the fraction of signal energy in the leading eigendirection.
The interaction parameter is
$\Delta_{\bm\theta}=s_0\{h\sqrt p+(1-h)\}/(2p-1)$; it tends to zero for $h=0$
but to $h$ for any fixed $h>0$. For the isotropic model,
$\Delta_{\bm\theta}=\sqrt{2/p}$.

The CDFs and central-reference critical values are calculated from the exact Gaussian laws. We compare CDFs on a standardized grid from $-6$ to $12$ with step $0.025$; the grid discrepancy is not a certified Kolmogorov supremum. For each base covariance, 10,000 independent Gaussian contrasts provide a separate Monte Carlo check at $\alpha\in\{0.01,0.05,0.10\}$. Directions and levels share the base draws. These are oracle experiments; Appendix~\ref{app:boundary-design} gives numerical accuracy checks and additional isotropic and singular examples.

\begin{figure}[t]
\centering
\includegraphics[width=\textwidth]{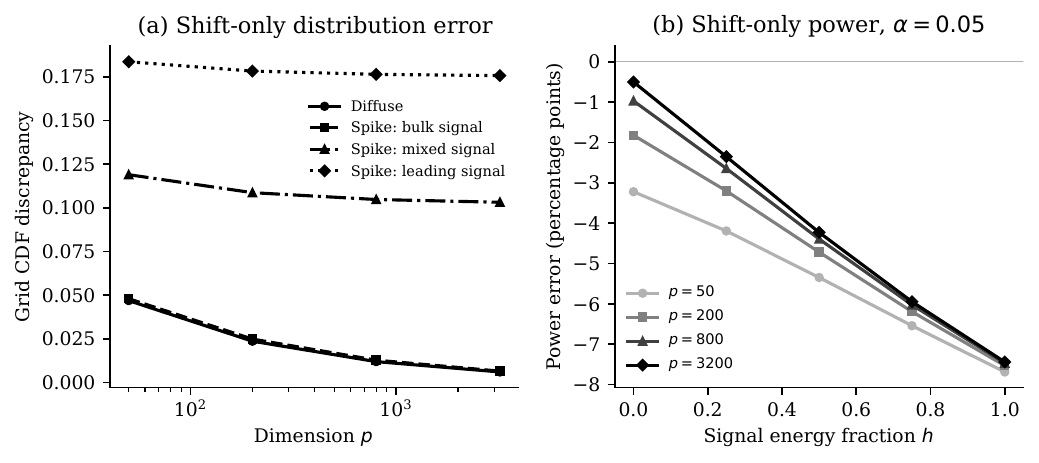}
\caption{Gaussian-reference approximation boundaries at $\kappa=1$.
Panel (a) gives finite-grid CDF discrepancies for diffuse and spiked
covariances; the mixed signal uses $h=0.5$. Panel (b) gives the
shifted-null power prediction minus the exact noncentral power in the
spiked model, in percentage points. Lines are deterministic distributional
calculations. Independent 10,000-replication checks are reported separately.}
\label{fig:boundary-diagnostics}
\end{figure}

Figure~\ref{fig:boundary-diagnostics} illustrates the distinction between
vanishing and persistent interaction. At $p=3200$, the grid discrepancies
are $0.0060$ for the isotropic model, $0.0065$ for a bulk signal in the
spiked model, and $0.1756$ for a leading-direction signal. For the last
two signals, the shifted-null approximation gives the same $5\%$-level
power of $18.56\%$, whereas the noncentral powers are $19.06\%$ and
$26.00\%$, respectively. The Monte Carlo estimates are $19.16\%$ and
$25.68\%$, with standard errors $0.39$ and $0.44$ percentage points.
The persistent discrepancy is therefore visible at this specified
significance level as well as over the CDF grid; the theoretical
all-level statement does not assert such a gap at every fixed level.

Consider $\bm{\Omega}_p=\operatorname{diag}(\sqrt p,1,\ldots,1)$. Its standardized
null reference converges to $\zeta=(Z_1^2-1)/2+G/\sqrt2$, where $Z_1$ and
$G$ are independent standard normal variables. Numerical integration gives
the upper $0.05$ critical value $1.751431$. The limiting two- and
three-cumulant critical values are $1.644854$ and $1.876828$, respectively.
Their rejection probabilities under the limiting reference are $0.057234$
and $0.042744$. Thus, knowing the population cumulants does not remove the
approximation error in this example.

\begin{figure}[t]
\centering
\includegraphics[width=.72\textwidth]{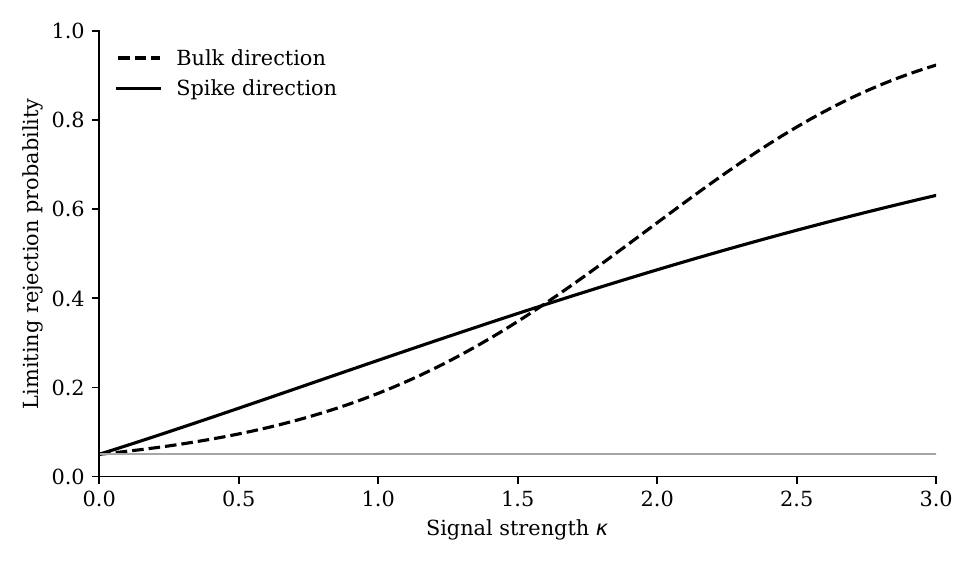}
\caption{Oracle limiting power at nominal level $0.05$ for equal-norm
signals in the single-spike model. The bulk and spike curves are
$\Pr(\zeta+\kappa>c_{.95})$ and
$\Pr(\zeta+\sqrt{2\kappa}Z_1+\kappa>c_{.95})$, respectively.}
\label{fig:directional-power}
\end{figure}

Figure~\ref{fig:directional-power} gives the limiting powers for signals
with $\|\bm{\theta}_p\|^2/s_{0,p}\to\kappa$. A signal in the spike direction
retains the linear term $\sqrt{2\kappa}Z_1$, whereas a signal in the bulk
direction does not. The resulting power curves differ and may cross.
These curves are evaluated by adaptive one-dimensional integration, with
absolute and relative tolerances $2\times10^{-11}$.

We next examine how accurately the reference predicts finite-sample rejection probabilities.
For the BF statistic, put $n=n_1+n_2$ and $v=(n_2+2n_1)/n$, and take
$\bm{\Sigma}_1=\bm{\Omega}_p/v$ and $\bm{\Sigma}_2=2\bm{\Sigma}_1$. The normalized mean
contrast then has covariance $\bm{\Omega}_p$. We consider $\bm{\Omega}_p=\bm{I}_p$ and
the single-spike covariance above, $p\in\{200,800\}$, and
$(n_1,n_2)\in\{(20,30),(50,80)\}$. Independent standardized coordinates
have a normal distribution, a variance-one $t_8$ distribution, or the
distribution of $(\chi_1^2-1)/\sqrt2$. Signals satisfy
$\|\bm{\theta}_p\|^2=\kappa s_{0,p}$, with $\kappa\in\{0,0.5,1,2\}$.
In the spiked model, they lie in the first or second eigendirection.
Only one direction is used in the diffuse model, where the Gaussian
reference is rotationally invariant.

Each of the 132 configurations uses 10,000 independent Monte Carlo
replications. Within each base error and sample-size setting, the same
noise is used across signal strengths, directions and calibration rules.
The finite Gaussian reference is evaluated by 256-node Gauss--Hermite
integration over its leading coordinate, with the remaining coordinates
handled by a noncentral chi-square distribution. Checks against adaptive
quadrature at the evaluated oracle critical points differed by less than
$9\times10^{-12}$. The integration rule is separate from the number of
Monte Carlo replications.

Table~\ref{tab:finite-power} reports selected results for $p=800$ and
$\kappa=1$. The complete Gaussian reference gives bulk and spike powers
$19.42\%$ and $25.91\%$, whereas the shift-only approximation gives $18.45\%$
in both directions. The empirical powers show the same directional
difference. The bulk power still has a finite-dimensional linear
contribution, which vanishes only along the asymptotic sequence.
\begin{table}[!htbp]
\centering\small\setlength{\tabcolsep}{4pt}
\begin{threeparttable}
\caption{Empirical sizes and powers (in \%) in the single-spike model.}
\label{tab:finite-power}\label{tab:oracle-selected}
\begin{tabular}{rrlrrrr}
\toprule
$n_1$ & $n_2$ & Error law & Size & Bulk power & Spike power & Max. MC SE \\
\midrule
20 & 30 & Normal & 4.83 & 19.07 & 26.12 & 0.44 \\
 &  & $t_8$ & 5.07 & 19.80 & 25.78 & 0.44 \\
 &  & $(\chi_1^2-1)/\sqrt2$ & 5.51 & 19.72 & 26.16 & 0.44 \\
\addlinespace
50 & 80 & Normal & 4.90 & 18.95 & 25.91 & 0.44 \\
 &  & $t_8$ & 5.46 & 19.82 & 26.38 & 0.44 \\
 &  & $(\chi_1^2-1)/\sqrt2$ & 5.24 & 20.05 & 26.38 & 0.44 \\
\midrule
\multicolumn{3}{l}{Gaussian reference} & 5.00 & 19.42 & 25.91 & --- \\
\multicolumn{3}{l}{Shift-only approximation} & 5.00 & 18.45 & 18.45 & --- \\
\bottomrule
\end{tabular}

\begin{tablenotes}[flushleft]
\footnotesize
\item[] \textit{Note.} $p=800$. Powers use $\kappa=1$ and sizes use $\kappa=0$. Each empirical entry uses 10,000 replications. The last column gives the largest Monte Carlo standard error among the three empirical entries, in percentage points. Reference rows are population calculations.
\end{tablenotes}
\end{threeparttable}
\end{table}

Across the 132 configurations, the median and maximum absolute differences
between empirical rejection and the Gaussian-reference prediction are
$0.32$ and $1.36$ percentage points, respectively. These differences include
Monte Carlo error. The complete results give empirical CDF discrepancies
and rejection probabilities under all three population calibration rules.
Figure~\ref{fig:finite-power} in Appendix~\ref{app:numerical-details}
shows the selected empirical powers with their Monte Carlo uncertainty.

\begin{figure}[t]
\centering
\includegraphics[width=\textwidth]{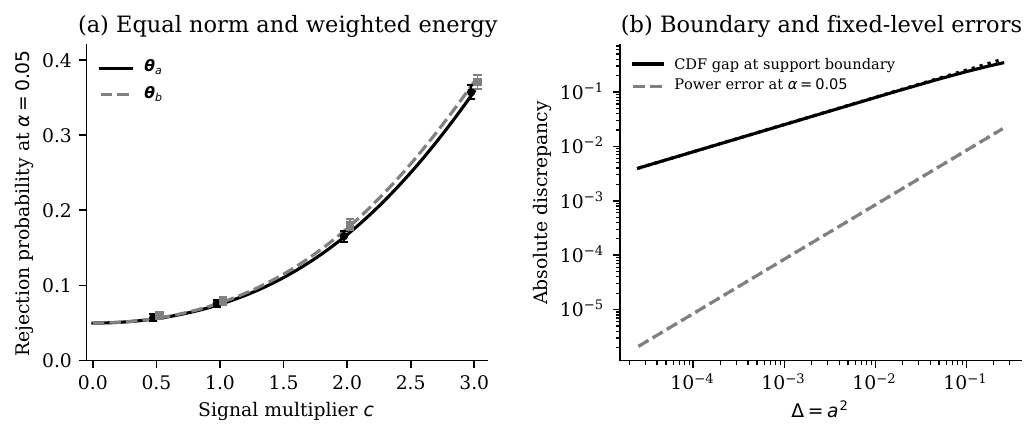}
\caption{Information about the alternative and approximation errors.
Panel (a) compares equal-norm signals with equal $\Delta$ under
$\bm\Omega=\diag(1,2,3)$. Lines are population powers from numerical
quadrature; points and error bars are Monte Carlo estimates and
$1.96$ standard errors from 10,000 contrasts. Points are displaced
horizontally for visibility. Panel (b) uses analytic rank-one formulas:
the CDF gap at the support boundary $a^2$ equals the Kolmogorov
distance by Proposition~\ref{prop:rankone-sharp-error}. Dotted lines give the
corresponding first-order expressions. Both panels use known
population covariances.}
\label{fig:information-rates}
\end{figure}

The code includes the seed, numerical environment and all 396
calibration-specific results.
To distinguish the weighted signal energy from the information needed for
the complete reference, we also consider $\bm\Omega=\diag(1,2,3)$ and
$\bm\theta_a=c(1/\sqrt2,0,1/\sqrt2)^\top$,
$\bm\theta_b=c(0,1,0)^\top$. Both signals have squared norm $c^2$ and
$\Delta=c^2/7$, but their third cumulants differ.
Figure~\ref{fig:information-rates}(a) shows their population powers and
Monte Carlo estimates from 10,000 independent Gaussian contrasts, shared
across the signal settings. At $c=3$ and $\alpha=0.05$, the population
powers are $35.82\%$ and $37.10\%$, respectively; the empirical powers are
$35.72\%$ and $37.05\%$. The paired standard error of their empirical
power difference is $0.62$ percentage points. Thus the distinction is
visible, although the power differences in this example are modest.

Panel (b) separates a support-boundary discrepancy from the error at a
fixed significance level. For $(Z+a)^2$ and $Z^2+a^2$, with
$\Delta=a^2$, the CDF discrepancy at $a^2$ is
$\Phi(2|a|)-1/2$, whereas the power difference at $\alpha=0.05$ has
first-order magnitude proportional to $\Delta$. At $a=0.1$, the former
is $0.07926$ and the latter is $0.0008472$. These are analytic population
calculations. They illustrate why a bound for the entire distribution
need not describe the error at a conventional testing level closely.
The numerical design, quadrature checks and paired Monte Carlo errors
are given in Appendix~\ref{app:information-numerics}.

We also examined the fully specified unknown-covariance bound in
Proposition~\ref{prop:finite-unknown-covariance-diagnostic}, using 10,000
independent Gaussian residual datasets in each of nine configurations.
The configurations include isotropic, rank-one and spiked covariances,
with dimensions up to 80,000; the diagnostic itself uses no population
covariance information. The joint sufficient confidence event had empirical
coverage between 99.41\% and 99.72\% at the nominal 95\% level.
The bounds were nevertheless conservative. For an isotropic covariance
with $(m,p)=(400,800)$, the median bound was 0.156, 0.255 and 0.449
at $K=0.001$, $0.01$ and $0.1$, respectively. At $(m,p)=(40,800)$,
the median was the universal cap $1/2$ at all four examined values of $K$;
the rank-one cases were also capped. The spacing refinement improved the
small-interval probability bound relative to a DKW-only construction
with the same total confidence budget. Thus the finite-sample result
supplies a computable certificate, but does not yet certify small
approximation error at modest residual sample sizes. These experiments
check the coverage and informativeness of the sufficient confidence
construction; they do not estimate the actual supremum shift error or
the rejection probability of a selected test. Full settings and results
are given in Appendix~\ref{app:unknown-covariance-experiment}.

We also compare the finite-sample approximation certificates on the same
90,000 Gaussian residual datasets used in the original confidence experiment,
with 10,000 independent datasets in each of nine configurations.
These are paired recalculations, not additional independent replications.
The linear bound under spectral spread changes the diagnostic substantially
when the spectral ratio can be certified. For isotropic covariance with
$m=400$ and $p=800$, at $K=0.25$ the median original interval bound is
0.500, compared with 0.0279 for the hybrid and 0.0270 for the smooth
certificate alone. Both new bounds are at most 0.05 in all 10,000 datasets.
At $K=0.5$, however, their medians are 0.0558 and 0.0539, respectively,
so neither certifies the 0.05 tolerance in this configuration.
The gain is not uniform: with $m=40$, or with the rank-one and spiked
spectra at $K=0.25$, the median remains 0.500.
The separate smooth certificate allocates the entire 0.05 error budget
to its four spectral events; the hybrid reserves two of six equal
allocations for the upper trace and interval bounds. Both have 95\%
nominal confidence, but the hybrid need not be smaller than this separately
allocated smooth certificate. Appendix~\ref{app:tightened-diagnostic-experiment}
reports the allocation, paired comparisons and observed coverage.

\subsection{Simulation 2: Feasible calibration}
We first examine Gaussian spherical calibration in five
configurations, including $p=20{,}000$ with sample sizes $(20,30)$.
The group covariance matrices are nonproportional; their weighted sum
is either $\bm I_p$ or $\diag(\sqrt p,1,\ldots,1)$.
Each configuration uses 10,000 independent Gaussian datasets and 999
fresh reference draws per feasible procedure. The joint contrast--residual
Gram matrices are generated through their exact Gaussian law, as detailed
in Appendix~\ref{app:rotation-numerics}; no population covariance is
supplied to the feasible calibration algorithms.

\begin{figure}[t]
\centering
\includegraphics[width=\textwidth]{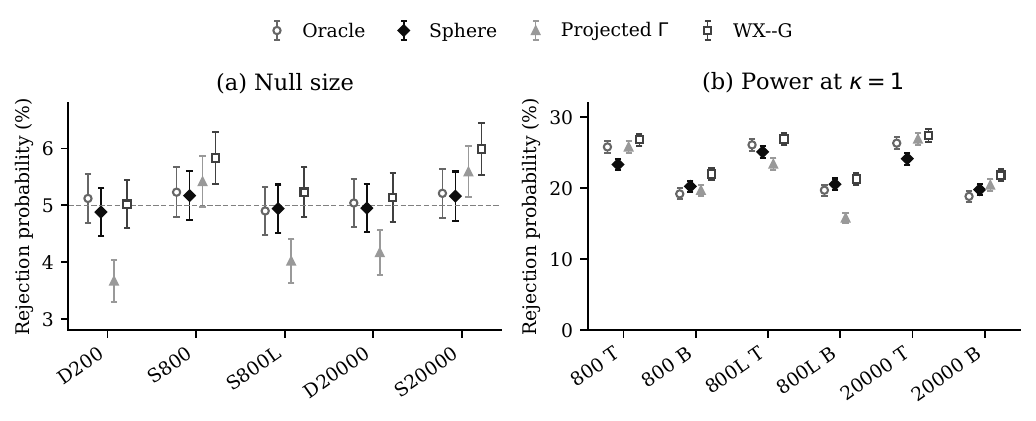}
\caption{Gaussian calibration at nominal level $0.05$. Panel (a) gives
null rejection probabilities; panel (b) gives power at $\kappa=1$ in the
spiked models, where T and B denote leading-axis and bulk-axis signals.
The suffix L denotes sample sizes $(40,60)$, and the other configurations
use $(20,30)$. Error bars are $1.96$ Monte Carlo standard errors from
10,000 datasets. PG and WX--G are defined in Table~\ref{tab:rotation-size}.}
\label{fig:rotation-calibration}
\end{figure}

Table~\ref{tab:rotation-size} presents the empirical sizes, and
Figure~\ref{fig:rotation-calibration} compares power. The competing
procedures are the population oracle, the projected gamma reference and
the Gaussian multiplier version of \citet{WangXu2022}.
The spherical sizes range from $4.88\%$ to $5.17\%$; the projected-gamma
sizes range from $3.67\%$ to $5.59\%$, and the Wang--Xu Gaussian sizes
from $5.02\%$ to $5.99\%$. Spherical calibration uses neither a spectral
threshold nor a variance-budget projection. The original threshold
equals $0.5$ in every replication of the comparison procedure.

The power comparison also shows a finite-sample cost. For a leading-axis
signal with $\kappa=1$, $p=800$ and $(n_1,n_2)=(20,30)$, spherical and
oracle powers are $23.34\%$ and $25.81\%$; with sizes $(40,60)$ they
are $25.12\%$ and $26.10\%$. At $p=20{,}000$ and sizes $(20,30)$,
the corresponding powers are $24.15\%$ and $26.34\%$.
These results support the stated Gaussian calibration guarantee, not a
uniform power ordering. Powers for all signal settings and paired
Monte Carlo errors are recorded in the separate computational archive.

\begin{table}[!htbp]
\centering\small
\begin{threeparttable}
\caption{Empirical sizes (in \%) in the Gaussian calibration study.}
\label{tab:rotation-size}
\begin{tabular}{lrr rrrr}
\toprule
& & & \multicolumn{4}{c}{Method}\\
\cmidrule(lr){4-7}
Configuration & $(n_1,n_2)$ & $p$ & Oracle & Sphere & PG & WX--G\\
\midrule
D200 & $(20,30)$ & 200 & 5.12 & 4.88 & 3.67 & 5.02 \\
S800 & $(20,30)$ & 800 & 5.23 & 5.17 & 5.42 & 5.83 \\
S800L & $(40,60)$ & 800 & 4.90 & 4.94 & 4.02 & 5.23 \\
D20000 & $(20,30)$ & 20,000 & 5.04 & 4.95 & 4.17 & 5.14 \\
S20000 & $(20,30)$ & 20,000 & 5.21 & 5.16 & 5.59 & 5.99 \\
\bottomrule
\end{tabular}
\begin{tablenotes}[flushleft]
\footnotesize
\item[] \textit{Note.} The nominal level is $0.05$. Each row uses 10,000 independent datasets and 999 fresh auxiliary reference draws per feasible procedure. The Monte Carlo standard error near $5\%$ is $0.22$ percentage points. PG denotes the projected gamma reference with all Gaussian residuals; WX--G denotes the Wang--Xu Gaussian multiplier implementation.
\end{tablenotes}
\end{threeparttable}
\end{table}

We further examine severe sample imbalance and covariances that do
not commute. Four fresh configurations include $(n_1,n_2)=(12,120)$,
nonaligned rank-one spikes in the two groups, and $p=20{,}000$.
Each configuration uses 10,000 independent Gaussian datasets and 999
fresh auxiliary draws per feasible procedure. The spherical procedure
uses all observations for the mean contrast, together with $11$ matched
residual columns in the imbalanced configurations and $23$ in the
balanced configuration. The exact Gaussian Gram simulation and all
covariance definitions are given in Appendix~\ref{app:noncommuting-numerics}.

Table~\ref{tab:noncommuting-size} reports the null sizes, and
Figure~\ref{fig:noncommuting-calibration} also displays power at $\kappa=1$.
The spherical sizes range from $4.77\%$ to $5.03\%$. At $p=200$, the two noncommuting configurations have the same
$\bm\Omega$, obtained by adjusting the group spike magnitudes as sample
allocation changes. Thus they compare calibration under a fixed
contrast covariance; the two group covariances themselves are not held
fixed. For a leading-axis signal with $\kappa=1$, the spherical powers are $21.90\%$ and $20.52\%$ in the balanced and imbalanced cases, respectively. These results examine the power changes accompanying this sample-allocation and covariance design; they do not isolate a universal residual-matching cost. Powers are reported at the nominal critical
values and should be interpreted with their corresponding sizes;
neither a uniform power ordering nor non-Gaussian finite-sample validity
is inferred from this experiment.

\begin{table}[!htbp]
\centering\small
\begin{threeparttable}
\caption{Empirical sizes (in \%) under sample imbalance and noncommuting covariances.}
\label{tab:noncommuting-size}
\begin{tabular}{lrr rrrr}
\toprule
& & & \multicolumn{4}{c}{Method}\\
\cmidrule(lr){4-7}
Configuration & $(n_1,n_2)$ & $p$ & Oracle & Sphere & PG & WX--G\\
\midrule
D200-I & $(12,120)$ & 200 & 5.38 & 5.03 & 4.64 & 6.02 \\
NC200-B & $(24,24)$ & 200 & 4.78 & 4.82 & 4.82 & 5.79 \\
NC200-I & $(12,120)$ & 200 & 4.99 & 4.77 & 5.82 & 5.87 \\
NC20000-I & $(12,120)$ & 20,000 & 4.97 & 4.82 & 5.77 & 5.83 \\
\bottomrule
\end{tabular}
\begin{tablenotes}[flushleft]
\footnotesize
\item[] \textit{Note.} The observations are Gaussian and the nominal level is $0.05$. Each configuration has 10,000 independent datasets and $B=999$ fresh reference draws. The Monte Carlo standard error near $5\%$ is $0.22$ percentage points. PG uses all Gaussian residuals; WX--G is the published Gaussian multiplier implementation.
\end{tablenotes}
\end{threeparttable}
\end{table}

We next examine the non-Gaussian approximation in
Theorem~\ref{thm:nongaussian-spherical-calibration}. The latent
coordinates are either Rademacher or standardized $t_8$, and the
noncommuting covariance matrices are held fixed while sample sizes
increase from $(20,30)$ to $(40,60)$. Each of the four configurations
uses 10,000 independent datasets and 999 fresh reference draws.
All raw observations are generated explicitly, so dependence between
the sample means and completed orthogonal residuals is preserved.
The data generation and residual basis are specified in
Appendix~\ref{app:non-gaussian-numerics}.

Table~\ref{tab:non-gaussian-size} reports the null sizes.
The spherical sizes range from $4.98\%$ to $5.51\%$. The Rademacher sizes are $5.46\%$ and $5.51\%$, indicating modest finite-sample inflation; this two-size grid does not demonstrate monotone improvement. The corresponding $t_8$ sizes are $5.02\%$ and $4.98\%$. Figure~\ref{fig:non-gaussian-calibration} additionally displays
leading-axis power at $\kappa=1$. These results examine the asymptotic
transfer to non-Gaussian data; they do not imply exact finite-sample
validity outside the Gaussian model. The Gaussian population reference
is included as an approximation benchmark and is not a non-Gaussian
oracle. The projected-gamma comparison uses independent disjoint-pair
residuals; its $t_8$ results are exploratory because the established
spectral sufficient conditions for that procedure are sub-Gaussian.
All reported powers are nominal and should be read together with the
sizes. Full power results and paired Monte Carlo standard errors are
recorded in the separate computational archive.

\begin{table}[!htbp]
\centering\small
\begin{threeparttable}
\caption{Empirical sizes (in \%) under non-Gaussian observations.}
\label{tab:non-gaussian-size}
\begin{tabular}{lrr rrrr}
\toprule
& & & \multicolumn{4}{c}{Method}\\
\cmidrule(lr){4-7}
Configuration & $(n_1,n_2)$ & $p$ & G-ref & Sphere & PG & WX--G\\
\midrule
Rad-S & $(20,30)$ & 200 & 5.00 & 5.46 & 6.04 & 6.12 \\
Rad-L & $(40,60)$ & 200 & 5.24 & 5.51 & 5.22 & 5.87 \\
t8-S & $(20,30)$ & 200 & 4.90 & 5.02 & 5.82 & 5.95 \\
t8-L & $(40,60)$ & 200 & 5.14 & 4.98 & 5.07 & 5.73 \\
\bottomrule
\end{tabular}
\begin{tablenotes}[flushleft]
\footnotesize
\item[] \textit{Note.} The latent variables are non-Gaussian and the covariance matrices are noncommuting. The nominal level is $0.05$. Each configuration has 10,000 independent datasets and $B=999$ fresh reference draws. The Monte Carlo standard error near $5\%$ is $0.22$ percentage points. G-ref is the Gaussian population reference, not a non-Gaussian oracle. PG uses disjoint-pair residuals; WX--G is the published Gaussian multiplier implementation.
\end{tablenotes}
\end{threeparttable}
\end{table}

\Needspace{5\baselineskip}\paragraph{Moment-based calibration.}

We consider six two-sample scenarios with $p\in\{100,400\}$ and sample
sizes $(40,60)$ or $(80,120)$. The two covariance matrices are
nonproportional and not simultaneously diagonalizable. Four scenarios
use Gaussian observations, covering both dimensions and both sample-size
pairs; two use transformed Rademacher observations at the larger sample
sizes. We compare the population Gaussian-mixture oracle, estimated
two-cumulant calibration, and raw and projected spectral references with
normal or gamma remainders. Gaussian Helmert residuals and non-Gaussian
same-data paired differences provide the corresponding feasible estimates.
Every scenario has 10,000 evaluation and 10,000 separate null-training
datasets; inner spectral references use 3,000 draws.
Appendix~\ref{app:expanded-design} gives the full design and further diagnostics.

\begin{table}[!htbp]
\centering\small\setlength{\tabcolsep}{4pt}
\begin{threeparttable}
\caption{Empirical sizes (in \%) for the feasible two-sample study.}
\label{tab:expanded-size}
\begin{tabular}{lrrrrrrr}
\toprule
& & & & \multicolumn{4}{c}{Method}\\
\cmidrule(lr){5-8}
Model & $p$ & $n_1$ & $n_2$ & Oracle & 2C & Proj. G & Proj. $\Gamma$ \\
\midrule
N & 100 & 40 & 60 & 5.16 & 6.10 & 5.50 & 5.36 \\
 & 100 & 80 & 120 & 4.84 & 5.57 & 5.39 & 5.14 \\
 & 400 & 40 & 60 & 4.69 & 5.90 & 4.60 & 4.59 \\
 & 400 & 80 & 120 & 4.48 & 5.52 & 4.86 & 4.84 \\
\midrule
R & 100 & 80 & 120 & 4.99 & 6.12 & 5.70 & 5.60 \\
 & 400 & 80 & 120 & 4.82 & 6.16 & 4.80 & 4.81 \\
\midrule
\multicolumn{4}{l}{ARE} & 4.47 & 17.90 & 7.77 & 6.20 \\
\bottomrule
\end{tabular}

\begin{tablenotes}[flushleft]
\footnotesize
\item[] \textit{Note.} The nominal level is $5\%$. N and R denote normal and Rademacher latent errors. Each entry uses 10,000 replications; the binomial Monte Carlo standard error near $5\%$ is $0.22$ percentage points. ARE is $100M^{-1}\sum_{j=1}^M|\widehat\alpha_j-\alpha|/\alpha$, with $\alpha=0.05$ and the $M=6$ configurations displayed here. It summarizes size distortion over these configurations.
\end{tablenotes}
\end{threeparttable}
\end{table}

Table~\ref{tab:expanded-size} reports the empirical sizes and their ARE
values. The estimated two-cumulant sizes range from $5.52\%$ to $6.16\%$, compared
with $4.60\%$--$5.70\%$ for the projected normal reference and
$4.59\%$--$5.60\%$ for the projected gamma reference. The corresponding
ARE values are $17.90$, $7.77$ and $6.20$. The oracle sizes range from
$4.48\%$ to $5.16\%$, including outer simulation variation and
Gaussian-reference error for non-Gaussian samples.
These finite designs do not establish a uniform advantage for either remainder.

The predetermined threshold is capped at $0.5$ in every replication:
$m_{\min}\le79$ makes this automatic. These results assess the capped rule
and do not establish finite-sample threshold adaptation. Projection occurs
in $25.30\%$ of Gaussian evaluations with $p=400$ and sizes $(40,60)$,
and $15.96\%$ of Rademacher evaluations with $p=400$ and sizes $(80,120)$.
It has a modest effect here and cannot repair an incorrectly estimated
spectral shape.

Independent null training gives a size-adjusted power diagnostic without changing the operational method. For Gaussian observations with $p=400$ and sizes $(40,60)$, raw spike powers of two-cumulant and projected gamma calibration are $26.01\%$ and $22.94\%$; adjusted powers are $22.55\%$ and $22.45\%$. Raw bulk powers are $20.93\%$ and $16.31\%$, with adjusted powers $15.83\%$ and $15.65\%$. The raw differences therefore need to be interpreted together with the size differences. Appendix~\ref{app:expanded-design} reports the corresponding Monte Carlo standard errors, including independent training uncertainty for adjusted power.

\begin{figure}[t]
\centering
\includegraphics[width=\textwidth]{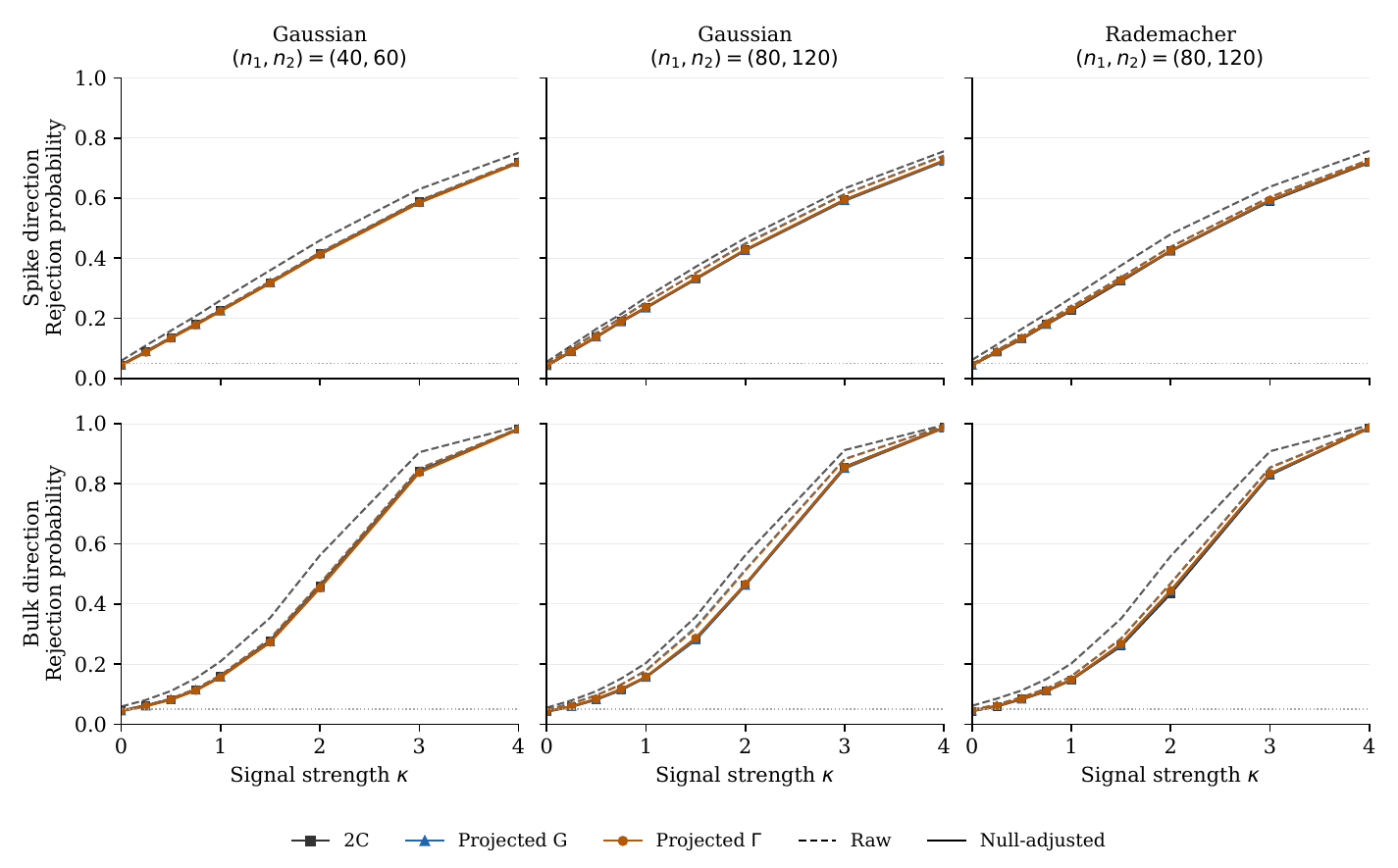}
\caption{Two-sample power curves at $p=400$. Columns show the three feasible-calibration designs and rows show the signal direction. Dashed and solid lines are raw and independently null-adjusted rejection probabilities. Each point uses 10,000 evaluation samples, shared across signal strengths; the adjustment uses 10,000 independent null samples. These are point estimates without error bars. Two-phase bootstrap standard errors at $\kappa=1$ are reported in Table~\ref{tab:expanded-power}.}
\label{fig:expanded-adjusted}
\end{figure}

Figure~\ref{fig:expanded-adjusted} examines how these powers change with signal strength. For each of the three $p=400$ designs, we use $\kappa\in\{0,0.25,0.5,0.75,1,1.5,2,3,4\}$ in the spike and bulk directions. The same 10,000 evaluation samples and the same estimated reference are used throughout each curve. Translation invariance permits exact reconstruction of the shifted rejection margins from the saved null and unit-signal margins; the identity is given in Appendix~\ref{app:expanded-design}. This does not add independent observations. The curves show that the directional ordering can change with signal strength and that most raw gaps between these calibrations narrow after the independent null adjustment. They do not establish uniform power superiority.

\subsection{Simulation 3: Comparison with existing tests}
We compare the augmented spherical test, the two-cumulant test, the four
spectral calibrations, the approximate randomization test of
\citet{WangXu2022}, and the split linear test based on
\citet{ShekharKimRamdas2022}. The separate unequal-covariance justification
for the latter is given in Proposition~\ref{prop:split-linear-validity}.
Wang--Xu uses the full-sample Chen--Qin statistic, whereas the spherical
and spectral tests use the squared sample mean contrast. Thus this
experiment compares complete testing procedures.

The six settings combine Gaussian or Rademacher observations with
diffuse or three-spike spectra and nonproportional diagonal group
covariances. The sample sizes are $(40,60)$ or $(80,120)$, and the
dimension is 200 or 800. Alternatives have $\kappa=1$ or 2 and lie
along a leading or a bulk coordinate. Every method receives the same
10,000 evaluation datasets per setting. The spherical reference uses
fixed Helmert residual contrasts and includes the observed mean contrast
under each alternative; it uses 1,999 auxiliary draws. The other spectral
estimates use paired residuals from the same observations and 1,999
reference draws when needed. Wang--Xu uses 999 randomizations.
These inner budgets are separate from the number of independent datasets.
A further 10,000 independent null-training datasets per setting are used
only for size-adjusted power diagnostics. Appendix~\ref{app:benchmark-design}
gives the construction, paired Monte Carlo errors and timing protocol.

\begin{table}[!htbp]
\centering\small\setlength{\tabcolsep}{3.5pt}
\begin{threeparttable}
\caption{Empirical sizes (in \%) for the external-method comparison.}
\label{tab:external-benchmarks}
\begin{tabular}{lrrrrrrrr}
\toprule
& & & \multicolumn{2}{c}{Unprojected} & \multicolumn{2}{c}{Projected} & \multicolumn{2}{c}{Benchmarks} \\
\cmidrule(lr){4-5}\cmidrule(lr){6-7}\cmidrule(lr){8-9}
Scenario & 2C & Sphere & G & $\Gamma$ & G & $\Gamma$ & WX & Split \\
\midrule
G1 & 5.05 & 4.85 & 0.21 & 0.21 & 3.66 & 3.66 & 4.68 & 5.93 \\
G2 & 6.52 & 5.02 & 0.21 & 0.21 & 5.49 & 5.49 & 5.35 & 5.84 \\
G3 & 6.19 & 5.07 & 3.15 & 3.16 & 4.33 & 4.34 & 5.21 & 5.20 \\
\midrule
R1 & 5.03 & 5.03 & 0.18 & 0.18 & 3.56 & 3.56 & 4.84 & 5.66 \\
R2 & 5.85 & 5.11 & 4.79 & 4.73 & 4.79 & 4.73 & 5.13 & 5.26 \\
R3 & 6.11 & 5.21 & 3.30 & 3.30 & 4.30 & 4.30 & 5.11 & 5.48 \\
\midrule
ARE & 15.83 & 1.97 & 60.53 & 60.70 & 16.17 & 16.33 & 4.27 & 11.23 \\
\bottomrule
\end{tabular}
\begin{tablenotes}[flushleft]\footnotesize
\item[] \textit{Note.} Each scenario uses the same 10,000 evaluation datasets for all methods. Sphere is the augmented spherical test with 1,999 auxiliary draws. WX uses 999 randomizations. The spectral Gaussian (G) and gamma ($\Gamma$) references use 1,999 auxiliary draws. Scenarios are defined in Appendix~\ref{app:benchmark-design}; ARE is defined as in Table~\ref{tab:expanded-size}.
\end{tablenotes}
\end{threeparttable}
\end{table}

Table~\ref{tab:external-benchmarks} reports the empirical sizes.
The spherical sizes are $4.85\%$--$5.21\%$, with ARE 1.97, compared
with $4.68\%$--$5.35\%$ and ARE 4.27 for Wang--Xu. These differences
are small relative to the Monte Carlo uncertainty in individual sizes.
The Gaussian results agree with finite-sample level control. The
Rademacher results are empirical evidence for the asymptotic extension;
they do not supply a non-Gaussian finite-sample guarantee.

The older spectral calibrations provide a useful comparison. In G1, G2
and R1, the retained eigenvalues exceed the estimated variance budget in
every replication. Projection changes sizes from $0.18\%$--$0.21\%$
to $3.56\%$--$5.49\%$, but some conservatism remains. The threshold is
$0.5$ throughout, and projection is active in almost every replication
in five settings. When it is active the remainder variance is zero, so
the projected Gaussian and gamma references coincide. These settings
mainly assess variance-budget correction, whereas the spherical rule
avoids that projection.

\begin{figure}[ht]
\centering
\includegraphics[width=\textwidth]{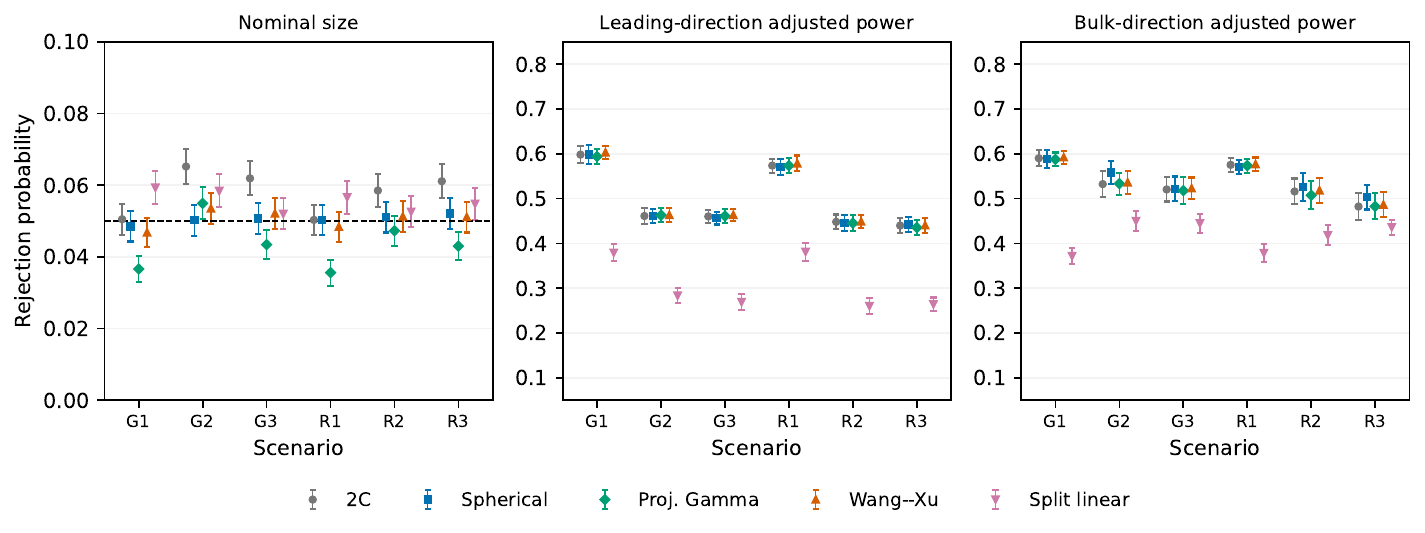}
\caption{External comparisons based on 10,000 evaluation datasets per
setting. Error bars are 1.96 Monte Carlo standard errors; those for
size-adjusted power include uncertainty from 10,000 independent
null-training datasets.}
\label{fig:external-benchmarks}
\end{figure}

Figure~\ref{fig:external-benchmarks} shows sizes and independently
size-adjusted powers. At $\kappa=2$, the spherical-minus-Wang--Xu
adjusted-power differences range from $-0.85$ to $0.18$ percentage points
in the leading direction and from $-0.59$ to $2.23$ in the bulk direction.
Their paired Monte Carlo errors incorporate both null-training and
evaluation uncertainty. Neither procedure dominates across these
directions and settings. Nominal powers, held-out adjusted sizes and
paired differences are reported in Appendix~\ref{app:benchmark-design}.
Size adjustment uses the generating null law and is a diagnostic, not
an additional implementable test. The fixed auxiliary draw counts also
distinguish these experiments from the limiting oracle-power statements.

\clearpage
\section{Discussion}\label{sec:discussion}
This paper studies noncentral normal-reference inference for the high-dimensional two-sample Behrens--Fisher problem. The covariance-weighted signal energy determines when the null reference can be shifted by the squared signal size to approximate the alternative distribution. The uniform fourth-root error order is sharp over unrestricted spectra. Three nonnegligible spectral blocks give a three-quarter-order upper bound, and stronger spectral spread gives a linear bound. A separate spectral-range criterion determines when a reference based only on signal magnitude is adequate. Otherwise, covariance-eigenspace signal energies are needed to describe directional power. The limitation on estimating power at an unknown local alternative explains why specified-effect power ranges have a sensitivity interpretation.

The augmented spherical reference gives a feasible test without estimating individual covariance eigenvalues. Its exact Gaussian validity follows from classical residual reduction and rotation invariance, while the comparison results justify its asymptotic use under the two stated moment conditions. The simulations examine size and power under unequal covariances, severe sample imbalance and non-Gaussian sampling. The modest Rademacher size inflation in Simulation 2 illustrates a finite-sample limitation of the asymptotic extension, while the paired comparison in Simulation 3 gives spherical sizes close to the nominal level in its six settings. The results do not establish a universal power advantage over existing tests.

For unknown Gaussian covariance, the confidence certificate can identify local effect ranges over which the shifted-null approximation is sufficiently accurate. Its improved informativeness in some diffuse-spectrum designs is useful, but concentrated spectra and small residual samples can still yield the universal bound $1/2$. Such an outcome does not establish that the actual approximation error is large. The certificate is simultaneous over effect ranges; choosing confidence budgets from the same data is not covered by its guarantee.

Further work is needed on sharper spherical-reference rates, informative finite-sample certificates under non-Gaussian sampling, and weaker conditions for diagonal fluctuations. The current empirical evidence is based on synthetic data, and minimax detection optimality remains outside the scope of the results.

\paragraph{Computational materials.}
Simulation settings and implementation details are described in
Appendix~\ref{app:numerical-details}. The computational archive is
maintained separately and is not included in this source-only submission.
References below to code, metadata and replication-level outputs refer
to that separate archive.

\clearpage

\appendix

\section{Gaussian comparison and proofs}\label{app:transfer}
\subsection{A finite-sample comparison}
\label{app:general-comparison}

To justify the reference under the alternative, we compare the two-sample
statistic with its Gaussian counterpart. Write
$\bm{\xi}_{ij}=\bm{Y}_{ij}-\bm{\mu}_i$ and $n=n_1+n_2$. The centered
normalized mean difference is the sum of the independent increments
\[
 \bm{X}_{1j}=\frac{\sqrt{n_1n_2/n}}{n_1}\bm{\xi}_{1j},
 \qquad
 \bm{X}_{2j}=-\frac{\sqrt{n_1n_2/n}}{n_2}\bm{\xi}_{2j}.
\]
For the comparison argument, relabel these increments as
$\bm{X}_1,\ldots,\bm{X}_N\in\mathbb R^p$, where $N=n_1+n_2$.
Let their covariance matrices be $\bm{\Gamma}_1,\ldots,\bm{\Gamma}_N$, and put
\[
 \bm{\Omega}=\sum_{t=1}^N\bm{\Gamma}_t,
 \qquad A=\operatorname{tr}(\bm{\Omega}^2)>0.
\]
On an enlarged probability space, let $\bm{G}_t\sim N(\bm{0},\bm{\Gamma}_t)$ be independent
Gaussian vectors that are also independent of the entire $\bm{X}$ array. For a deterministic
signal $\bm{\theta}\in\mathbb R^p$, define
\begin{equation}
 T_{\bm{\theta}}=\left\|\bm{\theta}+\sum_{t=1}^N\bm{X}_t\right\|^2,
 \qquad
 T_{\bm{\theta}}^{\mathrm G}=\left\|\bm{\theta}+\sum_{t=1}^N\bm{G}_t\right\|^2,
 \qquad
 s_{\bm{\theta}}^2=2A+4\bm{\theta}^\top\bm{\Omega}\bm{\theta}.
 \label{eq:uniform-transfer-statistics}
\end{equation}
Here $s_{\bm{\theta}}^2$ is the variance of the Gaussian reference, including
the linear contribution under the alternative. This population scale is
used for the comparison; it is not an estimated studentizing factor. Write
\[
 W_{\bm{\theta}}=\frac{T_{\bm{\theta}}-\operatorname{tr}(\bm{\Omega})-\|\bm{\theta}\|^2}{s_{\bm{\theta}}},
 \qquad
 W_{\bm{\theta}}^{\mathrm G}
 =\frac{T_{\bm{\theta}}^{\mathrm G}-\operatorname{tr}(\bm{\Omega})-\|\bm{\theta}\|^2}{s_{\bm{\theta}}}.
\]
For real random variables $V,V'$, denote their Kolmogorov distance by
\[
 d_{\mathrm K}(V,V')=\sup_{x\in\mathbb R}
 |\Pr(V\le x)-\Pr(V'\le x)|.
\]

\begin{assumption}[Moment conditions]
\label{ass:increment-moments}
\textnormal{(C1)} There is a constant $K<\infty$, common to all increments, such that
\begin{equation}
 \mathbb E(\bm{a}^\top \bm{X}_t)^4\le K(\bm{a}^\top\bm{\Gamma}_t\bm{a})^2
 \quad\text{for every }\bm{a}\in\mathbb R^p\text{ and }1\le t\le N.
 \label{eq:increment-projection-fourth}
\end{equation}
The variances $\operatorname{Var}(\|\bm{X}_t\|^2)$ are finite.

\textnormal{(C2)} For bounds stated only in terms of covariance matrices,
we additionally assume
\begin{equation}
 \operatorname{Var}(\|\bm{X}_t\|^2)
 \le K\operatorname{tr}(\bm{\Gamma}_t^2),\qquad 1\le t\le N.
 \label{eq:increment-norm-variance}
\end{equation}
\end{assumption}

The projection bound controls the off-diagonal products and the linear
term, whereas the norm-variance bound controls the diagonal terms.
The latter does not follow from bounded projection kurtosis.
Uniformity over $\bm{\theta}$ refers to location shifts of fixed centred
increment laws. Families of centred laws are also allowed when the
moment constants and influence bounds below hold uniformly.

Define the dimensionless quantities
\begin{align}
 D_{\bm{\theta}}
 &:=\frac{\sum_{t=1}^N\operatorname{Var}(\|\bm{X}_t\|^2)
              +2\sum_{t=1}^N\operatorname{tr}(\bm{\Gamma}_t^2)}{s_{\bm{\theta}}^2},
 \label{eq:transfer-D}\\
 B_{\bm{\theta}}
 &:=\sum_{t=1}^N
 \left\{
 \frac{\operatorname{tr}(\bm{\Omega}\bm{\Gamma}_t)+\bm{\theta}^\top\bm{\Gamma}_t\bm{\theta}}
      {s_{\bm{\theta}}^2}
 \right\}^{3/2}.
 \label{eq:transfer-B}
\end{align}
Here $D_{\bm{\theta}}$ measures the diagonal fluctuations in the two statistics,
and $B_{\bm{\theta}}$ bounds the observation influences in their quadratic
and linear terms.

\begin{theorem}[Finite-sample quadratic--linear Gaussian comparison]
\label{thm:finite-transfer}
Under \eqref{eq:increment-projection-fourth} and finite $D_{\bm{\theta}}$, there is a
finite constant $C_K$, depending only on $K$, such that
\begin{equation}
 d_{\mathrm K}(T_{\bm{\theta}},T_{\bm{\theta}}^{\mathrm G})
 =d_{\mathrm K}(W_{\bm{\theta}},W_{\bm{\theta}}^{\mathrm G})
 \le
 \min\left[1,\ C_K\inf_{0<\varepsilon\le1}
 \left\{\frac{D_{\bm{\theta}}}{\varepsilon^2}
       +\frac{B_{\bm{\theta}}}{\varepsilon^3}+\varepsilon^{1/2}\right\}\right].
 \label{eq:finite-transfer-unsimplified}
\end{equation}
In particular, after increasing $C_K$ if necessary,
\begin{equation}
 d_{\mathrm K}(T_{\bm{\theta}},T_{\bm{\theta}}^{\mathrm G})
 \le \min\{1,\ C_K(D_{\bm{\theta}}^{1/5}+B_{\bm{\theta}}^{1/7})\}.
 \label{eq:finite-transfer-rate}
\end{equation}
If \eqref{eq:increment-norm-variance} also holds, then
\begin{equation}
 D_{\bm{\theta}}\le
 (K+2)\frac{\sum_{t=1}^N\operatorname{tr}(\bm{\Gamma}_t^2)}{s_{\bm{\theta}}^2}.
 \label{eq:finite-transfer-D-covariance}
\end{equation}
The bounds hold for each finite $N$ and $p$, without a spectral diffuseness
condition or an assumption on the existence of a limiting distribution.
\end{theorem}

\begin{remark}[Scope of the rate]
The bound combines a third-order replacement argument, control of the
diagonal terms by Chebyshev's inequality, and Gaussian quadratic
anti-concentration. Its exponents are not asserted to be optimal.
Theorem~\ref{thm:finite-transfer} is a specialization of Gaussian
invariance results to the full squared mean statistic; it supplies the
comparison needed for the subsequent power and calibration results.
The relation to existing results is discussed below and summarized in
Appendix~\ref{app:literature-comparison}.
\end{remark}

\begin{cor}[Uniformity over signals and deterministic critical values]
\label{cor:uniform-signal-transfer}
For a set $\Theta\subseteq\mathbb R^p$, suppose
$D_\Theta:=\sup_{\bm{\theta}\in\Theta}D_{\bm{\theta}}<\infty$ and
$B_\Theta:=\sup_{\bm{\theta}\in\Theta}B_{\bm{\theta}}<\infty$. Then
\begin{equation}
 \sup_{\bm{\theta}\in\Theta}\sup_{x\in\mathbb R}
 |\Pr(T_{\bm{\theta}}\le x)-\Pr(T_{\bm{\theta}}^{\mathrm G}\le x)|
 \le \min\{1,C_K(D_\Theta^{1/5}+B_\Theta^{1/7})\}.
 \label{eq:uniform-signal-transfer}
\end{equation}
Consequently, the same right-hand side bounds
\[
 \sup_{\bm{\theta}\in\Theta}
 |\Pr(T_{\bm{\theta}}>c_{\bm{\theta}})-\Pr(T_{\bm{\theta}}^{\mathrm G}>c_{\bm{\theta}})|
\]
for any deterministic thresholds $c_{\bm{\theta}}$.
Random estimated critical values require the additional calibration
argument in Section~\ref{sec:spectral-calibration}.
\end{cor}

For the original observations, the increment conditions follow from \eqref{eq:raw-group-moments}. Scalar multiplication preserves both inequalities. The resulting bound depends on the smaller sample size, without requiring balanced samples.

\begin{cor}[Two-sample Gaussian comparison]
\label{cor:grouped-transfer}
Suppose the two independent samples satisfy
\eqref{eq:raw-group-moments}. Put
\begin{align*}
 \bm{\Omega}_1&=\frac{n_2}{n}\bm{\Sigma}_1,
 &\bm{\Omega}_2&=\frac{n_1}{n}\bm{\Sigma}_2,\\
 \bm{\Omega}&=\bm{\Omega}_1+\bm{\Omega}_2,
 &n_{\min}&=\min(n_1,n_2).
\end{align*}
Then $\operatorname{Cov}(\bm{X}_{ij})=\bm{\Omega}_i/n_i$, and, uniformly
for all $\bm{\theta}\in\mathbb R^p$,
\begin{equation}
 D_{\bm{\theta}}\le\frac{K+2}{2n_{\min}},
 \qquad
 B_{\bm{\theta}}\le\frac{1}{2^{3/2}\sqrt{n_{\min}}}.
 \label{eq:grouped-influence-bounds}
\end{equation}
Consequently,
\begin{equation}
 \sup_{\bm{\theta}\in\mathbb R^p}
 d_{\mathrm K}(T_{\bm{\theta}},T_{\bm{\theta}}^{\mathrm G})
 \le \min\{1,C_K(n_{\min}^{-1/5}+n_{\min}^{-1/14})\}.
 \label{eq:grouped-uniform-rate}
\end{equation}
The distance tends to zero as $\min(n_1,n_2)\to\infty$, with no restriction
on the dimension imposed by this bound. Neither $n_1/n$ nor $n_2/n$
needs to be bounded away from zero. Conditions for estimating the reference
distribution are stated separately in Section~\ref{sec:spectral-calibration}.
\end{cor}
\begin{proof}See Appendix~\ref{app:transfer}.\end{proof}

For example, \eqref{eq:raw-group-moments} holds with a uniform constant when
$\bm{\xi}_{ij}=\bm{L}_i \bm{z}_{ij}$,
$\bm{L}_i\bm{L}_i^\top=\bm{\Sigma}_i$, and the coordinates of each
$\bm{z}_{ij}$ are independent, centered, have unit variance, and have uniformly
bounded fourth moments. Indeed, with $\bm{H}_i=\bm{L}_i^\top \bm{L}_i$ and
$\kappa_{i\ell}=\mathbb E z_{ij,\ell}^4-3$,
\[
 \operatorname{Var}(\bm{z}_{ij}^\top \bm{H}_i \bm{z}_{ij})
 =2\operatorname{tr}(\bm{H}_i^2)+\sum_\ell\kappa_{i\ell}(\bm{H}_i)_{\ell\ell}^2
 \le C\operatorname{tr}(\bm{\Sigma}_i^2).
\]
The fourth-moment identity for a linear combination gives the projection
bound. These comparison conditions require no eighth moments.

The Gaussian reference is a quadratic polynomial and may have a
nonnormal limit. It is useful to relate the preceding comparison to existing
invariance results. Table~\ref{tab:closest-comparison} in
Appendix~\ref{app:literature-comparison} compares the relevant results. In particular, \eqref{eq:finite-transfer-rate} follows from the
degree-two specialization of existing invariance bounds.

For clarity, let $F_{\bm{\theta}}=T_{\bm{\theta}}-\tr(\bm{\Omega})-\|\bm{\theta}\|^2$ and
\[
 Q_{\bm{\theta}}=2\sum_{r<t}\bm{X}_r^\top \bm{X}_t+2\bm{\theta}^\top\sum_t\bm{X}_t.
\]
The polynomial $Q_{\bm{\theta}}$ has degree two and is multilinear in the independent
observation vectors. Applying the approximate-polynomial theorem of
\citet{HuangAusternOrbanz2025} with $f=F_{\bm{\theta}}$, $q_2=Q_{\bm{\theta}}$, and
$\nu=3$, then controlling the Gaussian diagonal remainder, recovers the
order $D_{\bm{\theta}}^{1/5}+B_{\bm{\theta}}^{1/7}$. The variance used in that mapping is
$\operatorname{Var}(Q_{\bm{\theta}})=s_{\bm{\theta}}^2-2\sum_t\tr(\bm{\Gamma}_t^2)$.
Their Proposition~3 and Example~3 also address the cost of retaining
diagonal terms. Our proof makes the moment and two-sample influence bounds
explicit; deterministic-shift uniformity and positive-semidefinite
aggregation are applications of this comparison principle.

The subsequent analysis concerns directional power and estimated spectral
calibration. The mixed null profile in
\citet[][Supplement, S.1]{WangXu2022} is used below to compare
signals of equal norm and to assess cumulant matching. Our additional
norm-variance condition controls the diagonal terms and is stronger than
bounded projection kurtosis. The refinement below is related to
\citet[][Theorems~1--2]{ChoudharyKuchibhotla2026}, whose
anti-concentration bounds adapt to the relative linear and quadratic
contributions.

A bounded reference density improves the comparison bound. Such a bound
requires additional spectral spread, since the centered rank-one reference
$(Z^2-1)/\sqrt2$ has square-root small-ball probabilities. Suppose the
eigenvalues can be partitioned into three sets $I_1,I_2,I_3$ such that,
for some fixed $c>0$,
\[
 \min_{1\le\ell\le3}\sum_{j\in I_\ell}\lambda_j^2\ge cA.
\]
The condition concerns the total mass of each spectral block and does not
require the third largest eigenvalue to be bounded below. It holds if
$\lambda_1^2\le A/4$, and in the single-spike example with
$c=1/6$ when $p\ge5$.

\begin{theorem}[Comparison under spectral spread]
\label{thm:spread-summary}
Under the displayed three-block condition, the standardized Gaussian
reference $W_{\bm{\theta}}^{\mathrm G}$ has a density bounded by a constant
$M_c$, uniformly over all deterministic $\bm{\theta}$ and the dimension.
Consequently, under the increment moment assumptions of
Theorem~\ref{thm:finite-transfer},
\[
 d_{\mathrm K}(T_{\bm{\theta}},T_{\bm{\theta}}^{\mathrm G})
 \le\min\{1,C_{K,c}(D_{\bm{\theta}}^{1/3}+B_{\bm{\theta}}^{1/4})\}.
\]
Under the two-sample conditions of
Corollary~\ref{cor:grouped-transfer}, this gives
\[
 \sup_{\bm{\theta}} d_{\mathrm K}(T_{\bm{\theta}},T_{\bm{\theta}}^{\mathrm G})
 \le\min\{1,C_{K,c}(n_{\min}^{-1/3}+n_{\min}^{-1/8})\}.
\]
No upper bound on the signal strength is imposed.
\end{theorem}
\begin{proof}
See Appendix~\ref{app:adaptive}, in particular
Theorems~\ref{adap:uniform-density} and \ref{adap:transfer}.
\end{proof}

The proof integrates the noncentral characteristic function and uses its
exponential factor to control large signals. This gives a sufficient
condition in the present positive-semidefinite setting; the more general
anti-concentration problem is treated by
\citet{ChoudharyKuchibhotla2026}. The resulting exponents remain
smoothing upper bounds. Appendix~\ref{app:adaptive} also treats
mean-dominated low-rank references and derives a linear
anti-concentration cost for random critical values under the
three-block condition.

\subsection{Proofs}

We first give the anti-concentration bound used in the Gaussian comparison.

\begin{lem}[Gaussian quadratic anti-concentration]
\label{lem:gaussian-quadratic-small-ball}
If $P(\bm{Z})$ is a polynomial of degree at most two in a standard Gaussian vector
$\bm{Z}$, and $\operatorname{Var}\{P(\bm{Z})\}=1$, then an absolute constant $C$ satisfies
\begin{equation}
 \sup_{x\in\mathbb R}\Pr\{|P(\bm{Z})-x|\le\varepsilon\}
 \le C\varepsilon^{1/2},\qquad \varepsilon>0.
 \label{eq:gaussian-quadratic-small-ball}
\end{equation}
The dimension of $\bm{Z}$ does not enter the bound.
\end{lem}

\begin{proof}
Theorem~8 of \citet{CarberyWright2001}, with polynomial degree $2$ and its
moment parameter $q=4$, gives
\[
 \Pr\{|P(\bm{Z})-x|\le\varepsilon\}
 \le C\varepsilon^{1/2}
       \{\mathbb E|P(\bm{Z})-x|^2\}^{-1/4}.
\]
The second moment is at least the variance, which equals one. This proves
\eqref{eq:gaussian-quadratic-small-ball}. Singular Gaussian covariance matrices
cause no difficulty, since their polynomial can be expressed in a standard
Gaussian vector of the appropriate rank.
\end{proof}

\begin{proof}[Proof of Theorem~\ref{thm:finite-transfer}]
Suppress the subscript $\bm{\theta}$ on $s_{\bm{\theta}}$. Decompose
\[
 W_{\bm{\theta}}=U_X+R_X,
 \qquad
 W_{\bm{\theta}}^{\mathrm G}=U_G+R_G,
\]
where
\begin{align*}
 U_X&=\frac{2\sum_{r<t}\bm{X}_r^\top \bm{X}_t+2\bm{\theta}^\top\sum_t\bm{X}_t}{s},
 &R_X&=\frac{\sum_t\{\|\bm{X}_t\|^2-\operatorname{tr}(\bm{\Gamma}_t)\}}{s},\\
 U_G&=\frac{2\sum_{r<t}\bm{G}_r^\top \bm{G}_t+2\bm{\theta}^\top\sum_t\bm{G}_t}{s},
 &R_G&=\frac{\sum_t\{\|\bm{G}_t\|^2-\operatorname{tr}(\bm{\Gamma}_t)\}}{s}.
\end{align*}
Independence across observations, without independence between $U_X$ and
$R_X$, yields
\begin{equation}
 \mathbb E R_X^2+\mathbb E R_G^2=D_{\bm{\theta}}.
 \label{eq:diagonal-second-moment}
\end{equation}

We next show that every bounded $C^3$ function $f$ with bounded third
derivative satisfies
\begin{equation}
 |\mathbb Ef(U_X)-\mathbb Ef(U_G)|
 \le C_K\|f'''\|_\infty B_{\bm{\theta}}.
 \label{eq:smooth-offdiagonal-replacement}
\end{equation}
Replace the increments successively by their Gaussian counterparts.
At step $m$, condition on all increments except the $m$th and put
\[
 \bm{\eta}_m=\bm{\theta}+\sum_{r<m}\bm{G}_r+\sum_{r>m}\bm{X}_r.
\]
The two arguments of $f$ differ only through
$D_m=2\bm{\eta}_m^\top \bm{X}_m/s$ and
$D_m^*=2\bm{\eta}_m^\top \bm{G}_m/s$; their remaining, common part is measurable with
respect to the conditioned variables. Their first two conditional moments
agree, so second-order Taylor expansion with its third-order remainder gives
a replacement error at most
\[
 \frac{\|f'''\|_\infty}{6}
 \mathbb E\{|D_m|^3+|D_m^*|^3\}.
\]
By the projection fourth-moment bound and Lyapunov's inequality, this is at
most
\begin{equation}
 C_K\|f'''\|_\infty s^{-3}
 \mathbb E(\bm{\eta}_m^\top\bm{\Gamma}_m\bm{\eta}_m)^{3/2}.
 \label{eq:replacement-one-step}
\end{equation}

To bound this expectation, write $\bm{V}_m=\bm{\eta}_m-\bm{\theta}$.
It is centred, is independent of $\bm{X}_m,\bm{G}_m$, and has covariance
$\bm{\Omega}-\bm{\Gamma}_m\preceq\bm{\Omega}$. Expanding the fourth moment of the
independent scalar summands gives
\[
 \mathbb E(\bm{a}^\top \bm{V}_m)^4\le K_0(\bm{a}^\top\bm{\Omega}\bm{a})^2,
 \qquad K_0=\max(K,3).
\]
Diagonalizing the nonnegative definite matrix $\bm{\Gamma}_m$ and applying
Cauchy--Schwarz term by term therefore yields
\[
 \mathbb E(\bm{V}_m^\top\bm{\Gamma}_m \bm{V}_m)^2
 \le K_0\{\operatorname{tr}(\bm{\Omega}\bm{\Gamma}_m)\}^2.
\]
Since
$\bm{\eta}_m^\top\bm{\Gamma}_m\bm{\eta}_m
\le2\bm{\theta}^\top\bm{\Gamma}_m\bm{\theta}+2\bm{V}_m^\top\bm{\Gamma}_m\bm{V}_m$,
we obtain
\[
 \mathbb E(\bm{\eta}_m^\top\bm{\Gamma}_m\bm{\eta}_m)^2
 \le C_K\{\bm{\theta}^\top\bm{\Gamma}_m\bm{\theta}+
                 \operatorname{tr}(\bm{\Omega}\bm{\Gamma}_m)\}^2.
\]
Lyapunov's inequality now bounds the expectation in
\eqref{eq:replacement-one-step} by
$C_K\{\bm{\theta}^\top\bm{\Gamma}_m\bm{\theta}+
\operatorname{tr}(\bm{\Omega}\bm{\Gamma}_m)\}^{3/2}$.
Summing over $m$ proves \eqref{eq:smooth-offdiagonal-replacement}.

To pass to distribution functions, fix $0<\varepsilon\le1$ and choose a
nonincreasing $C^3$ cutoff $f_{a,\varepsilon}$ such that
\[
 \mathbf1\{u\le a\}\le f_{a,\varepsilon}(u)
 \le\mathbf1\{u\le a+\varepsilon\},
 \qquad
 \|f_{a,\varepsilon}'''\|_\infty\le C\varepsilon^{-3}.
\]
Such a cutoff is obtained by integrating a fixed smooth nonnegative function
supported on $(0,1)$ and then translating and rescaling it. For every $x$,
Chebyshev's inequality, \eqref{eq:smooth-offdiagonal-replacement}, and
\eqref{eq:diagonal-second-moment} give
\begin{align*}
 \Pr(W_{\bm{\theta}}\le x)
 &\le\Pr(U_X\le x+\varepsilon)+\Pr(|R_X|>\varepsilon)\\
 &\le\mathbb E f_{x+\varepsilon,\varepsilon}(U_X)
                +\varepsilon^{-2}\mathbb E R_X^2\\
 &\le\mathbb E f_{x+\varepsilon,\varepsilon}(U_G)
                +C_KB_{\bm{\theta}}\varepsilon^{-3}
                +\varepsilon^{-2}\mathbb E R_X^2\\
 &\le\Pr(U_G\le x+2\varepsilon)
                +C_KB_{\bm{\theta}}\varepsilon^{-3}
                +\varepsilon^{-2}\mathbb E R_X^2\\
 &\le\Pr(W_{\bm{\theta}}^{\mathrm G}\le x+3\varepsilon)
                +C_KB_{\bm{\theta}}\varepsilon^{-3}
                +D_{\bm{\theta}}\varepsilon^{-2}.
\end{align*}
For the opposite direction, the same cutoff bounds give
\begin{align*}
 \Pr(W_{\bm{\theta}}\le x)
 &\ge\Pr(U_X\le x-\varepsilon)-\Pr(|R_X|>\varepsilon)\\
 &\ge\mathbb E f_{x-2\varepsilon,\varepsilon}(U_X)
                -\varepsilon^{-2}\mathbb E R_X^2\\
 &\ge\mathbb E f_{x-2\varepsilon,\varepsilon}(U_G)
                -C_KB_{\bm{\theta}}\varepsilon^{-3}
                -\varepsilon^{-2}\mathbb E R_X^2\\
 &\ge\Pr(U_G\le x-2\varepsilon)
                -C_KB_{\bm{\theta}}\varepsilon^{-3}
                -\varepsilon^{-2}\mathbb E R_X^2\\
 &\ge\Pr(W_{\bm{\theta}}^{\mathrm G}\le x-3\varepsilon)
                -C_KB_{\bm{\theta}}\varepsilon^{-3}
                -D_{\bm{\theta}}\varepsilon^{-2}.
\end{align*}
The Gaussian random variable $W_{\bm{\theta}}^{\mathrm G}$ is a degree-two
polynomial in standard Gaussian coordinates, with variance exactly one:
the centred quadratic and linear terms have zero covariance and variances
$2A/s^2$ and $4\bm{\theta}^\top\bm{\Omega}\bm{\theta}/s^2$, respectively. Applying
Lemma~\ref{lem:gaussian-quadratic-small-ball} to the two intervals of length
$3\varepsilon$ proves \eqref{eq:finite-transfer-unsimplified}.
The anti-concentration bound applies to the full Gaussian reference;
no such bound for the diagonal-deleted statistic is needed.

If $D_{\bm{\theta}}\le1$ and $B_{\bm{\theta}}\le1$, choose
$\varepsilon=\max\{D_{\bm{\theta}}^{2/5},B_{\bm{\theta}}^{2/7}\}$.
Then
\[
 D_{\bm{\theta}}\varepsilon^{-2}\le D_{\bm{\theta}}^{1/5},
 \qquad
 B_{\bm{\theta}}\varepsilon^{-3}\le B_{\bm{\theta}}^{1/7},
 \qquad
 \varepsilon^{1/2}\le D_{\bm{\theta}}^{1/5}+B_{\bm{\theta}}^{1/7}.
\]
The case where either quantity exceeds one follows from
$d_{\mathrm K}\le1$, after increasing $C_K$. The distance is unchanged by
applying the same positive affine transformation to both variables. Finally,
\eqref{eq:finite-transfer-D-covariance} follows directly from
\eqref{eq:increment-norm-variance}.
\end{proof}

\begin{remark}[Relation to the null fluctuation scale]
Let $s_0^2=2A$ and $\mathcal T_M=\{\bm{\theta}:
\bm{\theta}^\top\bm{\Omega}\bm{\theta}\le MA\}$ for a fixed $M<\infty$. On this set,
$1\le s_{\bm{\theta}}/s_0\le\sqrt{1+2M}$. The Kolmogorov comparison itself is
unchanged if both variables are centered and scaled using $s_0$, because the
transformation is common to the two laws. Thus the theorem also transfers
directional local limits expressed on the null fluctuation scale. In
particular, $\|\bm{\theta}\|^2/\sqrt{2A}\le K_1$ implies
$\bm{\theta}^\top\bm{\Omega}\bm{\theta}/A\le\sqrt2K_1$, since
$\lambda_{\max}(\bm{\Omega})\le\sqrt A$. This observation does not force the linear
term to vanish: that requires an additional directional or spectral condition.
\end{remark}

\begin{proof}[Proof of Corollary~\ref{cor:grouped-transfer}]
The increments displayed in Appendix~\ref{app:general-comparison} satisfy
Assumption~\ref{ass:increment-moments} by scalar multiplication of
\eqref{eq:raw-group-moments}. For $i=1,2$,
$\operatorname{Cov}(\bm{X}_{ij})=\bm{\Omega}_i/n_i$.
Since $\bm{\Omega}_1$ and $\bm{\Omega}_2$ are nonnegative definite,
$\operatorname{tr}(\bm{\Omega}_i\bm{\Omega}_j)\ge0$, and hence
\[
 \sum_{i=1}^2\operatorname{tr}(\bm{\Omega}_i^2)
 \le\operatorname{tr}\left\{\left(\sum_{i=1}^2\bm{\Omega}_i\right)^2\right\}=A.
\]
Using $s_{\bm{\theta}}^2\ge2A$ gives
\[
 D_{\bm{\theta}}\le\frac{K+2}{s_{\bm{\theta}}^2}
 \sum_{i=1}^2\frac{\operatorname{tr}(\bm{\Omega}_i^2)}{n_i}
 \le\frac{K+2}{2n_{\min}}.
\]
Put
\[
 q_i(\bm{\theta})=
 \frac{\operatorname{tr}(\bm{\Omega}\bm{\Omega}_i)+\bm{\theta}^\top\bm{\Omega}_i\bm{\theta}}
      {s_{\bm{\theta}}^2}\ge0.
\]
If $b=\bm{\theta}^\top\bm{\Omega}\bm{\theta}$, then
\[
 \sum_{i=1}^2q_i(\bm{\theta})=\frac{A+b}{2A+4b}\le\frac12.
\]
Therefore
\begin{align*}
 B_{\bm{\theta}}
 &=\sum_{i=1}^2 n_i\{q_i(\bm{\theta})/n_i\}^{3/2}
 \le\frac{\sum_{i=1}^2q_i(\bm{\theta})^{3/2}}{\sqrt{n_{\min}}}
 \\
 &\le\frac{\{\sum_{i=1}^2q_i(\bm{\theta})\}^{3/2}}{\sqrt{n_{\min}}}
 \le\frac{1}{2^{3/2}\sqrt{n_{\min}}}.
\end{align*}
The theorem proves the assertion with constants independent of $p,\bm{\theta}$,
and the covariance matrices.
\end{proof}

\section{Density bounds and sharper comparison}\label{app:adaptive}

\label{adap:section}

We prove the density bound used in Theorem~\ref{thm:spread-summary}
and then give refinements for low-rank and mean-dominated references.
The square-root bound is needed near rank one, whereas sufficient
spectral spread gives a bounded density uniformly over deterministic
signals. The results below provide explicit conditions and constants
for positive-semidefinite mean-testing forms. General quadratic--linear
Gaussian forms, including signed coefficients, are treated by
\citet{ChoudharyKuchibhotla2026}.

Let $\bm{G}\sim N(\bm{0},\bm{\Omega})$, where $\bm{\Omega}\succeq0$, and write
\[
 A=\operatorname{tr}(\bm{\Omega}^2)>0,\qquad
 b=\bm{\theta}^{\top}\bm{\Omega}\bm{\theta},\qquad
 s^2=2A+4b,\qquad
 W_{\bm{\theta}}^{\mathrm G}
 =\frac{\|\bm{G}+\bm{\theta}\|^2-\operatorname{tr}(\bm{\Omega})-\|\bm{\theta}\|^2}{s}.
\]
Let $\lambda_1\geq\cdots\geq\lambda_p\geq0$ be the eigenvalues
of $\bm{\Omega}$.  Put
\begin{equation}
 c_3=\frac{\Gamma(1/4)}{4\sqrt{\pi}\,\Gamma(3/4)},
 \qquad
 \mathcal C_{\bm{\theta}}(\varepsilon)
 =\sup_{x\in\mathbb R}
       \Pr\{|W_{\bm{\theta}}^{\mathrm G}-x|\leq\varepsilon\}.
 \label{adap:definitions}
\end{equation}
Here $\Gamma(\cdot)$ is the gamma function and $c_3$ is an absolute constant.

\subsection{Three spectral blocks and uniformity over the mean}

\begin{lem}[An integrable characteristic-function envelope]
\label{adap:fourier}
Suppose the eigenvalue indices can be partitioned into three disjoint
sets $I_1,I_2,I_3$ such that
\begin{equation}
 \sum_{j\in I_\ell}\lambda_j^2\geq m>0,
 \qquad \ell=1,2,3.
 \label{adap:three-blocks}
\end{equation}
Then $W_{\bm{\theta}}^{\mathrm G}$ has a continuous bounded density
$f_{\bm{\theta}}$, and, for every $\bm{\theta}$,
\begin{equation}
 \|f_{\bm{\theta}}\|_{\infty}\leq c_3\frac{s}{\sqrt m}.
 \label{adap:basic-density}
\end{equation}
If $b>0$, a second bound is
\begin{equation}
 \|f_{\bm{\theta}}\|_{\infty}
 \leq\frac{s}{2\sqrt{\pi b}}
       +c_3\frac{s}{\sqrt m}
                \exp\left(-\frac{b}{4\lambda_1^2}\right).
 \label{adap:noncentral-density}
\end{equation}
Both inequalities hold simultaneously.
\end{lem}

\begin{proof}
Choose an eigenbasis of $\bm{\Omega}$, and let $\theta_j$ be the signal
coordinates in that basis.  Gaussian integration gives the characteristic-function modulus of
$T_{\bm{\theta}}^{\mathrm G}=\|\bm{G}+\bm{\theta}\|^2$:
\begin{equation}
 |\varphi_{T_{\bm{\theta}}^{\mathrm G}}(t)|
 =\prod_j(1+4t^2\lambda_j^2)^{-1/4}
   \exp\left\{-\sum_j
        \frac{2t^2\lambda_j\theta_j^2}{1+4t^2\lambda_j^2}\right\}.
 \label{adap:exact-cf}
\end{equation}
Zero eigenvalues contribute only a deterministic shift and do not
affect this modulus.  Within each block,
\[
 \prod_{j\in I_\ell}(1+4t^2\lambda_j^2)
 \geq1+4t^2\sum_{j\in I_\ell}\lambda_j^2
 \geq1+4mt^2.
\]
Consequently the product in \eqref{adap:exact-cf}, denoted below by
$P(t)$, is bounded by $(1+4mt^2)^{-3/4}$.  It is integrable, and
\[
 \frac1{2\pi}\int_{\mathbb R}P(t)\,dt
 \leq\frac1{2\pi}\int_{\mathbb R}(1+4mt^2)^{-3/4}\,dt
 =\frac{c_3}{\sqrt m}.
\]
Fourier inversion therefore gives a continuous bounded density for
$T_{\bm{\theta}}^{\mathrm G}$, with the displayed upper bound.  Centering
and dividing by $s$ multiplies its density bound by $s$, proving
\eqref{adap:basic-density}.

For the second bound, observe that the exponent in
\eqref{adap:exact-cf} is bounded above by
\[
 -\frac{2t^2b}{1+4t^2\lambda_1^2}.
\]
On $|t|\leq(2\lambda_1)^{-1}$ this is at most $-bt^2$.
On the complementary region it is at most
$-b/(4\lambda_1^2)$.  Dropping $P(t)\leq1$ on the first region
and using its integrable envelope on the second gives
\begin{align*}
 \frac1{2\pi}\int_{\mathbb R}
       |\varphi_{T_{\bm{\theta}}^{\mathrm G}}(t)|\,dt
 &\leq\frac1{2\pi}\int_{\mathbb R}e^{-bt^2}\,dt
   +\frac{e^{-b/(4\lambda_1^2)}}{2\pi}
       \int_{\mathbb R}P(t)\,dt\\
 &\leq\frac1{2\sqrt{\pi b}}
       +\frac{c_3}{\sqrt m}e^{-b/(4\lambda_1^2)}.
\end{align*}
Scaling proves \eqref{adap:noncentral-density}.
\end{proof}

\begin{theorem}[Uniform bounded density from spectral mass]
\label{adap:uniform-density}
Suppose \eqref{adap:three-blocks} holds with $m\geq cA$ for some
$c\in(0,1/3]$.  Define the explicit constant
\begin{equation}
 M_c=\max\left\{
       \frac{\sqrt6\,c_3}{\sqrt c},\quad
       \frac{\sqrt6}{2\sqrt\pi}
       +\frac{2\sqrt2 e^{-3/8}c_3}{\sqrt c}\right\}.
 \label{adap:constant}
\end{equation}
Then, uniformly over all deterministic $\bm{\theta}\in\mathbb R^p$,
\begin{equation}
 \|f_{\bm{\theta}}\|_{\infty}\leq M_c,
 \qquad
 \mathcal C_{\bm{\theta}}(\varepsilon)
       \leq\min\{1,2M_c\varepsilon\},
       \quad\varepsilon>0.
 \label{adap:uniform-smallball}
\end{equation}
No upper bound on $b/A$ or on $\|\bm{\theta}\|$ is required.
\end{theorem}

\begin{proof}
If $b\leq A$, then $s\leq\sqrt{6A}$, and
\eqref{adap:basic-density} gives the first entry in the maximum
defining $M_c$.  If $b>A$, write $r=b/A>1$.  Since
$\lambda_1^2\leq A$, \eqref{adap:noncentral-density} gives
\[
 \|f_{\bm{\theta}}\|_{\infty}
 \leq\frac{\sqrt6}{2\sqrt\pi}
   +\frac{c_3}{\sqrt c}\sqrt{2+4r}\,e^{-r/4}.
\]
Elementary differentiation shows
\[
 \sup_{r\geq1}\sqrt{2+4r}\,e^{-r/4}
       =2\sqrt2 e^{-3/8},
\]
with the maximum attained at $r=3/2$.  This proves the density bound.
Integration over an interval of length $2\varepsilon$ proves the
concentration bound.
\end{proof}

\Needspace{8\baselineskip}
\begin{prop}[Simple ways to verify the block condition]
\label{adap:verify-blocks}
Let $m_3$ be the largest possible value of
$\min_{1\leq\ell\leq3}\sum_{j\in I_\ell}\lambda_j^2$ over
partitions into three sets; set $m_3=0$ if no partition has three
positive masses.  The following lower bounds are sufficient:
\begin{align}
 m_3&\geq\frac{A-2\lambda_1^2}{3}
       &&\text{if }\lambda_1^2<A/2,
       \label{adap:diffuse-blocks}\\
 m_3&\geq\min\left\{\lambda_1^2,
                 \frac{A-\lambda_1^2-\lambda_2^2}{2}\right\}
       &&\text{if }A>\lambda_1^2+\lambda_2^2.
       \label{adap:spike-blocks}
\end{align}
In particular, $\lambda_1^2\leq A/4$ implies $m_3\geq A/6$.
For $\bm{\Omega}_p=\operatorname{diag}(\sqrt p,1,\ldots,1)$ and
$p\geq3$, one has
\begin{equation}
 A_p=2p-1,\qquad
 m_3\geq\frac{p-2}{2},\qquad
 \liminf_{p\to\infty}\frac{m_3}{A_p}\geq\frac14.
 \label{adap:running-example}
\end{equation}
Thus the running spike-plus-bulk example belongs to the uniform
bounded-density regime, even though its null limit is nonnormal.
For every $p\geq5$, one may take $c=1/6$.
\end{prop}

\begin{proof}
For \eqref{adap:diffuse-blocks}, put
$a=(A-2\lambda_1^2)/3>0$.  Accumulate eigenvalue squares until their
sum first reaches $a$ to form the first block; its mass is at most
$a+\lambda_1^2$.  Repeat for the second block.  There is enough mass
for both steps, and the remaining mass is at least
$A-2(a+\lambda_1^2)=a$.  This proves the first bound.
For \eqref{adap:spike-blocks}, use the first eigenvalue as one block.
The tail has total mass $B=A-\lambda_1^2$ and maximum summand
$\lambda_2^2$.  Accumulate a tail block until it first reaches
$(B-\lambda_2^2)/2$.  That block has mass at most
$(B+\lambda_2^2)/2$, and the remaining tail mass is at least
$(B-\lambda_2^2)/2$.  These two blocks together with the first
eigenvalue prove the second bound.  Substitution gives the remaining
claims.
\end{proof}

\subsection{A mean-dominated bound without a spectral spread condition}

A large mean provides another route to anti-concentration, even for
low-rank covariance. The following result uses the linear contribution
while retaining its dependence on the quadratic term.

\begin{prop}[Conditional linearization of the Gaussian reference]
\label{adap:mean-dominated}
For every nonzero covariance $\bm{\Omega}\succeq0$ and every $\bm{\theta}$
such that $b=\bm{\theta}^{\top}\bm{\Omega}\bm{\theta}>0$,
\begin{equation}
 \mathcal C_{\bm{\theta}}(\varepsilon)
 \leq\min\left\{1,
       \frac4{\sqrt{2\pi}}\frac{s\varepsilon}{\sqrt b}
       +4\exp\left(-\frac{b}{32\lambda_1^2}\right)\right\}.
 \label{adap:mean-smallball}
\end{equation}
For example, if $b/A\geq r_0>0$, the coefficient of
$\varepsilon$ is bounded in terms of $r_0$ only.  The additive
remainder tends to zero if $b/\lambda_1^2\to\infty$, including the
mean-dominated regime $b/A\to\infty$.
\end{prop}

\begin{proof}
Represent $\bm{G}=\bm{\Omega}^{1/2}\bm{Z}$ with $\bm{Z}\sim N(\bm{0},\bm{I}_p)$ and put
$\bm{v}=\bm{\Omega}^{1/2}\bm{\theta}/\sqrt b$, so that $\|\bm{v}\|=1$.
Write $\bm{Z}=z\bm{v}+\bm{Z}_{\perp}$, where $z\sim N(0,1)$ is independent of
$\bm{Z}_{\perp}$.  Conditional on $\bm{Z}_{\perp}$, the unstandardized
centered reference is a quadratic polynomial in $z$:
\[
 a z^2+2(\sqrt b+c)z+q,
 \qquad
 a=\bm{v}^{\top}\bm{\Omega} \bm{v}\in[0,\lambda_1],\qquad
 c=\bm{v}^{\top}\bm{\Omega} \bm{Z}_{\perp}.
\]
Here $q$ is measurable with respect to $\bm{Z}_{\perp}$.  The variable
$c$ is centered Gaussian with variance
\[
 \bm{v}^{\top}\bm{\Omega}(\bm{I}-\bm{v}\bm{v}^{\top})\bm{\Omega} \bm{v}\leq\lambda_1^2.
\]
Thus the probability of $|c|>\sqrt b/2$ is at most
$2\exp\{-b/(8\lambda_1^2)\}$.  On the complementary event, and
for $|z|\leq\sqrt b/(4\lambda_1)$, the derivative of the
conditional polynomial is at least
\[
 2(\sqrt b+c)-2a|z|\geq\sqrt b/2.
\]
On this interval the polynomial is strictly increasing.  The inverse
image of any interval of length $2s\varepsilon$ therefore has length
at most $4s\varepsilon/\sqrt b$.  Since the standard normal density
is at most $(2\pi)^{-1/2}$, its conditional probability is bounded
by the first term in \eqref{adap:mean-smallball}.  The omitted event
$|z|>\sqrt b/(4\lambda_1)$ has probability at most
$2\exp\{-b/(32\lambda_1^2)\}$.  Adding the two exceptional probabilities and bounding the faster-decaying
term by the slower one proves the result. The conditioning argument
retains the dependence between the quadratic and linear components.
\end{proof}

\begin{remark}[A two-eigenvalue bound and the rank-one boundary]
\label{adap:low-rank}
If $\lambda_2>0$, an additional elementary bound is
\begin{equation}
 \mathcal C_{\bm{\theta}}(\varepsilon)
       \leq\min\{1,s\varepsilon/\sqrt{\lambda_1\lambda_2}\}.
 \label{adap:two-eigenvalues}
\end{equation}
Indeed, the Gaussian density in the first two eigen-coordinates is at
most $(2\pi\sqrt{\lambda_1\lambda_2})^{-1}$ for every mean.
An annulus whose squared-radius interval has length $h$ has area at
most $\pi h$.  Conditioning on all remaining coordinates proves
\eqref{adap:two-eigenvalues}, with $h=2s\varepsilon$.
This bound is useful for two comparable spikes and bounded $b/A$;
it is not asserted to be uniform in the mean after standardization.

At the other boundary, take $\bm{\Omega}$ of rank one and $\bm{\theta}=\bm{0}$.
Then $W_{0}^{\mathrm G}=(Z^2-1)/\sqrt2$, and
\[
 \Pr\left\{\left|W_{0}^{\mathrm G}+\frac1{\sqrt2}\right|
                    \leq\varepsilon\right\}
 =\Pr\{|Z|\leq2^{1/4}\varepsilon^{1/2}\}
 \sim\frac{2^{3/4}}{\sqrt\pi}\,\varepsilon^{1/2}.
\]
Thus a bound linear in $\varepsilon$ cannot hold uniformly over all
covariances and signals; an additional condition is needed for the
bounded-density result.
\end{remark}

\subsection{Consequences for the full-numerator comparison}

\begin{theorem}[A faster bound under spectral spread]
\label{adap:transfer}
Use the independent increments and moment conditions of
Theorem~\ref{thm:finite-transfer}, with $D_{\bm{\theta}}$ and $B_{\bm{\theta}}$
as defined there.  Suppose the aggregate covariance has a partition
satisfying \eqref{adap:three-blocks} with $m\geq cA$, for a fixed
$c>0$.  Then, for a constant $C_{K,c}$ depending only on the moment
constant $K$ and $c$,
\begin{align}
 d_{\mathrm K}(T_{\bm{\theta}},T_{\bm{\theta}}^{\mathrm G})
 &\leq\min\left[1,C_{K,c}\inf_{0<\varepsilon\leq1}
       \left\{D_{\bm{\theta}}\varepsilon^{-2}
                    +B_{\bm{\theta}}\varepsilon^{-3}+\varepsilon\right\}\right]
       \label{adap:transfer-epsilon}\\
 &\leq\min\{1,C_{K,c}(D_{\bm{\theta}}^{1/3}+B_{\bm{\theta}}^{1/4})\}.
       \label{adap:transfer-rate}
\end{align}
The constants are independent of the dimension and the mean signal.
Under the two-sample conditions of
Corollary~\ref{cor:grouped-transfer}, this implies
\begin{equation}
 \sup_{\bm{\theta}\in\mathbb R^p}
 d_{\mathrm K}(T_{\bm{\theta}},T_{\bm{\theta}}^{\mathrm G})
 \leq\min\{1,C_{K,c}(n_{\min}^{-1/3}+n_{\min}^{-1/8})\}.
 \label{adap:grouped-rate}
\end{equation}
\end{theorem}

\begin{proof}
The smoothing proof of Theorem~\ref{thm:finite-transfer} establishes,
before invoking any particular anti-concentration bound,
\[
 d_{\mathrm K}(T_{\bm{\theta}},T_{\bm{\theta}}^{\mathrm G})
 \leq C_KD_{\bm{\theta}}\varepsilon^{-2}
       +C_KB_{\bm{\theta}}\varepsilon^{-3}
       +\mathcal C_{\bm{\theta}}(3\varepsilon).
\]
The radius $3\varepsilon$ covers the intervals in the smoothing
argument. Applying
Theorem~\ref{adap:uniform-density} gives the first inequality.
If $D_{\bm{\theta}},B_{\bm{\theta}}\leq1$,
choose $\varepsilon=\max\{D_{\bm{\theta}}^{1/3},B_{\bm{\theta}}^{1/4}\}$;
each term is bounded by a constant times
$D_{\bm{\theta}}^{1/3}+B_{\bm{\theta}}^{1/4}$.  If one quantity exceeds one,
the trivial bound $d_{\mathrm K}\leq1$ suffices.  The two-sample
influence bounds already proved in
Corollary~\ref{cor:grouped-transfer} give
$D_{\bm{\theta}}=O_K(n_{\min}^{-1})$ and
$B_{\bm{\theta}}=O(n_{\min}^{-1/2})$, uniformly over $\bm{\theta}$.
Substitution proves \eqref{adap:grouped-rate}.
\end{proof}

\begin{cor}[Large-signal refinement without three spectral blocks]
\label{adap:mean-transfer}
Under the increment conditions of Theorem~\ref{thm:finite-transfer},
suppose $b/A\geq r_0>0$.  Then
\begin{equation}
 d_{\mathrm K}(T_{\bm{\theta}},T_{\bm{\theta}}^{\mathrm G})
 \leq\min\left\{1,
       C_{K,r_0}(D_{\bm{\theta}}^{1/3}+B_{\bm{\theta}}^{1/4})
       +4e^{-b/(32\lambda_1^2)}\right\}.
 \label{adap:mean-transfer-rate}
\end{equation}
This statement allows rank-one covariance.  A vanishing remainder
requires $b/\lambda_1^2\to\infty$; it does not follow merely from
$b/A$ being bounded away from zero.
\end{cor}

\begin{proof}
Here $s/\sqrt b=\sqrt{4+2A/b}\leq\sqrt{4+2/r_0}$.
Use Proposition~\ref{adap:mean-dominated} in the same smoothing
inequality and optimize $\varepsilon$ as in
Theorem~\ref{adap:transfer}.
\end{proof}

\begin{cor}[Random critical values on the null fluctuation scale]
\label{adap:random-critical}
Suppose the covariance satisfies the condition of
Theorem~\ref{adap:uniform-density}, and put $s_0=\sqrt{2A}$.
Let $q$ be any deterministic critical value and $\widehat q$ any
data-dependent critical value, which need not be independent of
$T_{\bm{\theta}}$.  For every $\delta>0$,
\begin{align}
 &\left|\Pr\{T_{\bm{\theta}}>\widehat q\}
                -\Pr\{T_{\bm{\theta}}^{\mathrm G}>q\}\right|\notag\\
 &\quad\leq d_{\mathrm K}(T_{\bm{\theta}},T_{\bm{\theta}}^{\mathrm G})
       +\Pr\{|\widehat q-q|>\delta s_0\}+2M_c\delta.
 \label{adap:random-power}
\end{align}
Thus threshold estimation has a linear anti-concentration cost on
this spectral class, uniformly over all two-sample mean differences.
\end{cor}

\begin{proof}
On $|\widehat q-q|\leq\delta s_0$, bracket the rejection event by
$\{T_{\bm{\theta}}>q+\delta s_0\}$ and
$\{T_{\bm{\theta}}>q-\delta s_0\}$.  Replace each deterministic-threshold
probability by its Gaussian companion, paying the Kolmogorov distance
once in either direction.  The intervening Gaussian interval has
probability at most $2M_c\delta s_0/s\leq2M_c\delta$, since
$s\geq s_0$.  The complementary threshold-estimation event contributes
its displayed probability.
\end{proof}

\begin{remark}[Regimes covered by the bounds]
The general comparison gives $D^{1/5}+B^{1/7}$ for every covariance.
Under the three-block condition, the bound improves to
$D^{1/3}+B^{1/4}$ uniformly over signals, covering diffuse spectra
and persistent spikes with sufficiently spread remaining mass.
The mean-dominated result gives the same smoothing exponents for
low-rank spectra with an additional exponential remainder.
These are sufficient regimes for finite-sample upper bounds;
sharp rates and a complete classification are not established.
\end{remark}

\section{Additional power statements and approximation boundaries}\label{app:power-details}\label{app:profiles}
This appendix records the eigenspace interpretation and approximation boundaries used in Section~\ref{prof:section}. The corresponding proofs are given below.

\subsection{Eigenspaces and uniformity}
\begin{remark}[Null-space signal and multiplicities]
\label{prof:multiplicity}
Null-space signal contributes to $\kappa_n$ but not to $L_n^*$,
so the theorem allows singular covariance matrices. Signal energy
in a diffuse spectral tail likewise enters the limit only through
$\kappa$.

The eigenvectors themselves need not converge; the relevant quantities
are the signal projections on the moving covariance eigenspaces.
For a positive limiting eigenvalue $r$ of multiplicity $m_r$, let
$h_r$ be the limiting total normalized signal energy in its eigenspace
cluster. The corresponding part of $\zeta+L$ has distribution
\[
 \frac{r}{\sqrt{2}}
 \left\{\chi^2_{m_r}(h_r/r)-m_r-h_r/r\right\}.
\]
Thus only $h_r$, not an arbitrary basis within that cluster, matters.
Convergence of these cluster energies can replace convergence of the
individual signal coordinates.  To see this for a cluster whose
normalized eigenvalues all tend to $r$, replace its finitely many
coefficients by $r$; bounded local signal energy makes this replacement
$L^2$-negligible.  Rotational invariance of the resulting Gaussian block
then gives the displayed law.  Each such positive cluster is finite
because $\sum_j\rho_j^2\leq1$.
\end{remark}

\begin{remark}[Sequential and uniform statements]
\label{prof:uniformity}
Theorem~\ref{prof:theorem} gives a limit along a specified covariance
and signal sequence. Both the normalized eigenvalues and the signal
energies in persistent eigenspaces determine local power.
The uniform comparison in Theorem~\ref{thm:finite-transfer} instead
uses the Gaussian reference for each member of the sequence and
requires no common limiting profile. A uniform limiting power formula
would additionally require uniform approximation of the spectral and
signal profiles and control of the quantile step. Although
\eqref{prof:linear-tail} is uniform over $\kappa_n\leq K$,
convergence of the finite leading eigenspace energies is still needed.
\end{remark}

\subsection{Approximation boundaries}
\label{prof:matching-section}

After Gaussian-reference transfer, a separate approximation is needed
if the reference is replaced by a parametric distribution. The following
example shows that two- or three-cumulant matching may retain a
nonvanishing error, even for Gaussian observations. These approximations are useful in existing normal-reference procedures \citep{Zhang2005,ZhangZhuBS2022,ZhangZhuCQ2022}; the example identifies a spectrum for which matching a fixed number of cumulants does not recover the full limiting reference.

\begin{prop}[Persistent error of two- and three-cumulant matching]
\label{prof:matching}
Use the covariance sequence of Proposition~\ref{prof:spike}.  Put
\[
 B_p=\operatorname{tr}(\bm{\Omega}_p)=p-1+\sqrt p,
 \qquad C_p=\operatorname{tr}(\bm{\Omega}_p^3)=p-1+p^{3/2}.
\]
Let $T_{0,p}^*=\|\bm{G}_p\|^2$, and define the population two- and
three-cumulant matched references by
\begin{align*}
 R_{2,p}&=\beta_{2,p}\chi^2_{\nu_{2,p}},
 &\beta_{2,p}&=A_p/B_p,
 &\nu_{2,p}&=B_p^2/A_p,\\
 R_{3,p}&=a_{3,p}+\beta_{3,p}\chi^2_{\nu_{3,p}},
 &\beta_{3,p}&=C_p/A_p,
 &\nu_{3,p}&=A_p^3/C_p^2,
 &a_{3,p}&=B_p-A_p^2/C_p.
\end{align*}
Then
\begin{align}
 \frac{T_{0,p}^*-B_p}{s_{0,p}}
 &\xrightarrow{d}\frac{Z_1^2-1}{2}+\frac{Z_0}{\sqrt{2}},
 \label{prof:true-null}\\
 \frac{R_{2,p}-B_p}{s_{0,p}}&\xrightarrow{d}N(0,1),
 \label{prof:2c-null}\\
 \frac{R_{3,p}-B_p}{s_{0,p}}&\xrightarrow{d}
                  \frac{\chi^2_8-8}{4}.
 \label{prof:3c-null}
\end{align}
All three limits are different.  Consequently the Kolmogorov distances
between the exact central Gaussian reference and either matched
reference converge to strictly positive constants.  In particular,
consistency of estimators of the matching parameters cannot, by itself,
remove this population approximation error.
\end{prop}

\begin{proof}See Appendix~\ref{app:profiles}.\end{proof}

\begin{remark}
Divergence of the two-cumulant degrees of freedom does not imply
asymptotic normality of the exact reference. Three-cumulant matching
also need not yield asymptotically exact level under persistent spectra.
These approximations may nevertheless be numerically accurate in
particular settings. Their error must therefore be assessed separately
from the Gaussian-reference transfer error.
\end{remark}

A second limitation concerns Gaussian-reference transfer itself.

The full squared mean statistic also contains diagonal terms.
The following example shows that bounded projection fourth moments
do not control their fluctuations. Thus a comparison result for a
diagonal-deleted statistic requires an additional argument before
it can be applied to the full numerator.

This construction also occurs within the two-sample model. Take $n_1=n_2=m$, so $n=2m$, and let both centered sample distributions equal that of $R\bm{Z}$. The normalized sample-mean difference is a signed sum of $2m$ independent such vectors divided by $\sqrt{2m}$. By symmetry, it has the same distribution as $\bm{W}_n$ in the next proposition. Thus the example is a two-sample obstruction to reference transfer, even when the two population covariance matrices coincide.

\begin{prop}[Projection fourth moments do not control the diagonal]
\label{prof:radial}
Let $\bm{X}_i=R_i \bm{Z}_i\in\mathbb{R}^{p}$, $i=1,\ldots,n$, where
$\bm{Z}_i\sim N(\bm{0},\bm{I}_p)$ are independent, the $R_i$ are independent and
identically distributed, and all radii are independent of all Gaussian
vectors.  Suppose
\[
 \operatorname{E}R_i^2=1,
 \qquad v=\operatorname{Var}(R_i^2)\in(0,\infty).
\]
For $\bm{W}_n=n^{-1/2}\sum_i\bm{X}_i$, $T_n=\|\bm{W}_n\|^2$, and
$S_n=n^{-1}\sum_iR_i^2$, one has
\begin{equation}
 T_n\ \stackrel{d}{=}\ S_n U_p,
 \qquad U_p\sim\chi^2_p\text{ independent of }S_n,
 \label{prof:radial-mixture}
\end{equation}
and hence
\begin{equation}
 \operatorname{E}T_n=p,
 \qquad
 \operatorname{Var}(T_n)=2p+\frac{v(p^2+2p)}{n}.
 \label{prof:radial-variance}
\end{equation}
If $n,p\to\infty$ with $p/n\to c\in(0,\infty)$, then
\begin{equation}
 \frac{T_n-p}{\sqrt{2p}}
 \xrightarrow{d}N\left(0,1+\frac{vc}{2}\right),
 \label{prof:radial-limit}
\end{equation}
whereas its covariance-matched Gaussian reference has limit $N(0,1)$.
Nevertheless, for every deterministic $\bm{a}\in\mathbb{R}^p$,
\begin{equation}
 \operatorname{E}(\bm{a}^{\top}\bm{X}_i)^4
       =3(1+v)\|\bm{a}\|^4.
 \label{prof:projection-fourth}
\end{equation}
Thus a projection fourth-moment bound with a constant independent of
dimension does not, by itself, validate the Gaussian reference for the
full numerator.
\end{prop}

\begin{proof}See Appendix~\ref{app:profiles}.\end{proof}

\begin{prop}[Exact diagonal contribution in the radial example]
\label{prof:diagonal}
In the setting of Proposition~\ref{prof:radial}, decompose
\[
 T_n-p=D_n+U_n,\qquad
 D_n=\frac1n\sum_{i=1}^n(\|\bm{X}_i\|^2-p),\qquad
 U_n=\frac2n\sum_{i<j}\bm{X}_i^{\top}\bm{X}_j.
\]
Then $\operatorname{Cov}(D_n,U_n)=0$ and
\begin{equation}
 \frac{\operatorname{Var}(D_n)}{2p}
 =\frac{1+v}{n}+\frac{vp}{2n},
 \qquad
 \operatorname{Var}(U_n)=2p\frac{n-1}{n}.
 \label{prof:diagonal-size}
\end{equation}
If $\bm{u}_i=\bm{X}_i/\sqrt n$ and $\bm{\Gamma}_i=\operatorname{Cov}(\bm{u}_i)=\bm{I}_p/n$,
then the projection condition
$\operatorname{E}(\bm{a}^{\top}\bm{u}_i)^4
\leq C(\bm{a}^{\top}\bm{\Gamma}_i \bm{a})^2$ holds with $C=3(1+v)$, but
\begin{equation}
 \frac{\operatorname{Var}(\|\bm{u}_i\|^2)}{
              \operatorname{tr}(\bm{\Gamma}_i^2)}
 =2(1+v)+vp\longrightarrow\infty.
 \label{prof:norm-condition-failure}
\end{equation}
\end{prop}

\begin{proof}See Appendix~\ref{app:profiles}.\end{proof}

\begin{remark}[Scope of the boundary example]
\label{prof:radial-scope}
The radii may be bounded: for example, take $R_i^2=1/2$ or $3/2$
with equal probabilities, so $v=1/4$.  The obstruction therefore does
not require heavy tails.  In this particular model, the additional
radial term is negligible when $vp/n\to0$; a uniform squared-norm
variance assumption is one sufficient route, but is stronger than the
aggregate negligibility needed for this example.  When $p/n\to c>0$
and $v>0$, the extra fluctuation remains even under the null hypothesis.
For the full numerator, one must therefore control the diagonal
fluctuation on the reference scale, modify the statistic, or include
the additional radial term in the reference. A result for $U_n$ alone
is insufficient. Common-radial-scale elliptical data illustrate this
boundary unless the required radial concentration is verified.
\end{remark}

\begin{cor}[A power-error bound for an estimated critical value]
\label{cor:shift-feasible-bound}
Fix $\alpha\in(0,1)$. Let $T_{\bm\theta}$ denote the original-data
statistic, set
$e_{\bm\theta}=d_{\mathrm K}(T_{\bm\theta},T_{\bm\theta}^*)$, and let
$\widehat q_{\alpha}$ be any data-dependent critical value. For every
$\varepsilon>0$,
\begin{align}
 &\left|\Pr\{T_{\bm\theta}>\widehat q_{\alpha}\}
   -\Pr\{T_0^*+\|\bm\theta\|^2>q_{\alpha}^{\mathrm{NR}}\}\right|
 \notag\\
 &\qquad\leq e_{\bm\theta}
       +\min\{1/2,C\Delta_{\bm\theta}^{1/4}\}
       +\Pr\{|\widehat q_{\alpha}-q_{\alpha}^{\mathrm{NR}}|
                        >\varepsilon\sqrt{2A}\}
       +C\varepsilon^{1/2}.
 \label{eq:shift-feasible-bound}
\end{align}
Under the three-block condition in Theorem~\ref{thm:shift-boundary},
the second and fourth terms can respectively be replaced by
$\min\{1/2,C_c\Delta_{\bm\theta}^{3/4}\}$ and $C_c\varepsilon$.
No independence between $T_{\bm\theta}$ and $\widehat q_{\alpha}$ is assumed.
\end{cor}
\begin{proof}See Appendix~\ref{app:shift-boundary-proofs}.\end{proof}

\begin{prop}[A fixed cumulant summary can lose signal information]
\label{prop:finite-cumulant-information}
For every fixed integer $J\geq1$, there exist a positive definite
$(J+1)\times(J+1)$ covariance matrix with eigenvalues in $[1,2]$ and
two unit-norm signals such that their Gaussian-reference distributions
have equal first $J$ cumulants and unequal $(J+1)$st cumulants. The
covariance matrix is the same for the two signals. By padding the
covariance and signals with zero coordinates, the two distinct laws
and their positive Kolmogorov distance remain unchanged as $p\to\infty$.
Consequently, a reference determined only by the covariance and those
$J$ cumulants cannot approximate both sequences with vanishing
Kolmogorov error.
\end{prop}
\begin{proof}See Appendix~\ref{app:spectral-information-proof}.\end{proof}

This is an information limitation of a fixed cumulant summary over
unrestricted spectra. It does not preclude accurate cumulant
approximations in specific spectral regimes or under additional
restrictions on the alternatives.

\subsection{Signal information and cumulant summaries}\label{app:spectral-information-proof}
\begin{proof}[Proof of Proposition~\ref{prop:spectral-signal-information}]
For $\lambda>0$, let $m_\lambda$ denote its multiplicity. Orthogonal
Gaussian projections are independent, and rotational invariance within
each eigenspace gives
\begin{equation}
 T_{\bm\theta}^*\ \stackrel{d}{=}\
 e_0(\bm\theta)+\sum_{\lambda\in\mathcal L:\lambda>0}
     \lambda\chi^2_{m_\lambda}
       \{e_\lambda(\bm\theta)/\lambda\},
 \label{eq:eigenspace-noncentral-decomposition}
\end{equation}
where the chi-squared variables are independent and $e_0=0$ when
zero is not an eigenvalue. This proves that (ii) implies (i), including
the case $\bm\Omega=\bm0$.

For $t$ in a neighborhood of zero, the cumulant-generating function is
\begin{equation}
 K_{\bm\theta}(t)
 =-\frac12\sum_{\lambda>0}m_\lambda\log(1-2t\lambda)
   +\sum_{\lambda\in\mathcal L}
        \frac{t e_\lambda(\bm\theta)}{1-2t\lambda}.
 \label{eq:eigenspace-cgf}
\end{equation}
The expression for $\lambda=0$ is simply $t e_0(\bm\theta)$.
Expanding the logarithms and geometric series proves
\eqref{eq:noncentral-cumulants-information}. If (i) holds, these
cumulants agree. Since $\bm\Omega$ is fixed, equality of cumulants
of orders $1,\ldots,|\mathcal L|$ implies (iii).

Finally, write $s=|\mathcal L|$ and enumerate its distinct values as
$\lambda_1,\ldots,\lambda_s$. By the spectral theorem,
\[
 \bm\theta^\top\bm\Omega^j\bm\theta
   -\bm\vartheta^\top\bm\Omega^j\bm\vartheta
 =\sum_{\ell=1}^s\lambda_\ell^j
   \{e_{\lambda_\ell}(\bm\theta)
      -e_{\lambda_\ell}(\bm\vartheta)\},
 \qquad 0\leq j\leq s-1.
\]
For $j=0$, every coefficient is $1$, including that corresponding to
zero. The coefficient matrix is a Vandermonde matrix with determinant
$\prod_{k<\ell}(\lambda_\ell-\lambda_k)\ne0$. Hence (iii) implies
(ii), completing the equivalence.
\end{proof}

\begin{proof}[Proof of Proposition~\ref{prop:finite-cumulant-information}]
Set $\lambda_\ell=1+\ell/J$ for $\ell=0,\ldots,J$ and take
$\bm\Omega=\diag(\lambda_0,\ldots,\lambda_J)$. Define
\[
 h_\ell=(-1)^{J-\ell}\binom{J}{\ell},\qquad
 h_\ell^+=\max(h_\ell,0),\qquad
 h_\ell^-=\max(-h_\ell,0),\qquad H=2^{J-1}.
\]
The finite-difference identities give
\begin{equation}
 \sum_{\ell=0}^J h_\ell\lambda_\ell^q=0
 \quad(0\leq q<J),\qquad
 \sum_{\ell=0}^J h_\ell\lambda_\ell^J=J!/J^J>0.
 \label{eq:cumulant-finite-difference}
\end{equation}
These identities follow by taking the $J$th forward difference,
with step $1/J$, of the polynomial $x^q$ at $x=1$.
Also $\sum_\ell h_\ell^+=\sum_\ell h_\ell^-=H$, because
$\sum_\ell h_\ell=0$ and $\sum_\ell|h_\ell|=2^J$.
Choose signals with coordinates
\[
 \theta_{+,\ell}=\sqrt{h_\ell^+/H},\qquad
 \theta_{-,\ell}=\sqrt{h_\ell^-/H},
 \qquad 0\leq\ell\leq J.
\]
Both have unit norm. Equations~\eqref{eq:cumulant-finite-difference}
and~\eqref{eq:noncentral-cumulants-information} show equality of
the first $J$ cumulants and a strictly positive difference in the
$(J+1)$st cumulant. Thus their CDFs $F_+$ and $F_-$ are distinct, and
$\delta_J=\|F_+-F_-\|_\infty>0$.

Appending zero covariance coordinates and zero signal coordinates
does not change either quadratic form. All nonzero eigenvalues and
the signal norm remain fixed, so the normalized local signal strength
is a fixed positive constant along this sequence. Any reference rule
that uses only the covariance and the first $J$ cumulants must return
the same CDF $H_p$ for the two alternatives at each $p$. The triangle
inequality gives
\[
 \max\{\|F_+-H_p\|_\infty,\|F_--H_p\|_\infty\}
 \geq\tfrac12\|F_+-F_-\|_\infty=\delta_J/2.
\]
The bound holds for every $p\geq J+1$ and proves the assertion.
\end{proof}

\subsection{Rank-one approximation error and power}
\label{app:rankone-sharp-proof}

\begin{proof}[Proof of Proposition~\ref{prop:rankone-sharp-error}]
The distributions depend only on $|a|$, so assume $a\ge0$.
The case $a=0$ is immediate. For $x\ge0$, the distribution functions are
\[
 F_a(x)=\Phi(\sqrt x-a)+\Phi(\sqrt x+a)-1,
 \qquad
 H_a(x)=
 \begin{cases}
  0,&0\le x<a^2,\\
  2\Phi(\sqrt{x-a^2})-1,&x\ge a^2.
 \end{cases}
\]
Both functions vanish for $x<0$. At the boundary $x=a^2$,
\[
 F_a(a^2)-H_a(a^2)=\Phi(2a)-\frac12=:b_a.
\]
For $0\le x<a^2$, monotonicity gives
$0\le F_a(x)-H_a(x)\le b_a$.
For $x\ge a^2$, put
\[
 s=\sqrt{x-a^2},\qquad t=\sqrt x,
 \qquad l=t-a,\qquad u=t+a.
\]
Then $0\le l\le s\le u$ and
\[
 0\le u-s\le2a,\qquad 0\le s-l\le a.
\]
Indeed, $t=\sqrt{s^2+a^2}$ implies
$\max(s,a)\le t\le s+a$, which gives these inequalities.
The CDF difference has the decomposition
\[
 F_a(x)-H_a(x)
 =\{\Phi(u)-\Phi(s)\}-\{\Phi(s)-\Phi(l)\}.
\]
Since the standard normal density decreases on $[0,\infty)$,
\[
 0\le\Phi(u)-\Phi(s)\le\Phi(2a)-\frac12=b_a,
 \qquad
 0\le\Phi(s)-\Phi(l)\le\Phi(a)-\frac12\le b_a.
\]
Thus $|F_a(x)-H_a(x)|\le b_a$ for all $x$, with equality at $a^2$.
This proves the exact identity. Taylor expansion of $\Phi(2a)$ at zero
gives the stated asymptotic expression.

For the fixed-level assertion, write $z=z_\alpha$.
On $\alpha\in[\alpha_0,\alpha_1]$, $z$ lies in a compact subset of
$(0,\infty)$, so $a^2<z^2$ uniformly for sufficiently small $|a|$.
The exact powers are then
\[
 \pi_\alpha(a)=2-\Phi(z-a)-\Phi(z+a),\qquad
 \pi_\alpha^{\mathrm{shift}}(a)
       =2\{1-\Phi(\sqrt{z^2-a^2})\}.
\]
The even Taylor expansion in $a$ and the identity
$\phi'(z)=-z\phi(z)$ yield
\[
 \pi_\alpha(a)=\alpha+a^2z\phi(z)+O(a^4).
\]
Also,
$\sqrt{z^2-a^2}=z-a^2/(2z)+O(a^4)$, whence
\[
 \pi_\alpha^{\mathrm{shift}}(a)
       =\alpha+a^2\phi(z)/z+O(a^4).
\]
The derivatives involved are bounded uniformly on the indicated compact
sets, and $z$ is bounded away from zero. Both remainder bounds are
therefore uniform over $\alpha\in[\alpha_0,\alpha_1]$.
Subtracting proves \eqref{eq:rankone-fixed-level}.
The leading coefficient is positive for $z>1$ and negative for $z<1$.
At the exceptional level $\alpha_*=2\{1-\Phi(1)\}$, one further even
Taylor term gives
\begin{equation}
 \pi_{\alpha_*}(a)-\pi_{\alpha_*}^{\mathrm{shift}}(a)
 =-\frac23\phi(1)a^4+O(a^6).
 \label{eq:rankone-exceptional-level}
\end{equation}
For completeness, the fourth-order terms at a general fixed $z>0$ are
$a^4\phi(z)(z^3-3z)/12$ in $\pi_\alpha(a)$ and
$a^4\phi(z)(z^{-3}+z^{-1})/4$ in
$\pi_\alpha^{\mathrm{shift}}(a)$; their difference at $z=1$ is
$-2\phi(1)a^4/3$.
\end{proof}

\subsection{Proofs for local-power estimation}
\label{app:power-estimation-proofs}

\begin{proof}[Proof of Proposition~\ref{prop:local-power-estimation}]
Put $\bm d=\bm\theta_1-\bm\theta_0$. Since both signals belong to
$\operatorname{range}(\bm\Omega)$, the two joint Gaussian laws have
the same affine support. The squared Mahalanobis distance between their
joint mean vectors is
\[
 \sum_{i=1}^2\sum_{j=1}^m
 \left(\frac{\bm d}{\sqrt{2m}}\right)^\top
 \bm\Omega^+
 \left(\frac{\bm d}{\sqrt{2m}}\right)=h^2.
\]
The signs of the two group shifts disappear after squaring. Equivalently,
an orthogonal transformation of the observations isolates the normalized
contrast
$\sqrt{m/2}(\bar{\bm Y}_1-\bar{\bm Y}_2)$, which has law
$N(\bm\theta,\bm\Omega)$; the remaining Gaussian coordinates have
a law independent of $\bm\theta$ and are independent of the contrast.
On the common support, whitening and rotating the joint mean difference
reduces the comparison to $N(0,1)$ and $N(h,1)$. Their likelihood ratio
exceeds one exactly when the scalar coordinate exceeds $h/2$, so
\begin{equation}
 d_{\mathrm{TV}}(P_{\bm\theta_0,m},P_{\bm\theta_1,m})
   =2\Phi(h/2)-1.
 \label{eq:local-experiment-total-variation}
\end{equation}
Here $d_{\mathrm{TV}}(P,Q)=\sup_B|P(B)-Q(B)|$.
Independent auxiliary randomization leaves this distance unchanged.

Relabel the two signals if necessary so that their powers are
$b_0<b_1$, and let
$\varphi=\boldsymbol1\{\widehat\beta\geq(b_0+b_1)/2\}$.
On $\{\varphi=1\}$ the estimation error at the first signal is at
least $\delta/2$, and on $\{\varphi=0\}$ the error at the second
signal is greater than $\delta/2$. For every test $\varphi$,
\[
 P_{\bm\theta_0,m}(\varphi=1)
 +P_{\bm\theta_1,m}(\varphi=0)
 \geq1-d_{\mathrm{TV}}(P_{\bm\theta_0,m},P_{\bm\theta_1,m})
 =2\Phi(-h/2).
\]
Taking half this sum proves
\eqref{eq:power-estimation-probability-bound}.

For the absolute-error bound, let $p_0,p_1$ be densities of the joint
laws, including any auxiliary randomization, with respect to a common
dominating measure. Pointwise,
\[
 p_0|\widehat\beta-b_0|+p_1|\widehat\beta-b_1|
 \geq\min(p_0,p_1)
      \{|\widehat\beta-b_0|+|\widehat\beta-b_1|\}
 \geq\delta\min(p_0,p_1).
\]
Integration and
$\int\min(p_0,p_1)=1-d_{\mathrm{TV}}(P_{\bm\theta_0,m},
P_{\bm\theta_1,m})$ prove
\eqref{eq:power-estimation-absolute-bound}. The sequence assertion
follows because the lower bound on the error probability stays positive
at an error tolerance bounded away from zero.
\end{proof}

\begin{proof}[Proof of Corollary~\ref{cor:persistent-power-estimation}]
Suppress the index $m$. In an eigenbasis containing $\bm e_m$, write
\[
 T_0^*=\lambda Z^2+R,\qquad
 T_a^*=\lambda(Z+a)^2+R,
\]
where $Z\sim N(0,1)$ is independent of $R$, and
$\Var(R)=2(A-\lambda^2)$. Let $q=q_\alpha^{\mathrm{NR}}$ and
$K=(q-R)/\lambda$. Extend the CDF $F_1$ of $\chi_1^2$ by zero on
$(-\infty,0]$. Continuity of the null distribution gives
\begin{equation}
 \E F_1(K)=1-\alpha.
 \label{eq:persistent-conditional-null-probability}
\end{equation}
Chebyshev's inequality yields
$q\leq\tr(\bm\Omega)+\sqrt{2A/\alpha}$, and hence
\[
 \E K\leq M_0:=1+\frac{\sqrt2}{c\sqrt\alpha},
 \qquad
 \Var(K)=\frac{2(A-\lambda^2)}{\lambda^2}\leq\frac2{c^2}.
\]
Choose
\[
 \tau=\left\{\Phi^{-1}\!\left(\frac12+
                         \frac{1-\alpha}{8}\right)\right\}^2,
 \qquad
 M=M_0+\sqrt{\frac8{c^2(1-\alpha)}}.
\]
Then $0<\tau<M$, $F_1(\tau)=(1-\alpha)/4$ and
$\Pr(K>M)\leq(1-\alpha)/4$. Splitting
\eqref{eq:persistent-conditional-null-probability} over
$\{K\leq\tau\}$, $\{\tau<K\leq M\}$ and $\{K>M\}$ therefore gives
\[
 \Pr(\tau<K\leq M)\geq\frac{1-\alpha}{2}.
\]
For $t>0$, put
$g_a(t)=1-\Phi(t-a)+\Phi(-t-a)$.
The derivative with respect to $a>0$ is
$\phi(t-a)-\phi(t+a)>0$. Thus the continuous function
$g_{a_1}(t)-g_{a_0}(t)$ has a strictly positive minimum on
$[\sqrt\tau,\sqrt M]$. Conditional on $R$, the power difference is
zero when $K\leq0$, and is
$g_{a_1}(\sqrt K)-g_{a_0}(\sqrt K)>0$ when $K>0$. Consequently,
\[
 \beta_\alpha(\bm\theta_1)-\beta_\alpha(\bm\theta_0)
 \geq b:=\frac{1-\alpha}{2}
 \min_{\sqrt\tau\leq t\leq\sqrt M}
          \{g_{a_1}(t)-g_{a_0}(t)\}>0.
\]
The Mahalanobis distance between these signals is
$h=a_1-a_0$. Proposition~\ref{prop:local-power-estimation} gives
\eqref{eq:persistent-power-estimation-bound}, since an error of at
least half the actual power gap is also at least $b/2$.
Finally,
$\kappa_\ell=a_\ell^2\lambda/\sqrt{2A}$ and
$c\leq\lambda/\sqrt A\leq1$, proving the stated local-scale bounds.
\end{proof}

For a concrete growing-dimensional submodel, take
$\bm\Omega_m=\lambda_m\bm e_1\bm e_1^\top$, $\lambda_m>0$, in any
dimension $p_m\to\infty$, and
$\bm\theta_{a,m}=a\sqrt{\lambda_m}\bm e_1$. With
$z_\alpha=\Phi^{-1}(1-\alpha/2)$,
\begin{equation}
 \beta_\alpha(a)
   =1-\Phi(z_\alpha-a)+\Phi(-z_\alpha-a),
 \qquad \kappa_a=\frac{a^2}{\sqrt2}.
 \label{eq:rank-one-local-power}
\end{equation}
This power is strictly increasing for $a>0$. Any fixed
$0<a_0<a_1$ therefore gives a positive power gap and
$h=a_1-a_0$, independently of $m$ and $p_m$. Both signals are on the
local scale. Increasing the number of observations does not remove the
estimation bound because the corresponding population mean differences,
standardized by the marginal standard deviation in that eigendirection,
shrink at rate $m^{-1/2}$.

\begin{proof}[Details for the rank-one submodel]
The statistic divided by $\lambda_m$ has law $(Z+a)^2$, and its null
critical value divided by $\lambda_m$ is $z_\alpha^2$.
This gives \eqref{eq:rank-one-local-power}. Its derivative for $a>0$ is
$\phi(z_\alpha-a)-\phi(z_\alpha+a)>0$, since $z_\alpha>0$ and
$(z_\alpha+a)^2>(z_\alpha-a)^2$. The formulas for $\kappa_a$ and
$h$ follow from $\tr(\bm\Omega_m^2)=\lambda_m^2$ and
$\bm\Omega_m^+=\lambda_m^{-1}\bm e_1\bm e_1^\top$.
\end{proof}

\begin{proof}[A counterexample with known signal norm]
Fix $\alpha\in(0,1)$ and abbreviate $z=z_\alpha>0$.
For $\varepsilon\in(0,1)$, let $q_\varepsilon$ be the null
$(1-\alpha)$ quantile of $Z_1^2+\varepsilon Z_2^2$, where
$Z_1,Z_2$ are independent standard normal variables. As
$\varepsilon\downarrow0$, this variable converges in distribution
to $Z_1^2$. The latter has a continuous, strictly increasing CDF on
$(0,\infty)$, so $q_\varepsilon\to z^2$.

Under the equal-norm signals $z\bm e_1$ and $z\bm e_2$, respectively,
the statistic has the distributions
\[
 (Z_1+z)^2+\varepsilon Z_2^2,
 \qquad
 Z_1^2+(\sqrt\varepsilon Z_2+z)^2.
\]
Neither limiting variable has an atom at $z^2$. Their powers therefore
converge, respectively, to
\[
 \Pr\{(Z_1+z)^2>z^2\}
       =\frac12+\Phi(-2z),
 \qquad
 \Pr\{Z_1^2+z^2>z^2\}=1.
\]
Let $d_\alpha=1/2-\Phi(-2z)>0$. There exists a fixed
$\varepsilon_\alpha\in(0,1)$ for which the power gap exceeds
$d_\alpha/2$. Fix this value for all $m$ and embed the covariance in
any sequence of dimensions $p_m\geq2$. The signals lie in its range,
have the same local strength
$\kappa=z^2/\sqrt{2(1+\varepsilon_\alpha^2)}$, and satisfy
\[
 h=z\sqrt{1+\varepsilon_\alpha^{-1}}<\infty.
\]
Proposition~\ref{prop:local-power-estimation} gives a positive,
$m$-independent lower bound. The construction proves existence of such
a fixed anisotropy; it does not assert a nonzero power gap for every
pair of equal-norm directions or every covariance at the chosen level.

The phenomenon also occurs with full-rank covariance matrices.
If $p_m\to\infty$, append $p_m-2$ eigenvalues equal to $p_m^{-1}$
instead of zero. Their total quadratic contribution has mean
$(p_m-2)/p_m\to1$ and variance
$2(p_m-2)/p_m^2\to0$. Null quantiles and alternative statistics
acquire the same limiting deterministic shift. The power gap has the
same positive limit, and $h$ is unchanged. If an exactly constant
local strength is desired, multiply both signals by
$\{A_m/(1+\varepsilon_\alpha^2)\}^{1/4}\to1$, where
$A_m=\tr(\bm\Omega_m^2)$; this leaves the limiting conclusions unchanged.
\end{proof}

\begin{proof}[Existence of the specified-effect power extrema]
If $\bm\theta_k\to\bm\theta$, then
$\|\bm G+\bm\theta_k\|^2\to\|\bm G+\bm\theta\|^2$
almost surely under the same Gaussian coupling. At least one covariance
eigenvalue is positive. Conditional on the other Gaussian coordinates,
the squared norm is a positive multiple of the square of a continuously
distributed nondegenerate Gaussian variable plus a constant. Thus
$\Pr\{\|\bm G+\bm\theta\|^2=q_\alpha^{\mathrm{NR}}\}=0$.
Bounded convergence shows that $\beta_\alpha$ is continuous.
The asserted extrema follow from compactness.
\end{proof}

\subsection{Assessing a shifted reference over an effect range}
\label{app:diagnostic}

Let $A=\tr(\bm\Omega^2)>0$, $s_0=\sqrt{2A}$, and let $H$ denote
the distribution function of
\[
 Q=\frac{T_0^*-\tr(\bm\Omega)}{s_0}.
\]
For a distribution function $F$, define its concentration function by
\[
 \omega_F(u)=\sup_{x\in\mathbb R}\{F(x+u)-F(x)\},\qquad u>0,
\]
and put $\bar\Phi=1-\Phi$. A term $\bar\Phi(u/\sqrt{2d})$ is
interpreted as zero when $d=0$ and $u>0$.

\begin{prop}[A concentration bound for shifted-null error]
\label{prop:shift-diagnostic}
For every deterministic signal,
\begin{equation}
 d_{\bm\theta}\le
 \min\left\{\frac12,\inf_{u>0}
  \left[\bar\Phi\!\left(\frac{u}{\sqrt{2\Delta_{\bm\theta}}}\right)
       +\omega_H(u)\right]\right\}.
 \label{eq:diagnostic-oracle}
\end{equation}
Let $\widehat{\bm\Omega}\succeq0$, $\widehat A$, and a random
distribution function $\widehat H$ be computed from the data. Suppose
nonnegative, possibly data-dependent bounds $e_{\mathrm{op}},e_A,e_H$
are available such that the event
\begin{equation}
 \mathcal E=\left\{
  \|\widehat{\bm\Omega}-\bm\Omega\|_{\mathrm{op}}\le e_{\mathrm{op}},\quad
  |\widehat A-A|\le e_A,\quad
  \|\widehat H-H\|_\infty\le e_H
 \right\}
 \label{eq:diagnostic-event}
\end{equation}
has probability at least $1-\gamma$, for $0<\gamma<1$. Write
$\widehat\lambda_1=\lambda_{\max}(\widehat{\bm\Omega})$ and define
\[
 r_U=\begin{cases}
  \min\{1,(\widehat\lambda_1+e_{\mathrm{op}})
                 /\sqrt{\widehat A-e_A}\},&\widehat A>e_A,\\
  1,&\widehat A\le e_A,
 \end{cases}
\]
\[
 r_L=\begin{cases}
  \min\{1,(\widehat\lambda_1-e_{\mathrm{op}})_+
                 /\sqrt{\widehat A+e_A}\},&\widehat A+e_A>0,\\
  0,&\widehat A+e_A\le0.
 \end{cases}
\]
For $K>0$, set
\[
 \mathcal D(K)=\sup_{\|\bm\theta\|^2\le K\sqrt{2A}}d_{\bm\theta},
 \qquad d_U=\sqrt2Kr_U,\qquad d_L=\sqrt2Kr_L,
\]
and
\begin{equation}
 B_K=\min\left\{\frac12,\inf_{u>0}
   \left[\bar\Phi\!\left(\frac{u}{\sqrt{2d_U}}\right)
          +\omega_{\widehat H}(u)+2e_H\right]\right\}.
 \label{eq:diagnostic-estimated}
\end{equation}
Then, with probability at least $1-\gamma$,
\begin{equation}
 \frac{d_L^2}{60+120d_L+24d_L^2}
       \le \mathcal D(K)\le B_K.
 \label{eq:diagnostic-confidence}
\end{equation}
The probability statement concerns the estimation of a population
approximation error. It does not assert that a test using an estimated
critical value has level $\gamma$ or that its size error is $B_K$.
\end{prop}

\begin{proof}
Let $L=2\bm\theta^\top\bm G/s_0$. Then
$L\sim N(0,2\Delta_{\bm\theta})$ and
$d_{\bm\theta}=d_{\mathrm K}(Q+L,Q)$.
For any $u>0$ and $x\in\mathbb R$, the inclusions
\[
 \{Q+L\le x\}\subseteq\{Q\le x+u\}\cup\{L<-u\},
 \qquad
 \{Q\le x-u\}\subseteq\{Q+L\le x\}\cup\{L>u\}
\]
give
\[
 H(x-u)-\bar\Phi\!\left(\frac{u}{\sqrt{2\Delta_{\bm\theta}}}\right)
 \le \Pr(Q+L\le x)
 \le H(x+u)+\bar\Phi\!\left(\frac{u}{\sqrt{2\Delta_{\bm\theta}}}\right).
\]
If $\Delta_{\bm\theta}=0$, then $L=0$ almost surely and the
conclusion also holds. Taking the supremum over $x$, followed by the
infimum over $u$, proves the concentration bound. The two events
$\{Q+L\le x<Q\}$ and $\{Q\le x<Q+L\}$ imply $L<0$ and
$L>0$, respectively. Each sign has probability $1/2$ when
$\Delta_{\bm\theta}>0$, while $L=0$ if
$\Delta_{\bm\theta}=0$. Hence $d_{\bm\theta}\le1/2$, proving
the cap in \eqref{eq:diagnostic-oracle}. No independence
between $Q$ and $L$ is used.

On $\mathcal E$, Weyl's inequality and $0<r\le1$ give
$r_L\le r\le r_U$, where $r=\lambda_{\max}(\bm\Omega)/\sqrt A$.
For every signal in the displayed ball,
\[
 \Delta_{\bm\theta}
 \le\frac{\lambda_{\max}(\bm\Omega)\|\bm\theta\|^2}{A}
 \le\sqrt2Kr\le d_U.
\]
Also, $\omega_H(u)\le\omega_{\widehat H}(u)+2e_H$ simultaneously
for all $u>0$. Since $d\mapsto\bar\Phi(u/\sqrt{2d})$ is
nondecreasing, \eqref{eq:diagnostic-oracle} proves the upper bound
on $\mathcal E$.
For the lower bound, choose a leading unit eigenvector $\bm v_1$ and
the signal $\bm\theta=(K\sqrt{2A})^{1/2}\bm v_1$. Its
$\Delta_{\bm\theta}$ equals $\sqrt2Kr$. The lower bound in
Theorem~\ref{thm:shift-boundary} and the monotonicity on $[0,\infty)$
of $d\mapsto d^2/(60+120d+24d^2)$ give the stated lower bound.
Both inequalities therefore hold on the same event $\mathcal E$.
\end{proof}

The infimum in \eqref{eq:diagnostic-estimated} may be replaced by a
minimum over any nonempty finite set of positive values; this produces
a conservative upper bound and permits direct numerical evaluation.
The set may depend on the observed estimates because the inequalities
hold simultaneously in $u$. Computing the concentration function on
a finite grid of $x$ gives only a lower approximation to its supremum
and cannot by itself certify the displayed upper bound. For an empirical
distribution function it can instead be evaluated by counting the
largest number of sorted observations in a half-open interval
$(x,x+u]$. Counting closed intervals of length $u$ is conservative
if ties or pairs separated by exactly $u$ occur.

The assumptions in \eqref{eq:diagnostic-event} are substantive.
For instance, if $\bm\Omega$ is known, $r_L=r_U=r$ and
$\widehat H$ may be the empirical distribution function of $B$
independent null-reference draws. The Dvoretzky--Kiefer--Wolfowitz
inequality with the sharp constant \citep{Massart1990} then permits
\[
 e_H=\sqrt{\frac{\log(2/\gamma)}{2B}}.
\]
This gives a computable statement for a specified population reference.
For unknown covariance, the conditional consistency of an estimated
reference does not supply the finite-sample event
\eqref{eq:diagnostic-event}; a valid simultaneous error bound must be
proved or supplied separately. Here $H$ uses the true centering
$\tr(\bm\Omega)$ and scale $\sqrt{2A}$. If $\widehat H$ is constructed
using estimated centering and scale, their errors must also be included
in $e_H$; a bound on the unstandardized reference alone cannot silently
be substituted for the last condition in \eqref{eq:diagnostic-event}.
The present abstract proposition alone therefore does not supply all
nuisance guarantees. Appendix~\ref{app:finite-unknown-covariance-diagnostic}
supplies an explicit alternative for Gaussian observations, working on
the raw squared-norm scale so that no estimated centering or scale is
substituted into the present standardized event.
Under the operator-norm and trace consistency conditions already used
for spectral calibration, the scalar ratio
$\widehat\lambda_1/\sqrt{\widehat A}$ consistently estimates $r$,
but a finite-sample tolerance assessment requires more than that fact.

Finally, if
$e_K=\sup_{\|\bm\theta\|^2\le K\sqrt{2A}}
d_{\mathrm K}(T_{\bm\theta},T_{\bm\theta}^*)$ is an available
uniform original-data comparison bound, then the same event yields
\[
 \sup_{\|\bm\theta\|^2\le K\sqrt{2A}}
 d_{\mathrm K}(T_{\bm\theta},T_0^*+\|\bm\theta\|^2)
 \le \min\{1,e_K+B_K\}.
\]
This last statement still uses the exact population null reference;
estimated critical values require their own error control.

\subsection{A finite-sample bound from a hollow Gram matrix}
\label{app:gaussian-hollow-ratio}

Let $\bm X_1,\ldots,\bm X_m$ be independent $N(\bm0,\bm\Omega)$
vectors, with $m\ge2$ and $A=\tr(\bm\Omega^2)>0$. In the Gaussian
two-sample setting, such vectors are available by adding one independent
weighted residual from each sample:
$\bm X_r=\bm V_{1r}+\bm V_{2r}$,
$r=1,\ldots,m=\min(m_1,m_2)$, where
$\Cov(\bm V_{ir})=\bm\Omega_i$ and
$\bm\Omega_1+\bm\Omega_2=\bm\Omega$.
Using the smaller residual count is explicit; extra residuals are not
needed for the following guarantee.

Define the $m\times m$ hollow Gram matrix and trace estimator by
\begin{equation}
 H_{rs}=\begin{cases}
  m^{-1}\bm X_r^\top\bm X_s,&r\ne s,\\
  0,&r=s,
 \end{cases}
 \qquad
 \widehat A=\frac{1}{m(m-1)}
       \sum_{r\ne s}(\bm X_r^\top\bm X_s)^2.
 \label{eq:gaussian-ratio-hollow-estimates}
\end{equation}
The diagonal deletion removes the trace contribution that inflates
sample-Gram eigenvalues when the dimension is large. The bound below
uses the largest eigenvalue of $\bm H$, although $\bm H$ is generally
indefinite.

\begin{prop}[A Gaussian hollow-Gram confidence bound]
\label{prop:gaussian-hollow-ratio}
Choose deterministic tail probabilities
$\delta_A,\delta_N,\delta_V,\delta_R\in(0,1)$ with sum less than one.
Put
\[
 t_A=\sqrt{\frac{8(1-\delta_A)}{(m-1)\delta_A}},
 \qquad A_L=\frac{\widehat A}{1+t_A},
 \qquad q_m=F_{\chi_m^2}^{-1}(\delta_N),
\]
and define
\begin{align*}
 b_m^{(1)}&=
 F_{\operatorname{Beta}(1/2,(m-1)/2)}^{-1}
                    (1-\delta_V/m),\\
 b_m^{(2)}&=\frac{3}{m+2}
 +\sqrt{\frac{24(m-1)(1-\delta_V)}
  {(m+2)^2(m+4)(m+6)\delta_V}},\qquad
 b_m=\min\{b_m^{(1)},b_m^{(2)}\},\\
 c_m&=\frac{q_m}{m}(1-b_m)>0,
 \qquad t_R=\frac1m\sqrt{\frac{2(1-\delta_R)}{\delta_R}}.
\end{align*}
On $A_L>0$, put $h=\lambda_{\max}(\bm H)/\sqrt{A_L}$ and
\begin{equation}
 r_U=\begin{cases}
  1,&h\ge c_m,\\[2pt]
  \displaystyle
  \frac{c_m h+t_R\sqrt{c_m^2+t_R^2-h^2}}
       {c_m^2+t_R^2},&0\le h<c_m.
 \end{cases}
 \label{eq:gaussian-ratio-upper}
\end{equation}
Assign $r_U=1$ if $A_L=0$. Then, for every dimension and every
nonzero positive semidefinite $\bm\Omega$,
\begin{equation}
 \Pr\left\{
 A_L\le A,\quad
 \frac{\lambda_{\max}(\bm\Omega)}{\sqrt A}\le r_U
 \right\}
 \ge1-\delta_A-\delta_N-\delta_V-\delta_R.
 \label{eq:gaussian-ratio-coverage}
\end{equation}
No independence between the estimators or confidence events is required.
\end{prop}

\begin{proof}
Write $B=\tr(\bm\Omega^4)$. For the order-two kernel
$k(\bm X,\bm Y)=(\bm X^\top\bm Y)^2$, Gaussian conditioning gives
\[
 \E k=A,\qquad \Var(k)=2A^2+6B,
 \qquad
 \Var\{\E(k\mid\bm X)\}=2B.
\]
Counting pairs of kernel terms that share one or two observations yields
\begin{equation}
 \Var(\widehat A)
 =\frac{(8m-4)B+4A^2}{m(m-1)}
 \le\frac{8A^2}{m-1},
 \label{eq:gaussian-trace-variance-ratio}
\end{equation}
where $B\le A^2$. Cantelli's one-sided variance inequality therefore
gives
\[
 \Pr\{\widehat A>A(1+t_A)\}\le\delta_A.
\]
Thus $A_L\le A$ except on an event of probability at most $\delta_A$.
Also, $\widehat A>0$ almost surely because $A>0$ and $m\ge2$;
the stated zero convention covers exceptional or numerical cases.

Let $\lambda=\lambda_{\max}(\bm\Omega)>0$ and fix one corresponding
unit eigenvector $\bm e$. The Gaussian observations admit the decomposition
\[
 \bm X_r=\sqrt\lambda\,z_r\bm e+\bm R_r,
 \qquad z_r\stackrel{\mathrm{i.i.d.}}{\sim}N(0,1),
 \qquad
 \bm R_r\stackrel{\mathrm{i.i.d.}}{\sim}
 N(\bm0,\bm\Omega-\lambda\bm e\bm e^\top).
\]
The sequence $(z_r)$ is independent of all $\bm R_r$, and
$\bm e^\top\bm R_r=0$ almost surely. Set
$\bm z=(z_1,\ldots,z_m)^\top$, $\bm v=\bm z/\|\bm z\|$,
and $S_v=\sum_r v_r^4$. The Gaussian norm and direction are independent,
$\|\bm z\|^2\sim\chi_m^2$, and $\bm v$ is uniform on the unit sphere.
Each $v_r^2$ has distribution
$\operatorname{Beta}(1/2,(m-1)/2)$. Since
$S_v\le\max_r v_r^2$, a union bound gives
$\Pr(S_v>b_m^{(1)})\le\delta_V$.
The spherical moments also give
\[
 \E S_v=\frac3{m+2},\qquad
 \Var(S_v)=\frac{24(m-1)}{(m+2)^2(m+4)(m+6)}.
\]
Indeed,
$\E v_1^8=105/\{m(m+2)(m+4)(m+6)\}$ and
$\E(v_1^4v_2^4)=9/\{m(m+2)(m+4)(m+6)\}$.
Cantelli's inequality gives $\Pr(S_v>b_m^{(2)})\le\delta_V$.
The smaller of these two deterministic constants is itself a valid
bound with the same tail probability; this choice uses no data.
Consequently,
\begin{equation}
 \frac{\|\bm z\|^2}{m}(1-S_v)\ge c_m
 \label{eq:gaussian-leading-coefficient}
\end{equation}
except on an event of probability at most $\delta_N+\delta_V$.
The beta quantile is strictly below one, so $c_m>0$.

The Rayleigh quotient of the observed hollow matrix in direction
$\bm v$ is
\[
 \bm v^\top\bm H\bm v
 =\lambda\frac{\|\bm z\|^2}{m}(1-S_v)+W,
 \qquad
 W=\frac1m\sum_{r\ne s}v_rv_s\bm R_r^\top\bm R_s.
\]
Conditionally on $\bm z$, $W$ has mean zero and
\begin{equation}
 \Var(W\mid\bm z)
 =\frac{2(A-\lambda^2)}{m^2}(1-S_v)
 \le\frac{2(A-\lambda^2)}{m^2}.
 \label{eq:gaussian-hollow-noise-variance}
\end{equation}
To verify this identity, terms with different unordered pairs have zero
covariance, whereas
$\E(\bm R_r^\top\bm R_s)^2
 =\tr\{(\bm\Omega-\lambda\bm e\bm e^\top)^2\}
 =A-\lambda^2$. Conditional Cantelli gives
\[
 \Pr\{W<-t_R\sqrt{A-\lambda^2}\mid\bm z\}\le\delta_R.
\]
If $A=\lambda^2$, the residual covariance vanishes and $W=0$, so
the same assertion holds. Integrating the conditional bound and combining
it with \eqref{eq:gaussian-leading-coefficient} shows that, except on
probability $\delta_N+\delta_V+\delta_R$,
\[
 \lambda_{\max}(\bm H)
 \ge\lambda c_m-t_R\sqrt{A-\lambda^2}.
\]
Put $r=\lambda/\sqrt A\in(0,1]$. On the intersection with $A_L\le A$,
the trace-zero symmetric matrix $\bm H$ has nonnegative largest
eigenvalue, and hence
\[
 c_m r-t_R\sqrt{1-r^2}
 \le\frac{\lambda_{\max}(\bm H)}{\sqrt A}
 \le\frac{\lambda_{\max}(\bm H)}{\sqrt{A_L}}=h.
\]
The left-hand side is strictly increasing on $[0,1]$, with values
$-t_R$ and $c_m$ at its endpoints. If $h\ge c_m$, the universal
bound $r\le1$ applies. Otherwise its unique solution at equality is
the expression in \eqref{eq:gaussian-ratio-upper}, obtained by solving
the quadratic and selecting the root with $c_m r-h\ge0$.
A union bound proves \eqref{eq:gaussian-ratio-coverage}.
\end{proof}

The confidence bound is finite and computable without knowing a population
trace, effective rank or tail constant. It is valid under the stated
Gaussian residual model, including singular covariance matrices and
dimensions exceeding the residual sample size. These properties do not
imply a narrow bound at small $m$: $r_U$ can equal one, and the use of
one-sided variance inequalities may be conservative. A non-Gaussian
extension would require replacement guarantees for the Gaussian
norm--direction and residual-independence steps.

\subsection{A finite-sample diagnostic with unknown Gaussian covariance}
\label{app:finite-unknown-covariance-diagnostic}

Suppose both samples are Gaussian. Form the independent weighted Helmert
residuals $\bm V_{ir}$ in Appendix~\ref{app:implementation-details} and put
\[
 m=\min(n_1-1,n_2-1)\ge2,\qquad
 \bm X_r=\bm V_{1r}+\bm V_{2r},\quad r=1,\ldots,m.
\]
Then $\bm X_1,\ldots,\bm X_m$ are independent $N(\bm0,\bm\Omega)$
vectors, with $\bm\Omega=w_1\bm\Sigma_1+w_2\bm\Sigma_2$ as before.
They are also independent of the normalized sample-mean contrast, although
that independence is not needed in the confidence-bound argument below.
This construction uses matched residuals to obtain identically distributed
reference observations; excess residuals in the larger group are not used
for this diagnostic. The statistic continues to use both full sample means.

Write $A=\tr(\bm\Omega^2)>0$ and define
\begin{equation}
 U_m=\frac{1}{m(m-1)}\sum_{r\ne s}(\bm X_r^\top\bm X_s)^2,
 \qquad L=\lfloor m/2\rfloor.
 \label{eq:diagnostic-trace-U}
\end{equation}
The following lower-tail bound gives a finite upper confidence bound for
$A$ without requiring $m$ to exceed a threshold depending on a variance bound.

\begin{lem}[A Gaussian lower-tail bound for the trace of the squared covariance]
\label{lem:gaussian-trace-square-upper}
For $s>0$, let
\begin{equation}
 M(s)=\E\exp(-sZ_1^2Z_2^2)
 =\E(1+2sZ_1^2)^{-1/2}
 =\frac{\exp\{1/(8s)\}K_0\{1/(8s)\}}{\sqrt{4\pi s}},
 \label{eq:gaussian-product-laplace}
\end{equation}
where $Z_1,Z_2$ are independent standard normal variables and $K_0$ is
the modified Bessel function of the second kind of order zero.
For every $c>0$ and $s>0$,
\begin{equation}
 \Pr(U_m\le cA)\le \exp(sLc)M(s)^L.
 \label{eq:gaussian-U-lower-tail}
\end{equation}
Fix $0<\gamma_A<1$. Choose any deterministic $s>0$ for which
\[
 c=c_{m,\gamma_A}(s)
   :=\frac{L^{-1}\log\gamma_A-\log M(s)}{s}>0,
 \qquad A_U=U_m/c.
\]
Then
\begin{equation}
 \Pr(A\le A_U)\ge1-\gamma_A.
 \label{eq:finite-A-upper}
\end{equation}
Such an $s$ exists for every $m\ge2$ and $0<\gamma_A<1$.
The largest positive $c$ over any deterministic finite grid of $s$ values
may be used without a multiple-comparison penalty.
\end{lem}

\begin{proof}
Let $\lambda_1,\ldots,\lambda_p$ be the covariance eigenvalues and
$w_j=\lambda_j^2/A$. Conditional on $\bm X_1$, the inner product
$\bm X_1^\top\bm X_2$ is centered normal with variance
$\bm X_1^\top\bm\Omega\bm X_1$. Hence
\[
 \E\exp\{-s(\bm X_1^\top\bm X_2)^2/A\}
  =\E\left(1+2s\sum_jw_jZ_j^2\right)^{-1/2}.
\]
The function $x\mapsto(1+2sx)^{-1/2}$ is convex on $[0,\infty)$,
and $\sum_jw_j=1$. Pointwise Jensen's inequality, followed by expectation,
bounds the last expression by
$\sum_jw_j\E(1+2sZ_j^2)^{-1/2}=M(s)$.

Independently draw a uniform random permutation of $1,\ldots,m$ and
match its first $2L$ indices into $L$ disjoint unordered pairs. Let
$V$ be the average of $(\bm X_r^\top\bm X_s)^2$ over these pairs.
Conditional on the observations, symmetry gives $\E(V\mid\bm X)=U_m$.
Convexity of the exponential therefore yields
\[
 \E\exp(-sLU_m/A)
 \le\E\exp(-sLV/A)\le M(s)^L,
\]
because the $L$ pair kernels are independent. The same argument covers
odd $m$, when one uniformly selected observation is left unused by each
matching. Markov's inequality applied to the decreasing exponential
proves \eqref{eq:gaussian-U-lower-tail}. Substituting $c_{m,\gamma_A}(s)$
proves \eqref{eq:finite-A-upper}.

For the Bessel expression, use the Gaussian integral for
$\E(1+2sZ^2)^{-1/2}$ and substitute
$z=\sinh(t)/\sqrt{2s}$ on the positive half-line. The identity
$\sinh^2(t)=\{\cosh(2t)-1\}/2$ gives
\[
 M(s)=\frac{e^{1/(8s)}}{2\sqrt{\pi s}}
        \int_0^\infty e^{-(1/(8s))\cosh(u)}\,du,
\]
which is \eqref{eq:gaussian-product-laplace} by the integral definition
of $K_0$. Dominated convergence gives $M(s)\to0$ as $s\to\infty$,
so $M(s)<\gamma_A^{1/L}$ eventually, proving existence of a positive
$c$. The maximizing grid index is deterministic: it depends only on
$m$, $\gamma_A$ and the specified grid, not on the observations.
Thus it selects one valid lower-tail bound rather than combining
data-selected confidence events.
\end{proof}

Let $F_0$ be the unstandardized Gaussian null CDF,
$F_0(x)=\Pr(\|\bm G\|^2\le x)$ for
$\bm G\sim N(\bm0,\bm\Omega)$, and let
\[
 \widehat F_m(x)=\frac1m\sum_{r=1}^m
          \boldsymbol1\{\|\bm X_r\|^2\le x\}.
\]
For a CDF $F$, use the concentration function
$\omega_F(t)=\sup_x\{F(x+t)-F(x)\}$, $t>0$.

An order-statistic bound can improve the additive DKW bound when an
interval contains few observations. Let $q_{a,b}(v)$ denote the
$v$ quantile of a $\operatorname{Beta}(a,b)$ distribution. For a
fixed $0<\gamma_F<1$, put
\begin{align}
 \varepsilon_m&=\sqrt{\frac{\log(4/\gamma_F)}{2m}},\nonumber\\
 b_{m,k}&=q_{k+1,m-k}\!\left(1-
             \frac{\gamma_F}{2m(m-k+1)}\right),
                      &&0\le k<m,\nonumber\\
 \eta_{m,k}&=\max_{0\le\ell\le k}b_{m,\ell},
                      &&0\le k<m,\qquad \eta_{m,m}=1,\nonumber\\
 C_m(k)&=\min\{1,k/m+2\varepsilon_m,\eta_{m,k}\},
                      &&0\le k\le m.
 \label{eq:finite-interval-band}
\end{align}
These constants depend only on $m$ and the error allocation. Define
$k_m(t)$ to be the largest number of observed squared norms in a closed
interval of length $t$.

\begin{lem}[A simultaneous empirical concentration bound]
\label{lem:finite-interval-band}
With probability at least $1-\gamma_F$,
\begin{equation}
 \omega_{F_0}(t)\le C_m\{k_m(t)\}
 \qquad\hbox{for every }t>0.
 \label{eq:finite-interval-coverage}
\end{equation}
\end{lem}

\begin{proof}
The distribution $F_0$ is continuous because at least one covariance
eigenvalue is positive. For the ordered squared norms
$Q_{(1)}<\cdots<Q_{(m)}$, the variables $U_i=F_0(Q_{(i)})$ are
uniform order statistics. Put $U_0=0$ and $U_{m+1}=1$.
For each $0\le k<m$ and $0\le j\le m-k$, the spacing
$U_{j+k+1}-U_j$ has distribution
$\operatorname{Beta}(k+1,m-k)$. The probability that any of these
spacings exceeds $b_{m,k}$ is at most
\[
 \sum_{k=0}^{m-1}\sum_{j=0}^{m-k}
       \frac{\gamma_F}{2m(m-k+1)}=\gamma_F/2.
\]
On the complementary event, an interval $(a,b]$ containing $k<m$
observations has population probability at most $b_{m,k}$. To see this,
put $j=\#\{i:Q_i\le a\}$ and take $Q_{(0)}=-\infty$ and
$Q_{(m+1)}=+\infty$. Then
$Q_{(j)}\le a<b<Q_{(j+k+1)}$, allowing equality at the infinite
endpoints, and
$F_0(b)-F_0(a)\le U_{j+k+1}-U_j$.
Every interval of length $t$ contains at most $k_m(t)$
observations, so its probability is at most $\eta_{m,k_m(t)}$.
The assertion is trivial if $k_m(t)=m$.

Independently of this argument, the DKW inequality gives
$\|\widehat F_m-F_0\|_\infty\le\varepsilon_m$ with probability
at least $1-\gamma_F/2$. On that event, every interval of length $t$
has probability at most $k_m(t)/m+2\varepsilon_m$.
Intersecting the two events, using a union bound, and taking the smaller
of their bounds proves \eqref{eq:finite-interval-coverage}.
No independence of the two events is needed.
\end{proof}

\begin{prop}[Unknown-covariance tolerance bound]
\label{prop:finite-unknown-covariance-diagnostic}
Let $A_U>0$ and $r_U\in[0,1]$ be data-based bounds satisfying
$\Pr(A\le A_U)\ge1-\gamma_A$ and $\Pr(r\le r_U)\ge1-\gamma_r$.
The constructions in Lemma~\ref{lem:gaussian-trace-square-upper} and
Proposition~\ref{prop:gaussian-hollow-ratio} are one choice; the tightened
constructions in Propositions~\ref{prop:v104-trace-convex-order}
and~\ref{prop:v104-exponential-hollow-ratio} are another. One may also
use $r_U=1$ and $\gamma_r=0$. Fix $0<\gamma_F<1$ with
$\gamma_A+\gamma_r+\gamma_F<1$, and put
\[
 \sigma_U(K)=2\sqrt{\sqrt2K r_U A_U},\quad K>0.
\]
For any nonempty finite set $\mathcal T\subset(0,\infty)$, possibly
chosen from the observed residuals, define
\begin{equation}
 B_{\mathrm{interval}}(K):=\widehat B_m(K)=\min\left\{\frac12,
    \min_{t\in\mathcal T}\left[
       \bar\Phi\{t/\sigma_U(K)\}
       +C_m\{k_m(t)\}\right]\right\}.
 \label{eq:finite-unknown-diagnostic}
\end{equation}
Then, simultaneously for every $K>0$, with probability at least
$1-\gamma_A-\gamma_r-\gamma_F$,
\begin{equation}
 \sup_{\|\bm\theta\|^2\le K\sqrt{2A}}
 d_{\mathrm K}(T_{\bm\theta}^*,T_0^*+\|\bm\theta\|^2)
 \le \widehat B_m(K).
 \label{eq:finite-unknown-coverage}
\end{equation}
All quantities on the right are computed from the observations and
the specified error allocations. No population centering, scale,
covariance eigenvector or effective-rank bound is required.
\end{prop}

\begin{proof}
Use the event in Lemma~\ref{lem:finite-interval-band}, which has
probability at least $1-\gamma_F$. Intersect it with $A\le A_U$ and
$r:=\lambda_{\max}(\bm\Omega)/\sqrt A\le r_U$.
The union bound gives the stated probability; independence among these
three events is not assumed.

For a deterministic signal, write
$T_{\bm\theta}^*=T_0^*+L_{\bm\theta}+\|\bm\theta\|^2$ using
the same $\bm G$, where
$L_{\bm\theta}=2\bm\theta^\top\bm G$ is centered normal with
variance $4\bm\theta^\top\bm\Omega\bm\theta$.
The one-sided bracketing argument of
Proposition~\ref{prop:shift-diagnostic}, before standardization, gives
\[
 d_{\mathrm K}(T_{\bm\theta}^*,T_0^*+\|\bm\theta\|^2)
 \le \bar\Phi\!\left\{\frac{t}
                 {2\sqrt{\bm\theta^\top\bm\Omega\bm\theta}}\right\}
        +\omega_{F_0}(t),\qquad t>0.
\]
The Gaussian-tail term is defined as zero if its variance is zero.
On the joint confidence event, every signal in the specified ball
satisfies
\[
 \bm\theta^\top\bm\Omega\bm\theta
 \le\lambda_{\max}(\bm\Omega)K\sqrt{2A}
 =\sqrt2K r A\le\sqrt2K r_U A_U,
\]
and $\omega_{F_0}(t)\le C_m\{k_m(t)\}$
for all $t>0$. These inequalities hold simultaneously in $K$, $t$ and
the signal. The Gaussian-tail function is nondecreasing in its standard
deviation, so minimization over $\mathcal T$ preserves the bound.
The universal $1/2$ cap follows from the sign-event argument in
Proposition~\ref{prop:shift-diagnostic}. This proves
\eqref{eq:finite-unknown-coverage}.
\end{proof}

The count $k_m(t)$ is computed by sorting the
$m$ observed squared norms and finding the largest number in a closed
interval of length $t$. Closed-window counts are conservative for the
half-open definition. A finite grid in $t$ gives a conservative numerical
bound; a finite grid in the interval location alone does not certify the
required supremum. In floating-point calculations, if $U_m$ or an input
confidence bound is nonpositive or nonfinite, the diagnostic returns the
universal bound $1/2$.

Unlike the conditional formulation in
Proposition~\ref{prop:shift-diagnostic}, the present construction supplies
all nuisance guarantees from an explicit Gaussian experiment. The sample
size $m$ controls the interval-probability and spectral confidence bounds, so the resulting
tolerance statement can still be uninformative at moderate sample sizes.
Using only the additive DKW term would impose an error floor
$2\varepsilon_m$; the order-statistic bound can improve that floor for
small empirical interval masses. Neither construction guarantees a
nontrivial certificate for every covariance or effect range. A reported
bound of $1/2$ does not mean that the actual approximation error equals
$1/2$. The probability in
\eqref{eq:finite-unknown-coverage} concerns an error bound for a population
reference over specified alternatives, not the type-I error of a newly
selected test. Non-Gaussian observations require separate finite-sample
arguments; the present Gaussian proof does not follow from a projection
fourth-moment condition alone.
\subsection{Sharper Gaussian trace and spectral confidence bounds}
\label{app:v104-trace-spectral-bounds}

Let $\bm X_1,\ldots,\bm X_m$ be independent
$N(\bm0,\bm\Omega)$ vectors, where $m\ge2$ and
$A=\tr(\bm\Omega^2)>0$. Put
\[
 U_m=\frac1{m(m-1)}\sum_{i\ne j}(\bm X_i^\top\bm X_j)^2,
 \qquad
 H_{ij}=m^{-1}\bm X_i^\top\bm X_j\ (i\ne j),\quad H_{ii}=0.
\]
The following comparison retains the complete statistic $U_m$.
The random partition is a proof device and does not require the user to
split the observed sample. It uses the averaging principle underlying
classical bounds for complete $U$-statistics \citep{Hoeffding1963},
together with the standard convex-order and positive-part transforms
\citep{ShakedShanthikumar2007}. All details needed here are proved below.

\begin{prop}[An explicit Gaussian trace comparison]
\label{prop:v104-trace-convex-order}
Let $a=\lfloor m/2\rfloor$, $b=m-a$, and let
\[
 Z_m=\frac{Y_aY_b}{ab},\qquad
 Y_a\sim\chi_a^2,\quad Y_b\sim\chi_b^2,
\]
where the factors are independent. For every convex function $\phi$
for which the expectations exist,
\begin{equation}
 \E\phi(U_m/A)\le\E\phi(Z_m).
 \label{eq:v104-trace-convex-order}
\end{equation}
Write $q_m(v)$ for the $v$ quantile of $Z_m$. For $0<\delta<1$, define
\begin{align}
 \ell_m(\delta)&=\frac{\E[Z_m\boldsymbol1\{Z_m\le q_m(\delta)\}]}{\delta},
 \nonumber\\
 u_m(\delta)&=\frac{\E[Z_m\boldsymbol1\{Z_m\ge q_m(1-\delta)\}]}{\delta}.
 \label{eq:v104-trace-tail-means}
\end{align}
Then $0<\ell_m(\delta)<1<u_m(\delta)<\infty$. For deterministic
$\delta_L,\delta_U\in(0,1)$,
\begin{equation}
 A_L=U_m/u_m(\delta_L),\qquad A_U=U_m/\ell_m(\delta_U)
 \quad\Longrightarrow\quad
 \Pr\{A_L\le A\le A_U\}\ge1-\delta_L-\delta_U.
 \label{eq:v104-trace-confidence}
\end{equation}
The separate one-sided guarantees use respectively $\delta_L$ and
$\delta_U$. All constants depend only on $m$ and the specified tail
probabilities.
\end{prop}

\begin{proof}
First fix a partition $I,J$ of the observation indices with $|I|=a$,
$|J|=b$, and set
\[
 \bm S_I=\frac1a\sum_{i\in I}\bm X_i\bm X_i^\top,
 \qquad
 V_{I,J}=\frac1{ab}\sum_{i\in I,j\in J}(\bm X_i^\top\bm X_j)^2
          =\tr(\bm S_I\bm S_J).
\]
Conditional on $\bm S_I$, put $D=\tr(\bm S_I\bm\Omega)$.
This quantity is positive almost surely. If $\mu_j$ are the nonnegative
eigenvalues of $\bm\Omega^{1/2}\bm S_I\bm\Omega^{1/2}$, then
\[
 V_{I,J}/D\ \stackrel{\mathcal L}{=}\
 \sum_j\frac{\mu_j}{D}Q_{bj},
 \qquad Q_{bj}\stackrel{\mathrm{i.i.d.}}{\sim}\chi_b^2/b.
\]
The weights sum to one. Pointwise Jensen, followed by conditional
expectation, therefore gives
\[
 \E\{\phi(V_{I,J}/A)\mid\bm S_I\}
 \le \E_{Q_b}\phi\{(D/A)Q_b\},
\]
where $Q_b\sim\chi_b^2/b$ is independent of $\bm S_I$.
If $\lambda_j$ are the covariance eigenvalues, then
\[
 D/A\ \stackrel{\mathcal L}{=}\
 \sum_j\frac{\lambda_j^2}{A}Q_{aj},
 \qquad Q_{aj}\stackrel{\mathrm{i.i.d.}}{\sim}\chi_a^2/a.
\]
For each fixed $Q_b$, the function $x\mapsto\phi(xQ_b)$ is convex.
A second application of Jensen yields
$\E\phi(V_{I,J}/A)\le\E\phi(Q_aQ_b)=\E\phi(Z_m)$.
These arguments cover singular covariance matrices by omitting zero
weights. Positivity of $D$ follows from $A>0$ and a positive Gaussian
quadratic form.

Now choose $I$ uniformly among all subsets of size $a$, independently
of the observations, and let $J=I^c$. Each unordered pair is separated
with probability $2ab/\{m(m-1)\}$. Consequently
\[
 \E\{V_{I,J}\mid\bm X_1,\ldots,\bm X_m\}=U_m.
\]
Conditional Jensen proves \eqref{eq:v104-trace-convex-order}.

For the lower tail, put $q=q_m(\delta_U)$ and
$\ell=\ell_m(\delta_U)$. The distribution of $Z_m$ is continuous
and strictly positive, so $0<\ell<q$ and
\[
 \E(q-Z_m)_+=\delta_U(q-\ell).
\]
Markov's inequality and \eqref{eq:v104-trace-convex-order}, applied
to the convex function $(q-x)_+$, give
\[
 \Pr\{U_m/A\le\ell\}
 \le\frac{\E(q-U_m/A)_+}{q-\ell}
 \le\delta_U.
\]
Likewise, with $q=q_m(1-\delta_L)$ and $u=u_m(\delta_L)$,
\[
 \E(Z_m-q)_+=\delta_L(u-q),\qquad
 \Pr\{U_m/A\ge u\}\le\delta_L.
\]
These are the claimed one-sided bounds; a union bound proves
\eqref{eq:v104-trace-confidence}. Since $Z_m$ is positive, nondegenerate,
has mean one and has a finite first moment, its lower and upper tail
means satisfy the displayed strict inequalities.
\end{proof}

The constants require only one-dimensional integration and scalar
quantile inversion. For $x>0$, define
\[
 P_{c,d}(x)=\Pr(Y_cY_d\le x)
 =\int_0^\infty f_{\chi_c^2}(y)F_{\chi_d^2}(x/y)\,dy.
\]
Then $P_{a,b}\{abq_m(v)\}=v$, and the size-bias identity
$y f_{\chi_c^2}(y)=c f_{\chi_{c+2}^2}(y)$ gives
\begin{equation}
 \ell_m(\delta)=\frac{P_{a+2,b+2}\{abq_m(\delta)\}}{\delta},
 \qquad
 u_m(\delta)=\frac{1-P_{a+2,b+2}\{abq_m(1-\delta)\}}{\delta}.
 \label{eq:v104-trace-explicit-integrals}
\end{equation}
The upper-tail integral may be evaluated directly with the chi-square
survival function to avoid cancellation. The quantiles and constants
are fixed before seeing the data. Selecting tail budgets or combining
separately valid random bounds after inspecting the data is not justified
by this proposition.

The next elementary signed form of Gaussian quadratic concentration
makes explicit the classical moment-generating-function argument of
\citet[Lemma~1]{LaurentMassart2000}; the original lemma is stated for
nonnegative weights. We prove the signed version used here.

\begin{lem}[A lower tail for a signed Gaussian quadratic form]
\label{lem:v104-signed-gaussian-lower}
Let $\bm B$ be a deterministic real symmetric matrix,
$\bm g\sim N(\bm0,\bm I)$, $v=\|\bm B\|_F$ and
$b=\|\bm B_-\|_{\mathrm{op}}$, where $\bm B_-$ is the negative
part in the spectral decomposition. For every $x>0$,
\begin{equation}
 \Pr\{\bm g^\top\bm B\bm g-\tr\bm B
       <-2v\sqrt x-2bx\}\le e^{-x}.
 \label{eq:v104-signed-gaussian-lower}
\end{equation}
\end{lem}

\begin{proof}
Diagonalize $\bm B$ and write the centered quadratic form as
$W=\sum_j d_j(g_j^2-1)$. For $0<t<(2b)^{-1}$, with no restriction
when $b=0$,
\[
 \log\E e^{-tW}
 =\sum_j\{td_j-\tfrac12\log(1+2td_j)\}
 \le\frac{t^2\sum_jd_j^2}{1-2bt}.
\]
For $d_j\ge0$, the summand is at most $t^2d_j^2$ by
$\log(1+z)\ge z-z^2/2$ for $z\ge0$. For $d_j<0$, its power
series is bounded by $t^2d_j^2/(1-2t|d_j|)$, which is no larger
than the displayed bound. Exponential Markov with
$t=\sqrt x/(v+2b\sqrt x)$ proves
\eqref{eq:v104-signed-gaussian-lower}. If $v=0$, the quadratic
form is identically zero and the conclusion is immediate.
\end{proof}

\begin{prop}[An exponential hollow-Gram spectral bound]
\label{prop:v104-exponential-hollow-ratio}
Choose deterministic $\delta_L,\delta_N,\delta_V,\delta_R\in(0,1)$,
with sum less than one, and use $A_L=U_m/u_m(\delta_L)$ from
Proposition~\ref{prop:v104-trace-convex-order}. Put
\begin{align*}
 q&=F_{\chi_m^2}^{-1}(\delta_N),\qquad x=\log(1/\delta_R),\\
 d_0&=F_{\operatorname{Beta}(1/2,(m-1)/2)}^{-1}(1-\delta_V/m),
       \qquad b_0=d_0,\\
 d_1&=F_{\operatorname{Beta}(1/2,(m-1)/2)}^{-1}(1-\delta_V/(2m)),\\
 b_1&=\min\left\{d_1,\frac3{m+2}
 +\sqrt{\frac{24(m-1)(1-\delta_V/2)}
 {(m+2)^2(m+4)(m+6)(\delta_V/2)}}\right\},\\
 c_j&=\frac qm(1-b_j),\qquad \beta_j=\frac{2d_jx}{m},\qquad j=0,1.
\end{align*}
Select $j=1$ if $c_1-c_0\ge\beta_1-\beta_0$, and $j=0$ otherwise.
This choice is deterministic. Put $c=c_j$, $d=d_j$, $b_*=b_j$,
$\beta=\beta_j$, and
\[
 \alpha=\frac{2\sqrt x}{m},\qquad
 h_* =\frac{c-\alpha-\beta}{\sqrt2}.
\]
For $a>0$, $b>0$ and $0\le h<a$, define
\[
 R(a,b,h)=\frac{ah+b\sqrt{a^2+b^2-h^2}}{a^2+b^2}.
\]
On $A_L>0$, let $h=\lambda_{\max}(\bm H)/\sqrt{A_L}\ge0$ and
\begin{equation}
 r_U=\begin{cases}
 1,& h\ge c,\\
 R(c-\beta,\alpha,h),&0\le h<c\ \hbox{and}\ h\le h_*,\\
 R(c,\alpha+\beta,h),&0\le h<c\ \hbox{and}\ h>h_*.
 \end{cases}
 \label{eq:v104-hollow-piecewise-root}
\end{equation}
If $A_L=0$, set $r_U=1$. Then
\begin{equation}
 \Pr\{A_L\le A,\ \lambda_{\max}(\bm\Omega)/\sqrt A\le r_U\}
 \ge1-\delta_L-\delta_N-\delta_V-\delta_R.
 \label{eq:v104-hollow-coverage}
\end{equation}
In the second branch $c-\beta>0$ and $h<c-\beta$, so the stated
inverse is well defined. The construction allows arbitrary dimension
and singular covariance.
\end{prop}

\begin{proof}
Write $\lambda=\lambda_{\max}(\bm\Omega)$ and fix a corresponding
unit eigenvector $\bm e$. As in the Gaussian hollow-Gram decomposition,
\[
 \bm X_i=\sqrt\lambda\,z_i\bm e+\bm R_i,
 \qquad z_i\stackrel{\mathrm{i.i.d.}}{\sim}N(0,1),\quad
 \bm R_i\stackrel{\mathrm{i.i.d.}}{\sim}N(\bm0,\bm\Psi),
 \qquad \bm\Psi=\bm\Omega-\lambda\bm e\bm e^\top,
\]
where the $z_i$ and all residuals are independent. Put
$\bm v=\bm z/\|\bm z\|$, $S_v=\sum_i v_i^4$, and
$M_v=\max_i v_i^2$. The chi-square norm and spherical-coordinate
marginals give, with failure probability at most $\delta_N+\delta_V$,
\begin{equation}
 \|\bm z\|^2\ge q,\qquad M_v\le d,
 \qquad S_v\le b_*.
 \label{eq:v104-norm-direction-event}
\end{equation}
Indeed, each $v_i^2$ is $\operatorname{Beta}(1/2,(m-1)/2)$.
For configuration $j=0$, the union bound gives $M_v\le d_0$ with
failure probability at most $\delta_V$, and $S_v\le M_v$.
For configuration $j=1$, the same union bound gives $M_v\le d_1$
at cost $\delta_V/2$. The sphere moments give
\[
 \E S_v=\frac3{m+2},\qquad
 \Var(S_v)=\frac{24(m-1)}{(m+2)^2(m+4)(m+6)}.
\]
These follow by expanding $S_v^2$ and using
$\E v_i^8=105/\{m(m+2)(m+4)(m+6)\}$ and
$\E(v_i^4v_k^4)=9/\{m(m+2)(m+4)(m+6)\}$ for $i\ne k$.
Cantelli's inequality supplies the stated upper bound on $S_v$ at
cost $\delta_V/2$; intersecting it with $M_v\le d_1$ gives $S_v\le b_1$.
Thus either deterministic configuration costs at most $\delta_V$.
Moreover, $d_1\ge d_0$, and the stated selection condition makes
$c_1r-\beta_1\min\{r,\sqrt{1-r^2}\}$ no smaller than its
$j=0$ counterpart for all $0\le r\le1$. The deterministic selection
therefore improves or retains the scalar constraint without an
additional confidence cost.

The leading-direction Rayleigh quotient is
\[
 \bm v^\top\bm H\bm v
 =\lambda\frac{\|\bm z\|^2}{m}(1-S_v)+W,\qquad
 W=\frac1m\sum_{i\ne j}v_iv_j\bm R_i^\top\bm R_j.
\]
Let $\bm C_v=\bm v\bm v^\top-\operatorname{diag}(v_i^2)$.
Conditionally on $\bm z$, $W$ is a Gaussian quadratic form with
coefficient matrix $\bm\Psi\otimes\bm C_v/m$, which has trace zero.
Moreover,
\[
 \|\bm C_v\|_F^2=1-S_v\le1,
 \qquad
 \| (\bm C_v)_-\|_{\mathrm{op}}\le M_v,
\]
because $-\bm C_v\preceq\operatorname{diag}(v_i^2)$.
Also
\[
 \tr(\bm\Psi^2)=A-\lambda^2,
 \qquad
 \|\bm\Psi\|_{\mathrm{op}}
 \le\min\{\lambda,\sqrt{A-\lambda^2}\}.
\]
Lemma~\ref{lem:v104-signed-gaussian-lower}, applied conditionally, therefore
gives with conditional failure probability at most $\delta_R$,
\[
 W\ge-\frac{2\sqrt x}{m}\sqrt{A-\lambda^2}
       -\frac{2M_vx}{m}\min\{\lambda,\sqrt{A-\lambda^2}\}.
\]
This is also valid when $\bm\Psi=\bm0$. Integrating the conditional
bound and intersecting with \eqref{eq:v104-norm-direction-event} gives,
for $r=\lambda/\sqrt A$,
\begin{equation}
 \frac{\lambda_{\max}(\bm H)}{\sqrt A}
 \ge g(r):=cr-\alpha\sqrt{1-r^2}
             -\beta\min\{r,\sqrt{1-r^2}\}
 \label{eq:v104-hollow-scalar-constraint}
\end{equation}
except on probability $\delta_N+\delta_V+\delta_R$.
On $A_L\le A$, the trace-zero matrix $\bm H$ has a nonnegative
largest eigenvalue, so the left side is at most $h$.
This intersection costs at most one additional $\delta_L$.

For $r\le1/\sqrt2$, the function in
\eqref{eq:v104-hollow-scalar-constraint} is
$(c-\beta)r-\alpha\sqrt{1-r^2}$; for $r\ge1/\sqrt2$ it is
$cr-(\alpha+\beta)\sqrt{1-r^2}$. The two pieces agree at
$g(1/\sqrt2)=h_*$. The second piece is strictly increasing and
ends at $g(1)=c$. If $c-\beta>0$, the first piece is also strictly
increasing. If $c-\beta\le0$, the first piece is nonpositive,
so it cannot cross a nonnegative value $h$. Consequently, for every
$0\le h<c$, the inequality $g(r)\le h$ defines an initial interval
whose upper endpoint is obtained from the first piece when $h\le h_*$
and from the second piece otherwise. Solving
$ar-b\sqrt{1-r^2}=h$ for its root with $ar-h\ge0$ yields $R(a,b,h)$.
In the first branch, $h\le h_*=(c-\beta-\alpha)/\sqrt2$ implies
$c-\beta>0$ and $h<c-\beta$. This proves
\eqref{eq:v104-hollow-piecewise-root} and the union-bound coverage.
\end{proof}

The four tail allocations above may sum to $0.05$ when only a spectral
certificate is needed. If an upper trace bound and an empirical interval
band are also used, their error allocations must be added to that sum.
These choices are deterministic. No confidence guarantee is obtained by
selecting the smallest of several separately calibrated random certificates
without a joint argument. The numerical constants and the effects of this
refinement are reported separately from its proof.

\subsection{A linear bound under spectral spread}
\label{app:smooth-shift-certificate}

The following argument uses only finite-dimensional Gaussian variables
and allows singular covariance matrices. Let
$A=\tr(\bm\Omega^2)>0$, $s_0=\sqrt{2A}$, and
\[
 Q=\frac{\|\bm G\|^2-\tr(\bm\Omega)}{s_0},\qquad
 \Delta_{\bm\theta}=
       \frac{\bm\theta^\top\bm\Omega\bm\theta}{A}.
\]

\begin{lem}[An independent centered perturbation]
\label{lem:poisson-centered-perturbation}
There exist independent random variables $Q_0$ and $H_{\bm\theta}$ such
that $Q_0\overset d=Q$ and
\begin{equation}
 \frac{\|\bm G+\bm\theta\|^2
              -\tr(\bm\Omega)-\|\bm\theta\|^2}{s_0}
       \overset d=Q_0+H_{\bm\theta},\qquad
 \E H_{\bm\theta}=0,\qquad
 \E H_{\bm\theta}^2=2\Delta_{\bm\theta}.
 \label{eq:independent-poisson-perturbation}
\end{equation}
If $Q$ has a differentiable density $f_Q$ with
$\|f_Q'\|_\infty<\infty$, then
\begin{equation}
 d_{\bm\theta}\le
            \Delta_{\bm\theta}\|f_Q'\|_\infty.
 \label{eq:density-derivative-shift-bound}
\end{equation}
The same assertion holds if $f_Q$ is Lipschitz, with its Lipschitz constant
in place of $\|f_Q'\|_\infty$.
\end{lem}

\begin{proof}
Let $\lambda_1,\ldots,\lambda_q$ be the positive eigenvalues of
$\bm\Omega$, and let $a_j$ be the corresponding coordinates of
$\bm\theta$. For each $j$, independently of all other variables, take
\[
 Z_j\sim N(0,1),\qquad
 N_j\sim\operatorname{Poisson}\{a_j^2/(2\lambda_j)\},\qquad
 E_{jk}\sim\lambda_j\chi_2^2.
\]
Set $J_j=\sum_{k=1}^{N_j}E_{jk}$, with an empty sum equal to zero.
The usual noncentral chi-squared Poisson-mixture identity gives
\[
 \lambda_j\left(Z_j+\frac{a_j}{\sqrt{\lambda_j}}\right)^2
       \overset d=\lambda_j Z_j^2+J_j,
\]
where the two terms on the right are independent. For completeness,
the moment generating function of the right-hand side, for $t$ near zero,
is
\[
 (1-2t\lambda_j)^{-1/2}
 \exp\left[\frac{a_j^2}{2\lambda_j}
       \{(1-2t\lambda_j)^{-1}-1\}\right]
 = (1-2t\lambda_j)^{-1/2}
       \exp\left\{\frac{a_j^2t}{1-2t\lambda_j}\right\},
\]
which is the required noncentral chi-squared transform. Moreover,
\[
 \E J_j=a_j^2,\qquad
 \operatorname{Var}(J_j)
       =\frac{a_j^2}{2\lambda_j}\E E_{j1}^2
       =4a_j^2\lambda_j.
\]
Signal in the null space contributes exactly its squared norm and is
removed by subtracting $\|\bm\theta\|^2$. Consequently,
\[
 Q_0=\frac{\sum_{j=1}^q\lambda_j(Z_j^2-1)}{s_0},\qquad
 H_{\bm\theta}=\frac{\sum_{j=1}^q(J_j-a_j^2)}{s_0}
\]
are independent and satisfy \eqref{eq:independent-poisson-perturbation}.
In particular, the assertion also covers a signal entirely in the null
space, when $H_{\bm\theta}=0$ and $d_{\bm\theta}=0$.

Let $F_Q$ denote the CDF of $Q$. A Lipschitz bound $M$ for $f_Q$ yields,
for every real $x,h$,
\[
 |F_Q(x-h)-F_Q(x)+h f_Q(x)|\le Mh^2/2.
\]
Taking expectation with $h=H_{\bm\theta}$, using independence and
$\E H_{\bm\theta}=0$, gives
\[
 |\Pr(Q_0+H_{\bm\theta}\le x)-F_Q(x)|
       \le (M/2)\E H_{\bm\theta}^2
       =M\Delta_{\bm\theta}.
\]
Taking the supremum over $x$ proves the claim. All comparisons use the
same deterministic location and positive scale, which leave the
Kolmogorov distance unchanged.
\end{proof}

\begin{prop}[A computable linear bound]
\label{prop:smooth-shift-bound}
If $r=\lambda_{\max}(\bm\Omega)/\sqrt A<1/2$, then
$Q$ has a continuously differentiable density satisfying
\begin{equation}
 \|f_Q'\|_\infty\le\frac{1}{\pi(1-4r^2)}.
 \label{eq:spectral-density-derivative-bound}
\end{equation}
Consequently, \eqref{eq:smooth-shift-linear-bound} holds for every
$\bm\theta$. More generally, if
$\rho_j=\lambda_j/\sqrt A$ and
\begin{equation}
 M_1(\bm\rho)=\frac1\pi\int_0^\infty
             t\prod_{j=1}^q(1+2\rho_j^2t^2)^{-1/4}\,dt<\infty,
 \label{eq:spectrum-specific-density-derivative}
\end{equation}
then $d_{\bm\theta}\le\min\{1/2,M_1(\bm\rho)
\Delta_{\bm\theta}\}$. The integral is finite whenever there are at
least five positive eigenvalues; its value may become large as some of
those eigenvalues approach zero.
\end{prop}

\begin{proof}
The characteristic function of $Q$ has modulus
\[
 |\varphi_Q(t)|=
       \prod_{j=1}^q(1+2\rho_j^2t^2)^{-1/4},\qquad
       \sum_{j=1}^q\rho_j^2=1,\qquad 0<\rho_j\le r.
\]
Concavity of $a\mapsto\log(1+2at^2)$ and its value zero at $a=0$
imply, for $0\le a\le r^2$,
\[
 \log(1+2at^2)\ge
          \frac{a}{r^2}\log(1+2r^2t^2).
\]
After summing over $a=\rho_j^2$,
\begin{equation}
 |\varphi_Q(t)|\le(1+2r^2t^2)^{-1/(4r^2)}.
 \label{eq:effective-rank-fourier-envelope}
\end{equation}
If $r<1/2$, both $\varphi_Q$ and $t\varphi_Q(t)$ are integrable.
Fourier inversion and differentiation under the integral give
\begin{align*}
 \|f_Q'\|_\infty
 &\le\frac1{2\pi}\int_{\mathbb R}|t|\,|\varphi_Q(t)|\,dt\\
 &\le\frac1\pi\int_0^\infty
       t(1+2r^2t^2)^{-1/(4r^2)}\,dt
   =\frac1{\pi(1-4r^2)}.
\end{align*}
Lemma~\ref{lem:poisson-centered-perturbation} gives the linear bound;
the established universal $1/2$ bound gives its displayed truncation.
Using the exact characteristic-function modulus in the same argument
gives \eqref{eq:spectrum-specific-density-derivative}. With $q\ge5$,
the integrand is $O(t^{1-q/2})$ at infinity and is $O(t)$ near zero,
which proves integrability. The rank or ratio conditions stated here
are sufficient conditions for this proof, not necessary conditions for
a bounded derivative or a small shifted-null error.
\end{proof}

\begin{prop}[A finite-sample certificate using only the spectral ratio]
\label{prop:smooth-finite-diagnostic}
Suppose an observed statistic $r_U\in[0,1]$ satisfies
$\Pr(r\le r_U)\ge1-\gamma_r$. Then
\eqref{eq:smooth-finite-coverage} holds with
$B_{\mathrm{smooth}}$ given in \eqref{eq:smooth-finite-certificate}.
In particular, the explicit Gaussian hollow-Gram upper confidence bound
can be used for arbitrary $p$ and unknown, possibly singular covariance.
No independence between $r_U$ and any statistic of interest is required.
\end{prop}

\begin{proof}
On the single event $r\le r_U$, every $K>0$ and every signal satisfying
$\|\bm\theta\|^2\le K\sqrt{2A}$ obey
\[
 \Delta_{\bm\theta}\le
  \frac{\lambda_{\max}(\bm\Omega)\|\bm\theta\|^2}{A}
       \le\sqrt2Kr.
\]
If $r_U<1/2$, the function
$u\mapsto u/(1-4u^2)$ is increasing on $[0,1/2)$, and
Proposition~\ref{prop:smooth-shift-bound} proves
\eqref{eq:smooth-finite-coverage}. Otherwise the universal $1/2$ bound
applies. The same event works for all $K$ and all signals in their
respective balls, which proves the simultaneous assertion, including
for a data-selected $K$.
\end{proof}

\begin{cor}[A certified local effect range]
\label{cor:certified-effect-range}
Under the conditions of Proposition~\ref{prop:smooth-finite-diagnostic},
let $0<\varepsilon<1/2$ and define
\[
 \widehat K_{\varepsilon}=
 \begin{cases}
 \displaystyle\frac{\varepsilon\pi(1-4r_U^2)}{\sqrt2r_U},
       &0<r_U<1/2,\\[5pt]
 0,&1/2\le r_U\le1,\\
 +\infty,&r_U=0.
 \end{cases}
\]
With probability at least $1-\gamma_r$, simultaneously for every
$0<\varepsilon<1/2$,
\[
 \sup_{\|\bm\theta\|^2\le
             \widehat K_{\varepsilon}\sqrt{2A}}d_{\bm\theta}
       \le\varepsilon.
\]
This is the largest range certified by the smooth bound, not the largest
range on which the actual approximation error is below $\varepsilon$.
The hybrid certificate may certify a larger range. A zero value certifies
only the null signal and does not imply that the error elsewhere is large.
\end{cor}
\begin{proof}
On the single event $r\le r_U$, invert
\eqref{eq:smooth-finite-certificate}; the same event covers every
$\varepsilon$ and every signal in the resulting ball. For
$r_U\ge1/2$, the bound cannot certify any positive $K$ at a tolerance
below $1/2$, while $d_{\bm0}=0$. Since $A>0$ implies $r>0$,
the convention at $r_U=0$ is immaterial on the coverage event.
\end{proof}

With the explicit Gaussian residual construction, the same event also
gives $A_L\le A$. Consequently the entirely data-based ball
$\|\bm\theta\|^2\le\widehat K_{\varepsilon}\sqrt{2A_L}$
inherits the guarantee. The tail allocations remain prespecified;
inverting the simultaneous bound does not authorize selecting among
separately calibrated confidence procedures after observing the data.

The previous empirical-interval certificate and the new certificate
can be combined without an additional error allocation. If
$B_{\mathrm{interval}}(K)$ uses the same $r_U$, an upper scale bound
$A_U$ and an interval-probability event, with respective error budgets
$\gamma_r,\gamma_A,\gamma_F$, then
\begin{equation}
 \sup_{\|\bm\theta\|^2\le K\sqrt{2A}}d_{\bm\theta}
  \le\min\{B_{\mathrm{smooth}}(K;r_U),
                    B_{\mathrm{interval}}(K)\}
 \quad\hbox{for every }K>0
 \label{eq:hybrid-smooth-interval-certificate}
\end{equation}
with probability at least $1-\gamma_r-\gamma_A-\gamma_F$.
Indeed, both conclusions hold on the joint event already used for the
interval certificate. There is no second $r_U$ event to pay for.
If only the smooth certificate is reported, all of the desired total
error probability can instead be allocated to its spectral-ratio bound.
This distinction should be kept in numerical comparisons.

For illustration, $r_U=0.2$ and $K=0.25$ give
$B_{\mathrm{smooth}}=0.0268$, while $r_U=0.1$ and $K=0.5$ give
$B_{\mathrm{smooth}}=0.0234$. These are evaluations at specified
confidence-bound values, not claims that those values are attainable
at a particular sample size. The guarantee concerns the Gaussian
population-reference approximation. Its finite-sample confidence
coverage relies on the assumptions of the chosen $r_U$ construction;
using the Gaussian residual bound does not extend that coverage to
non-Gaussian samples.
\begin{remark}[The linear order in a smooth model]
\label{rem:smooth-linear-order}
Fix an integer $p\ge5$, take $\bm\Omega=\bm I_p$ and
$\bm\theta=\sqrt\eta\,\bm e_1$, and let $\eta\downarrow0$.
Then $r=1/\sqrt p<1/2$ and $\Delta=\eta/p$. If $F_p$ and $f_p$
are the CDF and density of $\chi_p^2$, the classical Poisson mixture
gives, for fixed $x>0$,
\[
 F_{\chi_p^2(\eta)}(x)
 =F_p(x)+\frac\eta2\{F_{p+2}(x)-F_p(x)\}+O(\eta^2).
\]
The gamma recurrence
$F_{p+2}(x)-F_p(x)=-(2x/p)f_p(x)$ and a Taylor expansion of
$F_p(x-\eta)$ therefore yield
\[
 F_{\chi_p^2(\eta)}(x)-F_p(x-\eta)
 =\eta(1-x/p)f_p(x)+O(\eta^2).
\]
Taking $x=p+2$ gives
\[
 \liminf_{\eta\downarrow0}\frac{d_{\bm\theta}}{\Delta}
       \ge2f_p(p+2)>0.
\]
Thus a uniform $o(\Delta)$ replacement of the linear bound is
impossible even within this fixed smooth isotropic model. Here
$\kappa=\eta/\sqrt{2p}\to0$; this observation does not assert
sharpness at a fixed positive effect magnitude.
\end{remark}

\subsection{Sharpness of the fourth-root order}
\label{app:sharp-quarter-proof}

\begin{proof}[Proof of Proposition~\ref{prop:shift-quarter-sharp}]
Put $c=(\sqrt2\kappa_0)^{1/2}$. For $0<\delta<c^2$, take
\[
 \bm\Omega_\delta=
 \operatorname{diag}\!\left\{
    \lambda_\delta,\frac{\delta}{c^2}\right\},\qquad
 \lambda_\delta=\left(1-\frac{\delta^2}{c^4}\right)^{1/2},
 \qquad \bm\theta_\delta=(0,c)^\top.
\]
The three displayed identities follow immediately. With independent
standard normals $Z_1,Z_2$, set
\[
 Q_\delta=\lambda_\delta Z_1^2+\frac{\delta}{c^2}Z_2^2,
 \qquad L_\delta=2\sqrt\delta Z_2.
\]
After subtracting the common deterministic shift $c^2$, the
Kolmogorov distance in \eqref{eq:shift-quarter-lower} is
$d_{\mathrm K}(Q_\delta+L_\delta,Q_\delta)$. Since
$\Pr(Q_\delta\le0)=0$, it is at least
\begin{equation}
 D_\delta:=\Pr(Q_\delta+L_\delta\le0)
 =\mathbb E H\!\left(
   \frac{\{-2\sqrt\delta Z_2-\delta Z_2^2/c^2\}_+}
        {\lambda_\delta}\right),
 \quad H(v)=2\Phi(\sqrt v)-1,\quad v\ge0.
 \label{eq:shift-quarter-witness}
\end{equation}
Here $H(0)=0$, so the positive-part notation causes no atom at zero.
The elementary relations
\[
 H(v)\sim\sqrt{\frac2\pi}\sqrt v\quad(v\downarrow0),
 \qquad 0\le H(v)\le\sqrt{\frac2\pi}\sqrt v
\]
show that the integrand in $\delta^{-1/4}D_\delta$ converges
pointwise to $2\pi^{-1/2}\sqrt{(-Z_2)_+}$. For all sufficiently
small $\delta$, it is bounded by
$2\pi^{-1/2}\lambda_\delta^{-1/2}\sqrt{(-Z_2)_+}$, which has
an integrable bound independent of $\delta$. Dominated convergence gives
\begin{equation}
 \lim_{\delta\downarrow0}\frac{D_\delta}{\delta^{1/4}}
 =\frac2{\sqrt\pi}\mathbb E\sqrt{(-Z_2)_+}
 =\frac{2^{1/4}\Gamma(3/4)}{\pi}.
 \label{eq:shift-quarter-witness-limit}
\end{equation}
This proves the lower bound. The existing uniform upper bound is
$C\Delta^{1/4}$, completing the uniform-order conclusion.
\end{proof}

The construction distinguishes a rank-one covariance from a sequence
approaching rank one through positive eigenvalues. The signal lies in
the shrinking-eigenvalue direction, where it induces a Gaussian linear
fluctuation of order $\sqrt\Delta$. The dominant squared-normal term
has probability of order $\Delta^{1/4}$ within that distance of its
support boundary. The limiting covariance places the signal in its
null space, but this limit does not preserve the approximation rate.
The rank-one parallel-signal expansion in
Proposition~\ref{prop:rankone-sharp-error} therefore does not yield a
uniform $\sqrt\Delta$ bound. Faster uniform rates require additional
restrictions on the spectrum or on the signal direction. The next
corollary shows that the same obstruction can also occur at a fixed
significance level, at a particular fixed effect magnitude.

\begin{cor}[A fixed-level power obstruction]
\label{cor:shift-quarter-fixed-level}
Fix $\alpha\in(0,1)$ and let $q_0$ be the $(1-\alpha)$ quantile
of $\chi_1^2$. In the preceding two-dimensional construction, take
$c^2=q_0$, so $\kappa_0=q_0/\sqrt2$. Let $q_\delta$ be the
$(1-\alpha)$ quantile of $Q_\delta$. Define the true and shifted-null
reference powers by
\[
 \beta_\delta=\Pr(Q_\delta+L_\delta+c^2>q_\delta),\qquad
 \widetilde\beta_\delta=\Pr(Q_\delta+c^2>q_\delta).
\]
Then
\begin{equation}
 \widetilde\beta_\delta-\beta_\delta
   =c_*\delta^{1/4}+o(\delta^{1/4}).
 \label{eq:shift-quarter-fixed-level}
\end{equation}
Thus no uniform improvement of the exponent is available for fixed-level
power over unrestricted covariance spectra and local effect magnitudes.
This is an existence statement at $\kappa_0=q_0/\sqrt2$; it does not
assert this rate at every fixed effect magnitude. Both powers tend to one
in this construction.
\end{cor}

\begin{proof}
First, $q_\delta=q_0+O(\delta)$. To verify this, put
$V_\delta=\delta Z_2^2/c^2$. For $x$ in a fixed sufficiently small
neighborhood of $q_0>0$,
\[
 \Pr(Q_\delta\le x)
 =\mathbb E H\{(x-V_\delta)_+/\lambda_\delta\}.
\]
On $V_\delta\le q_0/4$, the arguments of $H$ lie in a compact
subset of $(0,\infty)$, where its density is bounded. Consequently,
the difference from $H(x)$ on this event is bounded by a constant times
$\mathbb E V_\delta+|1-\lambda_\delta|=O(\delta)$.
The complementary event has probability exponentially small in
$1/\delta$. The CDF difference is therefore $O(\delta)$ uniformly
near $q_0$. Since the $\chi_1^2$ density is continuous and positive
at $q_0$, inversion yields the asserted quantile bound.

Write $x_\delta=q_\delta-c^2=O(\delta)$. Conditional on $Z_2$,
\[
 \Pr(Q_\delta+L_\delta\le x_\delta)
 =\mathbb E H\!\left(
  \frac{\{x_\delta-2\sqrt\delta Z_2
                    -\delta Z_2^2/c^2\}_+}{\lambda_\delta}
                  \right).
\]
Because $x_\delta/\sqrt\delta\to0$, the dominated-convergence
argument in \eqref{eq:shift-quarter-witness-limit} applies with the
integrable bound $C\{1+\sqrt{|Z_2|}\}$. This probability equals
$c_*\delta^{1/4}+o(\delta^{1/4})$. In contrast,
\[
 \Pr(Q_\delta\le x_\delta)
 \le H\{(x_\delta)_+/\lambda_\delta\}=O(\delta^{1/2}).
\]
Subtracting the two CDFs proves
\eqref{eq:shift-quarter-fixed-level}.
\end{proof}

\begin{remark}[Increasing dimension with positive definite covariance]
\label{rem:shift-quarter-growing-dimension}
The lower bound is not a consequence of keeping the dimension fixed or
allowing singular covariance matrices. For any integer sequence
$p_\delta\ge3$, possibly tending to infinity, add $p_\delta-2$
eigenvalues equal to $\delta^2/(p_\delta-2)$ and replace
$\lambda_\delta$ by
\[
 \left(1-\frac{\delta^2}{c^4}
            -\frac{\delta^4}{p_\delta-2}\right)^{1/2}.
\]
Keep the signal $(0,c,0,\ldots,0)^\top$. For sufficiently small
$\delta$ the covariance is positive definite, while $A=1$,
$\kappa=\kappa_0$ and $\Delta=\delta$ still hold exactly.
The extra nonnegative quadratic term $W_\delta$ satisfies
$\mathbb E W_\delta=\delta^2$, hence
$W_\delta/\sqrt\delta\to0$ in probability. In the proof above,
the positive part acquires the additional term $-W_\delta$;
its normalized integrand converges in probability to the same limit
and has the same integrable bound depending only on $Z_2$.
Uniform integrability therefore preserves
\eqref{eq:shift-quarter-witness-limit}.
\end{remark}

\subsection{Proofs for the shifted-null boundary}
\label{app:shift-boundary-proofs}

\begin{proof}[Proof of Theorem~\ref{thm:shift-boundary}]
Set $s_0=\sqrt{2A}$ and $\kappa=\|\bm\theta\|^2/s_0$.
In an orthonormal eigenbasis of $\bm\Omega$, put
$\rho_j=\lambda_j/\sqrt A$ and
$\eta_j=\bm e_j^\top\bm\theta/A^{1/4}$, including zero eigenvalues.
For independent standard normal variables $Z_j$, define
\[
 Q=\frac1{\sqrt2}\sum_j\rho_j(Z_j^2-1),\qquad
 L=\sqrt2\sum_j\sqrt{\rho_j}\eta_jZ_j,\qquad
 \Delta=\sum_j\rho_j\eta_j^2.
\]
Then $\sum_j\rho_j^2=1$, $0\leq\rho_j\leq1$,
$\Delta=\Delta_{\bm\theta}$, and
\[
 \frac{T_{\bm\theta}^*-\tr(\bm\Omega)}{s_0}=Q+L+\kappa,
 \qquad
 \frac{T_0^*+\|\bm\theta\|^2-\tr(\bm\Omega)}{s_0}=Q+\kappa.
\]
A common deterministic translation and positive scaling leave
Kolmogorov distance unchanged, so $d_{\bm\theta}=d_{\mathrm K}(Q+L,Q)$.
Moreover, $\E Q=\E L=0$, $\E Q^2=1$, $\E L^2=2\Delta$ and
$\E QL=0$. No bound on $\kappa$ has been used.

For the upper bound, let $X=Q+L$ and $Y=Q$. The polynomial
comparison inequality of \citet[][Theorem 3.5]{Kosov2021}, applied with degree two
and $\Var(Y)=1$, gives
\[
 d_{\mathrm K}(X,Y)\leq d_{\mathrm{TV}}(X,Y)
 \leq C\|X-Y\|_2^{1/2}=C(2\Delta)^{1/4}.
\]
The possible factor of two between conventions for total variation is
absorbed in $C$. This imported upper bound applies to the joint Gaussian
polynomials, without requiring them to be independent.

The bound by $1/2$ follows directly from the sign of $L$.
For each $x$, the two crossing events
$\{Q+L\le x<Q\}$ and $\{Q\le x<Q+L\}$ are contained in
$\{L<0\}$ and $\{L>0\}$, respectively. When $\Delta>0$,
each has probability at most $1/2$, and their probability difference
is the CDF difference. If $\Delta=0$, the two laws coincide.

For the three-block refinement, use the independent representation in
Lemma~\ref{lem:poisson-centered-perturbation}: there exist independent
$Q_0\overset d=Q$ and $H$ such that
$X\overset d=Q_0+H$, $\E H=0$ and $\E H^2=2\Delta$.
The three-block condition and the characteristic-function formula in
Lemma~\ref{adap:fourier} give
\[
 |\varphi_Q(t)|\leq(1+2ct^2)^{-3/4}
       \leq(2c)^{-3/4}|t|^{-3/2}\qquad(t\ne0).
\]
In particular, $\varphi_Q$ is integrable and $Q$ has a continuous
density $f_Q$. Fourier inversion yields, for all real $x,h$,
\begin{align}
 |f_Q(x+h)-f_Q(x)|
 &\leq\frac1\pi\int_0^\infty
                \min\{2,t|h|\}|\varphi_Q(t)|\,dt\notag\\
 &\leq C_{c,H}|h|^{1/2},\qquad
 C_{c,H}=\frac{4\sqrt2}{\pi(2c)^{3/4}}.
 \label{eq:three-block-density-holder}
\end{align}
For $h\ne0$, the last integral is evaluated by splitting at
$t=2/|h|$; the assertion is immediate for $h=0$.
Writing $F_Q$ for the CDF of $Q$, integration of the density difference
gives
\[
 |F_Q(x-h)-F_Q(x)+h f_Q(x)|
       \leq\frac23C_{c,H}|h|^{3/2}.
\]
Independence permits conditioning on $H$, and its zero mean cancels
the linear term. Therefore
\begin{equation}
 d_{\bm\theta}
 \leq\frac23C_{c,H}\E|H|^{3/2}
 \leq\frac23C_{c,H}(2\Delta)^{3/4}
 =\frac{8\sqrt2}{3\pi c^{3/4}}\Delta^{3/4}.
 \label{eq:three-block-shift-explicit}
\end{equation}
The second inequality is the monotonicity of $L^q$ norms. Combining
this with the universal bound by $1/2$ proves
\eqref{eq:shift-boundary-spread}. This argument uses the independent
Poisson-mixture representation in distribution; it does not assume
independence of the original quadratic and linear terms $Q$ and $L$.

We next prove the lower bound. The fourth centered moment of $Q$ is
\[
 \E Q^4=3+12\sum_j\rho_j^4\leq15.
\]
Since $L$ is centered Gaussian, $\E L^4=12\Delta^2$.
The mixed moment is
\begin{equation}
 \E Q^2L^2=2\Delta+8\sum_j\rho_j^3\eta_j^2.
 \label{eq:shift-mixed-moment}
\end{equation}
For completeness, write $a_j=\rho_j/\sqrt2$ and
$b_j=\sqrt{2\rho_j}\eta_j$, so that
$Q=\sum_j a_j(Z_j^2-1)$ and $L=\sum_jb_jZ_j$.
Terms containing $Z_jZ_k$ with $j\ne k$ have expectation zero because
$Q^2$ is even in each coordinate. Independence and the standard
normal moments give
\[
 \E(Q^2Z_j^2)
 =\E Q^2+8a_j^2=1+4\rho_j^2.
\]
Indeed,
$\E[(Z_j^2-1)^2Z_j^2]=15-6+1=10$, whereas
$\E[(Z_j^2-1)^2]=2$; all cross terms involving the other
coordinates have zero expectation. Summing
$b_j^2(1+4\rho_j^2)$ proves \eqref{eq:shift-mixed-moment}.
The odd mixed moments $\E Q^3L$ and $\E QL^3$ vanish under
$\bm Z\mapsto-\bm Z$. Therefore
\begin{align*}
 \E X^4
 &=\E Q^4+6\E Q^2L^2+\E L^4\\
 &\leq15+12\Delta+48\sum_j\rho_j^3\eta_j^2+12\Delta^2\\
 &\leq15+60\Delta+12\Delta^2.
\end{align*}
Here $\rho_j^3\leq\rho_j$ was used in the last step.
Consequently
\[
 \E X^2-\E Y^2=2\Delta,\qquad
 \E X^4+\E Y^4\leq M_\Delta:=30+60\Delta+12\Delta^2.
\]
Let $d=d_{\mathrm K}(X,Y)$ and $f_R(x)=\min(x^2,R^2)$ for $R>0$.
The total variation of $f_R$ on $\mathbb R$ is $2R^2$.
The integration-by-parts bound for expectations of bounded-variation
functions gives
\[
 |\E f_R(X)-\E f_R(Y)|\leq2R^2d.
\]
Also,
$0\leq\E X^2-\E f_R(X)\leq\E X^4/R^2$, and the same
bound holds for $Y$. Hence
\begin{equation}
 2\Delta\leq2R^2d+\frac{M_\Delta}{R^2}.
 \label{eq:shift-truncation}
\end{equation}
If $d=0$, letting $R\to\infty$ gives $\Delta=0$.
If $d>0$, taking $R^2=\sqrt{M_\Delta/(2d)}$ yields
$2\Delta\leq2\sqrt{2M_\Delta d}$, or
\[
 d\geq\frac{\Delta^2}{2M_\Delta}
   =\frac{\Delta^2}{60+120\Delta+24\Delta^2}.
\]
This proves \eqref{eq:shift-boundary-bounds}.
The function on the right is strictly increasing on $(0,\infty)$,
is positive there, and tends to $1/24$ at infinity. Together with
the upper bound, it proves \eqref{eq:shift-boundary-equivalence}
even for unbounded sequences of $\Delta$.
\end{proof}

\begin{proof}[Proof of Corollary~\ref{cor:shift-uniform-directions}]
At least one covariance eigenvalue is positive. The distribution of
$T_0^*=\sum_j\lambda_jZ_j^2$ is continuous and strictly increasing
on $(0,\infty)$, and has support $[0,\infty)$. To see strict
increase, one positive summand has a density strictly positive on
$(0,\infty)$, and convolution with the independent nonnegative
remaining summands assigns positive probability to every positive
interval. Thus $q_{\alpha}^{\mathrm{NR}}$ sweeps $(0,\infty)$ as
$\alpha$ varies in $(0,1)$. Both $T_{\bm\theta}^*$ and
$T_0^*+\|\bm\theta\|^2$ are nonnegative and have no atom at zero.
Their distribution functions agree at every negative argument and
at zero. Taking the supremum over the positive quantiles therefore
gives their full Kolmogorov distance, which proves
\eqref{eq:shift-power-supremum}.

On $\Theta_n(K)$, the Rayleigh quotient gives
\[
 \Delta_{\bm\theta}
 \leq\frac{\lambda_{\max}(\bm\Omega_n)\|\bm\theta\|^2}{A_n}
 \leq r_n.
\]
Equality is attained by a signal along a leading eigenvector with
squared norm $K\sqrt{2A_n}$. Apply both bounds in
Theorem~\ref{thm:shift-boundary} to obtain
\eqref{eq:shift-uniform-bounds}. Since $K>0$ is fixed, the
monotonicity and positivity of its lower bound prove
\eqref{eq:shift-uniform-equivalence}.

Finally, suppose that
\[
 e_n=\sup_{\bm\theta\in\Theta_n(K)}
 d_{\mathrm K}(T_{\bm\theta},T_{\bm\theta}^*)\longrightarrow0.
\]
For each $\bm\theta$, the triangle inequality and its reverse give
\[
 \left|d_{\mathrm K}(T_{\bm\theta},T_0^*+\|\bm\theta\|^2)
      -d_{\bm\theta}\right|\leq e_n.
\]
Taking suprema proves the equivalence for the original-data
Kolmogorov distances. For power at central-reference thresholds,
\[
 \sup_{\bm\theta\in\Theta_n(K)}\sup_{0<\alpha<1}
 \left|\Pr\{T_{\bm\theta}>q_{\alpha}^{\mathrm{NR}}\}
       -\Pr\{T_{\bm\theta}^*>q_{\alpha}^{\mathrm{NR}}\}\right|
 \leq e_n.
\]
Combine this with \eqref{eq:shift-power-supremum} and the same
reverse-triangle argument. This proof does not require continuity
of the original-data law.
\end{proof}

\begin{proof}[Proof of Corollary~\ref{cor:shift-feasible-bound}]
Write $q=q_{\alpha}^{\mathrm{NR}}$, $s_0=\sqrt{2A}$ and
$r_\varepsilon=\Pr\{|\widehat q_{\alpha}-q|>\varepsilon s_0\}$.
Deterministic bracketing gives
\[
 \Pr\{T_{\bm\theta}>q+\varepsilon s_0\}-r_\varepsilon
 \leq\Pr\{T_{\bm\theta}>\widehat q_{\alpha}\}
 \leq\Pr\{T_{\bm\theta}>q-\varepsilon s_0\}+r_\varepsilon.
\]
At each endpoint, replace $T_{\bm\theta}$ by
$T_{\bm\theta}^*$ at a cost at most $e_{\bm\theta}$.
The centered Gaussian polynomial
$(T_{\bm\theta}^*-\tr\bm\Omega-\|\bm\theta\|^2)/s_0$
has variance $1+2\Delta_{\bm\theta}\geq1$.
After standardization, Lemma~\ref{lem:gaussian-quadratic-small-ball}
bounds its probability in any interval of length
$2\varepsilon$ by $C\varepsilon^{1/2}$. Thus
\[
 \left|\Pr\{T_{\bm\theta}>\widehat q_{\alpha}\}
        -\Pr\{T_{\bm\theta}^*>q\}\right|
 \leq e_{\bm\theta}+r_\varepsilon+C\varepsilon^{1/2}.
\]
Apply Theorem~\ref{thm:shift-boundary} at the fixed threshold $q$
to obtain \eqref{eq:shift-feasible-bound}.
Under the three-block condition, the density of
$T_{\bm\theta}^*/s_0$ is bounded by
$c_3s_0/\sqrt{cA}=c_3\sqrt{2/c}$, uniformly in the signal,
by Lemma~\ref{adap:fourier} and the change-of-scale formula.
The interval-probability cost is then $C_c\varepsilon$, and
\eqref{eq:shift-boundary-spread} gives the other stated replacement.
\end{proof}

\subsection{Proofs for directional approximation}\label{app:directional-boundary-proofs}

\begin{proof}[Proof of Theorem~\ref{thm:directional-diameter}]
Put $\bm B=\bm\Omega/\sqrt A$, and let $\rho_1\geq\cdots\geq\rho_p\geq0$
be its eigenvalues. Thus $\tr(\bm B^2)=1$ and
$d=\rho_1-\rho_p\in[0,1]$.
For $\bm\eta=\bm\theta/A^{1/4}$ and a standard Gaussian vector
$\bm Z\in\mathbb R^p$, write
\begin{equation}
 X_{\bm\eta}(\bm Z)
 =\frac{\bm Z^\top\bm B\bm Z-\tr(\bm B)}{\sqrt2}
    +\sqrt2\bm\eta^\top\bm B^{1/2}\bm Z,
 \qquad \|\bm\eta\|^2=h:=\sqrt2\kappa.
 \label{eq:directional-normalized-polynomial}
\end{equation}
Then $U_{\bm\theta}\stackrel d=X_{\bm\eta}+\kappa$ and
\[
 \E X_{\bm\eta}=0,\qquad
 \Var(X_{\bm\eta})=1+2\bm\eta^\top\bm B\bm\eta\geq1.
\]
Lemma~\ref{lem:gaussian-quadratic-small-ball}, after division by the
standard deviation, therefore gives
\begin{equation}
 \sup_{\bm\eta:\|\bm\eta\|^2=h}\sup_x
 \Pr\{|X_{\bm\eta}-x|\leq\varepsilon\}
 \leq C\varepsilon^{1/2}.
 \label{eq:directional-common-concentration}
\end{equation}

First consider the upper bound. For two distinct vectors
$\bm\eta,\bm\xi$ of squared norm $h$, the Householder reflection
\[
 \bm H=\bm I_p-2\bm u\bm u^\top,
 \qquad
 \bm u=\frac{\bm\eta-\bm\xi}{\|\bm\eta-\bm\xi\|}
\]
satisfies $\bm H^\top=\bm H$, $\bm H^2=\bm I_p$ and
$\bm H\bm\eta=\bm\xi$. For equal vectors take $\bm H=\bm I_p$.
Since $\bm H\bm Z\stackrel d=\bm Z$, a coupling of the two laws is
\begin{align*}
 X_{\bm\eta}(\bm Z)-X_{\bm\xi}(\bm H\bm Z)
  &=\frac{\bm Z^\top\bm D\bm Z}{\sqrt2}
     +\sqrt2\bm\eta^\top\bm E\bm Z,\\
 \bm D&=\bm B-\bm H\bm B\bm H,
 &\bm E&=\bm B^{1/2}-\bm H\bm B^{1/2}\bm H.
\end{align*}
Here $\tr(\bm D)=0$ and $\operatorname{rank}(\bm D)\leq2$:
the range of $\bm D$ is contained in the span of $\bm u$ and
$\bm B\bm u$. Subtracting the midpoint of the spectral interval before
conjugation gives
\[
 \|\bm D\|_{\mathrm{op}}\leq d,
 \qquad
 \|\bm E\|_{\mathrm{op}}
   \leq\sqrt{\rho_1}-\sqrt{\rho_p}\leq\sqrt d.
\]
The quadratic and linear differences are uncorrelated by Gaussian
symmetry, and hence
\begin{equation}
 \E|X_{\bm\eta}(\bm Z)-X_{\bm\xi}(\bm H\bm Z)|^2
 =\tr(\bm D^2)+2\bm\eta^\top\bm E^2\bm\eta
 \leq2d^2+2h d.
 \label{eq:householder-coupling-bound}
\end{equation}
The dimension-free stability bound of
\citet[][Theorem~3.5]{Kosov2021}, applied to degree-two Gaussian
polynomials $f,g$, gives
\[
 d_{\mathrm{TV}}(f,g)
 \leq C\{\Var(g)\}^{-1/4}\|f-g\|_2^{1/2}.
\]
Apply it to the two coupled polynomials above. Their variances are at
least one, and Kolmogorov distance is bounded by total variation.
Consequently,
\[
 d_{\mathrm K}(X_{\bm\eta},X_{\bm\xi})
 \leq C(2d^2+2h d)^{1/4}
 \leq C_\kappa d^{1/4},
\]
where the last step uses $0\leq d\leq1$. If $d=0$, the coupling
difference is exactly zero. Taking the supremum proves the upper bound.
The stability exponent is imported from the cited theorem; the
Householder coupling supplies the spectral-range dependence.

For the lower bound, first note the following elementary moment
comparison. If $X,Y$ are centered, $\E X^4+\E Y^4\leq M$, and
$e=d_{\mathrm K}(X,Y)$, then for every $r>0$,
\begin{equation}
 |\E X^2-\E Y^2|\leq2r^2e+\frac M{r^2}.
 \label{eq:kolmogorov-variance-comparison}
\end{equation}
Indeed, $f_r(x)=\min(x^2,r^2)$ has total variation $2r^2$, so the
difference between its expectations is at most $2r^2e$. Each second
moment differs from the corresponding truncated moment by at most its
fourth moment divided by $r^2$. Optimizing in $r$, or letting $r$ tend
to infinity when $e=0$, yields
\begin{equation}
 e\geq\frac{|\E X^2-\E Y^2|^2}{8M}.
 \label{eq:kolmogorov-variance-lower-bound}
\end{equation}

For the quadratic part $Q$ of
\eqref{eq:directional-normalized-polynomial},
\[
 \E Q^4=3+12\sum_j\rho_j^4\leq15.
\]
The linear part $L_{\bm\eta}$ has variance at most $2h$, and thus
$\E L_{\bm\eta}^4\leq12h^2=24\kappa^2$.
Consequently
\[
 \E X_{\bm\eta}^4\leq8(\E Q^4+\E L_{\bm\eta}^4)
 \leq120+192\kappa^2.
\]
Take $\bm\eta=\sqrt h\bm e_1$ and
$\bm\xi=\sqrt h\bm e_p$, using eigenvectors associated with the
largest and smallest eigenvalues. Their variance difference is
$2h d=2\sqrt2\kappa d$. Applying
\eqref{eq:kolmogorov-variance-lower-bound} with
$M=240+384\kappa^2$ gives
\[
 D_\kappa(\bm\Omega)
 \geq\frac{\kappa^2}{240+384\kappa^2}\,d^2.
\]
This also covers zero smallest eigenvalues. The equivalence and uniform
constants on $\kappa\in[\kappa_0,K]$ follow immediately.
Finally, $p\rho_p^2\leq\sum_j\rho_j^2=1$, and hence
\[
 0\leq\rho_1-d=\rho_p\leq p^{-1/2}.
\]
If $p\to\infty$, $d\to0$ is therefore equivalent to $\rho_1\to0$.
\end{proof}

\begin{proof}[Proof of Corollary~\ref{cor:direction-blind-reference}]
For any distribution function $H$, the triangle inequality gives
$D_\kappa\leq2\sup_{\bm\theta\in\mathcal S_\kappa}
\|F_{\bm\theta}-H\|_\infty$.
For the converse, define the pointwise envelopes
\[
 F_-(x)=\inf_{\bm\theta\in\mathcal S_\kappa}F_{\bm\theta}(x),
 \qquad
 F_+(x)=\sup_{\bm\theta\in\mathcal S_\kappa}F_{\bm\theta}(x).
\]
Both are nondecreasing. The common concentration bound
\eqref{eq:directional-common-concentration} makes both continuous.
Their limits at $-\infty$ and $+\infty$ are respectively zero and one,
since the family has common mean $\kappa$ and variance bounded by
$1+2\sqrt2\kappa$. Therefore $H_*=(F_-+F_+)/2$ is a distribution
function. Moreover,
\[
 \sup_{\bm\theta}\|F_{\bm\theta}-H_*\|_\infty
 \leq\frac12\sup_x\{F_+(x)-F_-(x)\}
 =\frac12D_\kappa.
\]
This proves \eqref{eq:direction-blind-radius}. The midpoint reference
is used only to characterize the best possible approximation error;
no claim about its computational implementation is needed.

The null law is a centered positive weighted chi-squared sum with
support $[a,\infty)$, where $a=-\tr(\bm\Omega)/\sqrt{2A}$.
Its density is positive on $(a,\infty)$, so its quantiles
$c_\alpha$, $0<\alpha<1$, cover this open interval.
Every alternative $U_{\bm\theta}$ is supported on $[a,\infty)$ as
well. The supremum of the difference of two alternative distribution
functions can thus be taken over these null quantiles, proving
\eqref{eq:directional-power-diameter}.
\end{proof}

\subsection{Proofs for noncentral power and boundary statements}
\begin{proof}[Proof of Theorem~\ref{prof:theorem}]
Fatou's lemma gives $\sum_j\rho_j^2\leq1$ and
$\sum_j\eta_j^2\leq\sqrt{2}\kappa$.  Since $\rho_1\leq1$, this also
gives $\sum_j\rho_j\eta_j^2<\infty$, proving the asserted $L^2$
convergence.  Choose a finite $K$ such that $\kappa_n\leq K$ for all
sufficiently large $n$.  The tail of the linear sum satisfies
\begin{equation}
 \operatorname{E}\left|
 \sqrt{2}\sum_{j>m}\sqrt{\rho_{n,j}}\eta_{n,j}Z_j
 \right|^2
 =2\sum_{j>m}\rho_{n,j}\eta_{n,j}^2
 \leq2\sqrt{2}K\rho_{n,m+1}.
 \label{prof:linear-tail}
\end{equation}
Because the $\rho_j$ are nonincreasing and square summable,
$\rho_j\to0$.  Hence the right-hand side in \eqref{prof:linear-tail} tends to zero
when first $n\to\infty$ and then $m\to\infty$.
This step uses the assumed bounded local-signal scale.

We next approximate the quadratic tail by a Gaussian remainder.
For fixed $t\in\mathbb R$, the characteristic function of
$r(Z^2-1)/\sqrt{2}$ has logarithm
\[
 -\frac{itr}{\sqrt{2}}-\frac12\log(1-i\sqrt{2}tr)
 =-\frac{t^2r^2}{2}+O_t(r^3)
\]
uniformly for sufficiently small nonnegative $r$.  Thus, for large $m$
and then large $n$,
\begin{equation}
 \log\operatorname{E}\exp\left\{
 \frac{it}{\sqrt{2}}\sum_{j>m}\rho_{n,j}(Z_j^2-1)\right\}
 =-\frac{t^2}{2}\sum_{j>m}\rho_{n,j}^2+R_{n,m}(t),
 \quad
 |R_{n,m}(t)|\leq C_t\rho_{n,m+1}\sum_{j>m}\rho_{n,j}^2.
 \label{prof:quadratic-tail}
\end{equation}
For any fixed $m$, the first $m$ coordinate pairs in
\eqref{prof:QL} converge jointly to the corresponding pairs in
\eqref{prof:limit}.  If $\rho_j=0$, their linear coefficients tend to
zero even without convergence of $\eta_{n,j}$, since
$|\eta_{n,j}|^2\leq\sqrt{2}K$.  The head and tail are independent
for each $n$.  In a joint characteristic function with arguments $(t,u)$,
deleting the linear tail changes its value by at most
$|u|\{2\sqrt{2}K\rho_{n,m+1}\}^{1/2}$, by
\eqref{prof:linear-tail}.  Equation~\eqref{prof:quadratic-tail} then
replaces the quadratic tail, with an error vanishing as $m\to\infty$,
by an independent normal variable with variance
$1-\sum_{j\leq m}\rho_j^2$.  Taking $n\to\infty$ and then
$m\to\infty$ proves \eqref{prof:joint-limit} by the characteristic
function continuity theorem.  Equation~\eqref{prof:alternative-limit}
follows from $\kappa_n\to\kappa$.  Finally,
$\operatorname{E}\{(Z_j^2-1)Z_j\}=0$ gives
$\operatorname{Cov}(\zeta,L)=0$, and the variance identities follow.
The joint representation, rather than independence of $\zeta$ and $L$,
is used throughout.
\end{proof}

\begin{proof}[Proof of Corollary~\ref{prof:power}]
We first establish continuity of the limiting laws and uniqueness
of the null quantile. If
$\sigma>0$, both $\zeta$ and $V$ have continuous distributions, and
the density of $\zeta$ is positive on $\mathbb{R}$, by convolution
with the independent Gaussian term.  If $\sigma=0$, some $\rho_j$
is positive.  Writing $a=\rho_j/\sqrt{2}>0$, the null limit has
the form $a(Z_j^2-1)+Y$, with $Y$ independent of $Z_j$.  The first
term is atomless and has a positive density throughout $(-a,\infty)$.
Consequently $\zeta$ is atomless and its distribution function is
strictly increasing above the lower endpoint of its support: every
open interval above that endpoint has positive conditional probability
for a set of $Y$ values of positive probability.  Its quantile at any
probability in $(0,1)$ is therefore unique and is a point of strict
increase.  The alternative $V$ is also atomless, because its dependence
on this same $Z_j$, conditional on all other terms, is a nonconstant
quadratic polynomial of a continuous Gaussian variable.

It follows that $c_{n,\alpha}\to c_{\alpha}$.  Apply
\eqref{prof:alternative-limit} and
$\Pr(V=c_{\alpha})=0$.  For an estimated critical value, sandwich the rejection
event between the events at $c_{\alpha}-\varepsilon$ and
$c_{\alpha}+\varepsilon$, together with the event that the estimated
critical value is within $\varepsilon$ of $c_{\alpha}$.  Convergence
in probability and then $\varepsilon\downarrow0$ prove the last claim.
The bracketing argument does not require independence between the
statistic and its critical value.
\end{proof}

\begin{proof}[Proof of Proposition~\ref{prof:spike}]
Here $\rho_{p,1}\to1/\sqrt{2}$, while every other fixed normalized
eigenvalue tends to zero.  The diffuse variance is
$1-\rho_1^2=1/2$.  For the parallel signal,
$\eta_{p,1}^2=\sqrt{2}\kappa$, so its limiting linear coefficient is
$\sqrt{2\rho_1}\eta_1=\sqrt{2\kappa}$.  The perpendicular signal
has no energy in the persistent eigendirection.  Apply
Theorem~\ref{prof:theorem}.  The variance calculations follow either
directly or from \eqref{prof:limit-variance}.  The central limit
$\zeta=(Z_1^2-1)/2+Z_0/\sqrt{2}$ has a strictly positive density on
$\mathbb{R}$ because of its independent Gaussian component.  Its
quantiles therefore range over $\mathbb{R}$.  If the two alternative
rejection probabilities agreed at all such quantiles, their continuous
distribution functions would agree everywhere, contradicting the
different variances.
\end{proof}

\begin{proof}[Proof of Proposition~\ref{prof:matching}]
The first assertion follows from Proposition~\ref{prof:spike} with
zero signal.  Since $\nu_{2,p}\sim p/2\to\infty$, the standardized
chi-squared central limit theorem proves \eqref{prof:2c-null}.
For three-cumulant matching,
$\nu_{3,p}\to8$ and $\beta_{3,p}/s_{0,p}\to1/4$.
Also $a_{3,p}+\beta_{3,p}\nu_{3,p}=B_p$, proving
\eqref{prof:3c-null}.  The third cumulant of the true limit is
$8(1/2)^3=1$, whereas that of $N(0,1)$ is zero.  The true and
three-cumulant limits both have third cumulant one, but their fourth
cumulants are respectively
\[
 48(1/2)^4=3,
 \qquad
 \frac{48\cdot8}{4^4}=\frac32.
\]
Thus the limiting laws differ. Since their distribution functions are
continuous, weak convergence gives uniform convergence of the
corresponding distribution functions. The Kolmogorov distances
therefore converge to positive distances between the limiting laws.
Invariance under a common positive affine transformation proves
the assertion for the unstandardized references.
\end{proof}

\begin{proof}[Proof of Proposition~\ref{prof:radial}]
Conditional on the radii, $\bm{W}_n\sim N(\bm{0},S_n\bm{I}_p)$, giving
\eqref{prof:radial-mixture}.  Since
$\operatorname{E}S_n=1$ and $\operatorname{Var}(S_n)=v/n$,
the conditional variance identity gives
\[
 \operatorname{Var}(S_nU_p)
 =2p\operatorname{E}S_n^2+p^2\operatorname{Var}(S_n)
 =2p+v(p^2+2p)/n.
\]
The ordinary central limit theorem applies to $R_i^2$, because it has
finite variance.  Independently,
$(U_p-p)/\sqrt{2p}\xrightarrow{d}N(0,1)$.  The decomposition
\[
 \frac{S_nU_p-p}{\sqrt{2p}}
 =S_n\frac{U_p-p}{\sqrt{2p}}
   +\sqrt{\frac{p}{2n}}\,\sqrt n(S_n-1)
\]
and $S_n\xrightarrow{p}1$ prove \eqref{prof:radial-limit}.
Finally $\bm{a}^{\top}\bm{Z}_i\sim N(0,\|\bm{a}\|^2)$ and
$\operatorname{E}R_i^4=1+v$, proving
\eqref{prof:projection-fourth}.
\end{proof}

\begin{proof}[Proof of Proposition~\ref{prof:diagonal}]
Since $\|\bm{X}_i\|^2\stackrel{d}{=}R_i^2\chi_p^2$,
\[
 \operatorname{Var}(\|\bm{X}_i\|^2)=2p(1+v)+p^2v.
\]
Independence across observations gives the first identity in
\eqref{prof:diagonal-size}.  Distinct off-diagonal products are
uncorrelated, including pairs sharing one index, because every $\bm{X}_i$
has mean zero.  Also
$\operatorname{E}(\bm{X}_i^{\top}\bm{X}_j)^2=p$ for $i\ne j$.  This proves
the variance formula for $U_n$.  The same zero-mean conditioning proves
$\operatorname{Cov}(D_n,U_n)=0$.  Scaling the preceding squared-norm
variance by $n^{-2}$ and dividing by
$\operatorname{tr}((\bm{I}_p/n)^2)=p/n^2$ proves
\eqref{prof:norm-condition-failure}.
\end{proof}

\begin{proof}[Proof of Proposition~\ref{prop:consistency}]
Since the standardized null reference has mean zero and variance one, its
\((1-\alpha)\)-quantile is bounded above and below by constants depending
only on \(\alpha\), by Chebyshev's inequality. Write its upper bound as
\(c_\alpha\). The mean of \(T_n^*\) is
\(\mu_n+\|\bm{\theta}_n\|^2\), and
\[
 \frac{s_{\bm{\theta},n}^2}{\|\bm{\theta}_n\|^4}
 =\frac{2A_n+4\bm{\theta}_n^\top\bm{\Omega}_n\bm{\theta}_n}{\|\bm{\theta}_n\|^4}
 \le\kappa_n^{-2}+\frac{2\sqrt2}{\kappa_n}\longrightarrow0,
\]
using \(\lambda_{\max}(\bm{\Omega}_n)\le\sqrt{A_n}\).
Chebyshev's inequality therefore makes the Gaussian probability of falling
below \(\mu_n+c_\alpha s_{0,n}\) tend to zero. Gaussian comparison
proves the oracle assertion. Bracketing the estimated threshold by the
oracle threshold plus \(\epsilon s_{0,n}\) proves the last assertion.
\end{proof}

\section{Spectral calibration: additional results and proofs}\label{app:calibration}

This appendix gives the Gaussian spherical-reference proof and its
non-Gaussian extension, followed by the alternative moment-based
spectral construction. The
latter estimates leading eigenvalues and a Gaussian or gamma remainder,
under stronger conditions than the original-data Gaussian comparison.

\subsection{Gaussian spherical calibration}
\label{app:spherical-proof}

\begin{proof}[Gaussian residual construction]
Within group $i$, the sample mean and the orthogonal residual contrasts
$\bm Y_i\bm H_i$ are jointly Gaussian and uncorrelated. They are
therefore independent. The residual columns are independent
$N(\bm0,\bm\Sigma_i)$ vectors because
$\bm H_i^\top\bm H_i=\bm I$. Multiplying by $\sqrt{w_i}$
gives independent $N(\bm0,\bm\Omega_i)$ columns, with
$\bm\Omega_i=w_i\bm\Sigma_i$. The two groups are independent.
Thus the sums in \eqref{eq:spherical-columns} are independent
$N(\bm0,\bm\Omega_1+\bm\Omega_2)$ vectors and are independent
of the normalized mean contrast. This argument remains valid for
singular covariance matrices because Gaussian independence requires
zero cross-covariance, not nonsingularity.
\end{proof}

\begin{proof}[Finite-sample assertion in Theorem~\ref{thm:spherical-calibration}]
Under $H_0$, the columns of $\bm X$ are independent
$N(\bm0,\bm\Omega)$. Hence $\bm X\bm O$ and $\bm X$
have the same distribution for every deterministic orthogonal
$N\times N$ matrix $\bm O$, and their left Gram matrices
$\bm X\bm X^\top$ are equal. Conditional on this left Gram
matrix, the distribution of $\bm X$ is therefore invariant under
right orthogonal multiplication. One may also see this directly by
writing $\bm X=\bm\Omega^{1/2}\bm Z$ with independent standard
Gaussian entries in $\bm Z$ and disintegrating its orthogonally
invariant law. The orthogonal group is compact, so averaging a
conditional version over Haar measure yields the invariant version.

If $\bm O$ is an independent Haar orthogonal matrix, its first
column is uniform on the unit sphere. Thus the conditional distribution
of $T=\|\bm X\bm e_0\|^2$, given $\bm X\bm X^\top$,
equals that of $\|\bm X\bm U\|^2$. The latter conditional
distribution depends on $\bm X$ only through its squared singular
values, or equivalently through $\bm X\bm X^\top$.
Given $\bm X$, the $B$ auxiliary draws are independent with that
same distribution. It follows that, conditional on
$\bm X\bm X^\top$, the $B+1$ quantities
$T,R_1^{\mathrm{sph}},\ldots,R_B^{\mathrm{sph}}$ are independent
and identically distributed. Their upper rank gives a super-uniform
$p$-value with the tie convention in
\eqref{eq:spherical-monte-carlo-p}. Without ties the rank is uniform
on $\{1,\ldots,B+1\}$, proving both assertions.
Finally, $\bm U\stackrel d=\bm Z/\|\bm Z\|$ for a standard
Gaussian vector $\bm Z\in\mathbb R^N$. Diagonalizing
$\bm X^\top\bm X$ gives \eqref{eq:spherical-reference}.
\end{proof}

\begin{proof}[Conditional reference approximation]
Write $\bm X_j=\bm G_j+\boldsymbol1\{j=0\}\bm\theta$,
where the $\bm G_j$ are independent $N(\bm0,\bm\Omega)$.
Let $\bm U,\bm V$ be independent uniform sphere vectors,
independent of the $\bm G_j$, and put $c=\bm U^\top\bm V$.
Conditional on $(\bm U,\bm V)$, the pair
\[
 \left(\sum_jU_j\bm G_j,\sum_jV_j\bm G_j\right)
\]
is jointly Gaussian with marginal covariance $\bm\Omega$ and
cross-covariance $c\bm\Omega$. Consequently it can be coupled as
\[
 (\bm G,\ c\bm G+\sqrt{1-c^2}\bm G'),
\]
where $\bm G,\bm G'$ are independent $N(\bm0,\bm\Omega)$
vectors, independent of $(\bm U,\bm V)$. This coupling also
applies to singular $\bm\Omega$. Since $N\ge2$, the event
$|c|=1$ has probability zero; defining the coupling there arbitrarily
does not affect the argument.

Put $\mu=\tr(\bm\Omega)$, $s_0=\sqrt{2A}$ and
$\Delta=\bm\theta^\top\bm\Omega\bm\theta/A$.
For the two auxiliary reference values, use the coupled versions
\[
 Y_1=\|\bm G+U_0\bm\theta\|^2,\qquad
 Y_2=\|c\bm G+\sqrt{1-c^2}\bm G'+V_0\bm\theta\|^2,
\]
and the independent central targets
$W_1=\|\bm G\|^2$, $W_2=\|\bm G'\|^2$.
The sphere moments are
\[
 \E U_0^2=\E V_0^2=\E c^2=1/N,\qquad
 \E U_0^4=\E V_0^4=\frac{3}{N(N+2)}.
\]
Expansion of the first square gives
\begin{equation}
 \E\left(\frac{Y_1-W_1}{s_0}\right)^2
 =\frac{2\Delta}{N}
   +\frac{3\kappa^2}{N(N+2)},\qquad
 \kappa=\frac{\|\bm\theta\|^2}{s_0}.
 \label{eq:spherical-first-coupling}
\end{equation}
For the second square, conditional on $c$,
\[
 \E\left\{\|c\bm G+\sqrt{1-c^2}\bm G'\|^2
                      -\|\bm G'\|^2\right\}^2=4Ac^2.
\]
Indeed, the difference is
$c^2(\|\bm G\|^2-\|\bm G'\|^2)
 +2c\sqrt{1-c^2}\bm G^\top\bm G'$; its two terms are
uncorrelated, with total variance
$4Ac^4+4Ac^2(1-c^2)=4Ac^2$.
The cross-products with the added signal-linear term vanish by odd
Gaussian moments; the cross-product with the deterministic squared
signal term vanishes because the preceding difference is centered.
Thus
\begin{equation}
 \E\left(\frac{Y_2-W_2}{s_0}\right)^2
 =\frac{2(1+\Delta)}{N}
   +\frac{3\kappa^2}{N(N+2)}.
 \label{eq:spherical-second-coupling}
\end{equation}
Since $\lambda_{\max}(\bm\Omega)\le\sqrt A$,
$\Delta\le\sqrt2\kappa$. Both expressions are bounded by
$C_K/N$ uniformly over $\kappa\le K$, all spectra and all dimensions.

Let $H$ be the CDF of $(W_1-\mu)/s_0$ and let
$\widehat H_{\bm X}$ be the conditional CDF of
$(R^{\mathrm{sph}}-\mu)/s_0$. The Gaussian quadratic small-ball
bound used in Lemma~\ref{lem:gaussian-quadratic-small-ball} gives
\[
 \sup_x\Pr\{x<(W_1-\mu)/s_0\le x+u\}\le C u^{1/2}.
\]
For any $u>0$, the coupling inequalities and Chebyshev's inequality
therefore give, uniformly in $x,y$,
\begin{align*}
 \left|\Pr\{(Y_1-\mu)/s_0\le x\}-H(x)\right|
 &\le C_K\{N^{-1}u^{-2}+u^{1/2}\},\\
 \left|\Pr\{(Y_1-\mu)/s_0\le x,
                   (Y_2-\mu)/s_0\le y\}-H(x)H(y)\right|
 &\le C_K\{N^{-1}u^{-2}+u^{1/2}\}.
\end{align*}
Choose $u=N^{-2/5}$. Given $\bm X$, the two reference draws are
independent with CDF $\widehat H_{\bm X}$. The preceding bounds imply
\begin{equation}
 \sup_x\E\{\widehat H_{\bm X}(x)-H(x)\}^2
 \le C_K N^{-1/5}.
 \label{eq:spherical-cdf-second-moment}
\end{equation}
No limiting spectrum has been used.

For $0<\eta<1$, choose $J=\lceil1/\eta\rceil$ and use the
quantile grid of the continuous CDF $H$ at probabilities $j/J$,
$1\le j<J$, together with endpoints $-\infty,+\infty$.
Monotonicity implies that the supremum CDF error is no larger than
the maximum grid error plus $1/J$. A union bound and
\eqref{eq:spherical-cdf-second-moment} give
\begin{equation}
 \Pr_{\bm\theta}\left\{
  \|\widehat H_{\bm X}-H\|_\infty>2\eta\right\}
 \le C_K N^{-1/5}\eta^{-3}.
 \label{eq:spherical-cdf-probability}
\end{equation}
Integrating the smaller of this bound and one, with a splitting point
of order $N^{-1/15}$, proves
\eqref{eq:spherical-reference-consistency}. Common centering and
positive scaling leave the supremum CDF error unchanged.
\end{proof}

\begin{proof}[Quantiles and local power]
Conditional reference consistency implies consistency of its quantile
on the $s_0$ scale. Here one may use the spectral-subsequence argument
in the proof of Theorem~\ref{thm:spectral-calibration}: from any
covariance sequence, extract coordinatewise limits of the normalized
eigenvalues. Every resulting central limit law is a Gaussian plus a
sum of centered nonnegative chi-squared components, with variance one;
it is continuous and strictly increasing at its interior quantiles.
Thus convergence of the conditional CDF entails
\[
 (\widehat q_\alpha^{\mathrm{sph}}-q_\alpha^{\mathrm{NR}})/s_0
 \longrightarrow0
\]
in probability. The bounds above hold uniformly for $\kappa\le K$,
so the same argument applied to a hypothetical violating sequence of
signals proves uniform quantile consistency over that signal ball.
An empirical reference with $B\to\infty$ has vanishing additional
CDF error by the conditional DKW inequality. The plus-one adjustment
changes the tail threshold by at most $1/(B+1)$ and hence has no effect
on this conclusion.

Finally, bracket the random critical value between
$q_\alpha^{\mathrm{NR}}-\varepsilon s_0$ and
$q_\alpha^{\mathrm{NR}}+\varepsilon s_0$. The noncentral Gaussian
quadratic small-ball bound is uniform in $\bm\theta$ on this scale,
because its variance is $2A+4\bm\theta^\top\bm\Omega\bm\theta
\ge2A$. Consequently the power difference is bounded by the
quantile-error probability plus $C\varepsilon^{1/2}$.
First let $\min(n_1,n_2)$ and $B$ tend to infinity and then let
$\varepsilon\downarrow0$. This proves
\eqref{eq:spherical-local-power}, without requiring independence of
the statistic and its critical value.
\end{proof}

\paragraph{Finite auxiliary draw counts.}
\label{par:spherical-finite-B-power}
The power limits above let the number of auxiliary draws increase.
For a fixed $B$, define the ideal conditional upper-tail probability
\[
 p_\infty=\Pr\{R^{\mathrm{sph}}\ge T\mid\bm X\},
 \qquad k_B=\lfloor(B+1)\alpha\rfloor-1.
\]
Conditional on $\bm X$, the exceedance count in
\eqref{eq:spherical-monte-carlo-p} has distribution
$\operatorname{Bin}(B,p_\infty)$. Hence the actual rejection probability is
\begin{equation}
 \Pr_{\bm\theta}(p_B\le\alpha)
 =\E_{\bm\theta}\!\left[
   \sum_{j=0}^{k_B}\binom Bj
       p_\infty^j(1-p_\infty)^{B-j}\right],
 \label{eq:spherical-finite-B-power}
\end{equation}
with the sum interpreted as zero when $k_B<0$.
This identity holds for any sampling distribution and also permits ties.
For every $t>0$, the conditional Hoeffding inequality
\citep{Hoeffding1963} gives
\begin{align}
 &\left|\Pr_{\bm\theta}(p_B\le\alpha)
              -\Pr_{\bm\theta}(p_\infty\le\alpha)\right|
 \nonumber\\
 &\quad\le
 \Pr_{\bm\theta}\!\left\{
   |p_\infty-\alpha|\le t+\frac1{B+1}\right\}
       +2\exp(-2Bt^2).
 \label{eq:spherical-finite-B-deviation}
\end{align}
Indeed, if $\widehat p=B^{-1}\sum_{b=1}^B
\boldsymbol1\{R_b^{\mathrm{sph}}\ge T\}$, then
$|p_B-p_\infty|\le|\widehat p-p_\infty|+1/(B+1)$.
The two decisions can therefore differ only in the displayed boundary
band or on the event $|\widehat p-p_\infty|>t$.
The bound separates auxiliary simulation error from reference
approximation. Under alternatives, it does not assert a uniform
$B$-only error rate without control of the boundary-band probability.

\paragraph{Relation to classical Behrens--Fisher reduction.}
Suppose $n_1\le n_2$. The transformed vectors
\[
 \bm y_j=\bm Y_{1j}-\sqrt{n_1/n_2}
  (\bm Y_{2j}-\bar{\bm Y}_{2,1:n_1})-\bar{\bm Y}_2,
 \qquad 1\le j\le n_1,
\]
are independent Gaussian with mean $\bm\mu_1-\bm\mu_2$ and
covariance $\bm\Sigma_1+(n_1/n_2)\bm\Sigma_2$; this is the
Bennett reduction used by \citet{NishiyamaEtAl2013}.
An orthogonal mean--residual transform, multiplied by $\sqrt{w_1}$,
has the column law in \eqref{eq:spherical-columns}. Its Dempster
criterion is
\[
 T_D=\frac{m\|\bm X_0\|^2}{\sum_{r=1}^m\|\bm X_r\|^2}
     =\frac{mT}{\tr(\bm X\bm X^\top)-T}.
\]
Conditionally on the left Gram matrix, this is strictly increasing in
$T$ wherever its denominator is positive. Thus spherical calibration
is also an exact conditional rotation calibration of this classical
criterion. Neither the reduction nor the rotation principle is claimed
as a new testing principle here.

\paragraph{Relation to the Wang--Xu multiplier test.}
Let $c=n_1n_2/(n_1+n_2)$ and let $T_{\mathrm{CQ}}$ be the usual
Chen--Qin statistic. For every dataset, without Gaussian assumptions,
\[
 cT_{\mathrm{CQ}}=T-\widehat\mu,\qquad
 \widehat\mu=w_1\tr(\bm S_1)+w_2\tr(\bm S_2).
\]
Consequently, adding $\widehat\mu$ to a scaled Wang--Xu critical
value gives exactly the same rejection region as the original CQ-based
test. Theorem~2 of \citet{WangXu2022} allows both Gaussian and
Rademacher multipliers. A full-norm formulation with their Gaussian
multiplier reference is therefore an existing comparator, not a new
calibration method.

\subsection{Non-Gaussian validity of spherical calibration}
\label{app:nongaussian-spherical-proof}

We first give the two-output comparison needed to control the random
conditional CDF. A marginal Gaussian comparison alone does not establish
conditional calibration.

\begin{lem}[Two quadratic outputs]
\label{lem:two-output-quadratic-comparison}
For $i=1,2$ and $1\le j\le n_i$, let $\bm Z_{ij}$ be independent centered
vectors with covariance $\bm\Omega_i$, satisfying
\[
 \E(\bm v^\top\bm Z_{ij})^4
 \le K_0(\bm v^\top\bm\Omega_i\bm v)^2,
 \qquad
 \Var(\|\bm Z_{ij}\|^2)\le K_0\tr(\bm\Omega_i^2).
\]
Let $\bm\Omega=\bm\Omega_1+\bm\Omega_2$,
$A=\tr(\bm\Omega^2)>0$, and $s_0=\sqrt{2A}$.
For $k=1,2$, let $a_{k,ij}$ be deterministic coefficients satisfying
$\sum_j a_{k,ij}^2=1$ for each group, and let deterministic shifts satisfy
$\|\bm s_k\|^2/s_0\le K_1$. Define
\[
 Q_k=\left\|\bm s_k+\sum_{i,j}a_{k,ij}\bm Z_{ij}\right\|^2,
 \qquad
 Q_k^{\mathrm G}
 =\left\|\bm s_k+\sum_{i,j}a_{k,ij}\bm Z_{ij}^{\mathrm G}\right\|^2,
\]
where the $\bm Z_{ij}^{\mathrm G}\sim N(\bm0,\bm\Omega_i)$ are independent,
and the same Gaussian vector is used in both outputs. Put
\[
 D_4=\sum_{k=1}^2\sum_{i,j}a_{k,ij}^4,
 \qquad
 D_3=\sum_{k=1}^2\sum_{i,j}|a_{k,ij}|^3.
\]
Then
\begin{align}
 &\sup_{x_1,x_2\in\R}
 \left|\Pr(Q_1\le x_1,Q_2\le x_2)
             -\Pr(Q_1^{\mathrm G}\le x_1,Q_2^{\mathrm G}\le x_2)\right|
 \nonumber\\
 &\hspace{12mm}\le
 C_{K_0,K_1}\inf_{0<\varepsilon\le1}
 \{D_4\varepsilon^{-2}+D_3\varepsilon^{-3}+\varepsilon^{1/2}\}.
 \label{eq:two-output-quadratic-comparison}
\end{align}
The same bound holds for each marginal. No nonsingularity or lower bound
on the correlation between the two Gaussian outputs is required.
\end{lem}

\begin{proof}
Suppress the pair of observation indices and write $t$ for an observation,
with covariance $\bm\Gamma_t=\bm\Omega_i$ when $t$ belongs to group $i$.
For each $k$, center $Q_k$ by
$\|\bm s_k\|^2+\tr(\bm\Omega)$ and divide by $s_0$ to obtain
$W_k=U_k+R_k$, where
\begin{align*}
 U_k&=s_0^{-1}\left\{
       2\sum_{r<t}a_{k,r}a_{k,t}\bm Z_r^\top\bm Z_t
       +2\bm s_k^\top\sum_ta_{k,t}\bm Z_t\right\},\\
 R_k&=s_0^{-1}\sum_ta_{k,t}^2
       \{\|\bm Z_t\|^2-\tr(\bm\Gamma_t)\}.
\end{align*}
Define $U_k^{\mathrm G},R_k^{\mathrm G}$ in the same way. Since
$\bm\Gamma_t\preceq\bm\Omega$, we have
$\tr(\bm\Gamma_t^2)\le A$. Independence across observations gives
\begin{equation}
 \sum_{k=1}^2\{\E R_k^2+\E(R_k^{\mathrm G})^2\}
 \le (K_0+2)D_4/2.
 \label{eq:two-output-diagonal}
\end{equation}

Let $f:\R^2\to\R$ be bounded and three times continuously differentiable,
with all third partial derivatives bounded by $M_3$. Replace one
observation at a time by its Gaussian counterpart. Conditional on the
other observations at step $t$, the two-dimensional argument changes
only through the vector with coordinates
\[
 L_{k,t}=\frac{2a_{k,t}}{s_0}\bm\eta_{k,t}^\top\bm Z_t,
 \qquad
 \bm\eta_{k,t}=\bm s_k+\sum_{r\ne t}a_{k,r}\widetilde{\bm Z}_r,
\]
where each $\widetilde{\bm Z}_r$ is the original or Gaussian vector
according to the replacement order. The vector formed with
$\bm Z_t^{\mathrm G}$ has the same first and second conditional moments,
including the cross moment. Taylor expansion to second order therefore
bounds the replacement error by a constant times
\[
 M_3\sum_{k=1}^2
 \E\{|L_{k,t}|^3+|L_{k,t}^{\mathrm G}|^3\}.
\]
Here the elementary inequality
$(|u_1|+|u_2|)^3\le4(|u_1|^3+|u_2|^3)$ bounds mixed remainder terms.

For completeness, put $\bm V_{k,t}=\bm\eta_{k,t}-\bm s_k$.
The scalar fourth-moment expansion for independent centered summands
gives, with $K_* = \max(K_0,3)$,
\[
 \E(\bm v^\top\bm V_{k,t})^4
 \le K_*(\bm v^\top\bm\Omega\bm v)^2,
\]
because $\Var(\bm V_{k,t})\preceq\bm\Omega$.
Diagonalizing $\bm\Gamma_t$ and applying Cauchy--Schwarz term by term
then yields
\[
 \E(\bm V_{k,t}^\top\bm\Gamma_t\bm V_{k,t})^2
 \le K_*\{\tr(\bm\Omega\bm\Gamma_t)\}^2.
\]
The projection fourth-moment condition, followed by Lyapunov's inequality,
therefore implies
\begin{align*}
 \E|L_{k,t}|^3+\E|L_{k,t}^{\mathrm G}|^3
 &\le C_{K_0}|a_{k,t}|^3s_0^{-3}
       \{\tr(\bm\Omega\bm\Gamma_t)
                    +\bm s_k^\top\bm\Gamma_t\bm s_k\}^{3/2}\\
 &\le C_{K_0,K_1}|a_{k,t}|^3.
\end{align*}
The last inequality uses
$\tr(\bm\Omega\bm\Gamma_t)\le A$ and
$\bm s_k^\top\bm\Gamma_t\bm s_k
 \le\sqrt A\|\bm s_k\|^2\le\sqrt2K_1A$.
Summing the replacement errors gives
\begin{equation}
 |\E f(U_1,U_2)-\E f(U_1^{\mathrm G},U_2^{\mathrm G})|
 \le C_{K_0,K_1}M_3D_3.
 \label{eq:two-output-smooth-replacement}
\end{equation}

Choose a fixed smooth decreasing cutoff of the half-line and form the
product of two cutoffs with transition width $\varepsilon$.
All its third partial derivatives are bounded by $C\varepsilon^{-3}$.
Chebyshev's inequality and \eqref{eq:two-output-diagonal} remove the
diagonal remainders at total cost $C_{K_0}D_4\varepsilon^{-2}$.
The product cutoff and \eqref{eq:two-output-smooth-replacement} compare
the resulting rectangles, shifted in each coordinate by at most
$3\varepsilon$, at cost $C_{K_0,K_1}D_3\varepsilon^{-3}$.
Each Gaussian marginal is a noncentral quadratic with covariance
$\bm\Omega$ and variance at least $s_0^2$. The Gaussian quadratic
small-ball bound of Lemma~\ref{lem:gaussian-quadratic-small-ball}
therefore bounds the probability of each shifted boundary interval by
$C\varepsilon^{1/2}$. A union bound handles the two boundaries even if
the Gaussian outputs are dependent or identical. The lower bound follows
using lower shifted cutoffs. This proves
\eqref{eq:two-output-quadratic-comparison}; sending one threshold to
$+\infty$ proves the marginal assertion.
\end{proof}

\begin{proof}[Proof of Theorem~\ref{thm:nongaussian-spherical-calibration}]
Define the $n_i\times N$ deterministic matrices
\[
 \bm C_1=(n_1^{-1/2}\bm1,\bm H_1\bm e_1,\ldots,\bm H_1\bm e_m),
 \qquad
 \bm C_2=(-n_2^{-1/2}\bm1,\bm H_2\bm e_1,\ldots,\bm H_2\bm e_m).
\]
Both satisfy $\bm C_i^\top\bm C_i=\bm I_N$. Writing
$\bm c_{ij}^\top$ for row $j$ of $\bm C_i$ and
$\bm Z_{ij}=\sqrt{w_i}\bm\xi_{ij}$, we have exactly
\begin{equation}
 \bm X=\bm\theta\bm e_0^\top
              +\sum_{i,j}\bm Z_{ij}\bm c_{ij}^\top,
 \qquad \Var(\bm Z_{ij})=\bm\Omega_i=w_i\bm\Sigma_i.
 \label{eq:spherical-raw-coefficients}
\end{equation}
The two moment inequalities are preserved by this scalar weighting.

Let $\bm U,\bm V$ be independent uniform unit vectors in $\R^N$,
independent of the observations. Conditional on these vectors,
the reference pair $(\|\bm X\bm U\|^2,\|\bm X\bm V\|^2)$
has the form in Lemma~\ref{lem:two-output-quadratic-comparison}, with
\[
 a_{1,ij}=\bm c_{ij}^\top\bm U,\quad
 a_{2,ij}=\bm c_{ij}^\top\bm V,\quad
 \bm s_1=U_0\bm\theta,\quad \bm s_2=V_0\bm\theta.
\]
The orthonormality of the columns gives
$\sum_j a_{k,ij}^2=1$ for each $i,k$. Also,
$\|\bm s_k\|^2/s_0\le K_1$ whenever
$\|\bm\theta\|^2/s_0\le K_1$.

To control the coefficient sums, put
$\ell_{ij}=\|\bm c_{ij}\|^2$.
These are diagonal entries of an orthogonal projection, so
$0\le\ell_{ij}\le1$ and $\sum_j\ell_{ij}=N$.
The fourth moment of a uniform spherical projection gives
\begin{equation}
 \E(\bm c_{ij}^\top\bm U)^4
   =\frac{3\ell_{ij}^2}{N(N+2)},\qquad
 \E\sum_j a_{1,ij}^4\le\frac{3}{N+2}.
 \label{eq:sphere-coefficient-fourth}
\end{equation}
Moreover, Cauchy--Schwarz and $\sum_j a_{1,ij}^2=1$ imply
$\sum_j|a_{1,ij}|^3\le(\sum_j a_{1,ij}^4)^{1/2}$.
Applying Jensen's inequality and treating both groups and both draws
gives
\begin{equation}
 \E D_4\le\frac{12}{N+2},\qquad
 \E D_3\le4\sqrt{\frac{3}{N+2}}.
 \label{eq:sphere-coefficient-budgets}
\end{equation}

Replace all $\bm Z_{ij}$ by independent Gaussian vectors with the same
covariances, and call the resulting matrix $\bm X^{\mathrm G}$.
It is precisely the Gaussian augmented matrix of
Theorem~\ref{thm:spherical-calibration}: its columns are independent
$N(\bm0,\bm\Omega)$ vectors, with $\bm\theta$ added to the first.
Integrating \eqref{eq:two-output-quadratic-comparison} over
$\bm U,\bm V$ and taking $\varepsilon=N^{-1/7}$ shows that the two
unconditional reference pairs differ in rectangular CDF distance by
at most $C_{K_0,K_1}N^{-1/14}$. The same bound holds for one reference draw.
In particular, conditional independence of $\bm U$ and $\bm V$ gives
\begin{align}
 \sup_x|\E\widehat F_{\bm X}(x)
                 -\E\widehat F_{\bm X^{\mathrm G}}(x)|
 &\le C_{K_0,K_1}N^{-1/14},\nonumber\\
 \sup_x|\E\widehat F_{\bm X}(x)^2
                 -\E\widehat F_{\bm X^{\mathrm G}}(x)^2|
 &\le C_{K_0,K_1}N^{-1/14}.
 \label{eq:spherical-two-moment-comparison}
\end{align}
This step uses conditional independence of the auxiliary draws, not
independence of the non-Gaussian residual columns.

The Gaussian bound \eqref{eq:spherical-cdf-second-moment}, together with
\eqref{eq:spherical-two-moment-comparison}, now yields
\begin{equation}
 \sup_x\E\{\widehat F_{\bm X}(x)-F_{\bm\Omega}(x)\}^2
 \le C_{K_0,K_1}(N^{-1/5}+N^{-1/14})
 \le C'_{K_0,K_1}N^{-1/14}.
 \label{eq:nongaussian-spherical-pointwise}
\end{equation}
Use the quantile grid of the continuous CDF $F_{\bm\Omega}$ at mesh
at most $\eta$. Monotonicity, a union bound over at most
$1+\eta^{-1}$ grid points, and Chebyshev's inequality imply
\[
 \Pr\{\|\widehat F_{\bm X}-F_{\bm\Omega}\|_\infty>2\eta\}
 \le C_{K_0,K_1}N^{-1/14}\eta^{-3},\qquad 0<\eta<1.
\]
Integrating the smaller of one and this bound, splitting at a constant
multiple of $N^{-1/42}$, proves \eqref{eq:nongaussian-spherical-cdf}.

For the power statement, the conditional CDF bound and the conditional
DKW inequality for $B\to\infty$ give reference consistency for the
empirical auxiliary distribution. The spectral-subsequence argument
in the proof of Theorem~\ref{thm:spherical-calibration} then gives
consistency of its $(1-\alpha)$-quantile on the $s_0$ scale, uniformly
over the stated signal ball. The plus-one adjustment is asymptotically
negligible. If an empirical reference has ties, their maximal conditional
mass tends to zero because its CDF converges to the continuous
$F_{\bm\Omega}$; hence the weak inequality in the definition of $p_B$
causes no additional limiting term.

For any fixed $\delta>0$, bracket the random critical value between
$q_\alpha^{\mathrm{NR}}-\delta s_0$ and
$q_\alpha^{\mathrm{NR}}+\delta s_0$.
Corollary~\ref{cor:grouped-transfer} compares the full-sample statistic
with $\|\bm G+\bm\theta\|^2$, uniformly in $\bm\theta$, at a vanishing
Kolmogorov error. The Gaussian quadratic small-ball bound controls
the probability between the two brackets by $C\delta^{1/2}$.
First let $N,B\to\infty$, and then let $\delta\downarrow0$.
This proves \eqref{eq:nongaussian-spherical-power}. No independence of
the original statistic and its estimated critical value is used.
\end{proof}

\subsection{Moment-based calibration and its conditions}
\label{app:moment-calibration-summary}
The statistic is calibrated under the null hypothesis. The noncentral
reference describes its power under specified alternatives. Write
$\mu_n=\tr(\bm\Omega_n)$, $A_n=\tr(\bm\Omega_n^2)>0$ and
$s_{0,n}=\sqrt{2A_n}$. We use an unbiased estimator $\widehat\mu_n$
of the trace, a nonnegative unbiased estimator $\widehat A_n$ of the
trace of the squared covariance matrix, and empirical covariance eigenvalues
$\widehat\lambda_{n,j}$. Gaussian orthogonal residuals provide
$m_i=n_i-1$ independent covariance observations within group $i$;
for non-Gaussian samples we use $m_i=\lfloor n_i/2\rfloor$ disjoint
paired differences. The test statistic still uses the full sample means.
The estimators and their Gram-matrix implementation are given in
Appendix~\ref{app:implementation-details}.

The threshold in \eqref{eq:main-threshold} selects $\widehat r_n$
leading empirical eigenvalues. The projection in
\eqref{eq:main-projection} rescales them to
$\widetilde\lambda_{n,j}$ when their squared sum exceeds
$\widehat A_n$. The remaining variance is
$2\widetilde B_n=2(\widehat A_n-
\sum_{j\le\widehat r_n}\widetilde\lambda_{n,j}^2)\ge0$.
The resulting reference has the form
\[
 \widehat\mu_n+
 \sum_{j\le\widehat r_n}\widetilde\lambda_{n,j}(Z_j^2-1)+H_n^M,
 \qquad M\in\{G,\Gamma\},
\]
where the independent remainder $H_n^M$ is either Gaussian or a
centered gamma variable with variance $2\widetilde B_n$.
The precise gamma definition, including its boundary cases, is
\eqref{eq:main-gamma-remainder}. The test rejects when $T$ exceeds
the conditional $(1-\alpha)$ quantile of this reference.
Projection restores the estimated first two moments but does not
correct every spectral-shape error. If $\widetilde B_n=0$, the
two references coincide. Reference draws are auxiliary computations;
their budgets are separate from the evaluation replications.

The sufficient conditions for reference estimation are stronger than
those for the Gaussian comparison. For the Gaussian construction, it is
enough that
\begin{equation}
 m_{\min}\longrightarrow\infty,\qquad
 \frac{\mu_n}{m_{\min}\sqrt{A_n}}\longrightarrow0.
 \label{eq:main-effective-dimension}
\end{equation}
For the paired construction, impose the same growth conditions and,
with $\bm{U}_{ij}=\sqrt{w_i}(\bm{Y}_{ij}-\bm{\mu}_i)$, assume
\[
 \|\bm{u}^\top\bm{U}_{ij}\|_{\psi_2}
       \le K(\bm{u}^\top\bm{\Omega}_i\bm{u})^{1/2},\qquad
 \operatorname{Var}(\|\bm{U}_{ij}\|^2)
       \le K\operatorname{tr}(\bm{\Omega}_i^2)
\]
uniformly in $n$, over both samples and all deterministic $\bm{u}$, where
$\|U\|_{\psi_2}=\inf\{c>0:E\exp(U^2/c^2)\le2\}$.
Proposition~\ref{prop:gaussian-calibration-verification} and
Corollary~\ref{cor:paired-calibration} verify the requisite trace and
operator-norm estimation conditions. Theorems~\ref{thm:spectral-calibration},
\ref{thm:gamma-quantitative-calibration} and
\ref{thm:projected-calibration} then establish consistency of the estimated
reference distributions and critical values without a fixed limiting
spectrum. For example, a spike-plus-bulk spectrum
$\bm{\Omega}_p=\operatorname{diag}(\sqrt p,1,\ldots,1)$ requires
$\sqrt p/m_{\min}\to0$ for this construction, although the comparison
result itself imposes no such dimension restriction.

The following conclusion applies to either projected critical value
above, denoted generically by $\widehat q_{n,\alpha}$, and also to the
unprojected references defined in Appendix~\ref{app:calibration}.
Here $T_{\bm{\theta}}=T$ and $T_{\bm{\theta}}^{\mathrm G}=T^*$ emphasize
the dependence on the normalized mean difference.
\begin{cor}[Uniform feasible power approximation]\label{cor:feasible-final}
Fix \(\alpha\in(0,1)\). Suppose the two-sample model satisfies the
conditions of Corollary~\ref{cor:grouped-transfer}, and use either the
Gaussian residual construction of
Proposition~\ref{prop:gaussian-calibration-verification} or the independent
pair construction of Corollary~\ref{cor:paired-calibration}, under its
stated conditions. For fixed centered sample distributions and covariance
matrices, let \(\bm{\theta}\) range over all deterministic normalized mean
differences. Then
\[
 \sup_{\bm{\theta}}\left|
 \Pr_{\bm{\theta}}\{T_{\bm{\theta}}>\widehat q_{n,\alpha}\}
 -\Pr\{T_{\bm{\theta}}^{\mathrm G}>q_{n,\alpha}^{\mathrm{NR}}\}
 \right|\longrightarrow0.
\]
In particular, the feasible test has asymptotic level \(\alpha\).
The assertion compares finite-sequence power functions; it does not
require either one to have a limit under every alternative sequence.
\end{cor}
\begin{proof}See Appendix~\ref{app:calibration}.\end{proof}

\subsection{Implementation details}\label{app:implementation-details}
The null Gaussian reference depends on the unknown covariance matrices.
We estimate its leading eigenvalues and approximate the remaining
components by a Gaussian or centered gamma variable. The critical value
is computed under the null; the noncentral reference will be used to
describe power in Section~\ref{prof:section}.

Write $\mu_n=\operatorname{tr}(\bm{\Omega}_n)$,
$A_n=\operatorname{tr}(\bm{\Omega}_n^2)>0$ and
$s_{0,n}=\sqrt{2A_n}$. Set $w_1=n_2/(n_1+n_2)$ and
$w_2=n_1/(n_1+n_2)$, so $\bm{\Omega}_i=w_i\bm{\Sigma}_i$.
For Gaussian samples, take $m_i=n_i-1\ge2$ and form
\[
 \bm{V}_{ir}=\sqrt{w_i}\,\bm{Y}_i\bm{H}_i\bm{e}_r,
 \qquad r=1,\ldots,m_i,\quad i=1,2,
\]
where $\bm{Y}_i=(\bm{Y}_{i1},\ldots,\bm{Y}_{in_i})$ and
$\bm{H}_i$ is an $n_i\times m_i$ matrix with
$\bm{H}_i^\top\bm{H}_i=\bm{I}_{m_i}$ and
$\bm{H}_i^\top\mathbf1=\bm{0}$. A Helmert matrix is one choice, and $\bm{e}_r$ here denotes the $r$th standard basis vector of $\mathbb R^{m_i}$.
These residuals are independent with covariance $\bm{\Omega}_i$.
For these projected moment estimators under non-Gaussian sampling,
independence is obtained from disjoint pairs. This requirement does not
apply to the augmented spherical reference. Take $m_i=\lfloor n_i/2\rfloor\ge2$ and use
\[
 \bm{V}_{ir}=\sqrt{w_i/2}
       (\bm{Y}_{i,2r}-\bm{Y}_{i,2r-1}),\qquad r=1,\ldots,m_i.
\]
Disjoint pairs are independent and have the same covariance
$\bm{\Omega}_i$. Pairing is used only for reference estimation; the
test statistic continues to use both full-sample means.

For either construction, compute
\begin{align}
 \widehat{\bm{\Omega}}_n
  &=\sum_{i=1}^2\frac1{m_i}
                  \sum_{r=1}^{m_i}\bm{V}_{ir}\bm{V}_{ir}^\top,
 &\widehat\mu_n&=\operatorname{tr}(\widehat{\bm{\Omega}}_n),
 \label{eq:main-reference-estimates}\\
 \widehat A_n
  &=\sum_{i=1}^2
       \frac{\sum_{r\ne s}(\bm{V}_{ir}^\top\bm{V}_{is})^2}
            {m_i(m_i-1)}
    +\frac{2}{m_1m_2}
       \sum_{r=1}^{m_1}\sum_{s=1}^{m_2}
                    (\bm{V}_{1r}^\top\bm{V}_{2s})^2.
 \label{eq:main-squared-trace}
\end{align}
The estimator $\widehat A_n$ is nonnegative and unbiased. In the Gaussian
case, $\widehat{\bm{\Omega}}_n=w_1\bm{S}_1+w_2\bm{S}_2$, where
$\bm{S}_i$ is the usual unbiased sample covariance matrix. A Gram matrix
of the columns $\bm{V}_{ir}/\sqrt{m_i}$ yields the nonzero eigenvalues
$\widehat\lambda_{n,j}$ without covariance inversion. Order them nonincreasingly and append zeros when needed.

Put $m_{\min}=\min(m_1,m_2)$. On $\widehat A_n>0$, use
\begin{equation}
 \widehat a_n=\frac{\widehat\mu_n}{m_{\min}\sqrt{\widehat A_n}},
 \qquad
 \tau_n=\min\left\{\frac12,
                 \widehat a_n^{1/4}+m_{\min}^{-1/8}\right\},
 \qquad
 \widehat r_n=\#\{j:\widehat\lambda_{n,j}>
                         \tau_n\sqrt{\widehat A_n}\}.
 \label{eq:main-threshold}
\end{equation}
Retaining empirical eigenvalues can assign more variance to the leading
components than the estimated total $2\widehat A_n$. We therefore put
$\widehat S_n=\sum_{j\le\widehat r_n}\widehat\lambda_{n,j}^2$ and define
\begin{align}
 \widehat c_n&=
 \begin{cases}
  \min\{1,\sqrt{\widehat A_n/\widehat S_n}\},&\widehat S_n>0,\\
  1,&\widehat S_n=0,
 \end{cases}
 &\widetilde\lambda_{n,j}&=\widehat c_n\widehat\lambda_{n,j},
 \quad j\le\widehat r_n,\nonumber\\
 \widetilde B_n&=\widehat A_n-
                     \sum_{j\le\widehat r_n}\widetilde\lambda_{n,j}^2,
 &\widetilde L_n&=\widehat\mu_n-
                     \sum_{j\le\widehat r_n}\widetilde\lambda_{n,j}.
 \label{eq:main-projection}
\end{align}
Here $\widetilde B_n\ge0$; in computation it can be truncated at zero to
remove roundoff error. The projected Gaussian and gamma references are
\begin{equation}
 \widetilde R_n^M=\widehat\mu_n+
  \sum_{j\le\widehat r_n}\widetilde\lambda_{n,j}(Z_j^2-1)
       +H_n^M,\qquad M\in\{G,\Gamma\},
 \label{eq:main-projected-reference}
\end{equation}
where $H_n^G=\sqrt{2\widetilde B_n}\,G$ and
\begin{equation}
 H_n^\Gamma=
 \begin{cases}
  0,&\widetilde B_n=0,\\[2pt]
  \displaystyle\frac{\widetilde B_n}{\widetilde L_n}
       (\chi^2_{\widetilde d_n}-\widetilde d_n),\quad
       \widetilde d_n=\frac{\widetilde L_n^2}{\widetilde B_n},
       &\widetilde B_n>0,\ \widetilde L_n>0,\\[5pt]
  \sqrt{2\widetilde B_n}\,G,
       &\widetilde B_n>0,\ \widetilde L_n\le0.
 \end{cases}
 \label{eq:main-gamma-remainder}
\end{equation}
Here the $Z_j$ and $G$ are standard normal variables. All auxiliary variables are mutually independent and independent of
the data; chi-squared degrees of freedom need not be integers.
Both references have conditional mean $\widehat\mu_n$ and variance
$2\widehat A_n$. If $\widehat S_n>\widehat A_n$, the remainder variance
is zero and the two references coincide. A fixed reference can be assigned
when $\widehat A_n\le0$, an asymptotically negligible event under the
conditions below.

Let $\widetilde q_{n,\alpha}^M$ be the conditional $(1-\alpha)$ quantile
of \eqref{eq:main-projected-reference}. The test rejects $H_0$ when
\[
 T>\widetilde q_{n,\alpha}^M.
\]
The quantile can be computed by independent draws from the reference.
An auxiliary sample size tending to infinity preserves its asymptotic
consistency; Section~\ref{sec:numerics} specifies the finite draw counts.
The projection corrects the estimated total variance but does not by
itself ensure accurate spectral shape. We retain the unprojected
references as diagnostics in the simulations.

\subsection{Reference estimation and calibration}
Write $A_n=\operatorname{tr}(\bm{\Omega}_n^2)>0$,
$\mu_n=\operatorname{tr}(\bm{\Omega}_n)$ and $s_{0,n}=\sqrt{2A_n}$.
Let $\lambda_{n,1}\ge\lambda_{n,2}\ge\cdots\ge0$ be the eigenvalues of
$\bm{\Omega}_n$, extended by zeros. Suppose that a positive semidefinite estimator
$\widehat{\bm{\Omega}}_n$, a trace estimator $\widehat\mu_n$, and an estimator
$\widehat A_n$ are available. The trace estimator need not equal
$\operatorname{tr}(\widehat{\bm{\Omega}}_n)$ in the generic result. Denote the ordered
eigenvalues of $\widehat{\bm{\Omega}}_n$ by $\widehat\lambda_{n,j}$.

Choose a positive threshold $\tau_n$, which may be data dependent, and set
\begin{align}
 \widehat r_n&=\#\{j:\widehat\lambda_{n,j}>\tau_n\sqrt{\widehat A_n}\},
 \label{eq:spectral-rank}\\
 \widehat B_n&=\left(\widehat A_n-
                  \sum_{j\le\widehat r_n}\widehat\lambda_{n,j}^2\right)_+.
 \label{eq:spectral-tail}
\end{align}
On the event $\widehat A_n>0$, the proposed conditional reference is
\begin{equation}
 \widehat R_n=\widehat\mu_n+
   \sum_{j\le\widehat r_n}\widehat\lambda_{n,j}(Z_j^2-1)
   +\sqrt{2\widehat B_n}\,G,
 \label{eq:spectral-reference}
\end{equation}
where the auxiliary variables $G,Z_1,Z_2,\ldots$ are independent standard
normals, independent of the data. Its upper $\alpha$ critical value is the
conditional $(1-\alpha)$ quantile, denoted by $\widehat q_{n,\alpha}$.
An arbitrary fixed reference may be assigned on $\{\widehat A_n\le0\}$;
the assumptions below make that event asymptotically negligible. The positive
part in \eqref{eq:spectral-tail} is a finite-sample safeguard and does not
replace the consistency conditions.

\begin{theorem}[Estimated leading-spectrum calibration]
\label{thm:spectral-calibration}
Suppose
\begin{equation}
 \frac{\widehat\mu_n-\mu_n}{\sqrt{A_n}}\xrightarrow{p}0,
 \qquad \frac{\widehat A_n}{A_n}\xrightarrow{p}1,
 \qquad \tau_n\xrightarrow{p}0,
 \qquad
 \frac{\|\widehat{\bm{\Omega}}_n-\bm{\Omega}_n\|_{\mathrm{op}}}
      {\tau_n\sqrt{A_n}}\xrightarrow{p}0.
 \label{eq:spectral-estimation}
\end{equation}
Then, with $T^*_{0,n}=\sum_j\lambda_{n,j}Z_j^2$,
\begin{equation}
 \sup_x\left|
   \Pr(\widehat R_n\le x\mid\mathcal Y_n)
      -\Pr(T^*_{0,n}\le x)\right|\xrightarrow{p}0.
 \label{eq:spectral-conditional-consistency}
\end{equation}
No eigenvalue gap and no fixed number of spikes are required.
For every fixed $\alpha\in(0,1)$, without assuming a limiting spectrum,
\begin{equation}
 \frac{\widehat q_{n,\alpha}-q^{\mathrm{NR}}_{n,\alpha}}
      {s_{0,n}}\xrightarrow{p}0.
 \label{eq:spectral-quantile-consistency}
\end{equation}
If additionally $(T^*_{0,n}-\mu_n)/s_{0,n}\Rightarrow\zeta$, then
\begin{equation}
 \frac{\widehat q_{n,\alpha}-\mu_n}{s_{0,n}}
      \xrightarrow{p}c_{\zeta,\alpha}.
 \label{eq:spectral-quantile-profile}
\end{equation}
Here $q^{\mathrm{NR}}_{n,\alpha}$ is the exact Gaussian null quantile.
\end{theorem}

Theorem~\ref{thm:spectral-calibration} requires the operator-norm error
to be smaller than the spectral threshold. It does not require a gap
between retained and omitted population eigenvalues. In particular,
covariance consistency stated in another metric does not by itself verify
\eqref{eq:spectral-estimation}. The theorem also does not justify using
all empirical eigenvalues with an uncorrected squared trace.

\begin{prop}[Random critical values and power]
\label{prop:spectral-random-power}
Let $T_n=\|\bm{\theta}_n+\bm{W}_n\|^2$ and
$T^*_n=\|\bm{\theta}_n+\bm{G}_n\|^2$, with $\bm{G}_n\sim N(\bm{0},\bm{\Omega}_n)$.
For $0<\epsilon\le1$, any data-dependent critical value $\widehat q_n$ and
any deterministic $q_n$ satisfy
\begin{align}
 &\left|\Pr_{\bm{\theta}_n}(T_n>\widehat q_n)
                -\Pr(T^*_n>q_n)\right|\nonumber\\
 &\quad\le d_{\mathrm K}(T_n,T_n^*)
       +\Pr_{\bm{\theta}_n}(|\widehat q_n-q_n|>\epsilon s_{0,n})
       +C\epsilon^{1/2}.
 \label{eq:random-critical-power-bound}
\end{align}
Consequently, if Gaussian comparison and
\eqref{eq:spectral-quantile-consistency} hold under the alternative,
rejection with $\widehat q_{n,\alpha}$ has the same limiting power as
rejection with $q^{\mathrm{NR}}_{n,\alpha}$. The directional power limits
of Section~\ref{prof:section} therefore also apply to the feasible test. No independence
between the statistic and its estimated critical value is required.
Uniform comparison and quantile consistency over a signal class give
the corresponding uniform power conclusion.
\end{prop}

We first verify the estimation conditions for Gaussian observations.

Suppose that the two samples are Gaussian. Write
\[
 w_1=\frac{n_2}{n_1+n_2},\qquad
 w_2=\frac{n_1}{n_1+n_2},\qquad
 \bm{\Omega}_i=w_i\bm{\Sigma}_i,\quad i=1,2,
\]
so that $\bm{\Omega}_n=\bm{\Omega}_1+\bm{\Omega}_2$. Put
$\nu_i=n_i-1$ and $\nu_{\min}=\min(\nu_1,\nu_2)$.
The weights $w_i$ are nonnegative covariance weights; the minus sign in
the sample-mean difference does not enter its covariance estimator.
Choose a deterministic $n_i\times\nu_i$ matrix $\bm{H}_i$ such that
$\bm{H}_i^\top \bm{H}_i=\bm{I}_{\nu_i}$ and $\bm{H}_i^\top\mathbf1=\bm{0}$; a Helmert matrix is one
choice. With $\bm{Y}_i=(\bm{Y}_{i1},\ldots,\bm{Y}_{in_i})$, define its residual contrasts
\begin{equation}
 \bm{V}_{ir}=\sqrt{w_i}\,\bm{Y}_i\bm{H}_i\bm{e}_r,
 \qquad i=1,2,\quad 1\le r\le\nu_i.
 \label{eq:helmert-residuals}
\end{equation}
They are mutually independent with $\bm{V}_{ir}\sim N(\bm{0},\bm{\Omega}_i)$.
Set
\begin{align}
 \widehat{\bm{\Omega}}_i&=\frac1{\nu_i}\sum_{r=1}^{\nu_i}\bm{V}_{ir}\bm{V}_{ir}^\top,
 &\widehat{\bm{\Omega}}_n&=\widehat{\bm{\Omega}}_1+\widehat{\bm{\Omega}}_2,
 &\widehat\mu_n&=\operatorname{tr}(\widehat{\bm{\Omega}}_n),
 \label{eq:residual-covariance}\\
 \widehat A_n&=
   \sum_{i=1}^2\frac{\sum_{r\ne s}(\bm{V}_{ir}^\top \bm{V}_{is})^2}{\nu_i(\nu_i-1)}
   +\frac{2}{\nu_1\nu_2}
     \sum_{r=1}^{\nu_1}\sum_{s=1}^{\nu_2}(\bm{V}_{1r}^\top\bm{V}_{2s})^2.
 \label{eq:residual-unbiased-A}
\end{align}
For Gaussian samples, $\widehat{\bm{\Omega}}_i=w_i\bm{S}_i$, where
$\bm{S}_i$ is the usual unbiased sample covariance matrix.
The first sum requires $n_1,n_2\ge3$. The estimator
\eqref{eq:residual-unbiased-A} is nonnegative and unbiased. Its inner-product
form also allows both trace estimates and the nonzero eigenvalues of
$\widehat{\bm{\Omega}}_n$ to be computed from a residual Gram matrix.

\begin{prop}[Gaussian primitive verification]
\label{prop:gaussian-calibration-verification}
For the construction above,
\begin{equation}
 \operatorname{Var}(\widehat\mu_n)\le\frac{2A_n}{\nu_{\min}},
 \qquad E\widehat A_n=A_n,
 \qquad
 \operatorname{Var}(\widehat A_n)\le \frac{C A_n^2}{\nu_{\min}}.
 \label{eq:trace-estimator-bounds}
\end{equation}
For the operator-norm verification, suppose
\begin{equation}
 \nu_{\min}\longrightarrow\infty,
 \qquad
 a_n:=\frac{\mu_n}{\nu_{\min}\sqrt{A_n}}\longrightarrow0.
 \label{eq:gaussian-calibration-effective-dimension}
\end{equation}
Then
\begin{equation}
 \frac{\|\widehat{\bm{\Omega}}_n-\bm{\Omega}_n\|_{\mathrm{op}}}{\sqrt{A_n}}
       =O_p(\sqrt{a_n}+a_n).
 \label{eq:gaussian-calibration-op-rate}
\end{equation}
In particular, the following data-dependent threshold satisfies
\eqref{eq:spectral-estimation}:
\begin{equation}
 \widehat a_n=\frac{\widehat\mu_n}{\nu_{\min}\sqrt{\widehat A_n}},
 \qquad
 \tau_n=\min\left\{\frac12,
          \widehat a_n^{1/4}+\nu_{\min}^{-1/8}\right\},
 \label{eq:implemented-spectral-threshold}
\end{equation}
on $\widehat A_n>0$. Thus Theorem~\ref{thm:spectral-calibration} gives a
feasible calibration for the Gaussian two-sample problem, including spiked
and nonnormal reference profiles.
\end{prop}

For $\bm{\Omega}_p=\operatorname{diag}(\sqrt p,1,\ldots,1)$,
condition \eqref{eq:gaussian-calibration-effective-dimension} is equivalent
to $\sqrt p/\nu_{\min}\to0$. Thus the dimension may exceed the sample size,
but its growth is restricted by spectral estimation. This restriction is
not needed for the Gaussian comparison result.

For non-Gaussian samples, orthogonal residual contrasts are generally
dependent. For the projected spectral construction considered next, we
use disjoint pairs to estimate its parameters, retaining the full-sample
mean statistic. The spherical moment theorem does not require this
replacement.

\begin{cor}[Independent-pair calibration under sub-Gaussian sampling]
\label{cor:paired-calibration}
Suppose that the two samples are independent and each consists of
independent, identically distributed observations. Put
$\bm{U}_{ij}=\sqrt{w_i}(\bm{Y}_{ij}-\bm{\mu}_i)$, $i=1,2$, so
$\operatorname{Cov}(\bm{U}_{ij})=\bm{\Omega}_i$. Assume, uniformly in $i,n$ and
deterministic $\bm{u}$,
\begin{equation}
 \|\bm{u}^\top \bm{U}_{ij}\|_{\psi_2}
       \le K(\bm{u}^\top\bm{\Omega}_i \bm{u})^{1/2},
 \qquad
 \operatorname{Var}(\|\bm{U}_{ij}\|^2)
       \le K\operatorname{tr}(\bm{\Omega}_i^2).
 \label{eq:paired-primitive-conditions}
\end{equation}
Here $\|U\|_{\psi_2}=\inf\{c>0:E\exp(U^2/c^2)\le2\}$.
For $m_i=\lfloor n_i/2\rfloor\ge2$, use independent paired differences
\begin{equation}
 \bm{V}_{ir}=\sqrt{\frac{w_i}{2}}(\bm{Y}_{i,2r}-\bm{Y}_{i,2r-1}),
 \qquad i=1,2,\quad r=1,\ldots,m_i,
 \label{eq:paired-residuals}
\end{equation}
and use \eqref{eq:residual-covariance}--\eqref{eq:residual-unbiased-A} with
$\nu_i$ replaced by $m_i$. Let $m_{\min}=\min(m_1,m_2)$ and suppose
\[
 m_{\min}\to\infty,
 \qquad \frac{\mu_n}{m_{\min}\sqrt{A_n}}\to0.
\]
Then all conclusions of Proposition~\ref{prop:gaussian-calibration-verification}
hold with constants depending on $K$, with $m_i$ in place of $\nu_i$, and
with the coefficient $2$ in the trace-mean variance bound replaced by a
constant depending on $K$. The threshold in
\eqref{eq:implemented-spectral-threshold}, with $m_{\min}$ substituted,
therefore gives a feasible calibration for this non-Gaussian subclass.
\end{cor}

Conditions \eqref{eq:paired-primitive-conditions} hold for deterministic
linear transforms of independent, standardized components with uniformly
bounded sub-Gaussian norms. The projection bound follows from their moment
generating functions, and the norm-variance bound follows from the
quadratic-form variance identity. This class includes Gaussian and bounded
components. Student and centered chi-squared components can satisfy the
transfer conditions, but are not covered by this calibration corollary.
Pairing is used only for reference estimation; the dependence it induces
with the full-sample statistic is allowed by
Proposition~\ref{prop:spectral-random-power}.

Both projected moment constructions are invariant to separate translations of the two samples.
For fixed centered sample distributions, the law of the estimated reference
therefore does not depend on either population mean. Critical-value consistency
is uniform over these translations. Together with the two-sample comparison
bound, this gives uniform feasible power approximation through
\eqref{eq:random-critical-power-bound}. Table~\ref{tab:calibration-conditions}
summarizes the different sufficient conditions.

\begin{table}[!htbp]
\centering\small
\begin{threeparttable}
\caption{Sufficient conditions for comparison and reference estimation.}
\label{tab:calibration-conditions}
\begin{tabular}{@{}>{\raggedright\arraybackslash}p{.22\textwidth}
                  >{\raggedright\arraybackslash}p{.46\textwidth}
                  >{\raggedright\arraybackslash}p{.26\textwidth}@{}}
\toprule
Result & Sufficient conditions & Scope\\
\midrule
Two-sample comparison
& Common projection fourth-moment and norm-variance bounds;
  $n_{\min}\to\infty$.
& No fixed sample-size ratio required.\\[4pt]
Gaussian spherical calibration
& Gaussian samples; $n_i\ge2$, $A>0$; $B\ge1$.
& Finite-sample level; arbitrary dimension and covariances.\\[4pt]
Moment-class spherical calibration
& Both conditions in \eqref{eq:raw-group-moments}; $n_{\min},B\to\infty$.
& Asymptotic level and bounded-local reference power; no dimension ratio.\\[4pt]
Gaussian projected estimation
& Gaussian samples; $m_{\min}\to\infty$ and
  $\mu_n/(m_{\min}\sqrt{A_n})\to0$.
& Helmert residuals, $m_i=n_i-1$.\\[4pt]
Paired projected estimation
& Observations i.i.d. within each sample, satisfying
  \eqref{eq:paired-primitive-conditions}, and the same
  residual-count conditions.
& Disjoint pairs, $m_i=\lfloor n_i/2\rfloor$.\\
\bottomrule
\end{tabular}
\begin{tablenotes}[flushleft]
\footnotesize
\item[] \textit{Note.} The results concern the two-sample problem. The spherical moment row includes reference and local-power validity. For projected estimators, combine the estimation row with the comparison conditions. Here $m_{\min}$ is the smallest residual count.
\end{tablenotes}
\end{threeparttable}
\end{table}

A deterministic threshold such as $\tau_n=n_{\min}^{-1/8}$ is also allowed
by Theorem~\ref{thm:spectral-calibration}, provided its displayed estimation
condition is checked. It is not automatically justified by
\eqref{eq:gaussian-calibration-effective-dimension}. In the spike-plus-bulk
example, $p=o(n_{\min}^{3/2})$ is a sufficient stronger condition for this
particular deterministic threshold. The data-dependent rule
\eqref{eq:implemented-spectral-threshold} only requires the weaker
effective-dimension condition above.

The cap in \eqref{eq:implemented-spectral-threshold} matters in finite
samples. If the smallest residual count is at most $256$, its inverse
eighth root is at least $1/2$, and the rule gives $\tau_n=1/2$ for every
dataset with $\widehat A_n>0$. All estimated-reference settings reported
in Section~\ref{sec:numerics} have residual counts in this range. Their
results therefore evaluate the capped rule, rather than variation of the
threshold with the data. The rule is asymptotically valid under the stated
conditions; its finite-sample sensitivity remains a separate question.

The omitted chi-squared terms have positive skewness, whereas the
Gaussian remainder has zero skewness. We consider a centered gamma
remainder determined by the remaining first and second spectral moments.
It uses the same eigenvalue and trace estimates and does not require a
cubic-trace estimator. The following results establish its validity under
the same estimation conditions.

In addition to $\widehat B_n$ in \eqref{eq:spectral-tail}, define
\begin{equation}
 \widehat L_n=\widehat\mu_n-
                 \sum_{j\le\widehat r_n}\widehat\lambda_{n,j}.
 \label{eq:gamma-tail-trace}
\end{equation}
All auxiliary variables below are independent of the data and of the
retained chi-squared coordinates. On $\widehat A_n>0$, define
\begin{equation}
 \widehat H_n=
 \begin{cases}
  0,&\widehat B_n=0,\\[2pt]
  \displaystyle
  \frac{\widehat B_n}{\widehat L_n}
   \left(\chi^2_{\widehat d_n}-\widehat d_n\right),
    \quad \widehat d_n=\frac{\widehat L_n^2}{\widehat B_n},
      &\widehat B_n>0,\ \widehat L_n>0,\\[5pt]
  \sqrt{2\widehat B_n}\,G,
      &\widehat B_n>0,\ \widehat L_n\le0.
 \end{cases}
 \label{eq:gamma-tail-rule}
\end{equation}
The gamma reference and its critical value are
\begin{equation}
 \widehat R_n^{\Gamma}=\widehat\mu_n+
       \sum_{j\le\widehat r_n}\widehat\lambda_{n,j}(Z_j^2-1)
          +\widehat H_n,
 \qquad
 \widehat q_{n,\alpha}^{\Gamma}
       =\inf\{x:\Pr(\widehat R_n^{\Gamma}\le x\mid\mathcal Y_n)
                         \ge1-\alpha\}.
 \label{eq:gamma-reference}
\end{equation}
Assign an arbitrary reference if $\widehat A_n\le0$. The Gaussian
fallback in \eqref{eq:gamma-tail-rule} handles a nonpositive estimated tail
trace. When no eigenvalue is retained and $\widehat\mu_n>0$, the reference
reduces to the ordinary two-cumulant calibration:
\begin{equation}
 \widehat R_n^{\Gamma}
       =\frac{\widehat A_n}{\widehat\mu_n}
          \chi^2_{\widehat\mu_n^2/\widehat A_n}.
 \label{eq:gamma-no-spike-two-cumulant}
\end{equation}
For the primitive covariance constructions, $\widehat\mu_n$ is a trace of
a positive semidefinite matrix, so this reduction applies whenever
$\widehat A_n>0$ and no eigenvalue is retained.

\begin{prop}[Population moment comparison]
\label{prop:gamma-moment-comparison}
For any nonempty finite set of positive remaining eigenvalues, write
$L=\sum_j\lambda_j$, $B=\sum_j\lambda_j^2$ and $C_3=\sum_j\lambda_j^3$.
The centered gamma variable
\[
 H_\Gamma=\frac BL\left(\chi^2_{L^2/B}-\frac{L^2}{B}\right)
\]
has mean zero and variance $2B$. Its third cumulant satisfies
\begin{equation}
 0< K_3(H_\Gamma)=\frac{8B^2}{L}
       \le 8C_3
       =K_3\left\{\sum_j\lambda_j(Z_j^2-1)\right\},
 \label{eq:gamma-third-cumulant}
\end{equation}
and the third-cumulant discrepancy is
\begin{equation}
 8\left(C_3-\frac{B^2}{L}\right)
       =\frac4L\sum_{j,\ell}
                     \lambda_j\lambda_\ell(\lambda_j-\lambda_\ell)^2.
 \label{eq:gamma-skewness-discrepancy}
\end{equation}
If all the positive remaining eigenvalues are equal, the gamma law is the
exact remaining chi-squared mixture. More generally, for every integer
$r\ge3$, its $r$th cumulant is positive and no larger than the corresponding
cumulant of that mixture.
\end{prop}

Proposition~\ref{prop:gamma-moment-comparison} shows that the gamma
remainder is exact for a flat remaining spectrum and retains part of its
higher cumulants in general. This does not establish smaller distributional
error or better level and power than the Gaussian remainder. The next
theorem gives a common sufficient error bound for the two calibrations.

\begin{theorem}[Quantitative estimated-reference comparison]
\label{thm:gamma-quantitative-calibration}
Use the notation of Theorem~\ref{thm:spectral-calibration}, and abbreviate
the original Gaussian-remainder reference by $\widehat R_n^G$.
For a realization of the estimators, put
\[
 e_{\mu,n}=\frac{|\widehat\mu_n-\mu_n|}{\sqrt{A_n}},
 \qquad e_{A,n}=\left|\frac{\widehat A_n}{A_n}-1\right|,
 \qquad
 \delta_n=\frac{\|\widehat{\bm{\Omega}}_n-\bm{\Omega}_n\|_{\mathrm{op}}}{\sqrt{A_n}},
\]
and
\begin{equation}
 \mathcal E_n=e_{\mu,n}^2+e_{A,n}
                 +\frac{\delta_n}{\tau_n}
                 +\left(\frac{\delta_n}{\tau_n}\right)^2.
 \label{eq:calibration-error-summary}
\end{equation}
On the event
\begin{equation}
 \frac12\le\frac{\widehat A_n}{A_n}\le2,
 \qquad 0<\tau_n\le\frac12,
 \qquad \delta_n\le\frac{\tau_n}{4},
 \label{eq:calibration-good-event}
\end{equation}
an absolute constant $C$ gives, for $M\in\{G,\Gamma\}$,
\begin{equation}
 d_K\!\left\{
  \mathcal L(\widehat R_n^M\mid\mathcal Y_n),\mathcal L(T^*_{0,n})
       \right\}
 \le \min\left[1,
       C\left\{\tau_n^{1/7}+\mathcal E_n^{1/5}\right\}\right].
 \label{eq:gamma-estimated-finite-bound}
\end{equation}
Consequently, under \eqref{eq:spectral-estimation}, the gamma calibration
has the same conditional distribution consistency, oracle-quantile
matching, and feasible power conclusions as the Gaussian-remainder
calibration. In particular,
\begin{equation}
 \frac{\widehat q_{n,\alpha}^{\Gamma}
              -q_{n,\alpha}^{\mathrm{NR}}}{s_{0,n}}
       \xrightarrow{p}0
 \label{eq:gamma-quantile-consistency}
\end{equation}
for every fixed $\alpha\in(0,1)$, without a fixed limiting spectrum.
\end{theorem}

In \eqref{eq:gamma-estimated-finite-bound}, the first term bounds the
approximation of the omitted spectrum and the second bounds estimation
of the trace, squared trace and retained weights. The exponents are
sufficient for consistency and are not claimed to be optimal. A quantitative
quantile-error bound would additionally require lower control of the
reference density near the critical value.

Under their stated moment and effective-dimension conditions, both
primitive constructions, with threshold
\eqref{eq:implemented-spectral-threshold}, satisfy the theorem's
conditions. The gamma reference also remains invariant to separate sample
translations, so the corresponding feasible power comparison is uniform
over the same mean translations.

\subsection{Variance correction and plug-in comparison}
\label{sec:variance-budget}

The positive part in \eqref{eq:spectral-tail} prevents a negative
remaining variance, but it does not limit the variance assigned to the
retained components. Indeed, if
$\widehat S_n=\sum_{j\le\widehat r_n}\widehat\lambda_{n,j}^2$
exceeds $\widehat A_n$, both unprojected references have conditional
variance $2\widehat S_n$, rather than $2\widehat A_n$. Empirical
eigenvalue inflation can therefore undo the squared-trace correction.
We recommend the following variance-budget projection in finite samples.
It leaves the threshold and the estimated total variance unchanged.

On $\widehat A_n>0$, define
\begin{align}
 \widehat c_n&=
 \begin{cases}
   \min\{1,\sqrt{\widehat A_n/\widehat S_n}\},&\widehat S_n>0,\\
   1,&\widehat S_n=0,
 \end{cases}
 &\widetilde\lambda_{n,j}&=\widehat c_n\widehat\lambda_{n,j},
       \quad j\le\widehat r_n,
 \label{eq:variance-budget-projection}\\
 \widetilde B_n&=\widehat A_n-
                   \sum_{j\le\widehat r_n}\widetilde\lambda_{n,j}^2,
 &\widetilde L_n&=\widehat\mu_n-
                   \sum_{j\le\widehat r_n}\widetilde\lambda_{n,j}.
 \label{eq:projected-tail-moments}
\end{align}
The projected Gaussian reference is
\begin{equation}
 \widetilde R_n^G=\widehat\mu_n+
       \sum_{j\le\widehat r_n}\widetilde\lambda_{n,j}(Z_j^2-1)
       +\sqrt{2\widetilde B_n}\,G.
 \label{eq:projected-gaussian-reference}
\end{equation}
For the projected gamma reference $\widetilde R_n^\Gamma$, replace the
last term by \eqref{eq:gamma-tail-rule}, using
$(\widetilde B_n,\widetilde L_n)$ in place of
$(\widehat B_n,\widehat L_n)$. Thus a zero remaining variance gives
the point mass at zero, and a positive remaining variance with
$\widetilde L_n\le0$ uses the Gaussian fallback. The same arbitrary
reference convention applies on $\widehat A_n\le0$. Auxiliary
variables are independent as before. In numerical implementation,
roundoff can be handled by replacing $\widetilde B_n$ by its positive
part; this is not an additional statistical adjustment.

\begin{prop}[Projection and exact estimated moments]
\label{prop:variance-budget-moments}
On $\widehat A_n>0$, the vector
$(\widetilde\lambda_{n,j})_{j\le\widehat r_n}$ is the unique minimizer of
\[
 \sum_{j\le\widehat r_n}(w_j-\widehat\lambda_{n,j})^2
 \quad\text{subject to}\quad
 w_j\ge0,\qquad \sum_{j\le\widehat r_n}w_j^2\le\widehat A_n.
\]
Moreover,
\begin{equation}
 \widetilde B_n=(\widehat A_n-\widehat S_n)_+=\widehat B_n,
 \qquad
 E(\widetilde R_n^M\mid\mathcal Y_n)=\widehat\mu_n,
 \qquad
 \operatorname{Var}(\widetilde R_n^M\mid\mathcal Y_n)=2\widehat A_n,
 \quad M\in\{G,\Gamma\}.
 \label{eq:projected-exact-moments}
\end{equation}
If $\widehat S_n\le\widehat A_n$, projection leaves both references
unchanged. If $\widehat S_n>\widehat A_n$, the two projected references
coincide, since their remaining variance is zero.
\end{prop}

The projection enforces the estimated variance while changing the
retained weights as little as possible in Euclidean distance. It cannot
correct an inaccurate spectral shape. Thus the moment identities alone
do not guarantee level accuracy; distributional consistency still requires
the estimation conditions in the following theorem. The unprojected
references provide a useful diagnostic for variance inflation.

\begin{theorem}[Accuracy of variance-budget calibration]
\label{thm:projected-calibration}
Use the error quantities in \eqref{eq:calibration-error-summary}.
On the event \eqref{eq:calibration-good-event}, for an absolute constant
$C$ and each $M\in\{G,\Gamma\}$,
\begin{equation}
 d_K\!\left\{
  \mathcal L(\widetilde R_n^M\mid\mathcal Y_n),
                  \mathcal L(T^*_{0,n})\right\}
 \le\min\left[1,
       C\{\tau_n^{1/7}+\mathcal E_n^{1/5}\}\right].
 \label{eq:projected-calibration-bound}
\end{equation}
Consequently, under \eqref{eq:spectral-estimation}, both projected
references satisfy conditional distribution consistency and
\[
 \frac{\widetilde q_{n,\alpha}^M-q_{n,\alpha}^{\mathrm{NR}}}{s_{0,n}}
       \xrightarrow{p}0
 \qquad (0<\alpha<1),
\]
without a fixed limiting spectrum. They inherit the feasible-power
conclusions of Proposition~\ref{prop:spectral-random-power}, including
uniform conclusions under uniform estimation and transfer conditions.
The primitive constructions and threshold in
\eqref{eq:implemented-spectral-threshold} require no modification.
\end{theorem}

\begin{cor}[Sharper calibration with three substantial spectral blocks]
\label{cor:calibration-bounded-density}
Suppose the eigenvalues of $\bm{\Omega}_n$ satisfy
\eqref{adap:three-blocks} with $m_n\ge cA_n$ for some fixed
$c\in(0,1/3]$. On \eqref{eq:calibration-good-event}, all four references
$R\in\{\widehat R_n^G,\widehat R_n^\Gamma,
        \widetilde R_n^G,\widetilde R_n^\Gamma\}$ satisfy
\begin{equation}
 d_K\!\left\{\mathcal L(R\mid\mathcal Y_n),
                           \mathcal L(T^*_{0,n})\right\}
 \le\min\left[1,C_c\{\tau_n^{1/4}+\mathcal E_n^{1/3}\}\right],
 \label{eq:calibration-bounded-density-rate}
\end{equation}
where $C_c$ depends only on $c$. This remains a distribution-function
bound; it does not assert a quantitative inverse-quantile rate.
\end{cor}

The need to correct the estimated total variance is already visible in an isotropic example.
The following proposition compares the direct empirical-spectrum
reference with the corrected references. Its conclusion concerns this
particular plug-in construction; the proof uses classical Gaussian
sample-covariance identities.

\begin{prop}[An isotropic separation from direct spectral plug-in]
\label{plug:proposition}
Let $\bm{Z}\sim N(\bm{0},\bm{I}_p)$ be independent of
$\bm{V}_1,\ldots,\bm{V}_m\stackrel{\mathrm{iid}}{\sim}N(\bm{0},\bm{I}_p)$, and suppose
$p,m\to\infty$.  Define
\[
 T=\|\bm{Z}\|^2,\qquad
 \bm{S}=\frac1m\sum_{r=1}^m \bm{V}_r\bm{V}_r^{\top},\qquad
 \widehat Q=\sum_{j=1}^p\lambda_j(\bm{S})\xi_j^2,
\]
where the $\xi_j$ are auxiliary independent standard normal variables.
Let $\widehat q_{\alpha}$ be the conditional $(1-\alpha)$ quantile
of $\widehat Q$, for fixed $\alpha\in(0,1/2)$, and write
$z_{1-\alpha}=\Phi^{-1}(1-\alpha)>0$.

If $p/m\to c\in[0,\infty)$, then
\begin{equation}
 \Pr\{T>\widehat q_{\alpha}\}
 \longrightarrow
 1-\Phi\{z_{1-\alpha}\sqrt{1+c}\}.
 \label{plug:finite-ratio}
\end{equation}
If $p/m\to\infty$, then
\begin{equation}
 \Pr\{T>\widehat q_{\alpha}\}\longrightarrow0.
 \label{plug:infinite-ratio}
\end{equation}
By contrast, the Gaussian and gamma spectral calibrations with
\begin{equation}
 \widehat\mu=\operatorname{tr}(\bm{S}),\qquad
 \widehat A=\frac1{m(m-1)}
       \sum_{r\ne s}(\bm{V}_r^{\top}\bm{V}_s)^2
 \label{plug:corrected-moments}
\end{equation}
and threshold \eqref{eq:implemented-spectral-threshold} have limiting
rejection probability $\alpha$ whenever $p=o(m^2)$.
This assertion applies also to their variance-budget projections.
In particular, throughout $m\ll p\ll m^2$, the direct empirical
spectrum reference has limiting size zero while the corrected
references have limiting size $\alpha$.
\end{prop}

\begin{proof}See Appendix~\ref{app:calibration}.\end{proof}

\begin{remark}[The correction in the regime $m\ll p\ll m^2$]
\label{plug:gamma-reduction}
In the displayed separation regime, the retained spectrum is empty
with probability tending to one under the implemented threshold.
Indeed, $\widehat a/a\xrightarrow{p}1$,
$\widehat A/p\xrightarrow{p}1$, and eventually
$\tau\geq\widehat a^{1/4}$ with probability tending to one.
Since $\lambda_1(\bm{S})=O_p(p/m)$,
\[
 \frac{\lambda_1(\bm{S})}{\tau\sqrt{\widehat A}}
       =O_p(a^{3/4})\longrightarrow0.
\]
On this event the Gaussian reference is
$N(\widehat\mu,2\widehat A)$ and the gamma reference is
\[
 \frac{\widehat A}{\widehat\mu}
       \chi^2_{\widehat\mu^2/\widehat A}.
\]
Both use the corrected variance $2\widehat A\sim2p$, whereas the direct
empirical reference uses $2\operatorname{tr}(\bm{S}^2)\sim2p(1+p/m)$.
The separation in this example is therefore due to variance correction,
not the choice of remainder. The independence model is included in the
two-sample construction as a boundary case: take $n_1=n_2=m+1$,
$\bm{\Sigma}_1=2\bm{I}_p$ and $\bm{\Sigma}_2=\bm{0}$ under the null.
Then $\bm{\Omega}_1=\bm{I}_p$ and $\bm{\Omega}_2=\bm{0}$; the first sample's
Helmert residuals supply the $m$ independent $\bm{V}_r$, while the second
sample's residuals vanish. The mean statistic is independent of these
residuals by Gaussianity.
\end{remark}

For computation, form the matrix with columns \(\bm{V}_{ir}/\sqrt{m_i}\), where \(m_i\)
denotes the residual count used by the chosen construction. Its Gram
matrix has the same nonzero eigenvalues as \(\widehat{\bm{\Omega}}_n\).
The terms of \(\widehat A_n\) are sums of squared inner products
within and between the two residual samples, with the within-sample diagonals
removed. This avoids inversion of a sample covariance matrix.

The conditional quantile can then be computed by sampling the reference using independent standard-normal and chi-squared
variables, with noninteger degrees of freedom allowed for the gamma
component. With $M_n$ independent auxiliary draws, $M_n\to\infty$ gives
uniform convergence of the conditional empirical distribution. The
subsequence argument used for Theorem~\ref{thm:spectral-calibration} then
gives the same quantile conclusion. This additional Monte Carlo error is
separate from covariance-estimation error; the numerical study reports
its finite reference-sampling budgets.

\subsection{Proofs for spectral calibration}

\begin{proof}[Proof of Theorem~\ref{thm:spectral-calibration}]
We first control the retained and omitted eigenvalues, then compare the
conditional laws, and finally prove quantile consistency. Suppress $n$ and
put $\delta=\|\widehat{\bm{\Omega}}-\bm{\Omega}\|_{\mathrm{op}}/\sqrt A$.
On an event with probability tending to one,
$1/2\le\widehat A/A\le2$ and $\delta\le\tau/4$. Weyl's inequality gives
$|\widehat\lambda_j-\lambda_j|\le\delta\sqrt A$. For retained indices,
$\lambda_j\ge(\tau/\sqrt2-\delta)\sqrt A\ge c\tau\sqrt A$.
Since $\sum_j\lambda_j^2=A$, it follows that
\[
 \widehat r\le C\tau^{-2},\qquad
 \frac{\sum_{j\le\widehat r}(\widehat\lambda_j-\lambda_j)^2}{A}
       \le C\frac{\delta^2}{\tau^2}=o_p(1).
\]
For omitted indices, Weyl's inequality also gives
$\lambda_j\le(\sqrt2\tau+\delta)\sqrt A$. These two bounds do not require
a gap at the threshold. Write $B=\sum_{j>\widehat r}\lambda_j^2$.
By Cauchy--Schwarz,
\[
 \frac{\left|\sum_{j\le\widehat r}
                      (\widehat\lambda_j^2-\lambda_j^2)\right|}{A}
 \le 2\left\{\frac{\sum_{j\le\widehat r}
                  (\widehat\lambda_j-\lambda_j)^2}{A}\right\}^{1/2}
      +\frac{\sum_{j\le\widehat r}(\widehat\lambda_j-\lambda_j)^2}{A}
 =o_p(1).
\]
Because $B\ge0$ and the positive-part map is one-Lipschitz,
$|\widehat B-B|/A=o_p(1)$.

For each fixed realization of the data, introduce the intermediate
reference
\[
 R^\circ=\mu+\sum_{j\le\widehat r}\lambda_j(Z_j^2-1)+\sqrt{2B}\,G.
\]
We first compare it with $T^*_0$. Replace each omitted summand
$\lambda_j(Z_j^2-1)$ by an independent $N(0,2\lambda_j^2)$ variable.
The first two moments match. Third-order Taylor expansion gives, for any
$f$ with bounded third derivative,
\[
 \left|E f\left(\frac{R^\circ-\mu}{s_0}\right)
       -E f\left(\frac{T^*_0-\mu}{s_0}\right)\right|
 \le C\|f'''\|_\infty
       \frac{\sum_{j>\widehat r}\lambda_j^3}{A^{3/2}}
 \le C\|f'''\|_\infty\,\tau.
\]
The last inequality follows from the bound on the largest omitted
eigenvalue. A smooth indicator with transition width $\epsilon$ has third
derivative bounded by $C\epsilon^{-3}$. Degree-two Gaussian
anti-concentration gives
\[
 \sup_x\Pr\left\{x<\frac{T^*_0-\mu}{s_0}\le x+\epsilon\right\}
       \le C\epsilon^{1/2}.
\]
Smoothing with the full Gaussian target therefore yields
\[
 d_K\{\mathcal L(R^\circ),\mathcal L(T^*_0)\}
       \le C(\tau\epsilon^{-3}+\epsilon^{1/2})
       \le C\tau^{1/7}
\]
by taking $\epsilon=\tau^{2/7}$. This proves the required approximation
of the omitted spectrum.

Couple $\widehat R$ and $R^\circ$ using the same $Z_j,G$. Conditional on
the data,
\begin{align*}
 E\left[\left.\left(\frac{\widehat R-R^\circ}{s_0}\right)^2
                  \right|\mathcal Y\right]
 &=\frac{(\widehat\mu-\mu)^2}{2A}
   +\frac{\sum_{j\le\widehat r}(\widehat\lambda_j-\lambda_j)^2}{A}
   +\frac{(\sqrt{\widehat B}-\sqrt B)^2}{A}\\
 &\le\frac{(\widehat\mu-\mu)^2}{2A}
   +\frac{\sum_{j\le\widehat r}(\widehat\lambda_j-\lambda_j)^2}{A}
   +\frac{|\widehat B-B|}{A}=o_p(1).
\end{align*}
Apply conditional Chebyshev to this coupling and bracket the two
conditional distribution functions at points separated by a fixed positive
multiple of $s_0$. The preceding comparison and anti-concentration bound,
followed by a limit as the bracketing width decreases to zero, prove
\eqref{eq:spectral-conditional-consistency}.

It remains to establish quantile consistency when the full sequence has
no limiting spectrum. Take an arbitrary subsequence and put
$a_{n,j}=\lambda_{n,j}/\sqrt{A_n}$. Since the weights are ordered,
$0\le a_{n,j}\le j^{-1/2}$. A further subsequence has coordinatewise limits
$a_j$, with $\sum_j a_j^2\le1$. On this subsequence the standardized null
reference converges to
\[
 \zeta=\frac1{\sqrt2}\sum_{j\ge1}a_j(Z_j^2-1)
          +\sqrt{1-\sum_j a_j^2}\,G.
\]
To obtain this limit, retain the first $J$ terms and expand the logarithm
of the tail characteristic function to second order. The remainder is
bounded by
$C|t|^3\max_{j>J}a_{n,j}\sum_{j>J}a_{n,j}^2
\le C|t|^3(J+1)^{-1/2}$; let $n\to\infty$ and then $J\to\infty$.
The limit has variance one. If the Gaussian coefficient is positive,
its density is positive everywhere. Otherwise at least one $a_j$ is
positive. The corresponding scaled chi-squared variable has a density
positive on the interior of its half-line support. Conditioning on the
remaining independent terms shows that the convolution is continuous and
strictly increasing at every quantile with probability in $(0,1)$.
Consequently, for every $\eta>0$, the oracle quantile eventually lies
between $\mu_n+s_{0,n}(c_{\zeta,\alpha}-\eta)$ and
$\mu_n+s_{0,n}(c_{\zeta,\alpha}+\eta)$. By
\eqref{eq:spectral-conditional-consistency}, the estimated quantile lies
between the same endpoints with probability tending to one. Their
standardized difference therefore tends to zero along this further
subsequence. Since the initial subsequence was arbitrary,
\eqref{eq:spectral-quantile-consistency} follows. The same argument gives
\eqref{eq:spectral-quantile-profile} when a global profile exists.
\end{proof}

\begin{proof}[Proof of Proposition~\ref{prop:spectral-random-power}]
On $\{|\widehat q_n-q_n|\le\epsilon s_{0,n}\}$, rejection is bracketed
by the two deterministic thresholds $q_n\pm\epsilon s_{0,n}$.
Comparing each deterministic rejection probability with its Gaussian
counterpart contributes at most the Kolmogorov distance. The remaining
interval probability is bounded by Gaussian polynomial anti-concentration.
Indeed,
$\operatorname{Var}(T_n^*)=s_{0,n}^2+4\bm{\theta}_n^\top\bm{\Omega}_n\bm{\theta}_n
\ge s_{0,n}^2$, so the interval of length $2\epsilon s_{0,n}$ has
probability at most $C\epsilon^{1/2}$, uniformly in $\bm{\theta}_n$.
This proves \eqref{eq:random-critical-power-bound}. The consequences
follow by first letting $n\to\infty$ and then $\epsilon\downarrow0$.
The bracketing uses event inclusions and hence allows dependence between
the statistic and the estimated critical value.
\end{proof}

\begin{proof}[Proof of Proposition~\ref{prop:gaussian-calibration-verification}]
The identities $\bm{H}_i^\top\mathbf1=\bm{0}$ and $\bm{H}_i^\top \bm{H}_i=\bm{I}$ show that
the residual contrasts are centered, have covariance $\bm{\Omega}_i$, and are
uncorrelated. Joint Gaussianity makes them independent. We first verify
the trace estimates. Gaussian quadratic-form identities give
\[
 E\widehat\mu=\mu,\qquad
 \operatorname{Var}(\widehat\mu)
       =2\sum_{i=1}^2\frac{\operatorname{tr}(\bm{\Omega}_i^2)}{\nu_i}
       \le\frac{2A}{\nu_{\min}}.
\]
Unbiasedness of $\widehat A$ follows from
$E(\bm{V}_{ir}^\top \bm{V}_{js})^2=\operatorname{tr}(\bm{\Omega}_i\bm{\Omega}_j)$
for independent residual contrasts. Its Hoeffding decomposition is
\[
 \widehat A-A
 =2\sum_{i=1}^2\frac1{\nu_i}\sum_r
       \{\bm{V}_{ir}^\top\bm{\Omega} \bm{V}_{ir}-\operatorname{tr}(\bm{\Omega}_i\bm{\Omega})\}
       +D,
\]
where $D$ is the sum of degenerate pair kernels. The two terms are
orthogonal in $L^2$. The variance of the linear part is
\[
 8\sum_{i=1}^2\frac{\operatorname{tr}(\bm{\Omega}_i\bm{\Omega}\bm{\Omega}_i\bm{\Omega})}{\nu_i}
 \le\frac{8\operatorname{tr}(\bm{\Omega}^4)}{\nu_{\min}}
 \le\frac{8A^2}{\nu_{\min}}.
\]
For the first inequality put $\bm{M}_i=\bm{\Omega}^{1/2}\bm{\Omega}_i\bm{\Omega}^{1/2}$;
$\bm{M}_i$ are positive semidefinite, $\sum_{i=1}^2\bm{M}_i=\bm{\Omega}^2$, and
$\sum_{i=1}^2\operatorname{tr}(\bm{M}_i^2)\le\operatorname{tr}(\bm{\Omega}^4)$.
For independent $\bm{V}\sim N(\bm{0},\bm{\Omega}_i)$ and $\bm{V}'\sim N(\bm{0},\bm{\Omega}_j)$,
\[
 E(\bm{V}^\top \bm{V}')^4
 =3\{\operatorname{tr}(\bm{\Omega}_i\bm{\Omega}_j)\}^2
       +6\operatorname{tr}(\bm{\Omega}_i\bm{\Omega}_j\bm{\Omega}_i\bm{\Omega}_j)
 \le9\{\operatorname{tr}(\bm{\Omega}_i\bm{\Omega}_j)\}^2.
\]
Orthogonal projection onto the degenerate kernel subspace cannot increase
its $L^2$ norm. Distinct degenerate pairs are uncorrelated, including
pairs sharing one observation. Consequently,
\[
 \operatorname{Var}(D)
 \le18\sum_{i=1}^2\frac{\{\operatorname{tr}(\bm{\Omega}_i^2)\}^2}
                         {\nu_i(\nu_i-1)}
       +36\frac{\{\operatorname{tr}(\bm{\Omega}_1\bm{\Omega}_2)\}^2}
                              {\nu_1\nu_2}
 \le\frac{C A^2}{\nu_{\min}^2}.
\]
This proves \eqref{eq:trace-estimator-bounds} and the two trace
consistency statements.

We next control the covariance estimator. The expectation bound of
\citet{KoltchinskiiLounici2017}, applied separately to the two residual samples,
gives
\[
 E\|\widehat{\bm{\Omega}}_i-\bm{\Omega}_i\|_{\mathrm{op}}
 \le C\left\{
 \sqrt{\frac{\|\bm{\Omega}_i\|_{\mathrm{op}}\operatorname{tr}(\bm{\Omega}_i)}{\nu_i}}
       +\frac{\operatorname{tr}(\bm{\Omega}_i)}{\nu_i}\right\}.
\]
Since $\bm{\Omega}_i\preceq\bm{\Omega}$, summation and Cauchy--Schwarz imply
\[
 E\|\widehat{\bm{\Omega}}-\bm{\Omega}\|_{\mathrm{op}}
 \le C\left\{\sqrt{\frac{2\|\bm{\Omega}\|_{\mathrm{op}}\mu}{\nu_{\min}}}
                  +\frac{\mu}{\nu_{\min}}\right\}.
\]
Since $\|\bm{\Omega}\|_{\mathrm{op}}\le\sqrt A$, Markov's
inequality gives \eqref{eq:gaussian-calibration-op-rate}. Finally,
$\widehat\mu/\mu\to_p1$ follows from $\mu^2\ge A$ and the trace variance
bound. Hence $\widehat a/a\to_p1$. The cap in
\eqref{eq:implemented-spectral-threshold} is inactive with probability
tending to one, $\tau\to_p0$, and
\[
 \frac{\|\widehat{\bm{\Omega}}-\bm{\Omega}\|_{\mathrm{op}}}{\tau\sqrt A}
       =O_p\left(\frac{\sqrt a+a}{a^{1/4}+\nu_{\min}^{-1/8}}\right)
       =o_p(1).
\]
All of \eqref{eq:spectral-estimation} is now verified.
\end{proof}

\begin{proof}[Proof of Corollary~\ref{cor:paired-calibration}]
Each residual in \eqref{eq:paired-residuals} equals
$(\bm{U}_{i,2r}-\bm{U}_{i,2r-1})/\sqrt2$. Disjoint pairs are independent and
have mean zero and covariance $\bm{\Omega}_i$. The triangle inequality for
the $\psi_2$ norm preserves the covariance-relative sub-Gaussian bound,
up to a fixed constant. For two independent copies $\bm{U},\bm{U}'$ of $\bm{U}_{ij}$,
\[
 \operatorname{Var}\left(\left\|(\bm{U}-\bm{U}')/\sqrt2\right\|^2\right)
 =\frac12\operatorname{Var}(\|\bm{U}\|^2)+\operatorname{tr}(\bm{\Omega}_i^2).
\]
The mixed covariances vanish by centering and independence. This gives the
required trace-mean variance bound.

Sub-Gaussian projections imply
$E(\bm{u}^\top \bm{V}_{ir})^4\le C_K(\bm{u}^\top\bm{\Omega}_i \bm{u})^2$. For every positive
semidefinite $\bm{M}$, expansion in its eigenbasis and Cauchy--Schwarz yield
\[
 E(\bm{V}_{ir}^\top \bm{M} \bm{V}_{ir})^2
       \le C_K\{\operatorname{tr}(\bm{M}\bm{\Omega}_i)\}^2.
\]
Consequently, for independent residuals from samples $i,j\in\{1,2\}$,
\[
 E(\bm{V}_{ir}^\top \bm{V}_{js})^4
 \le C_K E(\bm{V}_{js}^\top\bm{\Omega}_i \bm{V}_{js})^2
 \le C_K\{\operatorname{tr}(\bm{\Omega}_i\bm{\Omega}_j)\}^2.
\]
The Hoeffding decomposition used in the Gaussian proof remains exact.
Its linear variance is now bounded by
\[
 C_K\sum_{i=1}^2\frac{\{\operatorname{tr}(\bm{\Omega}_i\bm{\Omega})\}^2}{m_i}
 \le \frac{C_K A^2}{m_{\min}},
\]
because the nonnegative numbers $\operatorname{tr}(\bm{\Omega}_i\bm{\Omega})$ sum
to $A$. The degenerate variance is at most $C_K A^2/m_{\min}^2$ by the
preceding fourth-moment bound. Thus $E\widehat A=A$ and
$\operatorname{Var}(\widehat A)\le C_K A^2/m_{\min}$.

For the covariance estimate, apply the upper bound in Theorem~4 of
\citet{KoltchinskiiLounici2017}. Its sub-Gaussian condition holds by
assumption. Its pre-Gaussian condition also holds, since a Gaussian
vector with the same covariance exists in finite-dimensional Euclidean
space. The resulting bound is
\[
 E\|\widehat{\bm{\Omega}}_i-\bm{\Omega}_i\|_{\mathrm{op}}
 \le C_K\left\{
 \sqrt{\frac{\|\bm{\Omega}_i\|_{\mathrm{op}}\operatorname{tr}(\bm{\Omega}_i)}{m_i}}
             +\frac{\operatorname{tr}(\bm{\Omega}_i)}{m_i}\right\}.
\]
Summing over the two samples gives
\eqref{eq:gaussian-calibration-op-rate}, with $m_{\min}$ in place of
$\nu_{\min}$. The threshold verification is then the same as in the
Gaussian case. Theorem~\ref{thm:spectral-calibration} completes the proof.
\end{proof}

\subsection{Gamma remainder}
\label{app:gamma-refinement}

\begin{proof}[Proof of Proposition~\ref{prop:gamma-moment-comparison}]
For $a=B/L$ and $d=L^2/B$, the centered variable $a(\chi_d^2-d)$ has
$r$th cumulant $2^{r-1}(r-1)!a^rd$ for $r\ge2$.
This gives the variance and third cumulant. Cauchy--Schwarz yields
$B^2\le LC_3$, and expansion of
$\sum_{j,\ell}\lambda_j\lambda_\ell(\lambda_j-\lambda_\ell)^2$
gives \eqref{eq:gamma-skewness-discrepancy}. For $r\ge3$, Jensen's
inequality with weights $\lambda_j/L$ gives
\[
 \sum_j\lambda_j^r
   =L\sum_j\frac{\lambda_j}{L}\lambda_j^{r-1}
   \ge L\left(\frac BL\right)^{r-1}
   =\frac{B^{r-1}}{L^{r-2}}.
\]
The right side is $a^rd$, which proves the cumulant comparison. If the
remaining positive eigenvalues equal $\lambda$ and there are $m$ of them,
then $a=\lambda$ and $d=m$, giving exactly
$\lambda(\chi_m^2-m)$.
\end{proof}

\begin{lem}[Smooth comparison of a centered gamma variable]
\label{lem:gamma-smooth-comparison}
For $L,B>0$, let $H=(B/L)(\chi^2_{L^2/B}-L^2/B)$ and let
$N_B\sim N(0,2B)$. If $f$ has bounded third derivative, then
\begin{equation}
 |E f(H)-E f(N_B)|
       \le C\|f'''\|_\infty\frac{B^2}{L}.
 \label{eq:gamma-smooth-bound}
\end{equation}
The same inequality holds after adding an independent common random
variable to both arguments. No lower bound on $L^2/B$ is required.
\end{lem}

\begin{proof}
Put $d=L^2/B$, $a=B/L$, $m=\max\{1,\lfloor d\rfloor\}$ and $v=d/m$.
Then $0<v<2$, and infinite divisibility of the gamma law represents
$H$ as the sum of $m$ independent copies of $a(\chi_v^2-v)$.
The identity $E(\chi_v^2)^3=v(v+2)(v+4)$ and
$|x-y|^3\le4(x^3+y^3)$ imply
$E|\chi_v^2-v|^3\le Cv$ for $0<v\le2$.
The centered Gaussian variable with variance $2v$ also has third absolute
moment bounded by $Cv$ on this interval. Replace the $m$ summands one at
a time by independent Gaussians of matching variance. The constant,
linear, and quadratic Taylor terms cancel, and the sum of the remainders
is bounded by $C\|f'''\|_\infty ma^3v
=C\|f'''\|_\infty B^2/L$. Conditional expansion proves the assertion
with an independent common term.
\end{proof}

\begin{remark}[An oracle bound]
\label{rem:gamma-oracle-bound}
For a deterministic retained set consisting of the largest eigenvalues,
let $C_3=\sum_{j>r}\lambda_j^3$, and use the exact remaining $L,B$.
The Gaussian replacement of the omitted chi-squared variables has smooth
error at most $C\|f'''\|_\infty C_3/A^{3/2}$ on the standardized scale.
Lemma~\ref{lem:gamma-smooth-comparison} and $B^2/L\le C_3$ give the
same order for a gamma remainder. Smoothing against the full Gaussian
quadratic-form target yields
\[
 d_K\{\mathcal L(R^{\Gamma,\mathrm{or}}),\mathcal L(T_0^*)\}
       \le C\min\left\{1,
                 \left(\frac{C_3}{A^{3/2}}\right)^{1/7}\right\}.
\]
Thus both remainder approximations are consistent when the largest
omitted normalized eigenvalue tends to zero. The cumulant comparison
alone does not improve the uniform exponent in this bound.
\end{remark}

\begin{proof}[Proof of Theorem~\ref{thm:gamma-quantitative-calibration}]
We first prove the bound for the Gaussian remainder and then compare
the gamma remainder with it. Suppress $n$ and work on
\eqref{eq:calibration-good-event}. The spectral estimates in the proof of
Theorem~\ref{thm:spectral-calibration} give
\begin{equation}
 \widehat r\le C\tau^{-2},\qquad
 \frac{\sum_{j\le\widehat r}(\widehat\lambda_j-\lambda_j)^2}{A}
        \le C(\delta/\tau)^2,
 \qquad \max_{j>\widehat r}\lambda_j\le C\tau\sqrt A.
 \label{eq:gamma-head-tail-bounds}
\end{equation}
Let $B=\sum_{j>\widehat r}\lambda_j^2$, $L=\sum_{j>\widehat r}\lambda_j$
and $b=B/A$. Positive-part contraction and Cauchy--Schwarz imply
\begin{equation}
 \frac{|\widehat B-B|}{A}
     \le C\left\{e_A+\delta/\tau+(\delta/\tau)^2\right\}
     =:u.
 \label{eq:gamma-tail-variance-error}
\end{equation}
Choose the absolute constant in $u$ large enough for this bound.

First consider the Gaussian-remainder calibration. Its oracle intermediate
law, with true retained eigenvalues, true mean, and remaining variance
$2B$, is within $C\tau^{1/7}$ of $T^*_0$ by third-order replacement and
degree-two anti-concentration. Coupling the estimated and intermediate
laws by the same auxiliary normals gives conditional mean-square error,
on the scale $s_0=\sqrt{2A}$, at most
\[
 C\left\{e_\mu^2+(\delta/\tau)^2+u\right\}\le C\mathcal E.
\]
Indeed, $(\sqrt{\widehat B}-\sqrt B)^2\le|\widehat B-B|$.
For any $\epsilon>0$, conditional Chebyshev and the anti-concentration
bound for $T^*_0$ therefore give
\[
 r_G:=d_K\{\mathcal L(\widehat R^G\mid\mathcal Y),\mathcal L(T^*_0)\}
 \le C\{\tau^{1/7}+\mathcal E\epsilon^{-2}+\epsilon^{1/2}\}.
\]
Taking $\epsilon=\mathcal E^{2/5}$ for $0<\mathcal E\le1$ and using
the trivial bound for $\mathcal E>1$ proves
\eqref{eq:gamma-estimated-finite-bound} for $M=G$. The case
$\mathcal E=0$ follows by letting $\epsilon\downarrow0$.

For the gamma remainder, consistency of the total trace does not by
itself control the estimated remaining trace on the scale $\sqrt A$.
The rank bound instead gives
\begin{equation}
 \frac{|\widehat L-L|}{\sqrt A}
       \le e_\mu+C\frac{\delta}{\tau^2}.
 \label{eq:gamma-tail-mean-error}
\end{equation}
Define $h=u+\tau e_\mu+C\delta/\tau$, enlarging constants if needed.
Then $h\le C(\mathcal E+\tau^2)$, because
$\tau e_\mu\le(e_\mu^2+\tau^2)/2$.
We now distinguish a substantial remaining variance from a negligible
one. Fix a sufficiently large absolute constant $C_0$.

If $b>C_0h$, then $B>0$, and
\eqref{eq:gamma-head-tail-bounds} implies
\[
 L\ge\frac{B}{\max_{j>\widehat r}\lambda_j}
       \ge\frac{b\sqrt A}{C\tau}.
\]
Combining this with \eqref{eq:gamma-tail-mean-error} and the definition
of $h$ gives $\widehat L\ge b\sqrt A/(C\tau)>0$.
Also $0<\widehat B\le CbA$, by
\eqref{eq:gamma-tail-variance-error}. Thus the gamma branch is used and
\[
 \frac{\widehat B^2}{\widehat L A^{3/2}}\le C\tau b\le C\tau.
\]
Lemma~\ref{lem:gamma-smooth-comparison}, with the retained sum added,
compares $\widehat R^\Gamma$ and $\widehat R^G$ with standardized smooth
error at most $C\tau\|f'''\|_\infty$. The interval probabilities of
$\widehat R^G$ are at most the corresponding degree-two
anti-concentration bound plus $2r_G$. Hence smoothing with
$\epsilon=\tau^{2/7}$ gives
\[
 d_K\{\mathcal L(\widehat R^\Gamma\mid\mathcal Y),
                     \mathcal L(\widehat R^G\mid\mathcal Y)\}
       \le C\tau^{1/7}+2r_G.
\]

If $b\le C_0h$, then $\widehat B/A\le Ch$.
In every branch of \eqref{eq:gamma-tail-rule}, the remaining variable is
centered and has variance $2\widehat B$. Couple it with an independent
Gaussian remainder and use the same retained coordinates. The conditional
mean-square difference divided by $s_0^2$ is at most $C h$.
Chebyshev and the same interval-probability bound give
\[
 d_K\{\mathcal L(\widehat R^\Gamma\mid\mathcal Y),
                     \mathcal L(\widehat R^G\mid\mathcal Y)\}
       \le C h^{1/5}+2r_G
\]
when $0<h\le1$, with the trivial bound covering $h>1$ and a limiting
argument covering $h=0$. Finally,
$h^{1/5}\le C(\mathcal E^{1/5}+\tau^{2/5})$ and
$\tau^{2/5}\le\tau^{1/7}$ for $0<\tau\le1$.
The triangle inequality and the bound for $r_G$ prove
\eqref{eq:gamma-estimated-finite-bound} in both cases.

Under \eqref{eq:spectral-estimation}, the good event has probability
tending to one and the bound tends to zero in probability. Apply the
subsequence quantile argument from
Theorem~\ref{thm:spectral-calibration}, followed by the bracketing bound
in Proposition~\ref{prop:spectral-random-power}, to obtain quantile
consistency and the feasible power conclusions. Uniform estimation and
comparison bounds give their uniform versions.
\end{proof}

\subsection{Variance-budget projection}
\label{app:variance-budget}

\begin{proof}[Proof of Proposition~\ref{prop:variance-budget-moments}]
Suppress $n$, and let $\bm{h}$ be the nonnegative vector of retained estimated
eigenvalues. If $\|\bm{h}\|_2\le\sqrt{\widehat A}$, the unique minimizer is
$\bm{h}$. Otherwise, for any feasible $\bm{w}$, Cauchy--Schwarz gives
\[
 \|\bm{w}-\bm{h}\|_2^2
  =\|\bm{w}\|_2^2+\|\bm{h}\|_2^2-2\bm{w}^\top \bm{h}
  \ge(\|\bm{h}\|_2-\|\bm{w}\|_2)^2
  \ge(\|\bm{h}\|_2-\sqrt{\widehat A})^2.
\]
Equality is attained only at
$\bm{w}=\sqrt{\widehat A}\,\bm{h}/\|\bm{h}\|_2$. This proves the projection
assertion, including the zero-dimensional case. Its squared norm is
$\min\{\widehat S,\widehat A\}$, giving the first identity in
\eqref{eq:projected-exact-moments}. Every remaining-variable branch
has mean zero and variance $2\widetilde B$. Independence of auxiliary
variables then gives the asserted mean and variance. The two cases
for $\widehat S$ follow directly from the definitions.
\end{proof}

\begin{proof}[Proof of Theorem~\ref{thm:projected-calibration}]
Suppress $n$ and work on \eqref{eq:calibration-good-event}. Let
$\bm{h}=(\widehat\lambda_j)_{j\le\widehat r}$,
$\widetilde{\bm{h}}=(\widetilde\lambda_j)_{j\le\widehat r}$, and
$\bm{\ell}=(\lambda_j)_{j\le\widehat r}$. By
\eqref{eq:gamma-head-tail-bounds},
$\|\bm{h}-\bm{\ell}\|_2^2/A\le C(\delta/\tau)^2$, whereas
$\|\bm{\ell}\|_2^2\le A$. Thus
\begin{align}
 \frac{\|\widetilde{\bm{h}}-\bm{h}\|_2^2}{A}
 &=\frac{(\sqrt{\widehat S}-\sqrt{\widehat A})_+^2}{A}
 \le\frac{(\widehat S-\widehat A)_+}{A}\nonumber\\
 &\le e_A+\frac{|\|\bm{h}\|_2^2-\|\bm{\ell}\|_2^2|}{A}
 \le e_A+C\frac\delta\tau+C\left(\frac\delta\tau\right)^2.
 \label{eq:projection-displacement}
\end{align}
Consequently $\|\widetilde{\bm{h}}-\bm{\ell}\|_2^2/A\le C\mathcal E$.
Since $\widetilde B=\widehat B$, the same variance-error bound
\eqref{eq:gamma-tail-variance-error} applies. Couple
$\widetilde R^G$ with the oracle intermediate law having mean $\mu$,
retained weights $\bm{\ell}$, and Gaussian remaining variance $2B$, using
the same auxiliary normals. Its mean-square error divided by
$s_0^2=2A$ is at most $C\mathcal E$. The oracle approximation and
Chebyshev--anti-concentration argument in the proof of
Theorem~\ref{thm:gamma-quantitative-calibration} therefore prove
\eqref{eq:projected-calibration-bound} for $M=G$.

For $M=\Gamma$, if $\widehat S\le\widehat A$, the projected reference
equals the unprojected gamma reference, to which
Theorem~\ref{thm:gamma-quantitative-calibration} applies. If
$\widehat S>\widehat A$, then $\widetilde B=0$, so
$\widetilde R^\Gamma=\widetilde R^G$, whose bound was just proved.
This proves the finite comparison in every case. Under
\eqref{eq:spectral-estimation}, its right side tends to zero in
probability. The spectral-subsequence quantile argument and random
critical-value bracketing used previously prove the remaining
conclusions. Translation invariance of both primitive constructions
is preserved by the projection, giving their uniform statements.
\end{proof}

\begin{proof}[Proof of Corollary~\ref{cor:calibration-bounded-density}]
By Theorem~\ref{adap:uniform-density} at $\bm{\theta}=\bm{0}$, the standardized
target $(T_0^*-\mu)/\sqrt{2A}$ has a density bounded by $M_c$.
Its interval probabilities over intervals of length $\epsilon$ are
therefore at most $M_c\epsilon$. Repeat the Gaussian-reference proof
with this linear bound in place of the square-root bound. Third-order
replacement followed by smoothing gives
$C_c(\tau\epsilon^{-3}+\epsilon)$, optimized at
$\epsilon=\tau^{1/4}$. The coupling error gives
$C_c(\mathcal E\epsilon^{-2}+\epsilon)$, optimized at
$\epsilon=\mathcal E^{1/3}$. This proves
\eqref{eq:calibration-bounded-density-rate} for both Gaussian references,
using \eqref{eq:projection-displacement} for the projected one.

For the unprojected gamma reference, retain the two cases in the proof
of Theorem~\ref{thm:gamma-quantitative-calibration}, writing
$r_G=d_K\{\mathcal L(\widehat R^G\mid\mathcal Y),\mathcal L(T_0^*)\}$.
In the case $b>C_0h$, smoothing the gamma-to-Gaussian comparison now
gives $C_c\tau^{1/4}+2r_G$. In the case $b\le C_0h$, the
mean-square coupling gives $C_c h^{1/3}+2r_G$. Since
$h\le C(\mathcal E+\tau^2)$ and
$\tau^{2/3}\le\tau^{1/4}$ for $\tau\le1$, both cases give the stated
bound. Zero errors follow by limits, and the trivial bound covers
errors larger than one. Finally, the projected gamma reference equals
the unprojected gamma reference on $\widehat S\le\widehat A$ and the
projected Gaussian reference on its complement, completing the proof.
\end{proof}

\begin{proof}[Proof of Proposition~\ref{plug:proposition}]
We first obtain the trace estimates and a conditional normal limit for
the empirical reference. Gaussian quadratic-form identities give
\begin{equation}
 \operatorname{E}\operatorname{tr}(\bm{S})=p,
 \qquad
 \operatorname{Var}\{\operatorname{tr}(\bm{S})\}=2p/m,
 \qquad
 \frac{\operatorname{tr}(\bm{S})-p}{\sqrt{2p}}
                         \xrightarrow{p}0.
 \label{plug:first-trace}
\end{equation}
For the second trace, the exact decomposition is
\begin{equation}
 \operatorname{tr}(\bm{S}^2)
 =\frac{m-1}{m}\widehat A
       +\frac1{m^2}\sum_{r=1}^m\|\bm{V}_r\|^4.
 \label{plug:second-trace-decomposition}
\end{equation}
The unbiased squared-trace estimator satisfies
\begin{equation}
 \operatorname{E}\widehat A=p,
 \qquad
 \operatorname{Var}(\widehat A)
   =\frac{8p}{m}+\frac{4p(p+1)}{m(m-1)}.
 \label{plug:unbiased-trace-variance}
\end{equation}
Indeed, for the kernel $h(\bm{v},\bm{w})=(\bm{v}^{\top}\bm{w})^2$, its first centered
projection is $\|\bm{v}\|^2-p$, with variance $2p$.
Also $\operatorname{Var}\{h(\bm{V}_1,\bm{V}_2)\}=2p^2+6p$, so its
degenerate second projection has variance $2p(p+1)$.
Expanding the variance of the average over unordered pairs gives
\eqref{plug:unbiased-trace-variance}.  Consequently
$\widehat A/p\xrightarrow{p}1$.

The chi-squared moments give
\[
 \operatorname{E}\|\bm{V}_r\|^4=p(p+2),\qquad
 \operatorname{Var}(\|\bm{V}_r\|^4)=8p(p+2)(p+3).
\]
Let $H=p(1+p/m)$.  The variance of the second term in
\eqref{plug:second-trace-decomposition}, divided by $H^2$, equals
\[
 \frac{8(p+2)(p+3)}{mp(m+p)^2}\longrightarrow0.
\]
Its mean divided by $H$ is $(p+2)/(m+p)$, while the first term
has leading ratio $(m-1)/(m+p)$ by
\eqref{plug:unbiased-trace-variance}.  Since these two ratios sum
to $1+1/(m+p)$, we obtain
\begin{equation}
 \frac{\operatorname{tr}(\bm{S}^2)}{p(1+p/m)}\xrightarrow{p}1.
 \label{plug:second-trace-limit}
\end{equation}

We next verify the conditional normal limit.
If $\bm{\mathsf V}$ is the $p\times m$ matrix with columns $\bm{V}_r$, an
elementary Gaussian net bound yields
$\|\bm{\mathsf V}\|_{\mathrm{op}}=O_p(\sqrt{p+m})$.  One proof uses
$1/4$-nets of the two unit spheres, each with at most $9^p$ and
$9^m$ points.  The operator norm is at most twice the maximum of
$|\bm{u}^{\top}\bm{\mathsf V}\bm{v}|$ over these nets.  Each such scalar is standard
normal, so a union bound at $3\sqrt{p+m}$ proves the claimed order.
It follows that
\[
 \lambda_1(\bm{S})=\|\bm{\mathsf V}\|_{\mathrm{op}}^2/m
                         =O_p(1+p/m).
\]
Together with \eqref{plug:second-trace-limit}, this gives
\begin{equation}
 \frac{\lambda_1(\bm{S})}{\sqrt{\operatorname{tr}(\bm{S}^2)}}
       =O_p\left(\sqrt{\frac1p+\frac1m}\right)
       \xrightarrow{p}0.
 \label{plug:conditional-diffuseness}
\end{equation}
Conditional on $\bm{S}$, the standardized reference is
\[
 \frac{\widehat Q-\operatorname{tr}(\bm{S})}
      {\sqrt{2\operatorname{tr}(\bm{S}^2)}}
 =\frac1{\sqrt2}\sum_j
      \frac{\lambda_j(\bm{S})}{\sqrt{\operatorname{tr}(\bm{S}^2)}}
                       (\xi_j^2-1).
\]
The weights have sum of squares one and their maximum tends to zero
in probability. Expanding the logarithm of the conditional characteristic
function gives $-t^2/2+o_p(1)$ for each fixed $t$; the remainder is bounded
by a constant depending on $t$ times the maximum weight. To make the
conditional convergence explicit, take any subsequence and then a further
subsequence along which the maximum weight tends to zero almost surely.
The deterministic-array characteristic-function argument gives a normal
limit along that further subsequence. Since this limit is continuous,
conditional convergence holds in Kolmogorov distance, in probability.
Its quantiles consequently satisfy
\begin{equation}
 \widehat q_{\alpha}
 =\operatorname{tr}(\bm{S})
       +\sqrt{2\operatorname{tr}(\bm{S}^2)}
                         \{z_{1-\alpha}+o_p(1)\}.
 \label{plug:quantile-expansion}
\end{equation}

If $p/m\to c<\infty$, equations
\eqref{plug:first-trace}, \eqref{plug:second-trace-limit}, and
\eqref{plug:quantile-expansion} give
\[
 \frac{\widehat q_{\alpha}-p}{\sqrt{2p}}
                  \xrightarrow{p}z_{1-\alpha}\sqrt{1+c}.
\]
Since $(T-p)/\sqrt{2p}\xrightarrow{d}N(0,1)$, deterministic
threshold bracketing proves \eqref{plug:finite-ratio}.
If $p/m\to\infty$, the same display with its diverging scale and
$z_{1-\alpha}>0$ shows that the standardized critical value tends
to $+\infty$ in probability.  Tightness of $(T-p)/\sqrt{2p}$
proves \eqref{plug:infinite-ratio}.

For the corrected references, the population summaries are
$\mu=p$ and $A=p$.  The Gaussian primitive condition
\eqref{eq:gaussian-calibration-effective-dimension}, in the two-sample
boundary realization described in Remark~\ref{plug:gamma-reduction}, becomes
\[
 a=\frac{\mu}{m\sqrt A}=\frac{\sqrt p}{m}\longrightarrow0.
\]
It is exactly $p=o(m^2)$.  The estimators in
\eqref{plug:corrected-moments} coincide with the two-sample construction of
Proposition~\ref{prop:gaussian-calibration-verification} in that realization,
since the second sample has zero residuals and both residual counts equal $m$; its
implemented threshold therefore satisfies
\eqref{eq:spectral-estimation}.  The Gaussian, gamma, and projected
calibration theorems imply that their critical values differ from
the exact $\chi_p^2$ critical value by $o_p(\sqrt p)$.  The normal
limit of $\chi_p^2$ and the same threshold bracketing prove their
limiting level $\alpha$.
\end{proof}

\begin{proof}[Proof of Corollary~\ref{cor:feasible-final}]
The residual-based critical value is invariant to a separate mean shift
in each of the two samples. Hence the probability of
\( |\widehat q_{n,\alpha}-q_{n,\alpha}^{\mathrm{NR}}|>
\epsilon s_{0,n}\) is identical for all mean translations and tends
to zero. Apply Proposition~\ref{prop:spectral-random-power} with the
uniform comparison of Corollary~\ref{cor:grouped-transfer}, take the
supremum over normalized mean differences, and let \(n\to\infty\), followed by
\(\epsilon\downarrow0\). Under the null the continuous Gaussian
reference rejects with probability exactly \(\alpha\).
\end{proof}

\section{Validity of the split-linear benchmark}\label{app:split-validity}

The split linear statistic used in the numerical comparisons is the
linear-kernel specialization of the cross-MMD construction of
\citet{ShekharKimRamdas2022}, equations (2) and (4), which develops the
cross-statistic approach of \citet{KimRamdas2024}. The published null theorem assumes equality of distributions. The following
proposition gives the equal-mean, unequal-covariance justification used
in our simulations.

Divide each sample deterministically into an evaluation half and a training
half. Denote the evaluation observations by
$\bm{Y}_{1j}^{(a)}$, $j=1,\ldots,N_1$, and
$\bm{Y}_{2j}^{(a)}$, $j=1,\ldots,N_2$, and let
$\bar{\bm{Y}}_1^{(b)}$ and $\bar{\bm{Y}}_2^{(b)}$ be the corresponding
training means. Put $\bm{v}=\bar{\bm{Y}}_1^{(b)}-\bar{\bm{Y}}_2^{(b)}$ and set
\[
 U_{1i}=\bm{v}^\top \bm{Y}_{1i}^{(a)},\quad U_{2j}=\bm{v}^\top \bm{Y}_{2j}^{(a)},\qquad
 \widehat V=\frac1{N_1^2}\sum_{i=1}^{N_1}(U_{1i}-\bar U_1)^2
             +\frac1{N_2^2}\sum_{j=1}^{N_2}(U_{2j}-\bar U_2)^2.
\]
The statistic is $Z_{\rm split}=(\bar U_1-\bar U_2)/\sqrt{\widehat V}$,
and the nominal rule rejects when $Z_{\rm split}>z_{1-\alpha}$.
The within-half variances use divisor $N_i$, as in the cited construction;
the two group variances are estimated separately.

\begin{prop}[Validity of the split linear benchmark]
\label{prop:split-linear-validity}
Suppose the two samples are independent, observations are identically
distributed within each group, and $\bm{\mu}_1=\bm{\mu}_2$. Assume that for a fixed
$K<\infty$, writing $\bm{Y}_{ir}$ for an observation from group $i$, every
deterministic $\bm{u}$ and $i\in\{1,2\}$ satisfy
\[
 E\{\bm{u}^\top(\bm{Y}_{ir}-\bm{\mu}_i)\}^4
 \le K(\bm{u}^\top\bm{\Sigma}_i \bm{u})^2.
\]
Let $N_{\min}=\min(N_1,N_2)\to\infty$, and suppose
$\Pr\{\bm{v}^\top(\bm{\Sigma}_1/N_1+\bm{\Sigma}_2/N_2)\bm{v}>0\}\to1$.
Then $Z_{\rm split}\Rightarrow N(0,1)$, without a restriction on the
dimension or equality of the two covariance matrices.
\end{prop}

\begin{proof}
Condition on the training observations and write
$a_i=\bm{v}^\top\bm{\Sigma}_i \bm{v}/N_i$ and $V=a_1+a_2$.
On $V>0$, the testing numerator is a sum of independent centered variables
with variance $V$. The uniform fourth-moment bound and Cauchy--Schwarz give
the conditional Berry--Esseen ratio at most
\[
 \frac{\sqrt K}{\sqrt{N_{\min}}}
 \frac{a_1^{3/2}+a_2^{3/2}}{(a_1+a_2)^{3/2}}
 \le\frac{\sqrt K}{\sqrt{N_{\min}}}.
\]
Expanding each sample variance into the average squared projection
minus the squared sample mean gives
\[
 E\{(\widehat V/V-1)^2\mid \bm{v}\}\le C_K/N_{\min}.
\]
These bounds are uniform in $\bm{v}$ and the dimension. Conditional
studentization, followed by integration over the independent training
sample and the assumed nondegeneracy event, proves the assertion.
\end{proof}

\section{Simulation designs and further diagnostics}\label{app:numerical-details}
\subsection{Population-reference calculations}
\label{app:boundary-design}
The boundary experiment uses ten base covariance configurations and
28 signal configurations, yielding 84 power entries at the three
significance levels. Within each base configuration, the 10,000 Gaussian
contrasts are independent; directions and levels use common observations.
The ten base configurations use independent random streams. Sampling the
normalized mean contrast directly is exact for Gaussian observations.
For example, equal sample sizes with covariance matrices
$(2/3)\bm\Omega$ and $(4/3)\bm\Omega$ give the stated contrast covariance.
Only two Gaussian coordinates and the independent chi-squared norm of the
remaining coordinates are needed for this experiment. No data from the
earlier simulations are reused. The seed is 2026110801.

Two further configurations separate directional invariance from
shift-only validity. If $p=5$ and $\bm\Omega=\bm I_5$, the alternative
law is invariant to signal orientation, but the shift-only grid discrepancy
at $\kappa=1$ is $0.1694$. If instead
$\bm\Omega=\operatorname{diag}(\bm I_5,\bm0_{195})$, support and
null-space signals of the same norm have different distributions. Their
grid discrepancy is $0.1694$, while the shift-only formula is exact for
the null-space signal. A flat positive spectrum is therefore different
from a flat ambient spectrum when zero eigenvalues are present.
Table~\ref{tab:boundary-diagnostics} records selected outcomes.
The saved files provide all 84 powers, 30 null-size checks, CDF-grid
summaries, integration checks and numerical environment information.

The CDFs and the central-reference critical values are calculated
numerically from the exact Gaussian laws. We evaluate the CDF discrepancy
on the fixed standardized grid $\{-6,-5.975,\ldots,12\}$, so the
reported value approximates the discrepancy over that grid and is not a
certified Kolmogorov supremum. In the spiked model, 512-node
Gauss--Hermite integration agrees with 1,024 nodes to within
$2.24\times10^{-10}$ over all grids. Independent adaptive-quadrature
checks at the selected critical and maximum-gap points differ by at most
$1.87\times10^{-11}$. We also generate 10,000 independent Gaussian
mean contrasts per base covariance and compare empirical rejection rates
with the calculated powers at $\alpha\in\{0.01,0.05,0.10\}$. Directions
and significance levels share the same base draws. These are oracle
experiments and do not assess covariance estimation.

\begin{table}[!htbp]
\centering\small\setlength{\tabcolsep}{3.5pt}
\begin{threeparttable}
\caption{Selected Gaussian-reference boundary diagnostics.}
\label{tab:boundary-diagnostics}
\begin{tabular}{lrlrrrrr}
\toprule
& & & & & \multicolumn{3}{c}{Power (in \%)}\\
\cmidrule(lr){6-8}
Model & $p$ & Direction & $\Delta_{\bm\theta}$ & CDF gap & Exact & Shifted & Empirical \\
\midrule
Diffuse & 3200 & First axis & 0.0250 & 0.0060 & 25.85 & 25.34 & 25.63 (0.44) \\
Single spike & 3200 & Bulk & 0.0177 & 0.0065 & 19.06 & 18.56 & 19.16 (0.39) \\
Single spike & 3200 & Leading & 1.0001 & 0.1756 & 26.00 & 18.56 & 25.68 (0.44) \\
Isotropic & 5 & First axis & 0.6325 & 0.1694 & 23.34 & 16.14 & 23.36 (0.42) \\
Rank five & 200 & Null space & 0.0000 & 0.0000 & 16.14 & 16.14 & 15.78 (0.36) \\
Rank five & 200 & Support & 0.6325 & 0.1694 & 23.34 & 16.14 & 23.35 (0.42) \\
\bottomrule
\end{tabular}
\begin{tablenotes}[flushleft]
\footnotesize
\item[] \textit{Note.} $\kappa=1$ and $\alpha=0.05$. The CDF discrepancy is maximized over the stated finite grid, not over all real thresholds. Exact, shifted-null and empirical powers are in percentages; parentheses give binomial Monte Carlo standard errors in percentage points, based on 10,000 replications per base covariance.
\end{tablenotes}
\end{threeparttable}
\end{table}

In the first study, population covariance matrices determine all critical
values. Standard normal contrast noise is generated from its exact law.
For the centered chi-squared errors, sample sums are generated using the
corresponding chi-squared distributions. Student sample means are computed
from individual independent observations. In every case the resulting
noise has exactly the sampling law specified in the design.

The 132 configurations comprise 12 base error/sample-size/dimension
combinations and 11 spectrum/mean configurations per base combination:
four diffuse-spectrum configurations, and seven single-spike
configurations consisting of the null and two directions at each of
three nonzero signal strengths. Every configuration contains 10,000
independent replications. Reusing its noise across shifts and reference
rules induces dependence between the reported estimates, not between
replications within an estimate. The code reports all 396
calibration-specific rows.

The limiting power curves use deterministic adaptive integration, and
the finite Gaussian reference uses 256-node Gauss--Hermite integration.
Thus the number of quadrature nodes is not a Monte Carlo replication
count. Empirical CDF discrepancies include sampling error. The full
reference is computed at the generating mean; its alternative dependence
is used to predict power, whereas every rejection threshold is computed
under the null.

\begin{figure}[t]
\centering
\includegraphics[width=.86\textwidth]{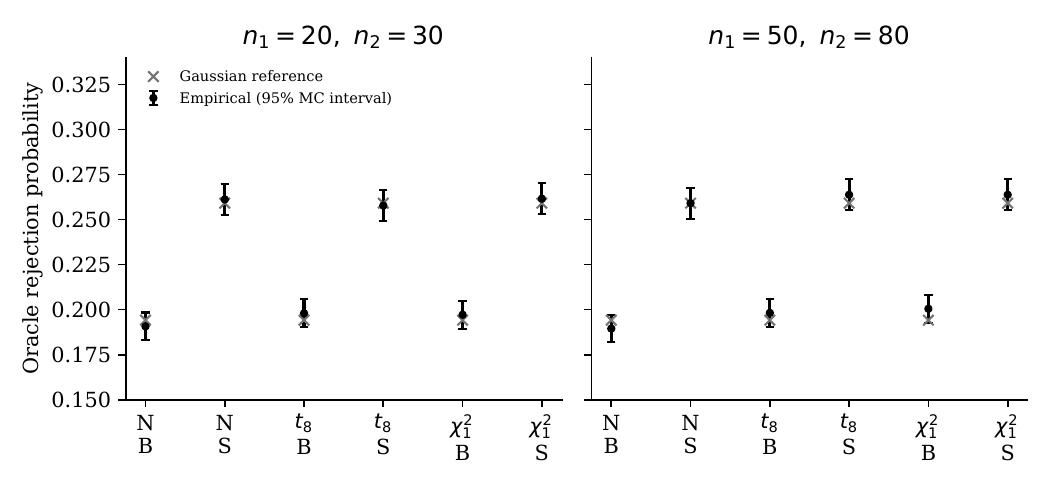}
\caption{Empirical oracle rejection probabilities and finite Gaussian
reference predictions at $p=800$ and $\kappa=1$. N denotes normal errors;
B and S denote bulk and spike directions. Error bars are $\pm1.96$ Monte
Carlo standard errors.}
\label{fig:finite-power}
\end{figure}
\subsection{Signal information and rank-one approximation errors}
\label{app:information-numerics}
This additional study uses two independent batches of 10,000 Gaussian
contrasts, one with covariance $\diag(1,2,3)$ and one with variance one.
The seed is 202611080102; independent child streams generate the two
batches. Within each batch, the contrasts are shared across signal
settings, and Monte Carlo standard errors for differences are calculated
from the paired rejection indicators. These experiments use known
population covariances and do not examine covariance estimation.

For the rank-three model, the signals are
$c(1/\sqrt2,0,1/\sqrt2)^\top$ and $c(0,1,0)^\top$, with
$c\in\{0.5,1,2,3\}$. Their squared norms and weighted energies are
$c^2$ and $\Delta=c^2/7$, while the third cumulants are
$288+120c^2$ and $288+96c^2$.
Population CDFs are evaluated by conditioning on the third Gaussian
coordinate and integrating the first two over a disk. The substitution
$(x_1,x_2)=\sqrt q\sin t(\cos\phi,\sin\phi)$,
$0\le t\le\pi/2$, $0\le\phi\le2\pi$, gives a smooth integrand and
Jacobian $q\sin t\cos t$. Gauss--Legendre orders $(96,192)$ are used.
Doubling both orders changes the tabulated probabilities by at most
$6.91\times10^{-14}$; conditioning instead on the first coordinate
and using orders $(160,320)$ changes them by at most
$5.40\times10^{-14}$. These are numerical checks, not certified
interval bounds. The saved CDF diagnostics use the grid
$\{0,0.25,\ldots,80\}$; their maxima are not certified Kolmogorov distances.

\begin{table}[!htbp]
\centering\small
\begin{threeparttable}
\caption{Powers for equal-norm, equal-$\Delta$ signal pairs.}
\label{tab:equal-information}
\begin{tabular}{rrrrrrrr}
\toprule
& & \multicolumn{2}{c}{Exact power} & \multicolumn{2}{c}{Empirical power} & \multicolumn{2}{c}{Paired difference}\\
\cmidrule(lr){3-4}\cmidrule(lr){5-6}\cmidrule(lr){7-8}
$c$ & $\Delta$ & Exact $a$ & Exact $b$ & MC $a$ & MC $b$
& MC diff. & SE(diff.)\\
\midrule
0.5 & 0.0357 & 5.59 & 5.63 & 5.72 & 5.96 & $-0.24$ & 0.18\\
1.0 & 0.1429 & 7.49 & 7.68 & 7.60 & 7.99 & $-0.39$ & 0.26\\
2.0 & 0.5714 & 16.76 & 17.59 & 16.57 & 18.08 & $-1.51$ & 0.45\\
3.0 & 1.2857 & 35.82 & 37.10 & 35.72 & 37.05 & $-1.33$ & 0.62\\
\bottomrule
\end{tabular}
\begin{tablenotes}[flushleft]
\footnotesize
\item[] \textit{Note.} $\alpha=0.05$. Powers and differences are percentages. The last two columns report the empirical difference $a-b$ and its paired Monte Carlo standard error, both in percentage points. Each empirical entry uses 10,000 contrasts.
\end{tablenotes}
\end{threeparttable}
\end{table}

For the rank-one model, $a\in\{0.01,0.02,0.05,0.1,0.2,0.4,0.8,1\}$,
$\Delta=a^2$, and $q_\alpha$ is the $(1-\alpha)$ quantile of $\chi_1^2$.
The exact power is
$1-\Phi(\sqrt{q_\alpha}-a)+\Phi(-\sqrt{q_\alpha}-a)$;
the shifted-null power is $\Pr(\chi_1^2>q_\alpha-a^2)$, taken as one
when $q_\alpha\le a^2$. Table~\ref{tab:rankone-information-rates}
uses these analytic formulas. Complementary paired Monte Carlo
calculations, null checks at $\alpha\in\{0.01,0.05,0.1\}$, and all
unrounded outputs are retained in the separate computational archive.
Monte Carlo estimates are not used to resolve the smallest displayed
population discrepancies.

\begin{table}[!htbp]
\centering\small
\begin{threeparttable}
\caption{Analytic rank-one discrepancies.}
\label{tab:rankone-information-rates}
\begin{tabular}{rrrrr}
\toprule
$a$ & $\Delta$ & CDF gap at $a^2$ & Power error & First order\\
\midrule
0.02 & 0.0004 & 0.0159534 & 0.0000338921 & 0.0000338923\\
0.10 & 0.0100 & 0.0792597 & 0.000847167 & 0.000847308\\
0.20 & 0.0400 & 0.155422 & 0.00338674 & 0.00338923\\
0.40 & 0.1600 & 0.288145 & 0.0135019 & 0.0135569\\
\bottomrule
\end{tabular}
\begin{tablenotes}[flushleft]
\footnotesize
\item[] \textit{Note.} $\alpha=0.05$. The CDF gap is evaluated at $a^2$. The power error is exact power minus shifted-null power; its first-order expression is $\Delta\phi(\sqrt{q_{0.05}})(q_{0.05}-1)/\sqrt{q_{0.05}}$. All entries are probabilities, not percentages.
\end{tablenotes}
\end{threeparttable}
\end{table}

For a numerical illustration of the power-estimation lower bound, take
$\alpha=0.05$, $a_0=0.5$ and $a_1=1.5$ in the rank-one Gaussian
experiment. The local strengths are $0.17678$ and $1.59099$, and the
powers are $0.07910$ and $0.32304$. Their gap is $0.24394$.
Consequently, for at least one of the two parameter values, every
power estimator has probability at least $\Phi(-0.5)=0.30854$ of an
absolute error of at least $0.12197$. The corresponding lower bound
on the maximum mean absolute error is $0.07527$.
These values evaluate the analytic two-point bound; they are not
Monte Carlo results or a comparison of testing procedures.

\paragraph{A known-covariance diagnostic illustration.}
To assess the numerical size of the upper bound in
\eqref{eq:diagnostic-estimated}, we use four specified population
covariances: $(1)$, $\diag(1,2,3)$, $\bm I_{800}$ and
$\diag(\sqrt{800},1,\ldots,1)$. For each covariance, we independently
generate $B=100{,}000$ draws of the standardized null reference
$Q=(T_0^*-\tr\bm\Omega)/\sqrt{2A}$. We compute
$\omega_{\widehat H}(u)$ by maximum sorted-window counts and minimize
the displayed expression over $121$ geometrically spaced values of
$u$ from $10^{-4}$ to $5$, and take the smaller of this result and
the universal bound $1/2$. Counts in closed intervals give a
conservative version of the empirical concentration function. With
$\gamma=0.01/4$ for each covariance, the DKW bound is
$e_H=0.0057813$. A union bound gives simultaneous coverage of at
least $99\%$ over all four covariances and all reported $K$.

\begin{table}[!htbp]
\centering\small
\begin{threeparttable}
\caption{Upper bounds for shifted-null CDF error with known covariance.}
\label{tab:known-covariance-diagnostic}
\begin{tabular}{lrrrr}
\toprule
Covariance & $K=0.1$ & $K=0.25$ & $K=0.5$ & $K=1$\\
\midrule
$(1)$ & 50.00 & 50.00 & 50.00 & 50.00\\
$\diag(1,2,3)$ & 47.35 & 50.00 & 50.00 & 50.00\\
$\bm I_{800}$ & 11.36 & 16.06 & 20.86 & 26.70\\
$\diag(\sqrt{800},1,\ldots,1)$ & 37.29 & 48.16 & 50.00 & 50.00\\
\bottomrule
\end{tabular}
\begin{tablenotes}[flushleft]
\footnotesize
\item[] \textit{Note.} The simulation-assisted bounds are for the supremum error over $\kappa\leq K$, with known population covariance. Entries are bounds multiplied by $100$, not rejection percentages or estimates of the actual error.
\end{tablenotes}
\end{threeparttable}
\end{table}

The bound is smaller for the diffuse covariance and weakens with the
allowed signal strength. It remains conservative, especially for the
concentrated spectra; even the smallest entry does not certify an
error tolerance of one percentage point. Several entries reduce to
the universal bound $1/2$, and several untruncated finite-grid minima
occur at the lower endpoint of the $u$ grid. These calculations illustrate
both the computability and the limitations of the diagnostic.
They use auxiliary draws from known null references, not outer
two-sample datasets, and do not validate a diagnostic with estimated
covariance. The probability statement refers to the auxiliary
simulation used to construct $\widehat H$.
\paragraph{Numerical check of the sharpness construction.}
For the family in Proposition~\ref{prop:shift-quarter-sharp}, the
discrepancy at the support boundary is the one-dimensional integral
\[
 D_\delta=\int_0^{2c^2/\sqrt\delta}
 \operatorname{erf}\!\left[
   \left\{\frac{2\sqrt\delta\,t-\delta t^2/c^2}
                    {2\lambda_\delta}\right\}^{1/2}\right]
 \phi(t)\,dt.
\]
We evaluate this expression for $\kappa_0\in\{0.25,1,4\}$ and
$\delta\in\{10^{-2},10^{-3},\ldots,10^{-12}\}$. Adaptive
quadrature is checked against 1024-node Gauss--Legendre quadrature after
the substitution $t=u^2$, with a separate 512-node check. The largest
discrepancy between the two methods is $2.53\times10^{-15}$; changing
the Gauss--Legendre order changes the results by at most
$2.03\times10^{-15}$. If the integral extends beyond $12$, the omitted
Gaussian tail is at most $1.78\times10^{-33}$. These numerical checks
support the analytic limit
$D_\delta/\delta^{1/4}\to c_*=0.4638648043$.
They evaluate a specified-threshold lower bound for the Kolmogorov
distance and do not numerically estimate its supremum.

\begin{figure}[t]
\centering
\includegraphics[width=\textwidth]{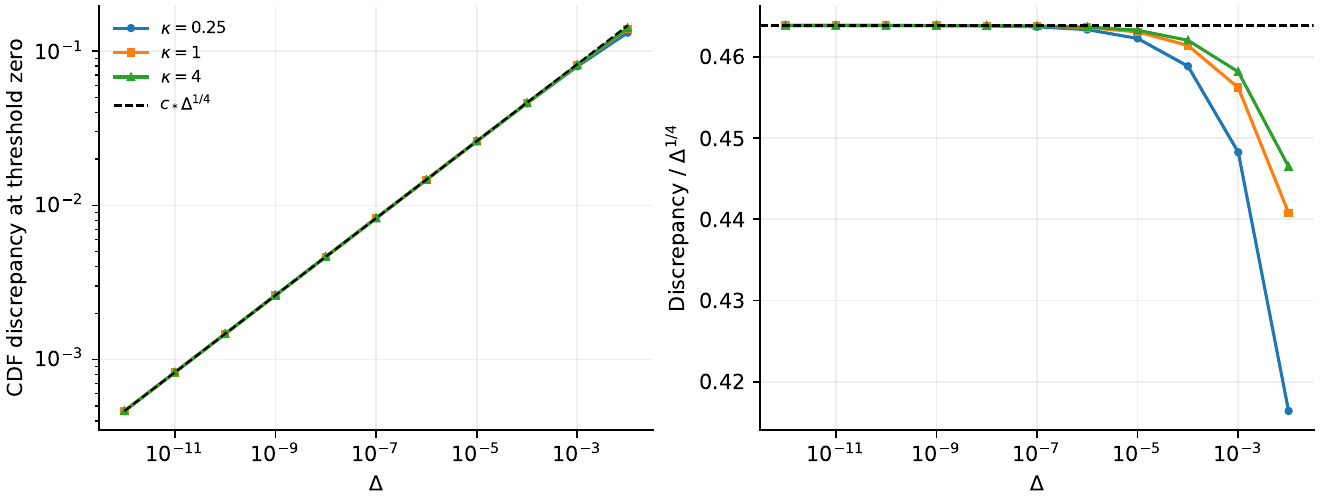}
\caption{Support-boundary discrepancy for the fixed-effect construction
in Proposition~\ref{prop:shift-quarter-sharp}. The right panel shows
$D_\delta/\delta^{1/4}$ approaching $c_*=0.4638648043$.
All curves use deterministic quadrature, and the plotted discrepancy
is a lower witness for the Kolmogorov distance, not its computed supremum.}
\label{fig:sharp-quarter}
\end{figure}

\subsection{Nonproportional covariance matrices}
\label{app:expanded-design}
We examine feasible two-sample calibration with several covariance spikes
and covariance matrices that are neither proportional nor simultaneously
diagonalizable. We take $p\in\{100,400\}$ and sample sizes $(40,60)$ or
$(80,120)$, using the normalized mean difference in Section~\ref{sec:setup}.

For an explicit reproducible covariance construction, the first three
coordinates in sample $i$ have covariance
\[
\sqrt p\,\bm{R}_i\operatorname{diag}(w_{i1},w_{i2},w_{i3})\bm{R}_i^\top,
\qquad
(w_{ij})=\begin{pmatrix}2&0.3&0.7\\0.4&1.5&0.9\end{pmatrix}.
\]
Here $\bm{R}_1=\bm{I}_3$ and
$\bm{R}_2=\bm{R}_{12}(0.61)\bm{R}_{23}(0.37)$, where
$\bm{R}_{ab}(t)$ is a coordinate-plane rotation through $t$ radians.
The remaining coordinates are divided at $\lfloor p/2\rfloor$ into two
diagonal variance blocks. Their variances are $(1,0.5)$ in sample 1 and
$(0.6,2)$ in sample 2. This construction changes both eigenvalues and
eigenspaces between the samples. The nonzero commutator
$\bm{\Sigma}_1\bm{\Sigma}_2-\bm{\Sigma}_2\bm{\Sigma}_1$ confirms that
the covariance matrices do not share an eigenbasis.

Four scenarios use Gaussian observations, covering both dimensions and
both sample-size pairs. Two further scenarios use independent Rademacher
latent coordinates transformed by the stated covariance factors, at the
larger sample sizes. Gaussian sample means and Helmert residual contrasts
are generated from their exact joint law. For Rademacher observations,
sample means and paired-difference covariance estimates use the same data,
preserving their finite-sample dependence. Alternatives have
$\|\bm{\theta}\|^2=s_0$, corresponding to $\kappa=1$. The spike signal is
aligned with a largest eigenvector of $\bm{\Omega}$; the bulk signal lies
along the fourth coordinate.

We compare the population Gaussian-mixture oracle, estimated two-cumulant
calibration, and spectral calibration with either a normal or a gamma
remainder. The two primary spectral methods use identical trace estimators,
retained eigenvalues and the predetermined capped threshold, followed by
the variance-budget projection. The original unprojected normal- and
gamma-tail calibrations are retained as two additional diagnostics.
Each conditional reference quantile uses 3,000 draws. Population covariance
parameters enter only the oracle critical value and the data-generating
mechanism. The oracle quantile uses 2,000,000 mixture draws, with its Monte
Carlo uncertainty recorded separately. Each scenario has 10,000 independent
evaluation replications, shared across methods and across the three signal
configurations by location shifts.

For adjusted power, 10,000 null-training samples are generated independently
of the 10,000 evaluation samples. Within each phase, observations are shared
across methods and alternatives by location shifts. The adjusted cutoff is
the empirical training $0.95$ quantile of
$(T-\widehat q_{.95})/\sqrt{2\widehat A}$. Both phases are independently
resampled in each of 1,000 bootstrap repetitions, retaining methodwise
pairing. This null-law adjustment is a diagnostic and is not used to select
the operational threshold.

\begin{table}[!htbp]
\centering\small\setlength{\tabcolsep}{4pt}
\begin{threeparttable}
\caption{Empirical powers (in \%) for the feasible two-sample study.}
\label{tab:expanded-power}
\begin{tabular}{lrrlrrr}
\toprule
\multicolumn{7}{l}{Panel A: Raw power} \\
\addlinespace
Model & $n_1$ & $n_2$ & Direction & 2C & Proj. G & Proj. $\Gamma$ \\
\midrule
N & 40 & 60 & Spike & 26.01 (0.44) & 22.89 (0.42) & 22.94 (0.42) \\
 &  &  & Bulk & 20.93 (0.41) & 16.30 (0.37) & 16.31 (0.37) \\
\addlinespace
N & 80 & 120 & Spike & 26.97 (0.44) & 25.30 (0.43) & 25.28 (0.43) \\
 &  &  & Bulk & 20.36 (0.40) & 17.93 (0.38) & 17.79 (0.38) \\
\addlinespace
R & 80 & 120 & Spike & 26.84 (0.44) & 23.87 (0.43) & 23.82 (0.43) \\
 &  &  & Bulk & 20.22 (0.40) & 15.97 (0.37) & 16.01 (0.37) \\
\midrule
\multicolumn{7}{l}{Panel B: Independently size-adjusted power} \\
\addlinespace
Model & $n_1$ & $n_2$ & Direction & 2C & Proj. G & Proj. $\Gamma$ \\
\midrule
N & 40 & 60 & Spike & 22.55 (0.63) & 22.46 (0.60) & 22.45 (0.67) \\
 &  &  & Bulk & 15.83 (0.77) & 15.72 (0.71) & 15.65 (0.80) \\
\addlinespace
N & 80 & 120 & Spike & 23.57 (0.75) & 23.46 (0.69) & 23.66 (0.68) \\
 &  &  & Bulk & 15.72 (0.78) & 15.61 (0.76) & 15.77 (0.74) \\
\addlinespace
R & 80 & 120 & Spike & 22.62 (0.57) & 22.92 (0.62) & 23.03 (0.68) \\
 &  &  & Bulk & 14.78 (0.59) & 14.77 (0.62) & 14.86 (0.73) \\
\bottomrule
\end{tabular}

\begin{tablenotes}[flushleft]
\footnotesize
\item[] \textit{Note.} $p=400$ and $\kappa=1$. Parentheses give Monte Carlo standard errors in percentage points. Each phase uses 10,000 independent samples. Raw standard errors are binomial; adjusted standard errors include independent null-training and evaluation uncertainty through 1,000 bootstrap resamples. The adjustment is a diagnostic based on the generating null law.
\end{tablenotes}
\end{threeparttable}
\end{table}

For the power curves, let $r_0$ and $r_1(\bm{u})$ be the saved rejection margins under the null and under a unit-strength signal in direction $\bm{u}$, standardized by $\sqrt{2\widehat A}$. Write $a=\widehat A/A>0$. Because residual estimates and their reference critical values are unchanged by a location shift, the margin at strength $\kappa$ is exactly
\[
 r_\kappa(\bm{u})=r_0+\sqrt\kappa\{r_1(\bm{u})-r_0\}
                         +\frac{\kappa-\sqrt\kappa}{\sqrt a}.
\]
The same independent null-training cutoff is used at every strength. The reconstructed rates at $\kappa=0$ and $1$ equal the original rates; a direct data-level identity check is recorded with the computational code. The curve CSV reports binomial standard errors for raw power and conditional standard errors for adjusted power. The latter condition on the saved null-training sample and do not include training uncertainty. The full two-phase bootstrap standard errors at $\kappa=1$ remain in Table~\ref{tab:expanded-power}.

The threshold cap rate is one in all six designs. This follows directly
from $m_{\min}\le79$ and the $m_{\min}^{-1/8}$ term, so increasing the outer
replication count does not evaluate an uncapped adaptive regime. Projection
rates are $25.30\%$ for Gaussian observations with $p=400$, $(40,60)$ and
$15.96\%$ for Rademacher observations with $p=400$, $(80,120)$; the total
number of projected evaluation replicates across the other four scenarios
is $1$. Projected normal and gamma reference laws coincide whenever the
retained weights exhaust the variance budget. The recorded mechanism
diagnostics report the positive-remainder rate and mean remaining variance
fraction, preventing this coincidence from being misread as evidence that
the two remainder approximations are generally interchangeable.

The population oracle uses 2,000,000 mixture draws for the quantile and
300,000 separate Gaussian draws for power. Repeated eigenvalues are grouped
after rounding at 11 decimal places solely for computation; the maximal
rounding discrepancy and oracle Monte Carlo errors are recorded. Feasible
references use 3,000 draws per dataset, and no population spectral quantity
enters their critical values. Calibration timings include trace estimation,
spectral calculation and conditional simulation where applicable, and
exclude data generation; they are descriptive measurements rather than
hardware-independent complexity claims.

The computational code specifies fixed seeds, generates replication scores and
estimated-spectrum diagnostics, and checks the variance-budget identities.
Current results use 10,000 replications in each phase; smaller runs were
used only as verification fixtures.
\subsection{External-method benchmarks under unequal covariances}
\label{app:benchmark-design}

We compare two-cumulant calibration, the unprojected and variance-budget
projected spectral references with Gaussian or Gamma remainders, the
approximate randomization test of \citet{WangXu2022}, and a split linear
test. The latter uses the linear-kernel cross-statistic construction of
\citet{ShekharKimRamdas2022}, with its heteroscedastic equal-mean
justification given in Proposition~\ref{prop:split-linear-validity}.

For Wang--Xu we implement the original procedure in Section~3, pages
1125--1126: the observed statistic is the full-sample Chen--Qin statistic,
while the reference statistics are Chen--Qin statistics of independently
signed pair differences $(\bm{X}_{i,2r}-\bm{X}_{i,2r-1})/2$, using the pair counts in
their denominators. There are $B=999$ Rademacher randomizations and the
plus-one rule $(1+\#\{T_b^*\ge T_{\rm CQ}\})/(B+1)\le0.05$.
This is their approximate randomization test, not an exact permutation
test. The split linear rule uses deterministic half splits, separate
group variance estimates with the divisors in equation~(4) of the cited
cross-MMD construction, and the upper standard-normal critical value.

All procedures receive the same full observations. Coordinates are
independent Gaussian or standardized Rademacher variables multiplied by
the square roots of their group variances. Put $a=n_2/(n_1+n_2)$,
$b=n_1/(n_1+n_2)$ and $r_j=1/2$ for odd $j$ and $r_j=2$ for even $j$.
The diagonal group variances are
\[
 (\bm{\Sigma}_1)_{jj}=\omega_j/(a+br_j),\qquad
 (\bm{\Sigma}_2)_{jj}=r_j(\bm{\Sigma}_1)_{jj},
\]
so that the covariance of
$\bm{d}=\sqrt{n_1n_2/(n_1+n_2)}(\bar{\bm{X}}_1-\bar{\bm{X}}_2)$ is
$\bm{\Omega}=\operatorname{diag}(\omega_j)$, although the two covariances
are not proportional. The diffuse spectrum has $\omega_j=1$; the
three-spike spectrum is
$(\sqrt p,0.7\sqrt p,0.4\sqrt p,1,\ldots,1)$.
The six scenarios are
\[
\begin{array}{ccll}
\text{ID}&(n_1,n_2),p&\text{distribution}&\text{spectrum}\\
\mathrm{G1}&(40,60),200&\text{Gaussian}&\text{diffuse}\\
\mathrm{G2}&(40,60),800&\text{Gaussian}&\text{three spikes}\\
\mathrm{G3}&(80,120),800&\text{Gaussian}&\text{three spikes}\\
\mathrm{R1}&(40,60),200&\text{Rademacher}&\text{diffuse}\\
\mathrm{R2}&(80,120),200&\text{Rademacher}&\text{three spikes}\\
\mathrm{R3}&(80,120),800&\text{Rademacher}&\text{three spikes}.
\end{array}
\]
Alternatives have $\|\bm{\theta}\|^2=\kappa\sqrt{2\operatorname{tr}(\bm{\Omega}^2)}$
for $\kappa\in\{1,2\}$, directed along $\bm{e}_1$ or $\bm{e}_4$. The corresponding
mean shift is applied to group~1. The names ``top'' and ``bulk'' designate
these coordinate choices; they are equivalent directions in the diffuse
population model, although their group variance ratios differ.

The spectral estimates use the same-data, disjoint pair construction for
both distributions, so the mean statistic and the reference estimates are computed from the
same data. The threshold is given by
\eqref{eq:implemented-spectral-threshold}. Nontrivial spectral quantiles use 1,999 reference
draws; if no eigenvalue is retained the Gaussian and two-cumulant
quantiles are evaluated analytically. Gaussian and Gamma variants share
their simulated retained component. When the variance-budget projection
is active the projected remainder variance is zero, so both projected
variants use the same scaled retained component. The two projected
references coincide whenever this projection is active. For the present
sample sizes, the threshold rule is identically $0.5$, since
$m_{\min}\le256$. We report the retained rank, projection frequency and
remaining-variance fraction to distinguish these effects.

Each scenario uses 10,000 evaluation replications, with the null and four
alternatives paired within a replication. Nominal results include
resampling and reference-quantile randomness. We additionally use 10,000
independent null replications to estimate the 95\% quantile of each
standardized critical-value margin. For the spectral rules the margin is
$(\|\bm{d}\|^2-\widehat q)/\sqrt{2\widehat A}$; for Wang--Xu it is
$(T_{\rm CQ}-\widehat q_{\rm WX})/\widehat{\operatorname{sd}}(T_b^*)$;
for the split linear rule it is $Z_{\rm split}-z_{.95}$.
For the spherical test, the margin is
$(T-\widehat q_B^{\mathrm{sph}})/\widehat\sigma_{\mathrm{sph}}$, where
$\widehat q_B^{\mathrm{sph}}$ is the order statistic equivalent to the
plus-one decision and
\[
 \widehat\sigma_{\mathrm{sph}}^2=
 \frac{2}{N(N+2)}\left\|
 \bm X^\top\bm X-\frac{\tr(\bm X^\top\bm X)}N\bm I_N
 \right\|_{\mathrm F}^2.
\]
This is the conditional spherical variance; it is recomputed from the
observed augmented Gram matrix under every signal shift. It uses no
unobserved null-data normalization.
The calibration order is $\lceil0.95(R_0+1)\rceil$ among $R_0=10\,000$
null margins. These size adjustments use knowledge of each simulation
null and serve only to compare power, not as additional feasible tests.
Held-out null sizes are reported with the adjusted power.
Binomial Monte Carlo standard errors are reported for the nominal
rejection rates. For size-adjusted power, a bootstrap with 1,000 resamples
accounts for both the null-training and evaluation samples. Paired differences use the same evaluation datasets. For comparisons
involving the spherical test, common bootstrap indices are also used
across methods in the null-training sample; the evaluation-bootstrap
variance is integrated analytically. These two-phase errors include
training uncertainty. Previously recorded individual-method bootstrap
errors are retained for the other procedures. Fixed seeds, scores and
timing records are retained in the separate computational archive. Adding the spherical procedure replays the
existing 60,000 training and 60,000 evaluation datasets, and reproduces
every archived two-cumulant score and trace diagnostic exactly; it does
not enlarge the independent sample used for each reported rate.

\begin{table}[!htbp]
\centering\small\setlength{\tabcolsep}{3.5pt}
\begin{threeparttable}
\caption{Nominal powers (in \%).}
\label{tab:benchmark-nominal-power}
\begin{tabular}{lrrrrrrrr}
\toprule
\multicolumn{9}{l}{Panel A: Leading direction} \\
& & & \multicolumn{2}{c}{Unprojected} & \multicolumn{2}{c}{Projected} & \multicolumn{2}{c}{Benchmarks} \\
\cmidrule(lr){4-5}\cmidrule(lr){6-7}\cmidrule(lr){8-9}
Scenario & 2C & Sphere & G & $\Gamma$ & G & $\Gamma$ & WX & Split \\
\midrule
G1 & 58.69 & 57.93 & 15.64 & 15.64 & 52.40 & 52.40 & 58.00 & 40.43 \\
G2 & 49.72 & 45.11 & 12.94 & 12.94 & 46.39 & 46.39 & 46.39 & 30.38 \\
G3 & 49.06 & 45.94 & 39.43 & 39.43 & 43.94 & 43.94 & 46.49 & 29.02 \\
\addlinespace
R1 & 58.75 & 58.42 & 16.59 & 16.59 & 52.26 & 52.26 & 58.04 & 40.47 \\
R2 & 47.47 & 44.61 & 44.05 & 44.07 & 44.05 & 44.07 & 45.06 & 28.43 \\
R3 & 48.57 & 45.64 & 38.99 & 38.99 & 42.80 & 42.80 & 46.08 & 28.39 \\
\midrule
\multicolumn{9}{l}{Panel B: Bulk direction} \\
& & & \multicolumn{2}{c}{Unprojected} & \multicolumn{2}{c}{Projected} & \multicolumn{2}{c}{Benchmarks} \\
\cmidrule(lr){4-5}\cmidrule(lr){6-7}\cmidrule(lr){8-9}
Scenario & 2C & Sphere & G & $\Gamma$ & G & $\Gamma$ & WX & Split \\
\midrule
G1 & 57.95 & 56.94 & 15.62 & 15.62 & 51.66 & 51.66 & 57.17 & 39.60 \\
G2 & 59.76 & 53.98 & 2.64 & 2.64 & 53.45 & 53.45 & 53.77 & 47.68 \\
G3 & 58.28 & 52.73 & 37.75 & 37.76 & 46.60 & 46.61 & 52.56 & 46.98 \\
\addlinespace
R1 & 58.79 & 58.27 & 16.34 & 16.34 & 52.55 & 52.55 & 57.78 & 40.45 \\
R2 & 56.54 & 52.80 & 49.99 & 49.83 & 49.99 & 49.83 & 52.08 & 44.71 \\
R3 & 58.42 & 53.55 & 38.00 & 38.00 & 46.66 & 46.66 & 52.74 & 45.58 \\
\bottomrule
\end{tabular}
\begin{tablenotes}[flushleft]\footnotesize
\item[] \textit{Note.} $\kappa=2$. Every method uses the same 10,000 evaluation datasets per scenario. Each method uses its operational nominal $5\%$ rule; no null adjustment is applied.
\end{tablenotes}
\end{threeparttable}
\end{table}

\begin{table}[!htbp]
\centering\small\setlength{\tabcolsep}{3.5pt}
\begin{threeparttable}
\caption{Independently size-adjusted powers (in \%).}
\label{tab:benchmark-adjusted-power}
\begin{tabular}{lrrrrrrrr}
\toprule
\multicolumn{9}{l}{Panel A: Leading direction} \\
& & & \multicolumn{2}{c}{Unprojected} & \multicolumn{2}{c}{Projected} & \multicolumn{2}{c}{Benchmarks} \\
\cmidrule(lr){4-5}\cmidrule(lr){6-7}\cmidrule(lr){8-9}
Scenario & 2C & Sphere & G & $\Gamma$ & G & $\Gamma$ & WX & Split \\
\midrule
G1 & 59.82 & 59.87 & 58.33 & 58.33 & 59.40 & 59.40 & 60.30 & 37.93 \\
G2 & 46.11 & 46.14 & 46.49 & 46.49 & 46.31 & 46.31 & 46.35 & 28.35 \\
G3 & 46.02 & 45.65 & 46.24 & 46.22 & 46.11 & 46.11 & 46.34 & 26.90 \\
\addlinespace
R1 & 57.37 & 57.01 & 56.24 & 56.24 & 57.35 & 57.35 & 57.86 & 38.14 \\
R2 & 44.81 & 44.52 & 44.27 & 44.52 & 44.27 & 44.52 & 44.88 & 26.02 \\
R3 & 43.98 & 44.22 & 43.67 & 43.68 & 43.53 & 43.53 & 44.04 & 26.42 \\
\midrule
\multicolumn{9}{l}{Panel B: Bulk direction} \\
& & & \multicolumn{2}{c}{Unprojected} & \multicolumn{2}{c}{Projected} & \multicolumn{2}{c}{Benchmarks} \\
\cmidrule(lr){4-5}\cmidrule(lr){6-7}\cmidrule(lr){8-9}
Scenario & 2C & Sphere & G & $\Gamma$ & G & $\Gamma$ & WX & Split \\
\midrule
G1 & 59.01 & 58.83 & 57.68 & 57.68 & 58.76 & 58.76 & 59.21 & 37.16 \\
G2 & 53.23 & 55.83 & 52.85 & 52.85 & 53.33 & 53.33 & 53.60 & 45.00 \\
G3 & 52.04 & 52.19 & 51.75 & 51.66 & 51.77 & 51.78 & 52.31 & 44.48 \\
\addlinespace
R1 & 57.53 & 57.07 & 56.26 & 56.26 & 57.36 & 57.36 & 57.66 & 37.86 \\
R2 & 51.59 & 52.59 & 50.56 & 50.74 & 50.56 & 50.74 & 51.83 & 41.82 \\
R3 & 48.23 & 50.30 & 47.97 & 47.97 & 48.26 & 48.27 & 48.63 & 43.62 \\
\bottomrule
\end{tabular}

\begin{tablenotes}[flushleft]
\footnotesize
\item[] \textit{Note.} $\kappa=2$. Each adjustment uses 10,000 independent null-training samples, followed by 10,000 independent evaluation samples. The adjustment uses the generating null law and is diagnostic. Method abbreviations are defined in Table~\ref{tab:external-benchmarks}.
\end{tablenotes}
\end{threeparttable}
\end{table}
\begin{table}[!htbp]
\centering\small\setlength{\tabcolsep}{3.5pt}
\begin{threeparttable}
\caption{Held-out sizes and Monte Carlo standard errors after null adjustment.}
\label{tab:benchmark-adjusted-size}
\begin{tabular}{lrrrrrrrr}
\toprule
\multicolumn{9}{l}{Panel A: Held-out size (in \%)} \\
& & & \multicolumn{2}{c}{Unprojected} & \multicolumn{2}{c}{Projected} & \multicolumn{2}{c}{Benchmarks} \\
\cmidrule(lr){4-5}\cmidrule(lr){6-7}\cmidrule(lr){8-9}
Scenario & 2C & Sphere & G & $\Gamma$ & G & $\Gamma$ & WX & Split \\
\midrule
G1 & 5.33 & 5.41 & 5.17 & 5.17 & 5.35 & 5.35 & 5.28 & 5.22 \\
G2 & 5.36 & 5.26 & 5.60 & 5.60 & 5.46 & 5.46 & 5.33 & 5.00 \\
G3 & 5.05 & 4.95 & 5.05 & 5.04 & 5.03 & 5.03 & 5.17 & 4.71 \\
\addlinespace
R1 & 4.68 & 4.62 & 4.64 & 4.64 & 4.79 & 4.79 & 4.81 & 5.05 \\
R2 & 4.92 & 5.09 & 4.85 & 4.88 & 4.85 & 4.88 & 5.06 & 4.48 \\
R3 & 4.48 & 4.65 & 4.53 & 4.53 & 4.42 & 4.42 & 4.41 & 4.97 \\
\midrule
\multicolumn{9}{l}{Panel B: Monte Carlo standard error (percentage points)} \\
& & & \multicolumn{2}{c}{Unprojected} & \multicolumn{2}{c}{Projected} & \multicolumn{2}{c}{Benchmarks} \\
\cmidrule(lr){4-5}\cmidrule(lr){6-7}\cmidrule(lr){8-9}
Scenario & 2C & Sphere & G & $\Gamma$ & G & $\Gamma$ & WX & Split \\
\midrule
G1 & 0.31 & 0.33 & 0.28 & 0.28 & 0.27 & 0.28 & 0.29 & 0.34 \\
G2 & 0.33 & 0.30 & 0.34 & 0.34 & 0.28 & 0.29 & 0.31 & 0.36 \\
G3 & 0.31 & 0.30 & 0.31 & 0.30 & 0.30 & 0.30 & 0.31 & 0.31 \\
\addlinespace
R1 & 0.28 & 0.28 & 0.30 & 0.31 & 0.30 & 0.30 & 0.30 & 0.35 \\
R2 & 0.31 & 0.35 & 0.29 & 0.31 & 0.29 & 0.33 & 0.29 & 0.33 \\
R3 & 0.28 & 0.28 & 0.32 & 0.32 & 0.28 & 0.28 & 0.26 & 0.28 \\
\bottomrule
\end{tabular}

\begin{tablenotes}[flushleft]
\footnotesize
\item[] \textit{Note.} The null adjustment is independent. Panel A reports percentages; Panel B reports percentage points. Standard errors use 1,000 bootstrap resamples and include both null-training and evaluation randomness. Method abbreviations are defined in Table~\ref{tab:external-benchmarks}.
\end{tablenotes}
\end{threeparttable}
\end{table}

The paired power differences in Table~\ref{tab:spherical-wx-paired}
use both the nominal decision rules and the independently adjusted
thresholds. Table~\ref{tab:benchmark-runtime-v105} reports a separate
profiling experiment on 200 common null datasets per setting, after ten
warm-up datasets and with randomized method order. Median spherical
runtimes are 1.34--3.56 ms, compared with 0.41--1.36 ms for Wang--Xu.
The respective auxiliary budgets are 1,999 and 999. These measurements
characterize the implementations used here on one host and do not suggest
a computational advantage for the spherical method.

\begin{table}[!htbp]
\centering\small\setlength{\tabcolsep}{4pt}
\begin{threeparttable}
\caption{Paired power differences between spherical and Wang--Xu tests.}
\label{tab:spherical-wx-paired}
\begin{tabular}{lrrrrrrrr}
\toprule
& \multicolumn{4}{c}{Leading direction} & \multicolumn{4}{c}{Bulk direction} \\
\cmidrule(lr){2-5}\cmidrule(lr){6-9}
& \multicolumn{2}{c}{Nominal} & \multicolumn{2}{c}{Adjusted} & \multicolumn{2}{c}{Nominal} & \multicolumn{2}{c}{Adjusted} \\
\cmidrule(lr){2-3}\cmidrule(lr){4-5}\cmidrule(lr){6-7}\cmidrule(lr){8-9}
Scenario & Diff. & MCSE & Diff. & MCSE & Diff. & MCSE & Diff. & MCSE \\
\midrule
G1 & -0.07 & 0.21 & -0.43 & 0.60 & -0.23 & 0.22 & -0.38 & 0.56 \\
G2 & -1.28 & 0.23 & -0.21 & 0.45 & 0.21 & 0.31 & 2.23 & 0.79 \\
G3 & -0.55 & 0.18 & -0.69 & 0.34 & 0.17 & 0.28 & -0.12 & 0.73 \\
\addlinespace
R1 & 0.38 & 0.22 & -0.85 & 0.49 & 0.49 & 0.21 & -0.59 & 0.40 \\
R2 & -0.45 & 0.19 & -0.36 & 0.42 & 0.72 & 0.27 & 0.76 & 0.79 \\
R3 & -0.44 & 0.19 & 0.18 & 0.40 & 0.81 & 0.27 & 1.67 & 0.72 \\
\bottomrule
\end{tabular}
\begin{tablenotes}[flushleft]\footnotesize
\item[] \textit{Note.} $\kappa=2$. Differences are spherical minus Wang--Xu, in percentage points. Nominal MCSE uses paired evaluation indicators. Adjusted MCSE includes the uncertainty in both 10,000 null-training and 10,000 evaluation datasets, preserving the pairing between methods in both stages. The 1,000-resample bootstrap integrates the conditional evaluation variance analytically. These are descriptive comparisons across the displayed configurations; no simultaneous superiority claim is made.
\end{tablenotes}
\end{threeparttable}
\end{table}

\begin{table}[ht]
\centering\small
\begin{threeparttable}
\caption{Median computation time (milliseconds).}
\label{tab:benchmark-runtime-v105}
\setlength{\tabcolsep}{4pt}
\begin{tabular}{lrrrrrrrr}
\toprule
Scenario & 2C & Raw G & Raw $\Gamma$ & Proj. G & Proj. $\Gamma$ & WX & Split & Sphere \\
\midrule
G1 & 0.181 & 0.546 & 0.521 & 0.522 & 0.527 & 0.414 & 0.144 & 1.410 \\
G2 & 0.280 & 0.846 & 0.824 & 0.819 & 0.822 & 0.541 & 0.172 & 1.706 \\
G3 & 0.548 & 0.876 & 0.847 & 0.851 & 0.849 & 1.276 & 0.226 & 3.404 \\
R1 & 0.170 & 0.522 & 0.502 & 0.505 & 0.506 & 0.405 & 0.137 & 1.344 \\
R2 & 0.246 & 0.553 & 0.568 & 0.551 & 0.569 & 0.918 & 0.154 & 2.630 \\
R3 & 0.560 & 0.901 & 0.867 & 0.868 & 0.855 & 1.357 & 0.238 & 3.564 \\
\bottomrule
\end{tabular}
\begin{tablenotes}[flushleft]\footnotesize
\item[] \textit{Note.} Each method is evaluated on the same 200 datasets per setting, following ten warm-up datasets. Method order is randomized within each dataset. Times include all required sample summaries and calibration, but exclude data generation. A single CPU thread is used. Spectral and spherical references use 1,999 draws; WX uses 999 randomizations. Projected methods are implemented separately and omit unused unprojected references. These implementation-specific timings are separate from the 10,000-replication size and power experiments.
\end{tablenotes}
\end{threeparttable}
\end{table}

\clearpage

\subsection{Gaussian spherical calibration}
\label{app:rotation-numerics}
We use five configurations: diffuse covariance at $p=200$ and
$p=20{,}000$ with $(n_1,n_2)=(20,30)$; single-spike covariance at
$p=800$ with $(20,30)$ and $(40,60)$; and single-spike covariance at
$p=20{,}000$ with $(20,30)$. Let $w_1=n_2/(n_1+n_2)$ and $w_2=1-w_1$.
For the prescribed diagonal entries $\omega_j$ of $\bm\Omega$, set
\[
 r_j=\begin{cases}1/2,&j\text{ odd},\\2,&j\text{ even},\end{cases}
 \qquad
 (\bm\Sigma_1)_{jj}=\frac{\omega_j}{w_1+w_2r_j},
 \qquad (\bm\Sigma_2)_{jj}=r_j(\bm\Sigma_1)_{jj}.
\]
All off-diagonal entries are zero. Thus the two covariances differ and
are nonproportional, while $w_1\bm\Sigma_1+w_2\bm\Sigma_2=\bm\Omega$.
The spiked model has $\omega_1=\sqrt p$ and $\omega_j=1$ for $j>1$;
the diffuse model has all $\omega_j=1$. Alternative signals are
$\bm\theta=\sqrt{\kappa\sqrt{2A}}\bm e_j$, with
$\kappa\in\{0.5,1\}$ and $j=1,2$. The term ``top'' denotes the first
coordinate in the diffuse case as well. Each configuration also includes
the null signal.

The Gaussian simulation uses the exact joint distribution of the
normalized mean contrast and independent within-group residuals, rather
than explicitly generating every raw coordinate. The augmented columns
are $\widehat{\bm\theta}$, the $n_1-1$ scaled residuals from group 1,
and the $n_2-1$ scaled residuals from group 2. The first two coordinate
rows are generated explicitly. Each remaining parity class has identical
coordinate variances and therefore contributes a scaled Wishart matrix
to their joint Gram matrix. A Bartlett factor generates that matrix when
the coordinate count is at least the Gram dimension; otherwise an
ordinary Gaussian rectangular matrix is used. This compression is exact
and changes only computation. The mean signal is added to the contrast
column; all feasible algorithms subsequently use the observed Gram
matrices, without access to the population covariance.

For the shared-data comparison, choose a fixed orthonormal residual
basis whose first $n_i/2$ columns are disjoint pair differences divided
by $\sqrt2$, completed inside the subspace orthogonal to the constant
vector. All sample sizes are even. This has the same Gaussian residual
law as a Helmert basis. The projected gamma procedure uses all
$n_i-1$ residuals and the estimators in Appendix~\ref{app:implementation-details}.
The Wang--Xu Gaussian multiplier procedure uses the paired residual
subset and the full CQ statistic, as in \citet{WangXu2022}; it is an
implementation of that published procedure, not a new proposed test.
The algebraic identity between the centered full norm and the scaled
CQ statistic is checked against raw observations in the separately maintained computational code.
The maximum error in these and the Gram-update checks is below
$1.1\times10^{-14}$.

All feasible methods use $B=999$ fresh auxiliary draws in each dataset.
Spherical weights are obtained by squaring independent standard normal
coordinates and dividing by their sum. For each mean shift, the observed
Gram spectrum is recomputed. Subtracting its average eigenvalue before
forming the reference removes a common deterministic shift and improves
numerical stability at large $p$. The spherical decision is checked
against the plus-one Monte Carlo p-value in every replication.
The projected-gamma and Wang--Xu reference quantiles use the same
plus-one order-statistic convention. Random draws are shared across
signal shifts for paired comparisons; all outer datasets are independent.
The base seeds are $202611080310+100000j$, $j=0,\ldots,4$, for
D200, S800, S800L, D20000 and S20000, respectively.

Population oracle critical values use the chi-squared CDF in the diffuse
case and adaptive one-dimensional integration in the single-spike case.
The latter integrates the bulk chi-squared CDF over the leading
Gaussian coordinate. A 512-node Gauss--Hermite calculation verifies
the resulting tail probabilities to within $2\times10^{-14}$.
Each row in Table~\ref{tab:rotation-power} uses the same 10,000 datasets
as its corresponding size entry. Paired differences and their standard
errors, projection frequencies and simulation seeds are included as CSV
files. No independent null-size adjustment is used in this study.

\begin{table}[!htbp]
\centering\small
\begin{threeparttable}
\caption{Empirical powers (in \%) in the Gaussian calibration study.}
\label{tab:rotation-power}
\begin{tabular}{lclrrrr}
\toprule
& & & \multicolumn{4}{c}{Method}\\
\cmidrule(lr){4-7}
Configuration & $\kappa$ & Direction & Oracle & Sphere & PG & WX--G\\
\midrule
D200 & 0.5 & Top & 12.80 (0.33) & 12.03 (0.33) & 10.31 (0.30) & 12.93 (0.34) \\
D200 & 0.5 & Bulk & 12.62 (0.33) & 12.04 (0.33) & 10.01 (0.30) & 12.73 (0.33) \\
D200 & 1 & Top & 25.28 (0.43) & 24.06 (0.43) & 21.28 (0.41) & 25.27 (0.43) \\
D200 & 1 & Bulk & 25.36 (0.44) & 24.16 (0.43) & 20.51 (0.40) & 25.39 (0.44) \\
\addlinespace[3pt]
S800 & 0.5 & Top & 14.96 (0.36) & 13.67 (0.34) & 14.75 (0.35) & 15.85 (0.37) \\
S800 & 0.5 & Bulk & 9.95 (0.30) & 10.50 (0.31) & 10.43 (0.31) & 11.66 (0.32) \\
S800 & 1 & Top & 25.81 (0.44) & 23.34 (0.42) & 25.81 (0.44) & 26.82 (0.44) \\
S800 & 1 & Bulk & 19.20 (0.39) & 20.26 (0.40) & 19.72 (0.40) & 22.01 (0.41) \\
\addlinespace[3pt]
S800L & 0.5 & Top & 15.04 (0.36) & 14.55 (0.35) & 13.12 (0.34) & 15.73 (0.36) \\
S800L & 0.5 & Bulk & 10.04 (0.30) & 10.17 (0.30) & 7.88 (0.27) & 10.88 (0.31) \\
S800L & 1 & Top & 26.10 (0.44) & 25.12 (0.43) & 23.41 (0.42) & 26.93 (0.44) \\
S800L & 1 & Bulk & 19.72 (0.40) & 20.57 (0.40) & 15.82 (0.36) & 21.29 (0.41) \\
\addlinespace[3pt]
D20000 & 0.5 & Top & 12.59 (0.33) & 12.21 (0.33) & 10.78 (0.31) & 12.19 (0.33) \\
D20000 & 0.5 & Bulk & 12.50 (0.33) & 11.92 (0.32) & 10.80 (0.31) & 12.23 (0.33) \\
D20000 & 1 & Top & 25.20 (0.43) & 23.89 (0.43) & 21.73 (0.41) & 24.84 (0.43) \\
D20000 & 1 & Bulk & 25.22 (0.43) & 23.83 (0.43) & 21.79 (0.41) & 24.71 (0.43) \\
\addlinespace[3pt]
S20000 & 0.5 & Top & 15.85 (0.37) & 14.20 (0.35) & 16.19 (0.37) & 16.65 (0.37) \\
S20000 & 0.5 & Bulk & 9.75 (0.30) & 10.42 (0.31) & 10.63 (0.31) & 11.31 (0.32) \\
S20000 & 1 & Top & 26.34 (0.44) & 24.15 (0.43) & 26.93 (0.44) & 27.41 (0.45) \\
S20000 & 1 & Bulk & 18.85 (0.39) & 19.79 (0.40) & 20.46 (0.40) & 21.84 (0.41) \\\bottomrule
\end{tabular}
\begin{tablenotes}[flushleft]
\footnotesize
\item[] \textit{Note.} Entries in parentheses are Monte Carlo standard errors in percentage points; all procedures are applied to the same 10,000 datasets within each configuration. Nominal powers should be read together with Table~\ref{tab:rotation-size}.
\end{tablenotes}
\end{threeparttable}
\end{table}

As an independent implementation check, we also generate 10,000
rank-one Gaussian contrast--residual datasets for each $m\in\{19,39\}$,
with $a\in\{0,1,2\}$ shared within each set. The conditional spherical
reference fraction is $\operatorname{Beta}(1/2,m/2)$, and its ideal
infinite-draw rejection probability is that of a noncentral
$F_{1,m}(a^2)$ statistic at the central $F_{1,m}$ critical value.
For the actual $B=999$ procedure, the population power is obtained by
integrating $\Pr\{\operatorname{Binomial}(999,p)\le49\}$ over the
alternative distribution of the exact spherical p-value $p$.
This explicitly accounts for auxiliary Monte Carlo variation.
The two integrations with tighter tolerances differ by at most
$4\times10^{-11}$; all six empirical rejection probabilities lie
within $1.17$ Monte Carlo standard errors of their finite-$B$
population values. These are 20,000 additional validation datasets,
separate from the 50,000 datasets in the main calibration comparison.

\subsection{Finite-sample unknown-covariance diagnostic experiment}
\label{app:unknown-covariance-experiment}

This experiment concerns the confidence construction in
Propositions~\ref{prop:gaussian-hollow-ratio} and
\ref{prop:finite-unknown-covariance-diagnostic}. Each configuration uses
10,000 independent datasets of $m$ Gaussian residual vectors with known
simulation covariance $\bm\Omega$. Population quantities are used only
to assess coverage, and are not inputs to the diagnostic. We consider
$\bm\Omega=\bm I_p$; a rank-one covariance with sole nonzero eigenvalue
one, represented in dimension one; and
$\bm\Omega=\operatorname{diag}(\sqrt p,1,\ldots,1)$ with $p=800$.
The last configuration has
$A=2p-1$ and $r=\sqrt{p/(2p-1)}$.
The count $m$ is the number of independent matched Gaussian residuals;
it is not the total size of the two original samples.

The six tail allocations for the lower bound on $A$, Gaussian norm,
spherical direction, hollow residual noise, upper bound on $A$, and
interval-probability band are each $0.05/6$. The last allocation is
split equally between DKW and the simultaneous spacing inequalities.
We record the joint sufficient event
$A_L\le A\le A_U$, $r\le r_U$, the DKW event, and all conservative
spacing inequalities using the cumulative upper constants $\eta_k$.
Its probability is at least $0.95$ by the proved bounds.
The upper-trace Chernoff constant is selected from the deterministic
721-point grid $s=\exp(u)$, $u\in[-12,24]$ equally spaced.
The diagnostic minimizes over
$t=\sqrt{2A_U}u$, with 81 logarithmically spaced values of
$u$ from $10^{-4}$ to $5$.
Maximum closed-window counts are computed from all sorted observed norms,
so the interval location is not discretized.
The DKW-only comparator uses the entire interval-probability allocation
$0.05/6$, whereas the combined band splits that same allocation equally.
Consequently both methods have the same overall confidence budget.

\begin{table}[!htbp]
\centering
\begin{threeparttable}
\caption{Coverage and median bounds in the Gaussian residual experiment.}
\label{tab:finite-unknown-coverage}
\small
\begin{tabular}{lrr rrrrrr}
\toprule
& & & & \multicolumn{5}{c}{Median}\\
\cmidrule(lr){5-9}
Spectrum & $m$ & $p$ & Coverage & $r_U$ & $A_U/A$ & $B_{.001}$ & $B_{.01}$ & $B_{.1}$\\
\midrule
Isotropic & 40 & 80 & 99.50 & 1.000 & 7.943 & 0.500 & 0.500 & 0.500 \\
Isotropic & 120 & 80 & 99.55 & 0.888 & 3.171 & 0.338 & 0.500 & 0.500 \\
Isotropic & 40 & 800 & 99.62 & 1.000 & 7.958 & 0.500 & 0.500 & 0.500 \\
Isotropic & 120 & 800 & 99.55 & 0.734 & 3.171 & 0.329 & 0.482 & 0.500 \\
Isotropic & 400 & 800 & 99.68 & 0.300 & 1.862 & 0.156 & 0.255 & 0.449 \\
Isotropic & 120 & 80,000 & 99.41 & 0.664 & 3.171 & 0.322 & 0.472 & 0.500 \\
\addlinespace[3pt]
Rank one & 40 & 1 & 99.62 & 1.000 & 7.341 & 0.500 & 0.500 & 0.500 \\
Rank one & 120 & 1 & 99.57 & 1.000 & 3.094 & 0.500 & 0.500 & 0.500 \\
\addlinespace[3pt]
Spiked & 120 & 800 & 99.72 & 1.000 & 3.127 & 0.358 & 0.500 & 0.500 \\\bottomrule
\end{tabular}
\begin{tablenotes}[flushleft]
\footnotesize
\item[] \textit{Note.} Coverage is the empirical percentage of the joint sufficient event; the nominal level is 95\%. The remaining entries are medians. Each row uses 10,000 independent residual datasets. $B_K$ denotes $\widehat B_m(K)$. All configurations have median $B_{0.5}=0.500$.
\end{tablenotes}
\end{threeparttable}
\end{table}

\begin{table}[!htbp]
\centering
\begin{threeparttable}
\caption{Median diagnostic bounds under a common confidence budget.}
\label{tab:finite-unknown-band-comparison}
\small
\begin{tabular}{lrrccc}
\toprule
Spectrum & $m$ & $p$ & $K=0.001$ & $K=0.01$ & $K=0.1$\\
\midrule
Isotropic & 40 & 80 & 0.500 / 0.500 & 0.500 / 0.500 & 0.500 / 0.500 \\
Isotropic & 120 & 80 & 0.338 / 0.440 & 0.500 / 0.500 & 0.500 / 0.500 \\
Isotropic & 40 & 800 & 0.500 / 0.500 & 0.500 / 0.500 & 0.500 / 0.500 \\
Isotropic & 120 & 800 & 0.329 / 0.433 & 0.482 / 0.500 & 0.500 / 0.500 \\
Isotropic & 400 & 800 & 0.156 / 0.233 & 0.255 / 0.313 & 0.449 / 0.488 \\
Isotropic & 120 & 80,000 & 0.322 / 0.430 & 0.472 / 0.500 & 0.500 / 0.500 \\
\addlinespace[3pt]
Rank one & 40 & 1 & 0.500 / 0.500 & 0.500 / 0.500 & 0.500 / 0.500 \\
Rank one & 120 & 1 & 0.500 / 0.500 & 0.500 / 0.500 & 0.500 / 0.500 \\
\addlinespace[3pt]
Spiked & 120 & 800 & 0.358 / 0.459 & 0.500 / 0.500 & 0.500 / 0.500 \\\bottomrule
\end{tabular}
\begin{tablenotes}[flushleft]
\footnotesize
\item[] \textit{Note.} Each pair reports the combined interval band / the DKW-only band. Both methods use the same $A_U$, $r_U$ and width grid. The value $0.500$ is the universal cap, not an estimate of approximation error.
\end{tablenotes}
\end{threeparttable}
\end{table}

The Monte Carlo standard errors of the reported coverage percentages
range from 0.053 to 0.077 percentage points.
The substantial overcoverage is compatible with the conservative
combination of one-sided inequalities. Increasing $p$ alone does not
produce a precise certificate: for isotropic data with $m=120$,
the median bound at $K=0.01$ is 0.482 for $p=800$ and 0.472 for
$p=80{,}000$. Increasing the residual count to $m=400$, at $p=800$,
reduces this median to 0.255. Even this value does not certify a small
error such as 0.01 or 0.05. Rank-one and spiked configurations emphasize
that the certificate need not be informative for concentrated spectra.

For isotropic cases, the Gram matrices are sampled exactly from the
appropriate Wishart law using the Bartlett construction when $p\ge m$,
and from a direct Gaussian matrix otherwise. The spiked Gram matrix
is the sum of an independent Wishart matrix with $p-1$ degrees of freedom
and $\sqrt p\,\bm z\bm z^\top$, with independent standard Gaussian
$\bm z$. Dense symmetric eigendecomposition supplies the largest
hollow-Gram eigenvalue. Isotropic and rank-one CDFs are evaluated as
chi-square CDFs. For the spiked coverage check only, the CDF of
$\chi^2_{p-1}+\sqrt p\,\chi^2_1$ is evaluated by one-dimensional
quadrature and interpolation, with numerical validation recorded in
the reproducibility metadata. This numerical step is not used to
compute the diagnostic bound. Full per-replicate results and random
seeds are retained. In total, the experiment comprises 90,000 independent
residual datasets; regenerating the same draws for the DKW comparison
does not add independent replicates. Neither the true supremum shift
error nor an outer two-sample rejection probability is estimated here.

\subsection{Noncommuting covariance and severe imbalance}
\label{app:noncommuting-numerics}
The four additional Gaussian configurations are D200-I, with $p=200$
and $(n_1,n_2)=(12,120)$; NC200-B, with $p=200$ and $(24,24)$;
NC200-I, with $p=200$ and $(12,120)$; and NC20000-I, with $p=20{,}000$
and $(12,120)$. D200-I has $\bm\Sigma_1=\bm\Sigma_2=\bm I_p$.
For the other three configurations, let
$\bm u=\bm e_1$, $\bm v=(\bm e_1+\bm e_2)/\sqrt2$,
$w_1=n_2/(n_1+n_2)$ and $w_2=1-w_1$, and set
\[
 \bm\Sigma_1=\bm I_p+\frac{\sqrt p}{w_1}\bm u\bm u^{\mathsf T},
 \qquad
 \bm\Sigma_2=\bm I_p+\frac{\sqrt p}{w_2}\bm v\bm v^{\mathsf T}.
\]
These two matrices do not commute: their commutator has Frobenius norm
$p/(\sqrt2w_1w_2)>0$. Their weighted sum is
$\bm\Omega=\bm I_p+\sqrt p(\bm u\bm u^{\mathsf T}+\bm v\bm v^{\mathsf T})$,
with leading eigenvalues
$\lambda_\pm=1+\sqrt p(1\pm1/\sqrt2)$ and all remaining eigenvalues equal
to one. The covariance of the normalized mean contrast is therefore
held fixed at $p=200$ while the sample allocation and group covariance
matrices change. This isolates an aspect of calibration under imbalance;
it is not a comparison with identical group covariances across sample
allocations.

Alternative signals are $\bm\theta=\sqrt{\kappa\sqrt{2A}}\bm d$,
where $\kappa\in\{0.5,1\}$ and $\bm d$ is either the leading eigenvector
$(\cos(\pi/8),\sin(\pi/8),0,\ldots,0)^{\mathsf T}$ or the bulk direction
$\bm e_3$. In D200-I, the first direction is $\bm e_1$.
Each configuration also includes the null signal. All signal settings
within a configuration share the same Gaussian errors.

The simulation generates the exact joint Gaussian Gram law of the
normalized mean contrast and all independent within-group residuals.
The first three coordinate rows are generated explicitly, with
covariance blocks $\bm\Omega$, $w_1\bm\Sigma_1$ and
$w_2\bm\Sigma_2$ for the contrast and the two residual groups,
respectively. Remaining coordinates are isotropic: their Gram
contribution is $\bm D\bm W\bm D$, where
$\bm W\sim W_{n_1+n_2-1}(p-3,\bm I)$ and $\bm D$ has diagonal
entries $1$, $n_1-1$ copies of $\sqrt{w_1}$ and $n_2-1$ copies of
$\sqrt{w_2}$. A Bartlett factor generates this independent Wishart
contribution. This is an exact computational compression, not an
approximation to a noncommuting covariance model. It does not explicitly
generate all raw observations; the feasible procedures receive only the
resulting observed contrast and residual Gram quantities.

For compatibility with the published Wang--Xu Gaussian multiplier
procedure, the fixed orthonormal residual basis begins with disjoint
pair differences divided by $\sqrt2$ and is completed in the subspace
orthogonal to the constant vector. All group sample sizes are even.
The spherical procedure uses the first $\min(n_1-1,n_2-1)$ residuals
from each group; the normalized mean contrast uses all observations.
The projected gamma comparison uses all $n_i-1$ residuals, whereas
WX--G uses the paired subset and the full CQ statistic, as in
\citet{WangXu2022}. Every procedure uses the same underlying dataset.
The spherical spectrum is recomputed after every mean shift, with a
common eigenvalue shift removed for stability. Fresh $B=999$ reference
draws are generated in every dataset, and the spherical decision is
checked against its plus-one p-value in every replication. All methods
use the same plus-one order-statistic convention. The original
projected-gamma threshold and variance-budget projection are retained
without modification; their realized frequencies are saved in the
reproducibility files. The main seeds are $202611090410+100000j$,
$j=0,\ldots,3$, in the configuration order given above.

The population oracle uses $\chi^2_p$ in the diffuse case. In the
noncommuting cases, its null distribution is
$\lambda_+\chi^2_1+\lambda_-\chi^2_1+\chi^2_{p-2}$, with independent
terms. In polar coordinates, the first two terms have distribution
$\Lambda(\phi)\chi^2_2$, where
$\Lambda(\phi)=\lambda_+\cos^2\phi+\lambda_-\sin^2\phi$ and
$\phi$ is uniform on $[0,\pi/2]$. Conditional on $\phi$, the CDF at
$q>0$ is
\[
 F_{p-2}(q)-\exp\!\left\{-\frac{q}{2\Lambda}
    -\frac{p-2}{2}\log(1-1/\Lambda)\right\}
 F_{p-2}\!\left(q(1-1/\Lambda)\right),
\]
where $F_d$ is the $\chi^2_d$ CDF. Adaptive angular integration gives
the critical value and is checked independently by Gauss--Legendre
integration and a tensor Gauss--Hermite calculation over the original
two Gaussian spike coordinates. Raw-data algebra checks separately
verify the spherical Gram construction, mean-shift update, full CQ
identity and the Wang--Xu matrix at both sample allocations.

Table~\ref{tab:noncommuting-power} reports all powers. Each configuration
contains 10,000 independent datasets, for 40,000 new main datasets in
this study. Mean shifts reuse those datasets and are not counted as
independent additional replications. CSV files include paired
differences and Monte Carlo standard errors. No null-size adjustment
is applied to the reported powers. These experiments concern Gaussian
validity; they do not establish a non-Gaussian rotation guarantee.

\begin{table}[!htbp]
\centering\small\setlength{\tabcolsep}{3.5pt}
\begin{threeparttable}
\caption{Empirical powers (in \%) under sample imbalance and noncommuting covariances.}
\label{tab:noncommuting-power}
\begin{tabular}{lclrrrr}
\toprule
& & & \multicolumn{4}{c}{Method}\\
\cmidrule(lr){4-7}
Configuration & $\kappa$ & Direction & Oracle & Sphere & PG & WX--G\\
\midrule
D200-I & 0.5 & Top & 12.71 (0.33) & 12.49 (0.33) & 11.52 (0.32) & 13.88 (0.35) \\
D200-I & 0.5 & Bulk & 12.92 (0.34) & 12.27 (0.33) & 11.30 (0.32) & 13.76 (0.34) \\
D200-I & 1 & Top & 25.63 (0.44) & 23.78 (0.43) & 22.45 (0.42) & 25.59 (0.44) \\
D200-I & 1 & Bulk & 25.35 (0.44) & 23.63 (0.42) & 22.44 (0.42) & 25.78 (0.44) \\
\addlinespace[3pt]
NC200-B & 0.5 & Top & 13.94 (0.35) & 13.28 (0.34) & 13.67 (0.34) & 15.33 (0.36) \\
NC200-B & 0.5 & Bulk & 8.43 (0.28) & 9.06 (0.29) & 8.65 (0.28) & 10.56 (0.31) \\
NC200-B & 1 & Top & 23.64 (0.42) & 21.90 (0.41) & 23.02 (0.42) & 25.35 (0.44) \\
NC200-B & 1 & Bulk & 15.19 (0.36) & 17.43 (0.38) & 15.66 (0.36) & 19.28 (0.39) \\
\addlinespace[3pt]
NC200-I & 0.5 & Top & 14.51 (0.35) & 12.81 (0.33) & 15.19 (0.36) & 15.84 (0.37) \\
NC200-I & 0.5 & Bulk & 8.92 (0.29) & 10.07 (0.30) & 10.27 (0.30) & 10.60 (0.31) \\
NC200-I & 1 & Top & 24.15 (0.43) & 20.52 (0.40) & 24.98 (0.43) & 25.59 (0.44) \\
NC200-I & 1 & Bulk & 15.88 (0.37) & 19.87 (0.40) & 18.52 (0.39) & 19.47 (0.40) \\
\addlinespace[3pt]
NC20000-I & 0.5 & Top & 13.88 (0.35) & 12.38 (0.33) & 15.10 (0.36) & 15.12 (0.36) \\
NC20000-I & 0.5 & Bulk & 8.46 (0.28) & 9.69 (0.30) & 9.92 (0.30) & 10.04 (0.30) \\
NC20000-I & 1 & Top & 23.26 (0.42) & 20.14 (0.40) & 24.61 (0.43) & 24.67 (0.43) \\
NC20000-I & 1 & Bulk & 14.80 (0.36) & 18.94 (0.39) & 17.69 (0.38) & 17.90 (0.38) \\\bottomrule
\end{tabular}
\begin{tablenotes}[flushleft]
\footnotesize
\item[] \textit{Note.} Parentheses give Monte Carlo standard errors in percentage points. Every configuration uses the same 10,000 datasets for all signal settings and methods. Top and bulk denote the leading eigenvector of $\bm\Omega$ and $\bm e_3$, respectively. Nominal powers should be read with the sizes in Table~\ref{tab:noncommuting-size}.
\end{tablenotes}
\end{threeparttable}
\end{table}

\begin{figure}[t]
\centering
\includegraphics[width=\textwidth]{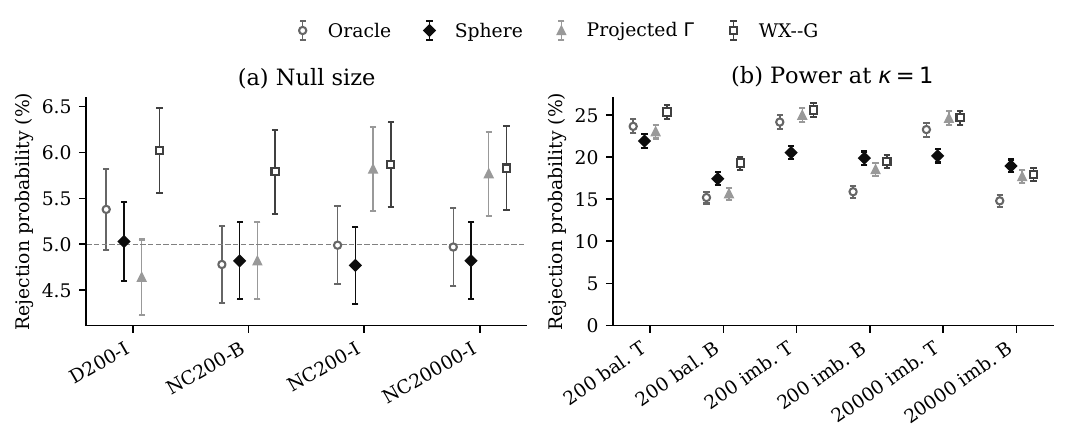}
\caption{Gaussian calibration with noncommuting covariances and severe
imbalance. Panel (a) gives null sizes; panel (b) gives power at $\kappa=1$
in the noncommuting models. T and B denote leading-axis and bulk-axis
signals. Balanced and imbalanced sample sizes are $(24,24)$ and
$(12,120)$, respectively. Error bars are $1.96$ Monte Carlo standard
errors based on 10,000 datasets.}
\label{fig:noncommuting-calibration}
\end{figure}
\subsection{Non-Gaussian latent models}
\label{app:non-gaussian-numerics}
We generate independent observations as
$\bm Y_{ij}=\bm\mu_i+\bm\Sigma_i^{1/2}\bm Z_{ij}$ using the symmetric
covariance square root. All coordinates of $\bm Z_{ij}$ are independent.
They are either Rademacher, taking $-1$ and $1$ with equal probability,
or $\sqrt{6/8}\,t_8$. Both laws have zero mean and unit variance; their
fourth moments are $1$ and $4.5$, respectively. The standardized $t_8$
law does not have a finite eighth moment. The independent latent
coordinates satisfy the stated fourth-moment and quadratic-form
conditions; no eighth-moment assumption is used for the spherical
approximation.

All four configurations have $p=200$. Rad-S and Rad-L use Rademacher
coordinates at sample sizes $(20,30)$ and $(40,60)$, respectively;
$t8$-S and $t8$-L use the corresponding standardized $t_8$ laws.
For both sample allocations, $w_1=3/5$ and $w_2=2/5$. Set
$\bm u=\bm e_1$, $\bm v=(\bm e_1+\bm e_2)/\sqrt2$ and
\[
 \bm\Sigma_1=\bm I_p+\frac{\sqrt p}{w_1}\bm u\bm u^{\mathsf T},
 \qquad
 \bm\Sigma_2=\bm I_p+\frac{\sqrt p}{w_2}\bm v\bm v^{\mathsf T}.
\]
These covariances do not commute. The two group covariances and
$\bm\Omega=\bm I_p+\sqrt p(\bm u\bm u^{\mathsf T}+\bm v\bm v^{\mathsf T})$
are held fixed when both sample sizes double. The leading direction
is $(\cos(\pi/8),\sin(\pi/8),0,\ldots,0)^{\mathsf T}$ and the bulk
direction is $\bm e_3$. We use
$\bm\theta=\sqrt{\kappa\sqrt{2A}}\bm d$ for $\kappa\in\{0.5,1\}$
and these two directions, together with the null signal. Signals
share the same raw errors. Adding $\bm\theta$ to the normalized
contrast is exactly the change induced by a constant group-1 mean
shift $\bm\theta/\sqrt{n_1n_2/(n_1+n_2)}$; residuals are unchanged.

Every raw latent coordinate and observation is generated explicitly.
No Gaussian Gram compression or replacement by independent Gaussian
mean and residual vectors is used. For each even sample size $n_i$,
write $k_i=n_i/2$. Let $\bm D_i$ and $\bm S_i$ be the $n_i\times k_i$
matrices of normalized disjoint-pair differences and sums, respectively.
Let $\bm K_i$ be the transpose of the standard $(k_i-1)\times k_i$
Helmert matrix, and use the fixed residual basis
$\bm H_i=(\bm D_i,\bm S_i\bm K_i)$.
Thus $\bm H_i^{\mathsf T}\bm H_i=\bm I_{n_i-1}$ and
$\bm H_i^{\mathsf T}\bm1=\bm0$. Residual columns are formed from the
actual observations using this basis and scaled by $\sqrt{w_i}$.
They are uncorrelated but generally dependent under the two latent
laws, and may also depend on the mean. Those dependencies are retained
in every replication. The spherical procedure uses all observations
for the normalized mean contrast and the first $\min(n_1-1,n_2-1)$
matched residual columns. The spectrum of the augmented Gram matrix
is recomputed after each mean shift.

The projected gamma implementation uses only the independent
disjoint-pair residual subset for both its covariance and trace
estimators. This preserves the independent-residual moment identities;
using all completed non-Gaussian residuals with those formulas would
not do so. Consequently this comparison has a different residual
choice from the all-residual Gaussian implementation in
Appendix~\ref{app:noncommuting-numerics}. Its $t_8$ results do not
establish the sub-Gaussian sufficient conditions used for the spectral
projection theory. WX--G uses the published Gaussian multiplier
reference based on the paired subset and the full CQ statistic,
as in \citet{WangXu2022}. The full CQ centering is computed from all
orthogonal residuals and is distinct from the paired trace used by
PG. Both identities are checked directly against raw observations.

Each feasible procedure uses $B=999$ fresh auxiliary draws per
dataset. Sphere weights are shared across mean shifts to support
paired comparisons, and every decision is checked against its
plus-one Monte Carlo p-value. The main seeds are
$202611090510+100000j$, $j=0,\ldots,3$, in the order Rad-S, Rad-L,
$t8$-S and $t8$-L. All outer datasets are independent.
The Gaussian population reference uses the known eigenvalues of
$\bm\Omega$ and the same angular integration as in
Appendix~\ref{app:noncommuting-numerics}; it is an approximation to
the non-Gaussian null law, not its exact population quantile.
No population covariance is supplied to any feasible calibration.

Table~\ref{tab:non-gaussian-power} reports all powers. Each of the
four configurations contains 10,000 independent datasets, for 40,000
new non-Gaussian datasets. The five mean settings reuse those datasets
and are not counted as additional independent replications. Results
include nominal sizes and powers, paired differences, Monte Carlo
standard errors, and projection frequencies. No null-size adjustment
is applied. Independent audit results are given below.

An independent raw-data audit verifies the covariance square roots, fixed
residual basis, separate paired-residual and full-CQ centerings, multiplier
matrix and signal-update identities; the maximum numerical discrepancy
is $8.3\times10^{-13}$. The null contrast provides a further distributional
check. For $T=\|\widehat{\bm\theta}\|^2$, latent fourth moment $\beta$ and
$c=n_1n_2/(n_1+n_2)$, the symmetric-root model gives
\[
 \mathbb E(T)=\tr(\bm\Omega),\qquad
 \operatorname{Var}(T)=2\tr(\bm\Omega^2)
 +(\beta-3)c^2\left\{
 \frac{\sum_{j=1}^p(\bm\Sigma_1)_{jj}^{2}}{n_1^3}
 +\frac{\sum_{j=1}^p(\bm\Sigma_2)_{jj}^{2}}{n_2^3}
 \right\}.
\]
All four simulated null means are within $0.80$ Monte Carlo standard
errors of this exact mean; empirical-to-theoretical variance ratios range
from $0.984$ to $1.026$. The variance comparison is descriptive: the
standardized $t_8$ latent variable has finite fourth moment $4.5$ but
infinite eighth moment, so no Gaussian standard error is assigned to its
empirical variance estimator.

\begin{table}[!htbp]
\centering\small\setlength{\tabcolsep}{3.5pt}
\begin{threeparttable}
\caption{Empirical powers (in \%) in the non-Gaussian study.}
\label{tab:non-gaussian-power}
\begin{tabular}{lclrrrr}
\toprule
& & & \multicolumn{4}{c}{Method}\\
\cmidrule(lr){4-7}
Configuration & $\kappa$ & Direction & G-ref & Sphere & PG & WX--G\\
\midrule
Rad-S & 0.5 & Top & 14.34 (0.35) & 13.00 (0.34) & 15.03 (0.36) & 15.33 (0.36) \\
Rad-S & 0.5 & Bulk & 8.54 (0.28) & 9.43 (0.29) & 10.55 (0.31) & 10.48 (0.31) \\
Rad-S & 1 & Top & 24.51 (0.43) & 21.85 (0.41) & 24.84 (0.43) & 24.96 (0.43) \\
Rad-S & 1 & Bulk & 15.95 (0.37) & 17.57 (0.38) & 18.53 (0.39) & 18.60 (0.39) \\
\addlinespace[3pt]
Rad-L & 0.5 & Top & 14.10 (0.35) & 13.53 (0.34) & 13.84 (0.35) & 14.65 (0.35) \\
Rad-L & 0.5 & Bulk & 8.98 (0.29) & 9.27 (0.29) & 8.97 (0.29) & 9.95 (0.30) \\
Rad-L & 1 & Top & 23.92 (0.43) & 22.33 (0.42) & 22.57 (0.42) & 24.00 (0.43) \\
Rad-L & 1 & Bulk & 15.79 (0.36) & 16.85 (0.37) & 15.56 (0.36) & 17.27 (0.38) \\
\addlinespace[3pt]
t8-S & 0.5 & Top & 14.38 (0.35) & 13.52 (0.34) & 15.92 (0.37) & 16.36 (0.37) \\
t8-S & 0.5 & Bulk & 8.52 (0.28) & 9.48 (0.29) & 10.71 (0.31) & 11.22 (0.32) \\
t8-S & 1 & Top & 24.43 (0.43) & 22.40 (0.42) & 25.86 (0.44) & 26.41 (0.44) \\
t8-S & 1 & Bulk & 15.26 (0.36) & 18.44 (0.39) & 19.63 (0.40) & 20.10 (0.40) \\
\addlinespace[3pt]
t8-L & 0.5 & Top & 14.40 (0.35) & 13.53 (0.34) & 13.78 (0.34) & 14.96 (0.36) \\
t8-L & 0.5 & Bulk & 8.67 (0.28) & 9.21 (0.29) & 8.80 (0.28) & 9.80 (0.30) \\
t8-L & 1 & Top & 23.73 (0.43) & 22.33 (0.42) & 23.04 (0.42) & 24.68 (0.43) \\
t8-L & 1 & Bulk & 15.64 (0.36) & 17.12 (0.38) & 16.09 (0.37) & 18.17 (0.39) \\\bottomrule
\end{tabular}
\begin{tablenotes}[flushleft]
\footnotesize
\item[] \textit{Note.} Parentheses give Monte Carlo standard errors in percentage points. Every configuration uses the same 10,000 datasets for all signal settings and methods. Top and bulk denote the leading eigenvector of $\bm\Omega$ and $\bm e_3$, respectively. Nominal powers should be read with the sizes in Table~\ref{tab:non-gaussian-size}.
\end{tablenotes}
\end{threeparttable}
\end{table}

\begin{figure}[t]
\centering
\includegraphics[width=\textwidth]{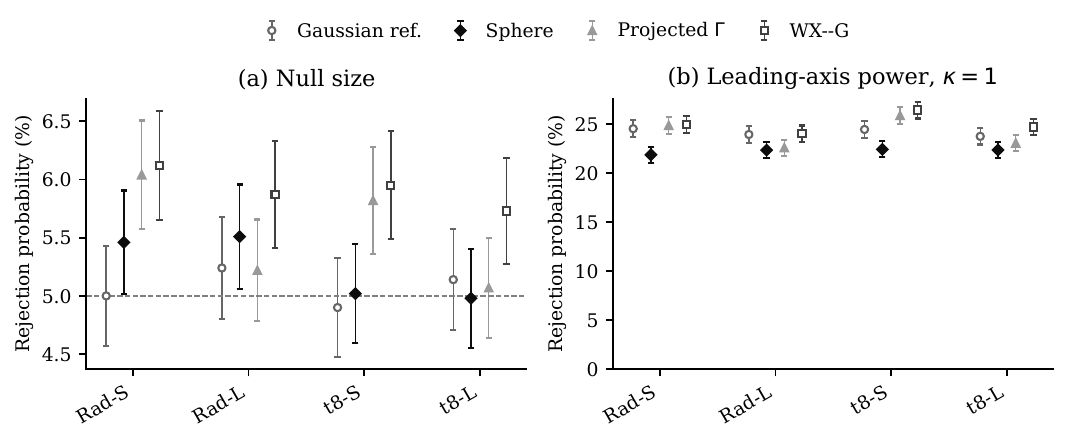}
\caption{Non-Gaussian latent models with fixed noncommuting covariance
matrices. Rad and $t8$ denote Rademacher and standardized $t_8$ latent
coordinates; S and L denote sample sizes $(20,30)$ and $(40,60)$.
Panel (a) gives null sizes and panel (b) leading-axis power at
$\kappa=1$. Error bars are $1.96$ Monte Carlo standard errors based
on 10,000 datasets per configuration. The Gaussian reference is an
approximation benchmark, not a non-Gaussian oracle.}
\label{fig:non-gaussian-calibration}
\end{figure}
\subsection{Paired comparison of tightened confidence certificates}
\label{app:tightened-diagnostic-experiment}

We reuse all nine Gaussian residual configurations of
Appendix~\ref{app:unknown-covariance-experiment}: six isotropic settings,
two rank-one settings, and the spiked covariance
$\operatorname{diag}(\sqrt p,1,\ldots,1)$ at $(m,p)=(120,800)$.
Every configuration contains 10,000 independent residual datasets.
The recorded seeds regenerate the original observations exactly; this
paired replay adds no independent datasets to the existing total of 90,000.
The original trace U-statistic and all four original diagnostic values,
at $K=0.001,0.01,0.1,0.5$, reproduce the saved arrays exactly in the
reported calculations. The additional values $K=0.05,0.25,1$ use the
same datasets. The sample size $m$ counts matched Gaussian residuals,
not all observations in the original two samples.

Four methods are compared. I0 is the original interval certificate.
I1 keeps its empirical concentration band and width grid but replaces
the trace and hollow-Gram bounds by the tightened constructions.
S uses the smooth certificate alone. H is the minimum of I1 and the
smooth bound calculated with I1's spectral upper bound.
I0, I1 and H use six allocations of $0.05/6$, for $A_L$, the Gaussian
norm, spherical direction, hollow residual noise, $A_U$, and the
empirical interval event. The last allocation is split equally between
DKW and the simultaneous beta-spacing inequalities. H uses the same
joint event as I1 and requires no extra allocation. S uses four
allocations of $0.05/4$ for its spectral construction and requires no
separate $A_U$ or interval event. Computing an unused upper-trace
constant in the software does not add a confidence event. All four
procedures therefore have the same 95\% nominal confidence level.
The differing allocation means that S and H are distinct prespecified
procedures: taking their minimum after viewing the data would need a
new joint-coverage argument and is not done here.

Coverage records sufficient population-parameter coverage events, not the
unknown true supremum approximation error. For S it records
$A_L\le A$ and $r\le r_U$; for I0, I1 and H it additionally records
$A\le A_U$ and the saved DKW and spacing event. Each recorded event
has probability at least 0.95 by the theoretical construction.
Population $A$, $r$ and the reference CDF enter only this assessment,
not calculation of the diagnostic. The width grid remains
$t=\sqrt{2A_U}u$ with 81 logarithmically spaced values of $u$ from
$10^{-4}$ to 5; closed-window locations are optimized over all observed
squared norms. The tightened trace constants use deterministic
one-dimensional integration. Recorded numerical error estimates and
conservative floating-point padding are implementation checks rather
than formal interval-arithmetic proofs.

An equivalent summary reports the largest effect range certified at a
fixed tolerance. At $\varepsilon=0.05$, the median
$\widehat K_{\varepsilon}$ from the standalone smooth certificate is
0.4635 for $(m,p)=(400,800)$ and 0.0906 for $(120,80{,}000)$;
their 10th--90th percentile intervals are $(0.4501,0.4754)$ and
$(0.0769,0.1025)$, respectively. For $(120,800)$ the median is
0.0332, and only 6.84\% of datasets certify $K\ge0.05$.
The other six configurations yield zero throughout their 10,000
replications. These are algebraic summaries of the same 90,000
residual datasets, not additional Monte Carlo experiments.

\begin{table}[!htbp]
\centering\small\setlength{\tabcolsep}{3.5pt}
\begin{threeparttable}
\caption{Coverage, spectral bounds and certified effect ranges.}
\label{tab:tightened-diagnostic-coverage}
\begin{tabular}{lrr rrrr cr}
\toprule
& & & \multicolumn{4}{c}{Coverage (in \%)} & \multicolumn{2}{c}{Median}\\
\cmidrule(lr){4-7}\cmidrule(lr){8-9}
Spectrum & $m$ & $p$ & I0 & I1 & S & H & $r_U$: S / H & $\widehat K_{.05}$: S\\
\midrule
Isotropic & 40 & 80 & 99.50 & 99.50 & 100.00 & 99.50 & 1.000 / 1.000 & 0.0000 \\
Isotropic & 120 & 80 & 99.55 & 99.55 & 100.00 & 99.55 & 0.596 / 0.627 & 0.0000 \\
Isotropic & 40 & 800 & 99.62 & 99.62 & 100.00 & 99.62 & 1.000 / 1.000 & 0.0000 \\
Isotropic & 120 & 800 & 99.55 & 99.55 & 100.00 & 99.55 & 0.464 / 0.489 & 0.0332 \\
Isotropic & 400 & 800 & 99.68 & 99.68 & 100.00 & 99.68 & 0.201 / 0.206 & 0.4635 \\
Isotropic & 120 & 80,000 & 99.41 & 99.41 & 100.00 & 99.41 & 0.408 / 0.431 & 0.0906 \\
\addlinespace[3pt]
Rank one & 40 & 1 & 99.62 & 99.16 & 99.55 & 99.16 & 1.000 / 1.000 & 0.0000 \\
Rank one & 120 & 1 & 99.57 & 98.89 & 99.49 & 98.89 & 1.000 / 1.000 & 0.0000 \\
\addlinespace[3pt]
Spiked & 120 & 800 & 99.72 & 99.72 & 100.00 & 99.72 & 1.000 / 1.000 & 0.0000 \\
\bottomrule
\end{tabular}
\begin{tablenotes}[flushleft]
\footnotesize
\item[] \textit{Note.} Coverage records sufficient confidence events at nominal level 95\%. I0 is the original interval certificate, I1 the tightened interval certificate, S the smooth certificate alone, and H the hybrid. Each row reuses the same 10,000 Gaussian residual datasets for all methods. The last column inverts S at error tolerance 0.05, as in Corollary~\ref{cor:certified-effect-range}; it is a sufficient effect range, not the largest true range of accurate approximation.
\end{tablenotes}
\end{threeparttable}
\end{table}

\begin{table}[!htbp]
\centering\small\setlength{\tabcolsep}{3.5pt}
\begin{threeparttable}
\caption{Median bounds and certification frequencies at $K=0.25$.}
\label{tab:tightened-diagnostic-usefulness}
\begin{tabular}{lrr rrrr rr}
\toprule
& & & \multicolumn{4}{c}{Median bound} & \multicolumn{2}{c}{Certified (in \%)}\\
\cmidrule(lr){4-7}\cmidrule(lr){8-9}
Spectrum & $m$ & $p$ & I0 & I1 & S & H & S $\le.05$ & H $\le.05$\\
\midrule
Isotropic & 40 & 80 & 0.500 & 0.500 & 0.500 & 0.500 & 0.00 & 0.00 \\
Isotropic & 120 & 80 & 0.500 & 0.500 & 0.500 & 0.500 & 0.00 & 0.00 \\
Isotropic & 40 & 800 & 0.500 & 0.500 & 0.500 & 0.500 & 0.00 & 0.00 \\
Isotropic & 120 & 800 & 0.500 & 0.500 & 0.376 & 0.500 & 0.00 & 0.00 \\
Isotropic & 400 & 800 & 0.500 & 0.485 & 0.027 & 0.028 & 100.00 & 100.00 \\
Isotropic & 120 & 80,000 & 0.500 & 0.500 & 0.138 & 0.188 & 0.00 & 0.00 \\
\addlinespace[3pt]
Rank one & 40 & 1 & 0.500 & 0.500 & 0.500 & 0.500 & 0.00 & 0.00 \\
Rank one & 120 & 1 & 0.500 & 0.500 & 0.500 & 0.500 & 0.00 & 0.00 \\
\addlinespace[3pt]
Spiked & 120 & 800 & 0.500 & 0.500 & 0.500 & 0.500 & 0.00 & 0.00 \\\bottomrule
\end{tabular}
\begin{tablenotes}[flushleft]
\footnotesize
\item[] \textit{Note.} The last two columns give the percentage of datasets on which the smooth or hybrid bound is at most 0.05. Definitions of I0, I1, S and H are as in Table~\ref{tab:tightened-diagnostic-coverage}. The value 0.500 is the universal bound, not an estimate of the true approximation error.
\end{tablenotes}
\end{threeparttable}
\end{table}

\begin{figure}[htbp]
\centering
\includegraphics[width=.98\textwidth]{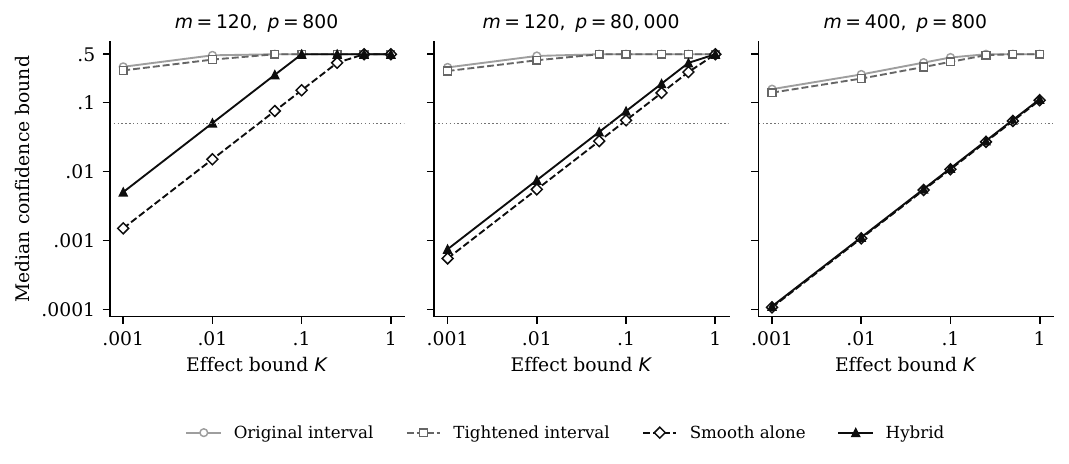}
\caption{Median 95\% confidence bounds over the effect range for three
isotropic configurations, using the same 10,000 residual datasets for
every method and effect value. Both axes are logarithmic. The dotted
line marks a 0.05 approximation tolerance. The curves are bounds, not
estimates of the actual approximation error. S and H use the different
valid budget allocations described in the text.}
\label{fig:tightened-diagnostic}
\end{figure}

For isotropic $(m,p)=(400,800)$, the median H bounds at
$K=0.05,0.25,0.5$ are 0.00558, 0.02791 and 0.05582;
the S medians are 0.00539, 0.02697 and 0.05394.
At $K=0.25$, all 10,000 datasets give bounds at most 0.05 for both
S and H, whereas at $K=0.5$ none do. Thus the certificate can identify
a useful effect range without asserting accuracy throughout a larger range.

The confidence budget matters near $r_U=1/2$. At isotropic
$(m,p)=(120,800)$, the median $r_U$ is 0.464 for S and 0.489 for H.
At $K=0.01$, the median bounds are respectively 0.0151 and 0.0503,
and 95.83\% versus 49.79\% of datasets certify the 0.05 tolerance.
At $K=0.1$, these medians increase to 0.1505 and 0.5000.
For $(m,p)=(120,80{,}000)$ at $K=0.05$, the S and H medians are
0.0276 and 0.0375, with 99.98\% and 94.37\% of datasets below
the 0.05 tolerance. Increasing $p$ therefore helps in this experiment,
although increasing $m$ remains important for a wider certified effect range.

Coverage of the recorded sufficient event ranges from 99.41\% to
99.72\% for I0, 98.89\% to 99.72\% for I1 and H, and
99.49\% to 100\% for S. The largest plug-in Monte Carlo standard
error among these coverage percentages is 0.105 percentage points.
An observed 100\% is a finite simulation count, not a claim of exact
coverage one. The reduced bounds remain conservative in these
configurations. Rank-one spectra have median 0.500 throughout
the evaluated effect grid; the spiked configuration also gives median
0.500 for $K\ge0.01$ under H and for all evaluated $K$ under S.
Failure to obtain a small upper bound does not establish a large true
approximation error.

All results retain the original Gaussian-model scope. This experiment
neither estimates the true supremum shifted-reference error nor measures
the rejection probability of an outer two-sample test. Paired per-dataset
bound differences and their Monte Carlo standard errors, all seven effect
values, seeds, exact replay checks and file hashes are retained in the
separate computational archive.

\section{Comparison with related approximation results}\label{app:literature-comparison}
Table~\ref{tab:closest-comparison} records the approximation targets and
conditions of the results discussed in Appendix~\ref{app:general-comparison}.
It distinguishes general Gaussian comparison from the spectral and
estimation statements used for the present two-sample problem.

The shifted-null comparison also differs from general perturbation of a Gaussian ball probability. The latter allows changes in means and covariances \citep{GoetzeEtAl2019}, whereas here the squared signal shift is retained in both laws. Smooth quadratic-form comparison can improve the dependence on an $L^2$ perturbation under eigenvalue conditions \citep{Kosov2022}. Our spectral refinements additionally use the zero mean of an independent perturbation in the Poisson-mixture representation. These are results for a particular positive quadratic-form family, rather than general improvements for signed kernels or maxima of quadratic forms.

\begin{table}[!htbp]
\centering
\begin{threeparttable}
\caption{Comparison with related approximation results.}
\label{tab:closest-comparison}
\small
\begin{tabular}{@{}>{\raggedright\arraybackslash}p{0.22\textwidth}>{\raggedright\arraybackslash}p{0.35\textwidth}>{\raggedright\arraybackslash}p{0.37\textwidth}@{}}
\toprule
Result & Statistic and approximation & Implication here\\
\midrule
\citet{WangXu2022}, Theorems~1--2, Corollary~2
& Heteroscedastic two-sample $U$-statistic; central Gaussian quadratic reference and conditional randomization.
& Arbitrary covariance spectrum and unbalanced samples are already covered. Their stated power result assumes a negligible linear cross term.\\[4pt]
\citet{HuangEtAl2023}, Theorem~2
& I.i.d. degree-two $U$-statistics; quadratic--linear Gaussian polynomial under kernel-moment conditions.
& Non-Gaussian alternatives and the generic $n^{-1/14}$ order are existing results.\\[4pt]
\citet{HuangAusternOrbanz2025}, Theorems~1--2
& Independent-vector multilinear polynomials and $L^2$ perturbations; Gaussian polynomial target.
& With degree $2$ and moment parameter $3$, the exponents are $1/7$ for influence and $1/5$ for a remainder variance ratio. Deterministic signals are allowed through polynomial coefficients.\\[4pt]
\citet{Imai2026}, Theorem~3.1, Corollary~3.2
& Vectors of independent, nonidentically distributed weighted/two-sample $U$- and $V$-statistics; multivariate normal target.
& Allows both Hoeffding components, subject to explicit error bounds. A persistent nonnormal scalar spectral limit is a different approximation target.\\[4pt]
\citet{CaiHu2026}, Theorems~2.5, 4.2
& I.i.d. canonical $U$-statistic maxima; signed Gaussian chaos, with feasible multiplier calibration under positive-semidefinite feature and covariance-estimation conditions.
& Spectral transitions and feasible chaos calibration already coexist. Their calibration uses Hilbert--Schmidt covariance conditions; the alternative projected calibration bounds use operator-norm error and two scalar trace estimates. The augmented spherical reference avoids those estimation conditions.\\
\bottomrule
\end{tabular}
\end{threeparttable}
\end{table}

\clearpage
\setlength{\bibsep}{3pt}

\end{document}